\documentclass[12pt,bezier]{article}
\usepackage{amssymb,amsfonts,amsthm}
\newcounter{ppp}
\newcounter{band}
\newcounter{annulus}
\newcounter{nonannulus}

\newcounter{pdtwo}

\newcommand{\pr}{\pageref}
\newcommand{\la}{\langle}
\newcommand{\ra}{\rangle}

\newcommand{\diff}{{\rm diff}}
\newcommand{\Lab}{\phi}
\newcommand{\gr}{{\cal G}}
\newcommand{\bgr}{{\bar{\cal G}}}
\newcommand{\bm}{{\bar m}}
\newcommand{\bee}{{\bar\eee}}
\newcommand{\hhh}{{\cal H}}
\newcommand{\tkk}{{\tilde\kkk}}
\newcommand{\tsigma}{{\tilde\Sigma}}
\newcommand{\tool}{\stackrel{\ell}{\too} }
\newcommand{\bsss}{{\bar\sss}}
\newcommand{\csss}{{\sss\cup\bsss}}
\newcommand{\aba}[1]{{\aaa(#1)\cup\bar\aaa(#1)}}
\newcommand{\ata}[1]{{\Theta(#1)\cup\bar\Theta(#1)}}

\newcommand{\mod}{{\hbox{mod }}}
\newcommand{\br}{{\hbox{br}}}

\newcommand{\ttt}{{\cal T}}
\newcommand{\kkk}{{\cal K}}
\newcommand{\eee}{{\cal E}}
\newcommand{\aaa}{{\cal A}}

\newcommand{\CC}{{\cal C}}

\newcommand{\LL}{{\cal L}}

\newcommand{\bb}{{\cal B}}
\newcommand{\topp}{{\bf top}}
\newcommand{\bott}{{\bf bot}}
\newcommand{\vk}{van Kampen }

\newcommand{\ccc}{{\cal C}}

\newcommand{\iv}{^{-1}}

\newcommand{\lel}{\!{\overleftarrow{L}\!}}
\newcommand{\rer}{\!{\overrightarrow{R}}\!}

\newcommand{\too}{\to }

\newcommand{\rr}{{\cal R} }
\newcommand{\rrr}{{\cal R} }
\newcommand{\xxx}{{\cal X} }

\newcommand{\kk}{{\cal K} }

\newcommand{\ww}{{\cal W} }

\newcommand{\sss}{{\cal S} }
\newcommand{\mm}{{\cal M} }

\begin{document}
\newtheorem{theo}{\quad Theorem}[section]
\newtheorem{lemma}{\quad Lemma}[section]
\newtheorem{cy}{\quad Corollary}[section]
\newtheorem{df}{\quad Definition}[section]
\newtheorem{rk}{\quad Remark}[section]
\newtheorem{prop}{\quad Proposition}[section]
\newtheorem{prob}{\quad Problem}
\renewcommand{\theequation}{\thesection.\arabic{equation}}

\title{The Conjugacy Problem and Higman Embeddings}
 \author{A.Yu. Ol'shanskii, M.V. Sapir\thanks{Both authors were supported in part by
the NSF grant DMS 0072307. In addition, the research of the first
author was supported in part by the Russian Fund for Basic
Research 99-01-00894 and by the INTAS grant 99-1224,  the research
of the second author was supported in part by the NSF grant DMS
9978802 and the US-Israeli BSF grant 1999298.}}
\date{}
\maketitle

\begin{abstract} For every finitely generated recursively presented group $\gr$ we construct a finitely presented group $\hhh$
containing $\gr$ such that $\gr$ is (Frattini) embedded into
$\hhh$ and the group $\hhh$ has solvable conjugacy problem if and
only if $\gr$ has solvable conjugacy problem. Moreover $\gr$ and
$\hhh$ have the same r.e. Turing degrees of the conjugacy problem.
This solves a problem by D. Collins.
\end{abstract}

\tableofcontents

\section{Introduction}

In 1961, G. Higman \cite{Hi} published the celebrated theorem that
a finitely generated group is recursively presented if and only if
it is a subgroup of a finitely presented group. Along with the
results of Novikov \cite{No} and Boone \cite{Bo} this result
showed that objects from logic (in that case, recursively
enumerable sets) have group theoretic characterizations (see Manin
\cite{Ma} for the philosophy of Higman embeddings).

Clapham \cite{Cla} (see also corrections in \cite{Va}) was
probably the first to investigate properties preserved under
Higman embeddings. In particular, he slightly modified the
original Higman construction and showed that his embedding
preserves solvability (and even the r.e. Turing degree) of the
word problem. Valiev \cite{Va} proved a much stronger result:
there exists a version of Higman embedding theorem that preserves
the truth-table degree of the word problem. In particular, if the
word problem of a finitely generated group is in P (i.e. can be
solved in polynomial time by a deterministic Turing machine) then
it can be embedded into a finitely presented group whose word
problem is also in P (a correction to Valiev's paper was published
in Mathematical reviews, review 54 \# 413, see also Lyndon and
Schupp \cite{LS} and Manin \cite{Ma}). A simplified version of
Valiev's proof was later published by Aanderaa and Cohen in
\cite{AC} (see also Kalorkoti \cite{Kal}). In \cite{BORS}, Birget,
Rips and the authors of this paper obtained a group theoretic
characterization of NP (non-deterministic polynomial time): a
finitely generated group $G$ has word problem in NP if and only if
$G$ is embedded into a finitely presented group with polynomial
Dehn function.

The conjugacy problem turned out to be much harder to preserve
under embeddings. Collins and Miller \cite{CM} and Gorjaga and
Kirkinski\u\i  \cite{GK} proved that even subgroups of index 2 of
finitely presented groups do not inherit solvability or
unsolvability of the conjugacy problem.

In 1976 D. Collins \cite{KT} posed the following question (problem
5.22): {\em Does there exist a version of the Higman embedding
theorem in which the degree of unsolvability of the conjugacy
problem is preserved?}

It was quickly realized that the main problem would be in
preserving the smallest Turing degree, that is the solvability of
conjugacy problem. In 1980, \cite{Collins}, Collins analyzed
existing proofs of Higman's theorem, and discovered that there are
essential difficulties. If a finitely generated group $C$ is
embedded into a finitely presented group $B$ using any existing at
that time constructions then the solvability of conjugacy problem
in $B$ implies certain properties of $C$ which are much stronger
than solvability of conjugacy problem. Collins even wrote the
following pessimistic comments: ``There seems at present to be no
hope to establishing the analogue of Clapham's theorem. ...
Furthermore these difficulties seem to be more or less inevitable
given the structure of the proof and probably a wholly new
strategy will be needed to avoid them. For the present the most
one can be hoped for is the isolation of conditions on $C$ that
are necessary and sufficient for the preservation of the
solvability of the conjugacy problem in the Higman embedding."

In this paper, we do find a ``new strategy" and prove the
following result. Recall (see R. Thompson \cite{Tho}) that a
subgroup $\gr$ of a group $\hhh$ is called \label{Frattinif}{\em
Frattini embedded} if any two elements of $\gr$ that are conjugate
in $\hhh$ are also conjugate in $\gr$. Clearly if $\gr$ is
Frattini embedded into $\hhh$ and $\hhh$ has solvable conjugacy
problem then $\gr$ has solvable conjugacy problem (by results from
\cite{CM} and \cite{GK}, cited above non-Frattini embedded
subgroups do not necessarily inherit solvability of the conjugacy
problem).

\begin{theo} \label{Col} A finitely generated group has solvable conjugacy
problem if and only if it is Frattini embedded into a finitely
presented group with solvable conjugacy problem.
\end{theo}

Moreover we prove the following stronger result solving the
original problem of Collins:

\begin{theo} \label{Col1}  A finitely generated group $\gr$ with
recursively enumerable set of relations has conjugacy problem of
r.e. Turing degree $\alpha$ if and only if $\gr$ can be Frattini
embedded into a finitely presented group $\hhh$ with conjugacy
problem of r.e. Turing degree $\alpha$.
\end{theo}

In the forthcoming paper \cite{OlSa}, we will present some
corollaries of these theorems (we did not include the proofs of
them here in order not to increase the difficulty level
unnecessarily). In particular, we will show that one can drop the
restriction that $\gr$ is finitely generated in Theorems
\ref{Col}, \ref{Col1}, replacing it with ``countably generated".
We shall also show that a finitely generated group with solvable
word, power and order problems can be embedded into a finitely
presented group with solvable word, order, power, and conjugacy
problems. Thus a Higman embedding can greatly improve algorithmic
properties of a group.

In this section, we explain the main ideas which lead us to our
construction. We try to keep notation as simple as possible in
this section of the paper, so the notation here does not
necessarily coincide with the notation in the rest of the paper.

We consider the case when \label{eee}$\gr=\langle A\ |\
\eee\rangle$ has decidable conjugacy problem and we want to embed
it into a finitely presented group with decidable conjugacy
problem (the case when the conjugacy problem in $\gr$ has
arbitrary r.e. Turing degree is similar). The standard idea is to
start with a machine which recognizes the set $\eee$ and then
interpret this machine in a group. Of course we would like to have
a machine which is easier to interpret.  The most suitable
machines for our purposes are the so called $S$-machines invented
by the second author in \cite{SBR} and successfully used in
\cite{SBR,BORS,OlSaBurns,talk,OSamen}. An $S$-machine works with
words which can have complicated structure determined by the
problem. Different parts of these words can be elements of
different groups (so we do not distinguish words whose
corresponding parts are equal in the corresponding groups). As in
Turing machines, every rule of an $S$-machine replaces subwords
containing state letters by other subwords. For an exact
definition of $S$-machines see Section \ref{smach}. For now the
reader can imagine just an ordinary Turing machine which works
with words from the free group instead of a free monoid.

Here we give just one example of an $S$-machine which essentially
goes back to C. Miller \cite{Mil} (although Miller did not use
$S$-machines). The real $S$-machine used in this paper is a
``descendant" of this $S$-machine (obtained from Miller's machine
by several ``mutations"). Assume that $A$ is a symmetric set that
is it contains the formal inverses of its elements, and let us
embed $\gr=\langle A\ |\ \eee\rangle$ into {\em some} finitely
presented group \label{bargr}$\bar \gr=\langle A\cup Y\ |\
\bee\rangle$ using, say, the standard Higman embedding. Now
consider the $S$-machine $\mm$ with one state letter $q$ and the
tape alphabet $A\cup Y$ considered as a generating set of the free
group, and the following commands (each command says which
subwords of a word should be replaced and the replacement word):

\begin{itemize}
\item $[q\to a\iv q a]$, for every $a\in A\cup Y$, \item $[q\to
rq]$, for every $r\in \bee$.
\end{itemize}
The first set of rules allows the state letter to move freely
along a word of the form $uqv$ where $u,v$ are words from the free
group on $A\cup Y$. Such words are called {\em admissible words}
for $\mm$.

The second set of rules allows us to insert any relation from
$\bee$ into the word (we are also allowed to insert and remove
subwords of the form $aa\iv$ because $A\cup Y$ is a generating set
of a free group).

It is easy to see that a word $qu$, where $u$ is a word over $A$,
is recognized by this machine (that is the machine takes it to the
word q) if and only if $u=1$ in the group $\bar\gr$, that is if
and only $u=1$ in $\gr$ (because $\gr$ is embedded into $\bar\gr
$).\footnote{Notice that by a result of Sapir \cite{talk}, for
every Turing machine $T$ there exists an $S$-machine $\mm$
corresponding a to a group $\bar\gr$ as above which is
polynomially equivalent to $T$. This shows that machines $\mm$ are
as powerful as ordinary Turing machines.}

Now let us present the ideas of a Higman embedding based on $\mm$.
Let $T_1, T_2,...,T_m=W_0$ be an accepting {\em computation} of
$\mm$, that is a sequence of admissible words such that $T_{i+1}$
is obtained from $T_i$ by applying an $S$-rule of $\mm$, the last
word in this sequence is equal to $q$ in the free group.

Suppose that modulo some set of group relations $Z$ we have for
every $i=1,...,m-1$

$$P_iT_i=T_{i+1}P_i'$$ for some words $P_i$
and $P_i'$ over a certain alphabet (it is included in the finite
generating set of the group we are constructing), that is we have
a \vk diagram with the boundary label $P_iT_i(P_i')\iv T_{i+1}\iv$
(Figure \theppp):

\bigskip

\unitlength=1mm \special{em:linewidth 0.5pt} \linethickness{0.5pt}
\begin{picture}(96.00,21.67)
\put(30.33,6.00){\framebox(62.67,12.33)[cc]{}}
\put(26.67,12.00){\makebox(0,0)[cc]{$P_i$}}
\put(60.00,2.67){\makebox(0,0)[cc]{$T_i$}}
\put(96.00,12.33){\makebox(0,0)[cc]{$P_i'$}}
\put(59.67,21.67){\makebox(0,0)[cc]{$T_{i+1}$}}
\put(30.33,14.67){\vector(0,-1){4.00}}
\put(93.00,14.67){\vector(0,-1){4.33}}
\put(52.33,6.00){\vector(1,0){11.33}}
\put(52.33,18.33){\vector(1,0){11.00}}
\end{picture}

\begin{center}
\nopagebreak[4] Figure \theppp.
\end{center}
\addtocounter{ppp}{1}

\noindent (the set $Z$ describes the method of tessellating this
rectangle into cells). Then we have that $PT_1=W_0P'$ where
$P=P_{m-1}...P_2P_1$, $P_i'=P_{m-1}'...P_2'P_1'$. The
corresponding \vk diagram $\Delta$ has the form of a {\em
trapezium}:

\bigskip

\unitlength=1mm \special{em:linewidth 0.5pt} \linethickness{0.5pt}
\begin{picture}(109.67,59.33)
\put(30.33,59.33){\line(2,-5){22.00}}
\put(52.33,4.67){\line(1,0){35.67}}
\put(88.00,4.67){\line(2,5){21.67}}
\put(109.67,59.00){\line(-1,0){79.33}}
\put(50.00,10.67){\line(1,0){40.33}}
\put(47.33,17.00){\line(1,0){45.33}}
\put(44.33,24.00){\line(1,0){51.33}}
\put(41.67,32.00){\line(1,0){57.33}}
\put(38.33,40.00){\line(1,0){63.33}}
\put(35.00,48.33){\line(1,0){70.33}}
\put(69.00,1.67){\makebox(0,0)[cc]{$T_1$}}
\put(69.00,8.33){\makebox(0,0)[cc]{$T_2$}}
\put(69.00,14.33){\makebox(0,0)[cc]{$T_3$}}
\put(48.00,7.00){\makebox(0,0)[cc]{$P_1$}}
\put(45.00,14.00){\makebox(0,0)[cc]{$P_2$}}
\put(95.00,14.00){\makebox(0,0)[cc]{$P_2'$}}
\put(92.00,7.33){\makebox(0,0)[cc]{$P_1'$}}
\put(42.67,29.33){\vector(1,-3){1.67}}
\put(97.67,29.00){\vector(-1,-3){1.33}}
\put(57.00,59.00){\vector(1,0){15.00}}
\put(64.67,4.67){\vector(1,0){7.67}}
\put(70.67,55.00){\makebox(0,0)[cc]{$T_m$}}
\end{picture}

\begin{center}
\nopagebreak[4] Figure \theppp.
\end{center}
\addtocounter{ppp}{1}

For example, if we use the $S$-machine $\mm$ described above then
the presentation $Z$ is easy to describe: for every rule
$\tau=[q\to uqv]$ we have the following relations in $Z$:
$\theta_\tau\iv q\theta_\tau=uqv$, $\theta_\tau a=a\theta_\tau$
for every $a\in A\cup Y$ (the letters $\theta_\tau$, corresponding
to the rules of $\mm$ are added to the generating set). In this
case the words $P_i, P'_i$ are of length 1, and are equal to
$\theta_\tau$.

The words $P$ and $P'$ contain the {\em history} of our
computation, because $P_i$ and $P_i'$ correspond to $S$-rules
applied during the computation. The words $T_1$ and $T_m$ are
labels of the top and the bottom paths of the trapezium.

If the set of relations $Z$ is chosen carefully (and there are
very many ways of doing this, see \cite{SBR,
BORS,talk,OlSaBurns,OSamen}), then the converse is also true:
\medskip

{\bf Condition 1.} If there exists a trapezium with the top
labelled by $T_1$, the bottom labelled by $W_0$ and the sides
labelled by words over $\Theta$ (the set of all $\theta_\tau$)
then $T_1$ is admissible and the labels of the sides correspond to
the history of an accepting computation.

\medskip

Thus the word $PT_1(P')\iv W_0\iv$ written on the boundary of a
trapezium if and only if the word $T_1$ is accepted.

The next step is to consider $N\ge 1$ copies of our trapezium with
labels taken from different alphabets. For technical reasons
(similar to the hyperbolic graphs argument in \cite{SBR},
\cite{Ol97}, \cite{BORS}) we need $N$ to be even and large enough
(say, $N\ge 8$). Let us choose disjoint alphabets in different
copies of the trapezium. Then we can glue two copies of the
trapezium by using a special letter $k$ which conjugates letters
from the right side of the first trapezium with letters from the
left side of the mirror image of another copy (Figure \theppp):

\bigskip
\unitlength=1mm \special{em:linewidth 0.5pt} \linethickness{0.5pt}
\begin{picture}(104.00,59.67)
\put(29.33,4.00){\framebox(33.67,53.33)[cc]{}}
\put(67.67,4.00){\framebox(33.67,53.33)[cc]{}}
\put(59.33,4.00){\line(1,0){11.67}}
\put(29.33,9.33){\line(1,0){72.00}}
\put(29.33,14.67){\line(1,0){72.00}}
\put(29.33,20.67){\line(1,0){72.00}}
\put(29.33,27.00){\line(1,0){72.00}}
\put(29.33,33.67){\line(1,0){72.00}}
\put(29.33,39.33){\line(1,0){72.00}}
\put(29.33,45.33){\line(1,0){72.00}}
\put(29.33,50.67){\line(1,0){72.00}}
\put(58.00,57.33){\line(1,0){15.33}}
\put(40.33,57.33){\vector(1,0){10.33}}
\put(91.33,57.33){\vector(-1,0){9.67}}
\put(88.67,4.00){\vector(-1,0){8.33}}
\put(41.00,4.00){\vector(1,0){9.67}}
\put(29.33,36.33){\vector(0,-1){7.67}}
\put(63.00,36.00){\vector(0,-1){7.33}}
\put(67.67,36.33){\vector(0,-1){7.67}}
\put(101.33,37.00){\vector(0,-1){8.33}}
\put(24.67,30.33){\makebox(0,0)[cc]{$P_i$}}
\put(60.00,30.00){\makebox(0,0)[cc]{$P_i'$}}
\put(70.00,30.00){\makebox(0,0)[cc]{$P_i'$}}
\put(104.00,30.00){\makebox(0,0)[cc]{$P_i$}}
\put(44.00,1.33){\makebox(0,0)[cc]{$T_1$}}
\put(83.00,1.33){\makebox(0,0)[cc]{$T_1$}}
\put(45.67,59.67){\makebox(0,0)[cc]{$T_m$}}
\put(86.00,59.67){\makebox(0,0)[cc]{$T_m$}}
\put(65.33,59.67){\makebox(0,0)[cc]{$k$}}
\put(65.33,1.33){\makebox(0,0)[cc]{$k$}}
\end{picture}

\begin{center}
\nopagebreak[4] Figure \theppp.
\end{center}
\addtocounter{ppp}{1}
(Notice that if $P$ and $P'$ are copies of
each other then we can also glue the right side of one trapezium
to the left side of another without taking mirror images.) Of
course the conjugacy relations involving letter $k$ should be
added into the set of group relations $Z$.

Suppose that $T_m$ is equal to $q$, i.e. the word $T_1$ is {\em
accepted} by the machine. If we connect the first copy of $\Delta$
with the second copy, then with the third copy, and finally
connect the $N$-th copy with the first copy, using different
letters $k_2,...,k_N,k_1$,  we get an {\em annular} diagram, a
{\em ring}. The outer boundary of this ring has label $\kkk(T_1)$
of the form $k_1T_1'k_2\iv (T_1'')\iv ...$ (here we use the fact
that $N$ is even). The words $T_1', ..., T_1^{(N)}$ are different
copies of $T_1$. The inner boundary has label $k_1q'k_2\iv
(q'')\iv ...$ (this word is called the {\em hub}) (Figure \theppp\
shows the diagram in the case $N=4$):

\bigskip

\unitlength=1mm \special{em:linewidth 0.5pt} \linethickness{0.5pt}
\begin{picture}(99.33,59.00)
\put(68.17,34.33){\oval(62.33,49.33)[]}
\put(69.67,33.50){\oval(15.33,13.67)[]}
\put(75.33,38.67){\line(1,1){20.00}}
\put(76.67,36.00){\line(1,1){21.00}}
\put(74.00,28.00){\line(5,-4){21.33}}
\put(97.67,12.33){\line(-6,5){21.33}}
\put(62.67,30.33){\line(-3,-2){25.33}}
\put(64.67,28.33){\line(-3,-2){25.33}}
\put(62.67,36.33){\line(-4,3){25.33}}
\put(64.00,38.67){\line(-4,3){24.67}}
\put(64.00,59.00){\vector(1,0){13.00}}
\put(99.33,28.33){\vector(0,1){13.00}}
\put(65.33,9.67){\vector(1,0){12.00}}
\put(37.00,26.00){\vector(0,1){15.33}}
\put(69.67,33.00){\makebox(0,0)[cc]{hub}}
\end{picture}

\begin{center}
\nopagebreak[4] Figure \theppp.
\end{center}
\addtocounter{ppp}{1}
We add the hub to the list of defining
relations $Z$. Now with every word $T$ of the form $uqv$ where $u,
v$ are words in $A\cup Y$ we have associated the word $\kkk(T)$,
and we see that if $T$ is accepted then $\kkk(T)$ is equal to 1
modulo $Z$. Again, if $Z$ is chosen carefully then the converse is
also true:

\medskip
{\bf Condition 2.} If $\kkk(T)$ is equal to 1 modulo $Z$ then $T$
is accepted by $\sss$. Thus we have an interpretation of $\sss$ in
the group generated by the tape alphabet and the state alphabet of
$\sss$, the set $\Theta$, and the set $\{k_1,...,k_N\}$, subject
to the relations from $Z$.

\medskip

The diagram on the previous Figure is called a {\em disk}. The
trapezia forming this disk are numbered in the natural order (from
1 to $N$).

We assume that the copy of the alphabet $A$ used in the first
subtrapezium of a disk coincides with $A$ itself (the generating
set of $\gr$).

Let us denote the group given by the presentation $Z$ we have got
so far by $G(\mm)$.

There are several ways to use the group $G(\mm)$ to embed $\gr$
into a finitely presented group. We cannot use the method used by
Higman \cite{Hi}, Aanderaa \cite{AC},  Rotman \cite{Rotman}, and
us in \cite{BORS} because it leads to problems discovered by
Collins in \cite{Collins}. We use another method, described in
\cite{talk} and used also in \cite{OSamen}.

Consider $N-1$ new copies of the trapezium $\Delta$. Number these
trapezia by $2',...,N'$. Identify the copy of $A$ (but not $Y$)
and the state letter $q$ in the trapezium number $i'$ with the
copy of $A$ and the copy of $q$ in the ``old" trapezium number
$i$. Now glue the trapezium number $2'$ with the trapezium number
$3'$ using letter $k_3$ (as before), then glue in the trapezium
number $4'$ and so on, but glue the trapezium number $N'$ with the
trapezium number $2'$ using $k_1qk_2\iv$ (that would be possible
because as it will turn out in that case the words $P, P'$ will be
copies of each other). That is conjugation of the letters on the
side of the trapezium number $N'$ by $k_1qk_2\iv$ gives the
letters of the side of the mirror image of the trapezium number
$2'$. Let us add the relations used in the new trapezia and the
conjugation relations (involving $k_i$) to $Z$.

\vskip 2 in

The resulting picture is an annular diagram. The inner boundary of
it is labelled by the hub. If we add the hub to the diagram, we
get a \vk diagram which also looks like a disk but has N-1
sectors. The outer boundary of this new disk is labelled by the
word $k_1qk_2\iv V$ where $V$ is the suffix of $\kkk(qu)$ staying
after $k_2\iv$. Thus the word $k_1qk_2\iv V$ and the word
$k_1quk_2\iv V$ are both equal to 1 modulo $Z$, so $Z$ implies
$u=1$. Let us call the group given by the set of relations $Z$
that we have constructed by $\hhh$.

We have proved that the identity map on $A$ can be extended to a
homomorphism from $\gr$ to the subgroup $\langle A\rangle$ of the
group $\hhh$ given by the set of relations $Z$. It is possible to
show that this homomorphism is injective, so $\gr$ is embedded
into $\hhh$.

Suppose that $\gr$ has solvable conjugacy problem. Is it true that
$\hhh$ has solvable conjugacy problem? The answer is ``not
always". Let us give two examples.

{\bf Example 1}. Consider two pairs of words $u_1, u_2$ and $v_1,
v_2$ over the alphabet $A$ (we can view these words as elements of
$\gr$). Let $q$ be the first state letter in $\kkk(qu)$. Consider
the words $W_1=q\iv u_1qu_2q\iv$ and $W_2=q\iv v_1qv_2q\iv$.
Suppose that there exists a word $W(A)$ in the alphabet $A$ such
that $W(A)u_1W(A)\iv =v_1, W(A)u_2W(A)\iv=v_2$. Let
$W(A)=a_1a_2...a_k$, $a_i\in A$. For every $a\in A$ let
$\theta(a)$ correspond to the $S$-rule of the form $[q\to a\iv
qa]$. Consider the word
$W(\Theta)=\theta(a_1)\theta(a_2)...\theta(a_k)$. Then it is easy
to see that $W(\theta)W_1W(\theta\iv)=W_2$. Thus if the pair
$(u_1,u_2)$ is conjugated to the pair $(v_1,v_2)$ in $\gr$ then
the words $W_1$ and $W_2$ are conjugated in $\gr$. One can also
prove that the converse statement is also true. In this example,
pairs of words can be obviously replaced by any $t$-tuples of
words, $t\ge 2$. Therefore if the conjugacy problem is solvable in
$\gr$ then the conjugacy of $t$-tuples of elements in $\gr$ is
solvable. It is known \cite{Collins} that there exists a finitely
generated group $\gr$ with solvable conjugacy problem and
unsolvable problem of conjugacy for sequences of elements. For
such a group $\gr$ the group $\hhh$ has undecidable conjugacy
problem.

{\bf Example 2.} Take two pairs of words $(u_1,u_2)$, $(v_1, v_2)$
over one of the copies of the alphabet $A$, say, $A'\ne A$, let
$q'$ be the corresponding copy of the letter $q$. Then it is
possible to show that the words $q'u_1(q')\iv u_2q'$ and
$q'v_1(q')\iv v_2q'$ (where $u_1, u_2, v_1, v_2$ are words in
$A'$) are conjugated in $\hhh$ if and only if for some words $P,Q$
$$Qu_1Q\iv=v_1, P\iv u_2P=v_2$$ in the free group and $PQ=1$ in
$\gr(A')$, a copy of $\gr$. This problem easily reduces
(exercise!) to the following problem about group $\gr$: {\em given
four words $s,t,p,r$ in the alphabet $A$ find two integers $m,n$
such that $st^n=p^mr$ in $\gr$.} It is quite possible that this
problem can be undecidable even if the conjugacy problem is
decidable.

These two examples suggest that we need to change the machine
$\mm$: we cannot allow trapezia with top path labelled by words of
the form $qu_1q\iv u_2q$ having long histories. In other words, if
a computation of the $S$-machine starts with $qu_1q\iv u_2q$, then
it should not be too long (more precisely, if it is too long, we
should be able to replace it with a shorter computation). If we
achieve that then we would have a recursive bound on the length of
the word $W$ in Example 1 and the words $P, Q$ in Example 2, so we
would be able to find these words by a simple search.

This lead us to the $S$-machine $\sss$ considered in this paper
(see Sections \ref{hardware},\ref{notation},\ref{rules}). This
machine has four sets of state letters: $K, L, P, R$-letters.
These letters are placed in an admissible word $\kkk(u)$ according
to the following pattern: $KLPRK\iv R\iv P\iv L\iv KLPR...$.

The $P$-letters play the role of the many copies of $q$ in
$\kkk(qu)$ above: they moves along an admissible word and ``find"
places where to insert relations from $\bee$. The $K$-letters are
the letters $k_i$ in $\kkk(qu)$ (they divide admissible words into
copies of the same word and do not move). The $L$- and $R$-letters
are ``support letters". While $P$ is moving, $L, R$ stay next to
$K$-letters. But when $P$ ``wants" to insert a relation, $L, R$
must move next to $P$ and they insert the relation ``together"
(executing the rule $LPR\to rLPR$, $r\in \bee$). After that $L, R$
must move back to their initial positions, and $P$ can move to a
new place in the admissible word (the work of the machine
resembles a complicated dance performed synchronously by several
groups of dancers).

It is important to mention that we cannot allow all parts of the
admissible words to be arbitrary group words. To avoid long
computations without inserting new relation from $\bee$, we had to
require that the words between neighbor $K$-letter and $L$-letter,
$L$-letter and $P$-letter, $R$-letter and $K$-letter (but not
between $P$-letter and $R$-letter!) to be positive (i.e. they
cannot contain letters $a\iv$).

In order to keep these subwords positive we use a trick that goes
back to Novikov, Boone and Higman (see Rotman \cite{Rotman}). They
used a special letter $x$ and Baumslag-Solitar relations of the
form $a\iv xa=x^2$ for every $a$ in the tape alphabet of the
machine. The idea is that if we have relations $a\iv xa=x^2$ for
all $a\in A$ and a word $u\iv xu$, where $u$ is a reduced word in
$A$, is equal modulo these relations to a power of $x$ then the
first letter of $u$ is positive.

We use a similar idea. But our task is more complicated than in
\cite{No,Bo,Hi} because we need to consider conjugacy of arbitrary
pairs of words, not only those that have ``nice" structure. So we
had to use many different letters $x$, replace $2$ by $4$ in the
Baumslag-Solitar relations, and make some other technical
modifications. Notice that one problem with using the $x$-letters
and Baumslag-Solitar relations is that since $N-1$ is odd we
cannot use $x$-letters in the second disk described above (indeed,
otherwise it would be impossible to glue the $N-1$ sectors
together since the words $P$ and $P'$ would not be copies of each
other). Thus we need to consider two $S$-machines (one responsible
for the first disk, another for the second disk) which are very
similar but have different ``hardware": the admissible words of
the first machine must have positive parts described above, the
admissible words of the second machine do not have to satisfy this
restriction.

Notice that Examples 1 and 2 are only the main obstacles that we
had to overcome in this paper. Other technical difficulties lead
to further fine tuning of our $S$-machine. Sections
\ref{convert1}, \ref{convert2}, \ref{convert3} contain precise
description of the presentation of our group $\hhh$.

Now let us briefly describe the strategy of the proof that the
conjugacy problem in $\hhh$ is Turing reducible to the conjugacy
problem in $\gr$. As in our previous papers the main tool in this
paper is bands (other people call them ``strips", ``corridors",
etc.).

In terms of annular (Schupp) diagrams, our task is the following:
given two words $u$ and $v$ in generators of $\hhh$, find out if
there exists an annular diagram over the presentation of $\hhh$
with boundary labels $u$ and $v$. It turns out that any annular
diagram over $\hhh$ (after removing some parts of recursively
bounded size) becomes a diagram of one of three main types: a {\em
ring} (where the boundary labels contain $K$-, $L$-, $P$-,
$R$-letters but do not contain $\theta$-letters, all maximal
$\theta$-bands are annuli surrounding the hole of the diagram), a
{\em roll} (where the boundary label does not contain $K$-, $L$-,
$P$-, $R$-letters but contain $\theta$-letters; all maximal $K$-,
$L$-, $P$-, $R$-bands are annuli surrounding the hole), and a {\em
spiral} (where the boundary labels contain both $K$-, $L$-, $P$-,
or $R$-letters and $\theta$-letters; both $K$-, $L$-, $P$-,
$R$-bands and the $\theta$-bands start and end on the different
boundary components, and each $\theta$-band crosses each $K$-,
$L$-, $P$-, and $R$-band many times\footnote{In that case the
$\theta$-band looks like a spiral on a plane that starts at the
origin and crosses many times straight lines (the $K$-, $L$-, $P$-
and $R$-bands.}). Different methods are used to treat different
cases. Roughly speaking, the study of rings amounts to study the
lengths of computations of our $S$-machines, the study of rolls
amounts to the study of the space complexity (how much space is
needed by the machines during a computation), and the study of
spirals amounts to the study of computations with periodic
history. The $x$-letters and the Baumslag-Solitar relations allow
us to treat the case of rings, but they cause the main technical
difficulties in the cases of rolls and spirals.

\section{List of relations}
\setcounter{equation}{0}
\medskip

\subsection{$S$-machines}
\label{smach}

Here we give the general definition of $S$-machines (see
\cite{talk}). The precise definition of the $S$-machines used in
this paper will be given later.

Let $T_1,...,T_k$ be a collection of groups, $T_i=\la A_i\ra$,
$i=1,...,k$, and the sets $A_i$ of generators pairwise disjoint.
Let $\kkk$ be a set of \label{State}{\em state letters}, disjoint
from $\cup A_i$. Let ${\bf A}$ be a set of \label{Admissible}{\em
admissible} words of the form $q_1w_1q_2...q_kw_kq_{k+1}$ where
$w_i$ is a word over some $A_j$, $q_i\in \kkk\cup\kkk\iv$.

We do not distinguish words over $A_i$ which are equal in the
group $T_i$. So we can view $w_i$ as elements in $T_i$.

Any \label{Srule}{\em $S$-rule} $\tau$ has the following form
$[k_1\too v_1k_1'u_1,...,k_n\too v_nk_n'u_n]$ where $u_j, v_j$ are
group words over $\cup A_i$, $k_j\in \kkk\cup\kkk\iv$. Every
$S$-rule is a partial transformation on the set of admissible
word. Every admissible word $W$ in the domain of the $S$-rule
$\tau$ must have $\kkk$-letters only from the set $\{k_1^{\pm
1},...,k_n^{\pm 1}\}$ (not all of these letters may occur, and
some of them may occur more than once). This transformation takes
an admissible word $W$, and replaces each $k_j^{\pm 1}$ by
$(v_jk_j'u_j)^{\pm 1}$, freely reduces the resulting word and then
trims maximal prefix and suffix of the resulting word consisting
of non-$\kkk$-letters. For example, if $W\equiv k_1p_1k_1p_2k_2\iv
p_3k_3$ for some words $p_1,p_2,p_3$ over $\cup A_i$ (we assume
$n\ge 3$) is in the domain of $\tau$, then the result of the
application of $\tau$ is the freely reduced form of
$k_1'u_1p_1v_1k_1'u_1p_2u_2\iv k_2\iv v_2\iv p_3v_3k_3'$. Notice
that the resulting word must be an admissible word for the
$S$-machine as well, otherwise $W$ is not in the domain of $\tau$.

The \label{Inverserule}{\em inverse} rule $\tau\iv$ is the inverse
partial transformation. It has the form  $[k_1'\too v_1\iv
k_1u_1\iv,...,k_n'\too v_n\iv k_nu_n\iv]$.

Thus we can view an $S$-machines as a semigroup of partial
transformations of the set of admissible words. If $W$ is an
admissible word of an $S$-machine $\sss$, and
$h=\tau_1\tau_2...\tau_n$ is a sequence (word) of rules of $\sss$,
applicable to $W$ (that is $W$ is in the domain of the partial map
$h$), then the result $W'$ of application of $h$ to $W$ is denoted
by $W\circ h$. The corresponding computation $W, W\circ \tau_1,
W\circ \tau_1\tau_2,...,W'$ will be denoted by $W\bullet h$ or
$W\bullet h=W'$.

\label{Smachines}The collection of groups $T_i$ and the set of
admissible words form the {\em hardware} of an $S$-machine. The
{\em software} of an $S$-machine is a finite collection of
$S$-rules closed under taking inverses.

\subsection{An auxiliary group} Let \label{grrr}$\gr$ be a group
given by the following
presentation
\begin{equation} \gr =\la a_1,\dots, a_m | w=1, w\in \eee\ra,
\label{3.1}\end{equation} where $\eee$  is  a   recursively
enumerable set of words in $a_1,\dots,a_m$. Adding if necessary,
new generators and relations of the form $a_ia'_i =1$ one can
assume that the set $\eee$ consists of positive words only, i.e.,
there are no occurrences of $a_i^{-1}$ in the words $w\in \eee$.

We first prove Theorem \ref{Col} and then show how to modify the
proof to obtain Theorem \ref{Col1}. So we shall assume that $\gr$
has solvable conjugacy problem.

By the Higman Embedding Theorem \cite{Rotman} there is a finitely
presented group \label{bargr1}$\bgr$ generated by
$a_1,...,a_m,...,a_\bm$ containing $\gr$ as a subgroup. By
Clapham-Valiev's result \cite{Cla,Va}, one can assume that the
word problem in $\bgr$ is decidable.

Besides, one may assume that again, for each generator $a$ of
$\bgr$, the inverse letter $a'$ is also included in the set of
generators $\{a_1,...,a_m,...,a_\bm\}$ of $\bgr$. For every $a\in
\{ a_1,...,a_\bm\}$, we assume that there are positive relators of
the form $aa'$ and $a'a$ in the finite set \label{beee}$\bee$ of
positive defining relators for $\bgr$. We will also suppose that
if $r\in \bee$ then $r'\in \bee$ where $r'$ is obtained from
$r^{-1}$ by replacing every occurrence of a letter $a\iv$ by the
letter $a'$. Finally, we will assume that $\bee$ contains the
empty word $\emptyset$. It is clear that $\bee$ is a presentation
of the group $\bgr$ in the class of all monoids, so we have the
following

\begin{lemma} \label{Positive}
Assume that $w_1=w_2$ for positive group words $w_1,w_2$ in the
generators of $\bgr$. Then there is a sequence of positive words
starting with $w_1$ and ending with $w_2$, where every word in the
sequence is obtained from the previous word by insertion or
deletion of a subword $w\in \bee$.
\end{lemma}

\medskip

The list of relations of the group $\hhh$ which we are going to
construct  will depend on the set $\bee$ of defining words for
$\bgr$. Lemma \ref{Emb} below claims that there is a natural
embedding of the subgroup $\gr\le \bgr$ into $\hhh$, but
(caution!) we will not embed the whole group $\bgr$ into $\hhh$.

\medskip

\subsection{The hardware of $\sss$}
\label{hardware}

The group $\hhh$ that we are going to construct is associated with
two very similar $S$-machines, $\sss$ and $\bsss$. We shall
describe the $S$-machines, then we will describe how to convert
these machines into a group presentation.

We fix an even number $N\ge 8$.

The set \label{kkk}$\kkk$ of state letters of $\sss$ consists of
letters $z_j(r,i)$ where  $z\in \{K,L,P,R\}$, $j=1,...,N$, $r\in
\bee$, $i\in\{1,2,3,4,5\}$.

We also define the set of \label{Basicletters}{\em basic} letters
\label{tkk}$\tkk$ which consists of letters $z_j$ where $z\in
\{K,L,P,R\}$, $j=1,...,N$. Letters from $\tkk\iv$ are also called
{\em basic} letters. There exists a natural map from
$\kkk\cup\kkk\iv$ to $\tkk\cup\tkk\iv$ which forgets indexes $r$
and $i$. If $z$ is a basic letter, $r\in \bee$, $i\in
\{1,2,3,4,5\}$ then $z(r,i)\in \kkk$. If $U$ is a word in $\tkk$
and other letters, $r\in \bee$, $i\in \{1,2,3,4,5\}$ then
\label{uri}$U(r,i)$ is a word obtained from $U$ by replacing every
letter $z\in\tkk$ by $z(r,i)$. The parameters $r$ and $i$ in the
letter $z(r,i)$ or in the word $U(r,i)$ will be called the
$\bee$-coordinate and the $\Omega$-coordinate of the word.

The set \label{aaa}$\aaa$ of \label{Tapeletters}{\em tape letters}
of the machine $\sss$ consists of letters $a_i(z)$ where
$i=1,...,\bm$, $z\in \tkk$. For every $z\in \tkk$ we define
\label{aaaz}$\aaa(z)$ as the set of all $a_i(z)$, $i=1,...,\bm$.

Let \label{tsigmaaaa}$\tsigma$ be the following word (considered
as a cyclic word):

\begin{equation}
\begin{array}{l}
K_1L_1P_1R_1K_2\iv R_2\iv P_2\iv L_2\iv K_3L_3P_3R_3 K_4\iv\\
... K_N\iv R_N\iv P_N\iv L_N\iv \end{array} \label{basehub}
\end{equation}
Notice that for every basic letter $z$ precisely one of $z$ and
$z\iv$ occurs in the word $\tsigma$. The word
\label{sigmaaaa}$\Sigma=\tsigma(\emptyset,1)$ will be called the
\label{hub}{\em hub}.

For every $z\in \tkk\cup\tkk\iv$ by \label{zmin}$z_-$ we denote
the letter immediately preceding $z$ in the cyclic word $\tsigma$
or in $\tsigma\iv$. This definition is correct because every basic
letter occurs exactly once in the word $\tsigma$ or its inverse
(see remark above). If $z'=z_-$ then we set \label{zplus}$z=z'_+$.
Similarly we define $z_-$ and $z_+$ for $z\in \kkk$. Notice that
for every $j=1,...,N$, $$(L_j)_-=\left\{\begin{array}{ll} K_j &
\hbox{ if } j \hbox{ is odd,}\\ K_{j+1}\iv &\hbox{ if } j \hbox{
is even.}\end{array}\right.$$ and
$$(R_j)_+=\left\{\begin{array}{ll} K_{j+1}\iv & \hbox{ if } j
\hbox{ is odd,}\\ K_j &\hbox{ if } j \hbox{ is
even.}\end{array}\right.$$ To simplify the notation and avoid
extra parentheses, we shall denote $(L_j)_-$ by
\label{lel}$\lel_j$ and $(R_j)_+$ by \label{rer}$\rer_j$. We also
define $\aaa(z\iv)$ for $z\in\tkk$ by setting
$\aaa(z\iv)=\aaa(z_-)$.

The language of \label{adm}{\em admissible words} of the machine
$\sss$ consists of all reduced words of the form $W\equiv
y_1u_1y_2u_2...y_tu_ty_{t+1}$ where $y_1,...,y_{t+1}\in \kkk^{\pm
1} $, $u_i$ are words in $\aaa(y_i)$, $i=1,2,...,t$, and for every
$i=1,2,...,t$, either $y_{i+1}\equiv (y_i)_+$ or $y_{i+1}\equiv
y_i\iv$. (Here \label{equiv}$\equiv$ is the letter-for-letter, or
graphical, equality of words.) The projection of $y_1...y_{t+1}$
onto $\tkk$ will be called the \label{baseadm}{\em base} of the
admissible word $W$. The subword $y_iu_iy_{i+1}$ is called the
\label{sector}$y_iy_{i+1}$-{\em sector} of the admissible word
$W$, $i=1,2....$. The word $u_i$ is called the
\label{innerpart}{\em inner part} of the $y_iy_{i+1}$-sector. We
assume that the inner part of any $zL_j$-, $zP_j$- or
$R_jz$-sector of $W$ and $W\iv$ is a positive word ($z\in\tkk$),
that is it does not contain $a\iv$ for any $a\in \aaa$.  Notice
that if $W$ is an admissible word for $\sss$ then $W\iv$ is an
admissible word as well.

By definition, all $\kkk$-letters in an admissible word have the
same \label{beeco}$\bee$-coordinates and
\label{omegaco}$\Omega$-coordinates. Notice that $\tsigma(r,i)$ is
an admissible word for every $r, i$.

Two admissible words $W\equiv y_1u_1y_2u_2...y_tu_ty_{t+1}$ and
$W'\equiv y_1'u_1'y_2'u_2'...y_{t'}'u'_{t'}y_{t'+1}$ are
considered {\em equal} if $W\equiv W'$.

\subsection{Notation}
\label{notation}

In this paper, the $S$-rules will have a specific form which
allows us to simplify the notation. In every rule $[k_1\too
v_1k_1'u_1,...,k_n\too v_nk_n'u_n]$, the projections $k_i$ and
$k_i'$ on $\tkk$ will be always the same, the $\bee$-coordinates
and the $\Omega$-coordinates of all state letters $k_1,...,k_n$
(resp. $k_1',...,k_n'$) will be the same. In addition, each $u_i$
(resp. $v_i$) will contain only letters from $\aaa(k_i)$ (resp.
$\aaa((k_i)_-)$), $i=1,...,n$.

For every word $W$ in the alphabet $\{a_1,...,a_{\bar m}\}$ and
every $k\in \tkk$, $W(k)$ will denote the word obtained from $W$
by substitution $a_i\to a_i(k)$, $i=1,...,\bar m$.

In any expressions like $vku$ where $v$ is a word in $\aaa(k_-)$,
$u$ is a word in $\aaa(k)$, $k\in\tkk$, we shall not mention what
alphabets these words are written in. For example the notation
$a_iP_ja_i\iv$ means $a_i(L_j)P_ja_i(P_j)\iv$, that is the two
$a_i$ are brothers taken from different alphabets.

Thus in order to simplify notation, we shall write a rule $\tau$
of the $S$-machines $\sss$ in the form $$[k_1\too v_1k_1u_1,...,
k_n\too v_nk_nu_n;r\too r', \omega\too\omega']$$ where $k_i\in
\tkk$, $v_i$ are words in $\aaa((k_i)_-)$, $u_i$ are words in
$\aaa(k_i)$. We shall denote $v_i$ by $v(\tau((k_i)_-))$, $u_i$ by
$u_\tau(k_i)$.

Some of the arrows in a rule $\tau$ can have the form
$k_i$\label{tool}$\tool$. This means that the rule $\tau$ can be
applied to an admissible word $W$ only if the inner part of
$k_ik_i'$-sectors and $k_{i''}(k_i)_+$-sectors in $W$ are empty
words and $k_i'=(k_i)_+$, $k_{i''}=k_i$ (i.e. $k_i'\ne k_i\iv$ and
$k_{i''}\ne (k_i)_+\iv$). We shall say that in this case the
$k_i(k_i)_+$-sectors are \label{locked}{\em locked} by the rule
and the rule {\em locks} these sectors. If $kk_+$-sectors are
locked by the rule, $u(k)$ and $v(k)$ must be empty. Thus if
$\tau$ locks $kk_+$-sectors then $\tau\iv$ locks these sectors as
well.

Let us expand the definition of $u_\tau(k), v_\tau(k_-)$ to the
base letters not occurring in the rule by setting that in this
case $u_\tau(k)$ and $v_\tau(k_-)$ are empty.

{\bf Thus with every rule $\tau$ and every $k\in\kk$ we associate
two words over $\aaa$: $u_\tau(k)$ and $v_\tau(k)$. We expand this
definition to $k\in \tkk\iv$ by saying that $v_\tau(k\iv)\equiv
u_\tau(k)\iv$, $u_\tau(k\iv)\equiv v_\tau(k_-)\iv$}.

By definition, to {\em apply} the rule $\tau=[z_1\too
v_1z_1u_1,..., z_s\too v_sz_su_s; r\too r', \omega\too\omega']$ to
an admissible word $W=y_1w_1y_2...w_ky_{k+1}$ means replacing each
$y_i(r,\omega)$ by $v_iy_i(r',\omega')u_i$, reducing the resulting
word, and trimming the prefix $v_\tau((y_1)_-)$ from the beginning
and the suffix $u_\tau(y_{k+1})$ from the end of the resulting
word. (Notice that $w_i$ is non-empty if $y_i\equiv y_{i+1}^{-1}$,
and therefore $y_i(r',\omega')$ will not cancel with
$y_{i+1}(r',\omega')$ in the resulting word, because the subword
between these two letters will not be empty.) The rule $\tau$ is
called \label{applicable}{\em applicable} to $W$ if:

(1) for every $z\in \tkk$ such that $\tau$ locks $zz_+$-sectors
the inner parts of $zz_+$-sectors in $W$ are empty words, and $W$
does not have $zz\iv$-sectors or $z_+\iv z_+$-sectors.

(2) the resulting word $W'$ is again admissible (that is the inner
parts of all sectors which are supposed to be negative or positive
are such).

\subsection{The rules of the $S$-machine $\sss$}
\label{rules}

The rules from $\sss^+\subset \sss$, described below, will be
called \label{pnr}{\em positive}, the inverses of these rules will
be called {\em negative}, each rule of $\sss$ will be either
positive or negative.

The set $\sss^+$ consists of 10 subsets
\label{somega}$\sss^+(\omega)$, $\omega\in \{1,12,2,23,3,34,4,
45,5, 51\}$.

The set $\sss^+(1)$ consists of rules $\tau(1,\emptyset,i)$, where
$i=1,...,\bm$. The rule $\tau(1,\emptyset,i)$ has the form

$$[\lel_j\tool \lel_j,P_j\too a_i P_ja_i\iv, R_j\tool R_j, j=1,...,N;
\emptyset\too \emptyset, 1\too 1].$$
The meaning of this set of rules is that the state letter $P_j$
can freely move (if the $(\bee, \Omega)$-coordinates are
$(\emptyset,1)$, and $L_j$-letters and $R_j$-letters stay next to
$K_j$-letters). Indeed, if in an admissible word, the
$\aaa$-letter next to the right of $P_j$ is $a_i$ (and rule
$\tau(1,\emptyset,i)$ is applicable to the word) then applying the
rule $\tau(1,\emptyset,i)$, we move all $P_j$ one letter to the
right. Similarly, if we want to move $P_j$ one letter to the left,
we need to apply $\tau(1,\emptyset,i)\iv$ for an appropriate $i$.

Notice that, for every $z\in \kkk$, all state letters $z_j$,
$j=1,\dots, N$ behave in the same way when we apply
$\tau(1,\emptyset,i)$. This important property will hold for all
other rules of $\sss$.

Recall also that the two $a_i$ in $a_iP_ja_i\iv$ are not the same:
they are brothers taken from different alphabets (see Section
\ref{notation}).

The set $\sss^+(12)$ consists of one rule $\tau(12,r)$ for each
$r\in\bee\backslash\{\emptyset\}$:

$$[\lel_j\tool \lel_j, R_j\tool R_j, j=1,...,N; \emptyset\too r, 1\too 2]$$
This rule does not insert any tape letters, it simply changes the
$\bee$- and $\Omega$-coordinates of state letters. This
\label{transrule1}{\em transition rule} prepares the machine for
step 2.

The set $\sss^+(2)$ consists of rules $\tau(2,r,i)$ where
$r\in\bee$, $i=1,...,\bm$:

$$\tau(2,r,i)=[L_j\too a_iL_ja_i\iv, R_j\tool R_j, j=1,...,N;
r\too r, 2\too 2].$$

The meaning of these rules is that they allow the state letters
$L_j$ to move freely while $R_j$ stays next to $\rer_j$.

The set $\sss^+(23)$ consists of one rule for each $r\in \bee$:

$$\tau(23,r)=[L_j\tool L_j, R_j\tool R_j; r\too r, 2\too 3].$$

The meaning of this rule is that the machine can start step 3 when
$L_j$ and $P_j$ meet (that is when the inner parts of
$L_jz$-sectors are empty).

The set $\sss^+(3)$ consists of one rule for each $r\in\bee$ and
each $i$ from 1 to $\bm$:

$$\tau(3,r,i)=[L_j\tool L_j, R_j\too a_i\iv R_ja_i,
j=1,...,N; r\too r, 3\too 3].$$

These rules allow the state letter $R_j$ to move freely between
$P_j$ and $\rer_j$.

The set $\sss^+(34)$ consists of one rule for each non-empty
$r\in\bee$:

$$\tau(34,r)=[L_j\tool rL_j, P_j\tool P_j, j=1,...,N; r\too r, 3\too 4].$$

This rule can be applied when the state letters $L_j, P_j, R_j$
meet together; it inserts $r$ to the left of the state letters
$L_j$, and prepares the machine for step 4.

The set $\sss^+(4)$ consists of rules $\tau(4,r,i)$, $r\in\bee,
i=1,...,\bm$:

$$\tau(4,r,i)=[L_j\too a_i L_ja_i\iv, P_j\tool P_j,
j=1,...,N; r\too r, 4\too 4].$$

These rules allow the state letter $L_j$ to move freely between
$\lel_j$ and $P_j$.

The set $\sss^+(45)$ consists of one rule for each $r\in\bee$:

$$\tau(45,r)=[\lel_j\tool \lel_j, P_j\tool P_j, j=1,...,N; r\too
r, 4\too 5].$$ The machine can start step $5$ when $L_j$ and
$\lel_j$ meet.

The set  $\sss^+(5)$ consists of one rule for each $r\in\bee,
i=1,...,\bm$:

$$\tau(5,r,i)=[\lel_j\tool \lel_j, R_j\too a_i\iv R_ja_i, j=1,...,N;
r\too r, 5\too 5].$$ These rules allow $R_j$ move freely between
$P_j$ and $\rer_j$.

Finally the set $\sss^+(51)$ consists of one  rule for each $r\in
\bee$:

$$\tau(51,r)=[\lel_j\tool \lel_j, R_j\tool R_j, j=1,...,N; r\too\emptyset, 5\too 1].$$

The cycle is complete, the machine can start step $1$ again when
$\lel_j$ meets $L_j$, and $R_j$ meets $\rer_j$.

Rules from the set $\sss(12)\cup \sss(23)\cup \sss(34)\cup
\sss(45)\cup\sss(51)$ will be called \label{trrules}{\em
transition rules}.

A graphical description of the work of the $S$-machine $\sss$ is
presented on Figure \theppp\  below.

\bigskip

\begin{center}
\unitlength=0.8mm \linethickness{0.5pt}
\begin{picture}(109.34,157.67)
\put(15.67,6.00){\framebox(4.67,151.67)[cc]{}}
\put(21.00,6.00){\framebox(4.00,26.33)[cc]{}}
\put(47.00,54.33){\line(-1,0){3.67}}
\put(43.00,54.33){\framebox(3.67,25.33)[cc]{}}
\put(43.00,79.67){\line(5,6){4.67}}
\put(47.67,85.33){\line(1,0){4.00}}
\put(51.67,85.33){\line(-5,-6){4.67}}
\put(47.34,85.33){\line(-5,6){26.00}}
\put(51.67,85.33){\line(-5,6){25.67}}
\put(26.00,116.33){\line(-1,0){5.00}}
\put(21.00,116.33){\framebox(5.00,41.33)[cc]{}}
\put(47.34,54.33){\framebox(3.67,25.33)[cc]{}}
\put(47.34,32.33){\framebox(3.67,22.00)[cc]{}}
\put(38.34,6.00){\line(1,3){8.67}}
\put(51.00,32.33){\line(-1,-3){8.67}}
\put(42.34,6.00){\line(-1,0){4.00}}
\put(47.34,79.67){\line(5,6){4.67}}
\put(51.00,79.67){\line(1,1){5.33}}
\put(52.00,85.33){\framebox(4.33,62.33)[cc]{}}
\put(52.00,147.67){\line(-2,5){4.00}}
\put(56.34,147.67){\line(-2,5){4.00}}
\put(52.34,157.33){\line(-1,0){4.33}}
\put(55.00,79.67){\line(1,1){5.67}}
\put(55.34,85.33){\line(1,0){5.33}}
\put(56.67,85.33){\framebox(4.00,31.00)[cc]{}}
\put(82.34,147.67){\framebox(4.00,9.67)[cc]{}}
\put(82.34,6.00){\framebox(4.00,48.33)[cc]{}}
\put(51.00,79.67){\line(5,-4){31.33}}
\put(55.00,79.67){\line(5,-4){31.33}}
\put(87.34,6.00){\framebox(4.67,151.33)[cc]{}}
\put(93.00,6.00){\framebox(4.00,48.33)[cc]{}}
\put(93.34,54.33){\line(3,4){8.33}}
\put(97.00,54.33){\line(4,5){8.67}}
\put(105.67,65.33){\line(-1,0){4.33}}
\put(11.34,6.00){\framebox(3.67,26.33)[cc]{}}
\put(11.34,32.33){\line(-1,2){5.33}}
\put(6.00,43.00){\line(1,0){3.67}}
\put(15.00,116.00){\line(-1,-3){5.33}}
\put(11.34,116.00){\line(-1,-3){5.00}}
\put(6.34,100.67){\line(1,0){3.33}}
\put(93.00,147.67){\framebox(4.00,9.67)[cc]{}}
\put(17.67,2.00){\makebox(0,0)[cc]{$L_{j-1}\iv K_jL_j$}}
\put(40.34,2.00){\makebox(0,0)[cc]{$P_j$}}
\put(92.67,2.00){\makebox(0,0)[cc]{$R_jK_{j+1}\iv R_{j+1}\iv$}}
\put(5.00,29.00){\line(1,0){104.33}}
\put(109.34,32.33){\line(-1,0){104.33}}
\put(5.00,54.33){\line(1,0){104.33}}
\put(5.34,50.67){\line(1,0){104.00}}
\put(43.00,50.67){\framebox(3.67,3.67)[cc]{}}
\put(21.00,32.33){\line(6,5){22.00}}
\put(25.34,32.33){\line(6,5){21.33}}
\put(5.00,79.67){\line(1,0){104.33}}
\put(5.34,85.33){\line(1,0){103.67}}
\put(5.00,116.33){\line(1,0){104.33}}
\put(5.34,120.67){\line(1,0){103.67}}
\put(56.34,116.33){\framebox(4.33,4.33)[cc]{}}
\put(56.67,116.33){\line(0,1){4.33}}
\put(5.00,147.67){\line(1,0){104.33}}
\put(5.34,143.33){\line(1,0){103.67}}
\put(93.00,143.33){\framebox(4.00,4.33)[cc]{}}
\put(93.00,143.33){\line(5,-6){8.00}}
\put(97.00,143.33){\line(5,-6){7.67}}
\put(104.67,134.00){\line(-1,0){4.00}}
\put(82.34,143.33){\framebox(4.00,4.33)[cc]{}}
\put(11.34,116.33){\framebox(3.67,41.33)[cc]{}}
\put(15.00,32.33){\line(-1,2){5.33}}
\put(82.34,143.33){\line(-1,-1){12.67}}
\put(69.67,130.67){\line(-4,-3){13.00}}
\put(86.34,143.33){\line(-1,-1){14.33}}
\put(72.00,129.00){\line(-4,-3){11.33}}
\put(-1.33,16.33){\makebox(0,0)[cc]{$\sss(1)$}}
\put(-1.33,30.00){\makebox(0,0)[cc]{$\sss(12)$}}
\put(-1.33,42.00){\makebox(0,0)[cc]{$\sss(2)$}}
\put(-1.33,52.33){\makebox(0,0)[cc]{$\sss(23)$}}
\put(-1.33,66.67){\makebox(0,0)[cc]{$\sss(3)$}}
\put(-1.33,82.33){\makebox(0,0)[cc]{$\sss(34)$}}
\put(-1.33,97.67){\makebox(0,0)[cc]{$\sss(4)$}}
\put(-1.33,118.00){\makebox(0,0)[cc]{$\sss(45)$}}
\put(-1.33,132.00){\makebox(0,0)[cc]{$\sss(5)$}}
\put(-1.33,145.33){\makebox(0,0)[cc]{$\sss(51)$}}
\put(-1.33,154.33){\makebox(0,0)[cc]{$\sss(1)$}}
\end{picture}
\end{center}
\begin{center}
\nopagebreak[4] Figure \theppp.
\end{center}
\addtocounter{ppp}{1}

Some of the main properties of $\sss$ are the following.

\begin{rk} \label{remark} Notice that if $W=y_1u_1y_2...y_{t+1}$ is an admissible
word for $\sss$, $y_1$ (resp. $y_{t+1}$) is $K_j(r,i)^{\pm 1}$ for
some $j,r,i$ then one does not need trimming the prefix (resp.
suffix) in order to produce $W\circ \tau$ because in this case
$v((y_1)_-)$ (resp. $u(y_{t+1})$) is empty for all rules $\tau$.
\end{rk}

For every word $w$ in $\{a_1,...,a_\bm\}$ let us denote the word
\begin{equation}\label{Sigma(w)}\begin{array}{l}
K_1(\emptyset,1)L_1(\emptyset,1)P_1(\emptyset,1)w(P_1)R_1(\emptyset,1)K_2(\emptyset,1)\iv...\\
K_N(\emptyset,1)\iv R_N(\emptyset,1)\iv w(P_N)\iv
P_N(\emptyset,1)\iv L_N(\emptyset,1)\iv\end{array}\end{equation}
by \label{sigmaw}$\Sigma(w)$. Notice that
$\Sigma(\emptyset)=\Sigma$. Clearly, $\Sigma(w)K_1(\emptyset,1)$
is admissible for $\sss$ for every positive $w$.

\begin{lemma} \label{noPivP} An admissible word for $\sss$ cannot have  $L_j\iv L_j$-, $P_j\iv
P_j$-, or $R_jR_j\iv$-sectors.
\end{lemma}

\proof Indeed if $W$ is an admissible word, then $W\iv$ is also an
admissible word. So the inner part of a  $L_j\iv L_j$-, $P_j\iv
P_j$-, or $R_jR_j\iv$-sector of $W$ must be both positive and
negative. Hence the inner part must be empty, and $W$ is not
reduced, a contradiction.
\endproof

\begin{lemma} \label{mainsss1}

Suppose that a positive word $w'$ in $\{a_1,...,a_\bm\}$ can be
obtained from a positive word $w$ by insertion or deletion of a
word $r\in\bee$. Then there exists a sequence $h$ of rules of
$\sss$ such that
\begin{equation}\label{1}\begin{array}{l}\Sigma(w) K_1(\emptyset,1)\circ
h=\Sigma(w')K_1(\emptyset,1).\end{array}\end{equation} The length
of $h$ is bounded by a linear function in $|w|+|w'|$.

Thus, in particular, for every positive word $w$ which is equal
to 1 in $\bgr$ there exists a word $h$ in $\sss$ such that
$$\Sigma K_1(\emptyset,1)\circ h= \Sigma(w) K_1(\emptyset,1).$$
\end{lemma}
\proof Without loss of generality we can assume that $w'$ is
obtained from $w$ by inserting $r$. Then we can first apply to
$\Sigma(w)K_1(\emptyset,1)$ rules from $\sss(1)$ to move the
letters $P_j^{\pm 1}$ to the places in copies of $w^{\pm 1}$ where
the insertion will occur ($j=1,...,N$). Then a rule from
$\sss(12)$ (which is applicable) will change the coordinates from
$(\emptyset, 1)$ to $(r,2)$. After that we can apply rules from
$\sss(2)$ to move $L_j^{\pm 1}$ next to $P_j^{\pm 1}$ and
($j=1,...,N$). Then a rule from $\sss(23)$ will change the
coordinates to $(r,3)$. After that we can apply rules from
$\sss(3)$ to move $R_j^{\pm 1}$ next to $P_j^{\pm 1}$. Then a rule
from $\sss(34)$ will insert $r$ (resp. $r\iv$) to the left of
$L_j$ (resp. to the right of $L_j\iv$). Then rules from $\sss(4)$
will move $L_j$ back to $K_j$, for odd $j$ and $L_j\iv$  to
$K_{j+1}$ for even $j$. Then a rule from $\sss(45)$ will change
the coordinates to $(r,5)$, and rules from $\sss(5)$ will move the
letters $R_j$ back to $K_{j+1}\iv$ ($j$ odd), and $R_j\iv$ to
$K_j\iv$ ($j$ even). Finally we can apply a rule from $\sss(51)$
and obtain $\Sigma(w') K_1(\emptyset,1)$ as required.

It is clear that the number of rules applied in the described
process is bounded by a linear function in $|w|+|w'|$.

The ``in particular" statement immediately follows from the first
statement and Lemma \ref{Positive}.
\endproof

\subsection{The $S$-machines $\bsss$ and $\csss$}

The $S$-machine $\bsss$ is similar to $\sss$ so we only present
the differences between $\sss$ and $\bsss$.

We introduce disjoint copies of sets $\kkk$, $\tkk$, $\aaa$,
denote them by \label{bkta}$\bar\kkk$, $\bar\tkk$, $\bar\aaa$,
respectively. In order to make $\sss$ and $\bsss$ ``communicate",
we identify $z_j(\emptyset,1)$ with $\bar z_j(\emptyset,1)$, $z\in
\{K, L, P, R\}$, $j=1,...,N$, and $a_i(P_j)$ with $\bar a_i(\bar
P_j)$ $i=1,...,m$,$j=1,...,N$. Notice that we do not identify
$a_{m+1},...,a_\bm$ with their ``bar"-brothers $\bar
a_{m+1},...,\bar a_\bm$, and we do not identify $a_i(z)$ with
$\bar a_i(z)$ for $z\ne P_j$. Notice also that this identification
of letters makes the word $\bar\Sigma$ coincide with the hub
$\Sigma$.

\label{admbar}{\em Admissible words} of $\bsss$ have the same form
as admissible words of $\sss$ except for the following
differences:

\begin{itemize}
\item All letters are replaced by their ``bar"-brothers. \item We drop
the restriction that certain parts of admissible words are
positive (negative). \item The $\bar K_1z$-, $\bar L_1z$-, $\bar
P_1z$- and $\bar R_1z$-sectors must be empty ($z\in\tkk$). Thus in
every admissible word of $\bsss$ the part between $\bar K_1$ and
$\bar K_2\iv$ contains no $\bar\aaa$-letters.
\end{itemize}

The rules of $\bsss$ are obtained from the corresponding rules of
$\sss$ by replacing every letter by its ``bar"-brother and by
removing letters from $\bar\aaa(\bar z_1)$, $z\in \{K,L,P,R\}$.
Thus rules of $\bsss$ do not insert any $\bar\aaa$-letters in the
interval between $\bar K_1$ and $\bar K_2\iv$, and work in other
parts of the admissible words just as $\sss$.

For every admissible word $W$ for $\sss$, we define an admissible
word $\bar W$ of $\bsss$ by adding {\em bar} to all letters, and
removing non-$\kkk$-letters from all $zz'$-sectors, such that $z$
or $(z')^{-1}$ belongs to $ \{K_1, L_1, P_1, R_1, K_2\}$.

\begin{lemma} \label{mainsss3}

Suppose that $w,u$ are reduced words in $\{a_1,...,a_\bm\}$, and
$w'$ is the freely reduced form of $wuru^{-1}$, $r\in\bee$. Then
there exists a sequence $h$ of rules of $\bsss$ such that
\begin{equation}\label{3}\begin{array}{l}\bar\Sigma(w) K_1(\emptyset,1)\circ
h=\bar\Sigma(w')K_1(\emptyset,1).\end{array}\end{equation}

Therefore, for every word $w$ which is equal to 1 in $\bgr$ there
exists a word $h$ in $\sss$ such that $$\Sigma
K_1(\emptyset,1)\circ h= \bar\Sigma(w) K_1(\emptyset,1).$$

If all letters in $w$ are with indices $\le m$, then one can
remove $\bar{ }$ from the right-hand side of the previous formula
(because we have identified these letters with their
``bar"-brothers).
\end{lemma}

\proof The proof is similar to the proof of Lemma \ref{mainsss1}.
The only difference is that at the beginning we transform every
subword of the form  $P_j(\emptyset, 1)w(P_j)R_j(\emptyset,1)$
into $w(L_j)q(L_j)P_jq(L_j)\iv R_j(\emptyset,1)$ using rules from
$\bar\sss(1)$.\endproof

We shall also need the combined $S$-machine $\csss$. An admissible
word for $\csss$ is either an admissible word for $\bsss$ or a
word satisfying all conditions from the definition of an
admissible word for $\sss$ except for the positivity condition.
The set of rules of $\csss$ is the union of the set of rules of
$\sss$ and $\bsss$.

For the sake of simplicity of notation, for every word $W$ we
denote the projection $W_{\aaa\cup\bar\aaa}$ of $W$ onto the
alphabet $\aaa\cup\bar\aaa$ by $W_a$.

We shall call a rule $\tau\in\sss(g)\cup\bsss(g)$
\label{activerule}{\em active with respect to $zz'$-sectors} if
$u_\tau(z)\ne \emptyset$ or $v_\tau(z'_-)\ne\emptyset$.

We define a map \label{beta}$\beta$ of the subset of generators
$\aaa\cup\bar\aaa$ of the group $\hhh$ onto $\{a_1,\dots,a_{\bar
m}\}\cup \{1\}$ by the rule $\beta(a_i(z))=\beta(\bar a_i(z))=a_i$
for all the $\aaa$- and $\bar \aaa$-letters. Hence the
$\beta$-image of arbitrary word in $\aaa\cup\bar\aaa$ is defined
as an element of the free group with basis $\{a_1,\dots,a_{\bar
m}\}$.

The following two properties of $\csss$ immediately follows from
its definition.

\begin{lemma}\label{csss1} If a rule $\sigma\in\sss(\omega)$ (resp.
$\sigma\in\bar\sss(\omega)$) is active with respect to
$zz'$-sectors then every rule in $\sss(\omega)$ (resp.,
$\bar\sss(\omega)$) is active with respect to $zz'$-sectors.
\end{lemma}

\begin{lemma}\label{csss}
Let $h$ be a sequence of rules from $\csss$ applicable to an
admissible (for $\csss$) word $W$. Then for every $zz'$-sector
$z(r,i)W'z'(r,i)$ of $W$ there exist two words $u=u(h)$ and
$v=v(h)$ such that $z(r,i)W'z'(r,i)\circ
h=z(r_1,i_1)uW'vz'(r_1,i_1)$ in the free group for some
coordinates $r_1, i_1$. The words $u=u(z,h)$ and $v=v(z',h)$
depend on $z,z'$ and $h$ only (i.e. they do not depend on the
inner part of $W'$). Moreover we have the following equalities:

$$u(K_j,h)\equiv v(K_j,h)\equiv\emptyset,
\beta(u(P_j,h))=\beta(v(P_j,h))^{-1}, \beta(u(R_j,h))\equiv
\beta(v(R_j,h))\iv.$$ We also have that
$\beta(u(L_j,h))=\beta(v(L_j,h))\iv$ provided $h$ contains no
rules from $\sss(34)$.
\end{lemma}

\begin{lemma}\label{bsss} Let $U$ be any word in $\tkk$,
$g\in\{1,2,3,4,5\}$. Then there exist two sequences of numbers
$\{\varepsilon_z, \delta_z\in \{-1,0,1\}\ |\ z\in \tkk\}$, such
that for every admissible word $W=W(r,g)$ with base $U$, and any
group word $w$ in the alphabet $\{a_1,...,a_{\bar m}\}$ there
exists a word $h$ in $\bsss(g)$ satisfying the following
properties:

\begin{enumerate}
\item $h$ is a copy\footnote{Here and below, we say that a word
$v$ is a {\em copy} of a word $u$ if $v$ is obtained from $u$ be
replacing equal letters by equal letters and different letters by
different letters.} of $w$ in the alphabet $\bsss(g)$; \item $W$
is in the domain of $h$; \item $W\circ h$ is obtained from $W$ by
all replacements of the form $$\bar z(r,g)\to
w(z_-)^{\varepsilon_z} \bar z(r,g)w(z)^{\delta_z},$$ where $z$
runs over all letters in the base of $W$. Here $w(z)$ and $w(z_-)$
are copies of $w$ in the alphabets $\bar\aaa(z)$ and
$\bar\aaa(z_-)$.
\end{enumerate}
\end{lemma}

Now we need to explain how to simulate the $S$-machines $\sss$ and
$\bsss$ in a group $\hhh$.

\subsection{Converting $\sss$ into a set of relations}
\label{convert1}

In order to convert the machine $\sss$ into a list of defining
relations, we need to introduce one more set of letters
\label{xxx}$\xxx=\{x(a,\tau)\ |\ a\in \aaa\backslash
\cup_j\aaa(P_j), \tau\in\sss^+\}$. We will subdivide this set into
subsets \label{xxxztau}$\xxx(z,\tau)=\{x(a,\tau) \ |\ a\in
\aaa(z)\}$, $z\in\tkk, z\ne P_j (j=1,...,N), \tau\in\sss$ and into
subsets \label{xxxz}$\xxx(z)=\cup_{\tau\in\sss}\xxx(z,\tau)$,
$z\in \tkk$, and
\label{xxxtau}$\xxx(\tau)=\cup_{z\in\tkk}\xxx(z,\tau)$,
$\tau\in\sss^+$. We also set $\xxx(z\iv)=\xxx(z_-)$ for every
$z\in\tkk$.

For every $\tau\in\sss$, we also consider map
\label{alphatau}$\alpha_\tau$ from $\aaa$ to the free group
generated by $\aaa$ and $\xxx\cup\xxx\iv$. Let $\tau\in\sss^+$,
$a\in \aaa(z), z\in \tkk$. Then:

$$\alpha_\tau(a)=\left\{ \begin{array}{ll} x(a,\tau) a & \hbox{ if
} z=K_j, j=1,...,N,\\ x(a,\tau) a & \hbox{ if } z=L_j,
j=1,...,N,\\ a & \hbox{ if } z=P_j, j=1,...,N,\\ ax(a,\tau)  &
\hbox{ if } z=R_j, j=1,...,N.
\end{array}\right.
$$ The definition of $\alpha_\tau^{-1}$, $\tau\in \sss^+$, is
obtained from the definition of $\alpha_\tau$ by replacing all
$\xxx$-letters by their inverses.

We expand these maps to all words in the alphabet
$\aaa\cup\aaa\iv\cup\kkk\cup\kkk\iv$ in the natural way (mapping
letters from $\kkk\cup\kkk\iv$ to themselves). Let $U$ be any word
in $\aaa$. Then $U$ is obviously equal to the {\em projection} of
$\alpha_\tau(U)$ onto the alphabet $\aaa$, i.e. the word obtained
from $\alpha_\tau(U)$ by removing all non-$\aaa$-letters. In
general, the projection of a word $U$ onto an alphabet $Y$ will be
denoted by $U_Y$.

Also with every $\tau\in\sss^+$, we associate a set of letters
\label{thetatau}
$$\Theta(\tau)=\{\theta(\tau,z)\ |\ z\in \{K_j, L_j, P_j, R_j\},
j=1,...,N\}.$$ Let \label{thetaz}$\Theta(z)=\{\theta(\tau,z)\ |\
\tau\in\sss^+\}$, $z\in \{K_j, L_j, P_j, R_j\}$. Finally the union
of all $\Theta(\tau)$ is denoted by \label{theta}$\Theta$.

Let $\tau=[U_1\too V_1,..., U_k\too V_k; r\too r', l\too l']$ be
one of the rules of $\sss^+$. Let $u(z), v(z)$, $z\in \tkk$, be
the words associated with this rule as in Section \ref{notation}.

Then we associate with $\tau$ the following list of
\label{mainrelations}{\em main relations} (or just
\label{thetak}$(\Theta,k)$-{\em relations}):

\begin{equation}
\theta(\tau,z_-)\iv z(r,l)\theta(\tau,z) =
\alpha_{\tau^{-1}}(v_i)z(r',l')\alpha_{\tau^{-1}}(u_i),
i=1,...,k,\label{mainrel}
\end{equation}
(Thus, the replacements associated with the rule $\tau$ are
simulated by ``conjugation". Figure \theppp\ shows a relation
involving a $P_j$-letter, corresponding to a rule
$\tau=\tau(1,\emptyset,i)$; there $P=P_j(1,0)$, $a=a_i(P_j)$,
$a'=a_i(L_j)$, $x'=x(a_i(L_j),\tau)$, $\theta=\theta(\tau,P_j)$,
$\theta'=\theta(\tau,L_j)$.)

\unitlength=1.00mm \special{em:linewidth 0.5pt}
\linethickness{0.5pt}
\begin{picture}(94.00,22.00)
\put(58.33,15.00){\line(2,-3){7.67}}
\put(66.00,3.67){\line(1,0){21.33}}
\put(87.33,3.67){\line(3,5){6.67}}
\put(94.00,14.67){\line(-1,0){35.67}}
\put(68.33,14.67){\vector(-1,0){9.67}}
\put(68.00,14.67){\vector(1,0){12.33}}
\put(69.33,14.67){\circle*{0.57}}
\put(66.00,3.67){\vector(-2,3){6.33}}
\put(73.00,3.67){\vector(1,0){8.33}}
\put(90.81,9.20){\vector(1,2){1.90}}
\put(60.81,7.92){\makebox(0,0)[cc]{$\theta'$}}
\put(61.23,17.00){\makebox(0,0)[cc]{$x'$}}
\put(75.00,17.00){\makebox(0,0)[cc]{$P$}}
\put(93.08,14.72){\vector(-1,0){8.49}}
\put(82.60,14.72){\circle*{0.57}}
\put(89.11,17.00){\makebox(0,0)[cc]{$a$}}
\put(92.09,7.92){\makebox(0,0)[cc]{$\theta$}}
\put(76.94,0.12){\makebox(0,0)[cc]{$P$}}
\put(64.07,14.72){\circle*{0.57}}
\put(64.92,14.72){\vector(1,0){3.68}}
\put(66.61,17.12){\makebox(0,0)[cc]{$a'$}}
\end{picture}

\begin{center}
\nopagebreak[4] Figure \theppp.
\end{center}
\addtocounter{ppp}{1}

 We also need the following \label{auxiliary}{\em auxiliary relations}. Let
$z\in \{K_j, L_j, P_j, R_j\}$, $j=1,...,N$, $a, b\in \aaa(z)$.
Then we add {\em auxiliary $\Theta$-relations} (or
\label{thetaa}$(\theta,a)$-{\em relations}):

\begin{equation}
\theta(\tau,z)\iv \alpha^+_\tau(a) \theta(\tau,z)=\alpha^-_\tau(a)
\hbox{ if } \tau \hbox{ does not lock } zz_+\hbox{-sector}
\label{auxtheta}
\end{equation}
(On Figure \theppp, $a=a_i(L_j)$, $x=x(a,\tau)$,
$\theta=\theta(\tau,L_j)$ where $\tau=\tau(2,r,i')$ for some
$i,i',r$.)

\bigskip

 \unitlength=1.00mm \special{em:linewidth 0.5pt}
\linethickness{0.5pt}
\begin{picture}(93.92,17.72)
\put(71.14,5.09){\framebox(21.37,7.78)[cc]{}}
\put(81.61,12.88){\circle*{0.57}}
\put(79.21,12.88){\vector(-1,0){5.52}}
\put(83.88,12.88){\vector(1,0){6.65}}
\put(81.61,5.09){\circle*{0.57}}
\put(73.83,5.09){\line(1,0){5.09}}
\put(73.97,5.09){\vector(1,0){4.81}}
\put(84.87,5.09){\vector(1,0){5.80}}
\put(71.14,7.78){\vector(0,1){3.25}}
\put(71.14,11.04){\vector(0,0){0.00}}
\put(92.51,7.50){\vector(0,1){3.54}}
\put(69.16,8.77){\makebox(0,0)[cc]{$\theta$}}
\put(93.92,8.77){\makebox(0,0)[cc]{$\theta$}}
\put(76.38,14.72){\makebox(0,0)[cc]{$x$}}
\put(86.28,14.72){\makebox(0,0)[cc]{$a$}}
\put(76.09,2.69){\makebox(0,0)[cc]{$x$}}
\put(86.28,2.69){\makebox(0,0)[cc]{$a$}}
\end{picture}
\begin{center}
\nopagebreak[4] Figure \theppp.
\end{center}
\addtocounter{ppp}{1} and {\em auxiliary}
\label{axrel}$(a,x)$-{\em relations}:

\begin{equation}
\begin{array}{ll}a x(b,\tau)a\iv = x(b,\tau)^4 & \hbox{ if } z=K_j \hbox{ or } z=L_j,\\
a\iv x(b,\tau)a=x(b,\tau)^4 & \hbox{ if } z=R_j.
\end{array}
\label{auxax}
\end{equation}

(On Figure \theppp, $a=a_i(L_j), x=x(a_{i'}(L_j),\tau)$ where
$\tau=\tau(2,r,i')$ for some $i, j, i', r$.)

\unitlength=1.00mm \special{em:linewidth 0.5pt}
\linethickness{0.5pt}
\begin{picture}(85.72,27.05)
\put(71.99,2.97){\line(4,3){12.03}}
\put(84.02,12.03){\line(0,1){3.96}}
\put(84.02,15.99){\line(-4,3){12.03}}
\put(71.99,25.05){\line(0,-1){22.08}}
\put(71.99,8.21){\circle*{0.57}} \put(71.99,13.58){\circle*{0.57}}
\put(71.99,19.10){\circle*{0.57}}
\put(71.99,20.66){\vector(0,1){2.97}}
\put(71.99,14.29){\vector(0,1){3.82}}
\put(71.99,9.06){\vector(0,1){3.11}}
\put(71.99,4.39){\vector(0,1){2.69}}
\put(84.02,12.17){\vector(0,1){3.11}}
\put(74.82,5.09){\vector(4,3){5.09}}
\put(74.96,22.92){\vector(4,-3){5.38}}
\put(85.72,13.73){\makebox(0,0)[cc]{$x$}}
\put(80.20,21.51){\makebox(0,0)[cc]{$a$}}
\put(80.20,5.94){\makebox(0,0)[cc]{$a$}}
\put(69.58,5.52){\makebox(0,0)[cc]{$x$}}
\put(69.58,11.04){\makebox(0,0)[cc]{$x$}}
\put(69.58,16.27){\makebox(0,0)[cc]{$x$}}
\put(69.58,21.65){\makebox(0,0)[cc]{$x$}}
\end{picture}
\begin{center}
\nopagebreak[4] Figure \theppp.
\end{center}
\addtocounter{ppp}{1}

Finally let $b'$ be the ``brother" of $b$ in the alphabet
$\aaa(z_-)$ if $z=K_j$ or $z=L_j$. Then we need auxiliary
\label{kxrel}$(k,x)$-{\em relations}:

\begin{equation}\label{auxkx}
\begin{array}{ll} z(r,i) x(b,\tau)=x(b',\tau)z(r,i) & \hbox{ if }
z=K_j,\\
z(r,i) x(b,\tau)=x(b',\tau)^4z(r,i) & \hbox{ if }
z=L_j.\end{array}
\end{equation}
(On Figure \theppp, $L=L_j(0,1)$, $x=x(a_i(L_j),\tau)$,
$x'=x(a_i(\lel_j),\tau)$, where $\tau=\tau(2,r,i')$ for some
$i,j,r,i'$.)

\unitlength=1.00mm \special{em:linewidth 0.5pt}
\linethickness{0.5pt}
\begin{picture}(85.72,27.05)
\put(71.99,2.97){\line(4,3){12.03}}
\put(84.02,12.03){\line(0,1){3.96}}
\put(84.02,15.99){\line(-4,3){12.03}}
\put(71.99,25.05){\line(0,-1){22.08}}
\put(71.99,8.21){\circle*{0.57}} \put(71.99,13.58){\circle*{0.57}}
\put(71.99,19.10){\circle*{0.57}}
\put(71.99,20.66){\vector(0,1){2.97}}
\put(71.99,14.29){\vector(0,1){3.82}}
\put(71.99,9.06){\vector(0,1){3.11}}
\put(71.99,4.39){\vector(0,1){2.69}}
\put(84.02,12.17){\vector(0,1){3.11}}
\put(74.82,5.09){\vector(4,3){5.09}}
\put(74.96,22.92){\vector(4,-3){5.38}}
\put(85.72,13.73){\makebox(0,0)[cc]{$x$}}
\put(80.20,21.99){\makebox(0,0)[cc]{$L$}}
\put(80.20,5.94){\makebox(0,0)[cc]{$L$}}
\put(69.58,5.52){\makebox(0,0)[cc]{$x'$}}
\put(69.58,11.04){\makebox(0,0)[cc]{$x'$}}
\put(69.58,16.27){\makebox(0,0)[cc]{$x'$}}
\put(69.58,21.65){\makebox(0,0)[cc]{$x'$}}
\end{picture}
\begin{center}
\nopagebreak[4] Figure \theppp.
\end{center}
\addtocounter{ppp}{1}

Let $\rr(\sss)$ be the set of all relations (\ref{mainrel}) -
(\ref{auxkx}) corresponding to the set of $S$-rules $\sss^+$. It
is easy to see that the relations corresponding to the negative
rules from $\sss$ follow from $\rr(\sss)$.

\subsection{Converting $\bsss$ into a set of relations}
\label{convert2}

To convert $\bsss$ into a set of relations we do not need
$\xxx$-letters, but we need ``bar"-brothers of the letters from
$\Theta$. Let this set be denoted by \label{btheta}$\bar\Theta$.
The easiest way to explain how to convert an $S$-rule of $\bsss$
to a set of relations is the following: convert it as in Section
\ref{convert1}, then add bars $\bar{ }$ to all letters, replace
all letters from $\xxx$ and from $\cup\bar A(z_1), z\in
\{K,L,P,R\})$ by 1. Let $\rr(\bsss)$ be the set of all relations
corresponding to $\bsss^+$.

\subsection{The presentation of the group $\hhh$}
\label{convert3}

The \label{hhhh}group $\hhh$ which is the main object studied in
this paper is given by the set of generators
$\kkk\cup\bar\kkk\cup\aaa\cup\bar\aaa\cup\Theta\cup\bar\Theta\cup\xxx$
and the set of relations
$\rr=\rr(\sss)\cup\rr(\bsss)\cup\{\Sigma\}$.

\section{The first properties of $\hhh$}

\subsection{Diagrams}
\medskip

Recall the well-known van Kampen-Lyndon topological interpretation
of the consequences of defining relations in groups \cite{LS},
\cite{Ol_book}. We define a \vk \label{vkd}{\it diagram} over some
group presentation $\LL=\langle x_1,\dots,x_k\mid {\cal P}\rangle$
(or briefly, by misuse of language, over the group $\LL$), where
the defining words from ${\cal P}$ are cyclically reduced. It is
an oriented, connected, and simply connected, finite planar
2-complex such that every oriented edge $e$ of it has a {\it
label} $\phi(e)\equiv x_i^{\pm1}$, $i=1,\dots,k,$ where
$\phi(e^{-1})\equiv\phi(e)^{-1}$. In addition, the word read on
the boundary $\partial\Pi$ of an arbitrary $2$-cell $\Pi$ (we call
it a {\it cell} latter on), coincides with one of defining word,
up to a cyclic permutation and an inversion.

According to the van Kampen lemma a word $w=w(x_1,\dots,x_k)$
equals $1$ in a group $\LL$ if and only if there exists a \vk
diagram $\Delta$ over $\LL$ with $w$ read clockwise on the
boundary $\partial(\Delta)$ starting from some vertex. All details
can be found in \cite{LS} or \cite{Ol_book}.

An \label{annular}{\it annular} diagram $\Delta$ is not simply
connected (unlike a \vk diagram), but its complement on the plane
has two components: the inner and the outer. So $\Delta$ has the
{\em inner} and {\em outer} contours. By the Schupp lemma
\cite{LS}, \cite{Ol_book} their labels, read clockwise, are
conjugate in $\LL$. Conversely, for any pair of words $u$, $v$,
which represent non-trivial conjugate elements of $\LL$, there
exists an annular diagram whose contours are labelled by $u$ and
$v$ respectively.

\medskip

\subsection{$\gr$ embeds into $\hhh$}
\label{embeds}

We define a map \label{delta}$\delta$ of the set of generators of
the group $\hhh$ into the group $\bgr$ by setting
$\delta(a_i(z))=\delta( a_i(z))=a_i$ for all the $\aaa$- and
$\bar\aaa$-letters and mapping other letters from the whole set of
generators to 1.

\begin{lemma}\label{alpha} The map $\delta$ extends to the epimorphism
$\delta:\;\hhh\rightarrow \bgr.$
\end{lemma}

\proof The map $\delta$ can be considered as an endomorphism of
the free group over the generators of $\hhh$. It immediately
follows from the list of relations of the group $\hhh$, that the
$\delta$-image of every defining word for $\hhh$ vanishes in
$\bgr$. Indeed, it is obvious for the hub word $\Sigma$ since it
does not contain $\aaa\cup\bar\aaa$-letters. It is also true for
all auxiliary relations (\ref{auxtheta}), (\ref{auxax}),
(\ref{auxkx}) because each of them contains exactly two
occurrences of $\aaa\cup\bar\aaa$-letters whose $\delta$-images
are mutually inverse. Same is true for all relations
(\ref{mainrel}) corresponding to $S$-rules $\tau$ and $z\in\tkk$
unless $\tau=\tau(34,r)$ and $z=L_j$ or $\tau=\bar\tau(34,r)$ and
$z=\bar L_j$ (for some $r$ and $j$). In the latter case the the
$\delta$-image of the relator (\ref{mainrel}) is equal to $r$
which belongs to $\bee$ and is equal to $1$ in $\bgr$.
\endproof

The following lemma is similar to Lemma \ref{alpha} and
immediately follows from the definition of $\csss$.

\begin{lemma}\label{alpha1} For every word $U$ in the domain of
$\tau\in\csss$, $\delta(U)=\delta(U\circ \tau)$ in $\bar\gr$.
\end{lemma}

\begin{lemma} \label{W1}
a) Let $W$ be an admissible word for the $S$-machine $\sss$, $z$
be its first $\kkk$-letter and $z'$ be its last $\kkk$-letter. Let
$\tau\in\sss$, and let $v=v(z_-)$ and $u=u(z')$ be the words
associated with $\tau$, defined in Section \ref{notation}. Suppose
that $\tau\in\sss$ is applicable to $W$. Then
$$\theta(\tau,z_-)\iv\alpha_\tau
(W)\theta(\tau,z')=\alpha_{\tau\iv} (v)\alpha_{\tau\iv}
(W\circ\tau)\alpha_{\tau\iv} (u)$$ in $\hhh$.

b) Let $W$ be an admissible word for the $S$-machine $\bsss$, $z$
is its first $\kkk$-letter and $z'$ is its last $\kkk$-letter. Let
$\tau\in\bsss$, and $v=v(z_-)$ and $u=u(z')$ be the words
associated with $\tau$, defined in Section \ref{notation}. Suppose
that $\tau$ is applicable to $W$. Then
$$\theta(\tau,z_-)\iv W\theta(\tau,z')=v \cdot(W\circ \tau)\cdot u$$
in $\hhh$.

\end{lemma}

\proof The proof follows from the definition of admissible words
for $\sss$ and $\bsss$ (we do not need the positivity of certain
segments of admissible words of $\sss$), relations
(\ref{mainrel}), (\ref{auxtheta}), and similar relations for
$\bsss$.
\endproof

\begin{df}\label{sigmari} {\rm Let $w$ be a word over $\{a_1,...,a_\bm\}$. For every
$z\in\tkk$, let $w(z)$ be the word $w$ rewritten in the alphabet
$\aaa(z)$. Let $w_1, w_2, w_3, w_4$ be words over
$\{a_1,...,a_\bm\}$, let $j=1,...,N$. Denote the word
$$w_1(\lel_j)L_jw_2(L_j)P_jw_3(P_j)R_jw_4(R_j)$$ by
$W_j=W_j(w_1,w_2,w_3,w_4)$. Denote the word $$K_1W_1K_2\iv
W_2\iv... K_N\iv W_N\iv$$ by \label{sigmawwww}
$\Sigma(w_1,w_2,w_3,w_4)$. Let the word
\label{barsigmawwww}$\bar\Sigma(w_1,w_2,w_3,w_4)$ be obtained from
$\Sigma(w_1,w_2,w_3,w_4)$ by replacing $W_1$ with $L_1P_1R_1$, and
adding $\bar{ }$ to all letters. Finally for every
$i\in\{1,2,3,4,5\}$ and $r\in\bee$ let
\label{sigmariwwww}$\Sigma_{r,i}(w_1,w_2,w_3,w_4)$ be the word $
\Sigma(w_1,w_2,w_3,w_4)(r,i)$ and
\label{barsigmariwwww}$\bar\Sigma_{r,i}(w_1,w_2,w_3,w_4)$ be the
word $\bar\Sigma(w_1,w_2,w_3,w_4)(r,i)$.}
\end{df}

Notice that if $w_1,w_2,w_4$ are positive words (we do not care
about $w_3$) then the word $\Sigma_{r,i}(w_1,w_2,w_3,w_4)K_1(r,i)$
(see the definition of $U(r,i)$ in \ref{hardware}) is admissible
for $\sss$ for every $r, i$ . The word
$\bar\Sigma_{r,i}(w_1,w_2,w_3,w_4)\bar K(r,i)$ is admissible for
$\bsss$ for every $w_1, w_2, w_3, w_4, r, i$.

\begin{lemma} \label{W}
Preserving the notation from Definition \ref{sigmari}, suppose
that $w_1,w_2,w_4$ are positive words. Let $j\in\{1,...,N\}$,
$r\in\bee$, $i\in\{1,2,3,4,5\}$, $W=W_j(r,i)$, and $\tau$ is
applicable to $W$. Then there exist words $X_1$, $X_1'$, $X_2$,
and $X_2'$ in the alphabet $\xxx(\tau)$, depending on $w_1, w_2,
w_4$ (but not on $w_3$) such that
\begin{itemize}
\item[a)] $X_1WX_1'=\alpha_\tau(W)$; \item[b)]
$X_2WX_2'=\alpha_{\tau\iv}(W)$; \item[c)] the words $X_1$ and
$X_2$ contain only letters from the $\xxx(\lel_j,\tau)$, the words
$X_1'$, $X_2'$ contain letters only from $\xxx(R_j,\tau)$.
\end{itemize}
The lengths of $X_i, X_i'$ ($i=1,2$) are bounded by a recursive
function in $|w_1|+|w_2|+|w_4|$ (in fact it does not exceed
$4^{|w_1|+|w_2|+|w_4|+1}$).
\end{lemma}

\proof The statement immediately follows from relations
(\ref{auxax}).
\endproof

\begin{lemma} \label{sigma} a) In the above notation, suppose that $w_1,w_2,
w_4$ are positive words, $r\in \bee, i\in \{1,2,3,4,5\}$. And
suppose that $\tau\in\sss$ is applicable to
$$\Sigma_{r,i}(w_1,w_2,w_3,w_4)K_1(r,i),$$ $\tau$ changes the
coordinates $(r,i)$ into $(r',i')$. Then, by (\ref{mainrel}) -
(\ref{auxkx}),
$$X\iv\Sigma_{r,i}(w_1,w_2,w_3,w_4)XK_1(r',i')=(\Sigma_{r,i}(w_1,w_2,w_3,w_4)K_1(r,i))\circ
\tau$$ for some word $X=X_1\theta X_2$ where $X_1$ and $X_2$ are
words in $\xxx(\tau)$, whose lengths are bounded by a recursive
function on $|w_1|+|w_2|+|w_4|$,  $\theta\in \Theta(\tau,
(K_1)_-)$.

b) Suppose that $w_1,w_2,w_3,w_4$ are arbitrary words, and $\tau$
is applicable to $$\bar\Sigma_{r,i}(w_1,w_2,w_3,w_4)\bar
K_1(r,i).$$ Then for some $\theta\in\bar\Theta(\tau, (K_1)_-)$ we
have $$\theta\iv\bar\Sigma_{r,i}(w_1,w_2,w_3,w_4)\theta
K_1(r',i')=(\bar\Sigma_{r,i}(w_1,w_2,w_3,w_4)K_1(r,i))\circ
\tau.$$
\end{lemma}

\proof a) Fix a number $j$ between $1$ and $N$.  Let
$K=\lel_j(r,i)$, $K'=\rer_j(r,i)$. Since $\tau$ is applicable to
$\Sigma_{r,i}(w_1,w_2,w_3,w_4)K_1(r,i)$, it is applicable to
$KW_j(w_1,w_2,w_3,w_4)K'$. Hence by Lemma \ref{W} a), there exist
words $X_j, X_j'$ in the alphabets $\xxx(K_j,\tau)$ and
$\xxx(R_j,\tau)$ respectively such that
\begin{equation}
X_jKW_jK'X_j'=K\alpha_\tau(W_j)K'. \label{9*}
\end{equation}
The words $X_j$ and $X'_j$ depend only on $w_1,w_2,w_4$ and
$\tau$. So all $X_j$ (resp. $X_j'$), $j=1,...,N$, are copies of
each other. The word $X_j$ contains letters from the set
$\xxx(K_-,\tau)$ only, the word $X_j'$ contains letters from
$\xxx(K',\tau)$ only. Hence by relations (\ref{auxkx}) we have

\begin{equation}\label{10*}
KX_j=X_{j-1}K, \quad (X'_j)\iv K' = K'(X_{j+1}')\iv.
\end{equation}

Combining (\ref{9*}) with (\ref{10*}) and using the definition of
the word $\Sigma_{r,i}(w_1,w_2,w_3,w_4)$, we get

\begin{equation}
X_N\iv
\Sigma_{r,i}(w_1,w_2,w_3,w_4)X_NK_1(r,i)=\alpha_\tau(\Sigma_{r,i}(w_1,w_2,w_3,w_4))K_1(r,i).
\label{12*}
\end{equation}

Now using Lemma \ref{W1}, relation (\ref{mainrel}), and the fact
that $u(K_1)=v(K_1)=\emptyset$ (Remark \ref{remark}) we deduce
that

\begin{equation}\label{11*}
\theta\iv\alpha_\tau(\Sigma_{r,i}(w_1,w_2,w_3,w_4))\theta
K_1(r',i')
=\alpha_{\tau\iv}(\Sigma_{r,i}(w_1,w_2,w_3,w_4)K_1(r,i)\circ \tau)
\end{equation}
where $\theta=\theta(\tau,(K_1)_-)$.

Now using Lemma \ref{W} b), c), we deduce (similarly to
(\ref{12*})) that for some word $Y_N$ in $\xxx((K_1)_-,\tau)$, we
have
\begin{equation}
Y_N\iv \alpha_{\tau\iv}(\Sigma_{r',i'}(w_1,w_2,w_3,w_4)\circ
\tau)Y_N K_1(r',i') =(\Sigma_{r,i}(w_1,w_2,w_3,w_4)K_1(r,i))\circ
\tau. \label{13*}
\end{equation}
Finally combining (\ref{12*}), (\ref{11*}) and (\ref{13*}), we
deduce statement a) of the lemma with $X=X_N\theta Y_N$.

\medskip
Statement b) immediately follows from Lemma \ref{W1} b).
\endproof

Let \label{hhh0}$\hhh_0$ be the group given by all defining
relations of $\hhh$ except for the hub $\Sigma$.

\begin{lemma} \label{Special}
Assume a positive word $w'=w'(a_1,\dots,a_{\bar m})$ can be
obtained from a positive word $w$ by an insertion (deletion) of a
word $r\in \bee$. Then the words $\Sigma(w)$ and $\Sigma(w')$ are
conjugate in the group $\hhh_0$. There is an annular diagram for
this conjugacy without $\bar\theta$-cells, in which the contours
$p$ and $p'$ can be connected by a path whose length is bounded by
a recursive function of $|w|+|w'|$.
\end{lemma}

\proof  Indeed, by Lemma \ref{mainsss1} there exists a word $h$ in
the alphabet $\sss$ which is applicable to $\Sigma(w)
K_1(\emptyset,1)$ and such that $$(\Sigma K_1(\emptyset,1))\circ
h=\Sigma(w') K_1(\emptyset,1).$$ The length of $h$ is bounded by a
recursive function of $|w|+|w'|$.

Notice that for every words $w_1,w_2, w_3, w_4$ where $w_1, w_2,
w_4$ are positive and for every rule $\tau$ applicable to
$\Sigma_{r,i}(w_1,w_2,w_3,w_4)K_1(r,i)$, the word
$\Sigma_{r,i}(w_1,w_2,w_3,w_4)K_1(r,i)\circ \tau$ also has the
form $\Sigma_{r',i'}(w'_1,w'_2,w'_3,w'_4)K_1(r',i')$ for some
$w_1', w_2', w_3', w_4', r', i'$. Now the statement of the lemma
follows by $|h|$ applications of Lemma \ref{sigma} a).
\endproof

\begin{lemma} \label{Special(bar)}
Assume a reduced word $w'=w'(a_1,\dots,a_{\bar m})$ is freely
equal to a product $wqrq^{-1}$, $r\in {\cal R}$. Then the words
$\bar\Sigma(w)$ and $\bar\Sigma(w')$ are conjugate in $\hhh_0$.
There is an annular diagram over $\hhh_0$ without $\theta$-cells
for this conjugacy, in which the contours $p$ and $p'$ can be
connected by a path whose length is bounded by a linear function
of $|w|+|w'|+|q|$.
\end{lemma}

\proof The proof is completely analogous to the proof of Lemma
\ref{Special}\endproof

\begin{lemma} \label{Sigma}
a) For every positive word $w$ in the alphabet $\{a_1,...,a_\bm\}$
which is equal to 1 in $\bgr$ the word $\Sigma(w)$ is a conjugate
of the hub $\Sigma$ in $\hhh_0$. The length of the conjugating
word is a bounded by a recursive function in $|w|$. The
corresponding annular diagram does not contain $\bar\theta$-cells.

b) For every word $w$ in the alphabet $\{a_1,...,a_\bm\}$ which is
equal to 1 in $\bgr$, the word $\bar\Sigma(w)$ is a conjugate of
$\Sigma$ in $\hhh_0$. The length of the conjugating word is a
bounded by a recursive function in $|w|$. The corresponding
annular diagram does not contain $\theta$-cells.
\end{lemma}

\proof It immediately follows from Lemmas \ref{Positive},
\ref{Special}, \ref{Special(bar)} and the solvability of the word
problem in $\bgr$.
\endproof

\begin{lemma} \label{Emb}
a) The map $a_1\to a_1(P_1)$, $a_2\to a_2(P_1)$, ..., $a_m\to
a_m(P_1)$ extends to an injective homomorphism from $\gr$ into
$\hhh$.

b) For every word $w$ over $a_1,...,a_m$ which is equal to 1 in
$\gr$, there exists a \vk diagram $\Delta$ over the presentation
$\rr$ of $\hhh$ with boundary label $w$ and the number of cells
bounded by a recursive function in $|w|$.
\end{lemma}

\proof a) Let $w\in \eee$. Then $w=\bar w$ (because we agreed to
identify $a_i$ with $\bar a_i$, $i=1,...,m$). By Lemma
\ref{Sigma}, both words $\Sigma(w)$ and $\bar\Sigma(w)$ are
conjugated to $\Sigma=1$ in $\hhh$. These words differ only by the
subword $w(P_1)$ (i.e. $\Sigma(w)$ is obtained from
$\bar\Sigma(w)$ by inserting $w(P_1)$ ). Hence $w(P_1)=1$ in
$\hhh$. Thus the map $a_i\to a_i(P_1)$ extends to a homomorphism
from $\gr$ to $\hhh$.

Since $\gr$ is a subgroup of $\bgr$, this homomorphism is
injective by Lemma \ref{alpha}.

b) This statement immediately follows from Lemmas \ref{Positive}
and \ref{Sigma}.
\endproof

\subsection{Bands and annuli} \setcounter{equation}{0}
\medskip

Let $s^{\pm 1}$ be a letter,  and $S$ be a set of letters. We call
a cell (edge) $s$-cell ($s$-edge) or $S$-cell ($S$-edge) provided
its label contains a letter $s^{\pm 1}$ or a letter from $S\cup
S\iv$, respectively. If a cell is both $S_1$-cell and $S_2$-cell,
we call it a $(S_1,S_2)$-{\em cell}.

For example any cell corresponding to one of the relations
(\ref{mainrel}), (\ref{auxtheta}) is a $\Theta(\tau)$-cell.

For $z\in \{K,L,P,R\}$, we shall call a cell a $z$-{\em cell} or
$z_j$-{\em cell} if it is a $z_j(r,i)$-cell for some $j,r,i$.
Similarly we define $\bar z$-cells.

For simplicity we shall call an $\aaa\cup\bar\aaa$-edge, (letter,
cell) as an \label{aedge}$a$-edge (letter, cell), a
$\Theta\cup\bar\Theta$-edge (letter, cell) as
$\theta$-edge(letter, cell), an $\xxx$-edge (letter, cell) as an
\label{xedge}$x$-edge (letter, cell), and $\kkk\cup\bar\kkk$-edge
(letter, cell) \label{kedge}$k$-edge (letter, cell). Thus we can
use such notions as \label{axcell}$(a,x)$-{\em cell},
\label{kxcell}$(k,x)$-{\em cell},
\label{thetaacell}$(\theta,a)$-{\em cell} and
\label{thetakcell}$(\theta,k)$-{\em cell}.

As in \cite{talk}, we enlarge the set of defining relations of the
group $\hhh$ by adding some consequences of defining relations.
Namely, by Lemma \ref{Emb}, we can include all cyclically reduced
relations of the copy of the group $\gr$ generated by the set
$\{a_1(P_1),...,a_m(P_1)\}$, i.e. all cyclically reduced words in
$\{a_1(P_1),...,a_m(P_1)\}$ which are equal to 1 in the copy of
$\gr$. These relations will be called \label{grrel}$\gr$-{\em
relations}. The cells corresponding to $\gr$-relations will be
called \label{grcell}$\gr$-{\em cells}.

We denote by \label{hhh1}$\hhh_1$ the group given by all
generators of $\hhh$, by all $\gr$-relations, and by all defining
relations of $\hhh$, except for the hub $\Sigma$.



The new, expanded, set of relations of $\hhh$ will be denoted by
$\rrr'$.

It is convenient to turn $\rrr'$ into a \label{gradedpres}{\em
graded presentation}.

A hub is ``higher" than $(\Theta, k)$-cells corresponding to
(\ref{mainrel}). (We also say that the \label{rankcell}{\em rank}
of a hub is greater than the rank of $(\Theta,k)$-cell.)  In turn,
$(\Theta,k)$-cells are ``higher" than $(\bar\Theta,\bar k)$-cells,
which are ``higher" than $\gr$-cells which are ``higher" than
$(\Theta,a)$-cells and $(\bar\Theta,\bar a)$-cells which are
``higher" than $x$-{\em cells} corresponding to relations
(\ref{auxkx}), (\ref{auxax}). The $(\Theta, k)$ and
$(\bar\Theta,\bar k)$-relations are also stratified: we set that
relations (\ref{mainrel}) corresponding to rules $\tau$ from
connecting steps (12), (34), (51) are ``higher" than others.

If $\Delta$ and $\Delta'$ are diagrams over $\rrr'$ then we say
that $\Delta$ has a higher \label{typed}{\em type} than $\Delta'$,
if $\Delta$ has more hubs, or the numbers of hubs are the same,
but $\Delta$ has more cells which are next in the hierarchy, etc.

It is obvious, that in this way we define a partial quasi-order (a
transitive reflective relation) on the set of diagrams over
$\rrr'$, and that this partial quasi-order satisfies the
descending chain condition, so we can consider type as an
inductive parameter.

A diagram $\Delta$ that has the smallest type among all diagrams
with the given boundary label (boundary labels, in the annular
case), then $\Delta$ is called a \label{minimald}{\it  minimal}
diagram. The following lemma is obvious (see Figure \theppp).

\begin{lemma} \label{Minimal}
Assume a diagram $\Delta$ over $\hhh$ has two cells $\Pi$ and
$\Pi'$ of the same rank $\rho$ such that words $V$ and $V'$ can be
read on the boundary paths $\partial\Pi$ and $\partial\Pi'$ when
going around the cells in opposite directions from vertices $o$
and $o'$. Assume further that there is a path $p=o-o'$ in $\Delta$
with no self-intersections, $\Lab(p)=U$, and there exists a
diagram with boundary label $UV'U\iv V\iv$ containing at most one
cell of rank $\rho$ and some cells of smaller ranks. Then the
number of cells of rank $\rho$ in $\Delta$ can be decreased, while
the boundary label (boundary labels, in the annular case) and
numbers of cells of higher ranks are preserved. In particular
$\Delta$ is not a minimal diagram.
\end{lemma}

\begin{center}
\unitlength=1.00mm \special{em:linewidth 0.5pt}
\linethickness{0.5pt}
\begin{picture}(76.83,25.33)
\put(12.17,12.83){\oval(17.33,16.33)[]}
\put(68.17,12.83){\oval(17.33,16.33)[]}
\bezier{132}(20.83,12.67)(30.83,25.33)(42.17,13.33)
\bezier{116}(41.83,13.33)(49.83,2.00)(59.50,13.00)
\put(40.17,19.67){\makebox(0,0)[cc]{$U$}}
\put(11.83,12.67){\makebox(0,0)[cc]{$\Pi$}}
\put(68.17,12.67){\makebox(0,0)[cc]{$\Pi'$}}
\put(19.83,4.00){\makebox(0,0)[cc]{$V$}}
\put(75.83,3.67){\makebox(0,0)[cc]{$V'$}}
\put(65.17,4.67){\vector(1,0){6.00}}
\put(15.17,4.67){\vector(-1,0){5.33}}
\end{picture}
\end{center}

\begin{center}
\nopagebreak[4] Figure \theppp.
\end{center}
\addtocounter{ppp}{1}

We say that an annular diagram $\Gamma$ over $\hhh_1$ is
\label{compressible}{\em compressible} if
\begin{itemize}
\item[(1)] It contains no hubs, and there are no $\theta$-edges on
its boundary, and

\item[(2)] There exists another annular diagram $\Gamma'$ which
\begin{itemize}
\item[(a)] has no hubs and has the same boundary labels as
$\Gamma$,

\item[(b)] the type of $\Gamma'$ is smaller than the type of
$\Gamma$, and these types differ already by the number of
$\theta$-cells, and

\item[(c)] $\Gamma'$ has no $(\bar\Theta,k)$-cells if $\Gamma$ has
no $(\bar\Theta,k)$-cells.
\end{itemize}
\end{itemize}
By definition, a diagram $\Delta$ is \label{reducedd}{\em not
reduced} if either
\begin{itemize}
\item[(i)] $\Delta$ contains a \label{reduciblec}{\it reducible}
pair of cells having a common edge $o-o'$ and boundary labels $w$
and $w\iv$ (read starting at the common vertex $o$) or


\item[(ii)] $\Delta$ contains
a pair of $(k,\Theta)$-cells like in
Lemma \ref{Minimal} where $V'\equiv V^{-1}$ and $U$ is a word in
$\aaa(P_1)$ which is equal to 1 modulo $\gr$-relations, or

\item[(iii)] $\Delta$ contains a pair of $\gr$-cells, and a
connecting path $p$ like in Lemma \ref{Minimal}, where the label
$U$ of $p$ is a word in $\theta(\tau,P_1)$-letters commuting with
letters from $\aaa(P_1)$ by relations (\ref{auxtheta}) or

\item[(iv)] $\Delta$ is an annular diagram and contains a
compressible annular subdiagram surrounding the hole.
\end{itemize}

It is easy to see that every non-reduced diagram is not minimal.
In case (iii), for example, the word $U$ commutes with all
$\aaa(P_1)$-letters by relations (\ref{auxtheta}), hence the
subdiagram consisting of the two $\gr$-cells and the connecting
path can be replaced by a subdiagram having at most one $\gr$-cell
and several commutativity $(\Theta, a)$-cells. In case (iv), the
annular subdiagram can be replaced by an annular diagram of a
smaller type with the same boundary labels.


Of course, a diagram is said to be {\em reduced} if it is not
non-reduced. Thus every minimal diagram is reduced (but the
converse is not necessarily true).

\medskip

As in our previous papers (see \cite{talk}), it is convenient to
study properties of diagrams using bands. Let us fix a point
$o_{\Pi}$ inside every cell $\Pi$ and some point $o_{\Delta}$ on
the plane outside the diagram (or two points $o$ and $o'$ in each
of two components of the complement of $\Delta$ on the plane, in
the annular case). Similarly we fix points $o_e$ inside every edge
$e$ of diagram $\Delta.$ In the case when $e$ is an edge on
$\partial(\Pi)$ (on $\partial(\Delta)$), we fix a simple Jordan
curve $l(\Pi,e)$ (respectively $l(\Delta,e)$) that connects the
vertex $o_{\Pi}$ ($o_{\Delta}$) with $o_e$ and has no other common
points (except for $o_e$) with the edges of the diagram $\Delta$.
We also require that $l(\Pi,e)$ and $l(\Pi,e')$ have the only
common point $o_{\Pi}$ for $e\neq e'$. The similar requirement
concerns the lines $l(\Delta,e)$.

Let $S$ be a set of letters and let $\Delta$ be a van Kampen
diagram.  Fix pairs of $S$-edges on the boundaries of some cells
from $\Delta$ (we assume that each of these cells has exactly two
$S$-edges).

Suppose that $\Delta$ contains a sequence of cells
($\pi_1,\ldots,\pi_n$) such that for each $i=1,\ldots,n-1$ the
cells $\pi_i$ and $\pi_{i+1}$ have a common $S$-edge $e_i$ and
this edge belongs to the pair of $S$-edges fixed in $\pi_i$ and
$\pi_{i+1}$. Let $e$ be the other $S$-edge ($\neq e_1$) of
$\pi_1$, and let $f$ be the other $S$-edge ($\neq e_{n-1}$) of
$\pi_n$. Consider the line which is the union of the lines
$\ell(\pi_i, e_i)$ and $\ell(\pi_{i+1}, e_i)$, $i=1,\ldots,n-1$,
together with $\ell(\pi_{1}, e)$ and $\ell(\pi_{n}, f)$. This
broken line is called the {\em median} of this sequence of cells.
Then our sequence of cells ($\pi_1,\ldots,\pi_n$) with common
$S$-edges $e_1, \ldots, e_{n-1}$, and distinguished edges $e, f$,
 is called an \label{Sband}{\em $S$-band} (with start edge $e$ and end edge $f$)
if the median is an absolutely simple curve or an absolutely
simple closed curve (\setcounter{band}{\value{ppp}}see Fig.
\theband).

\bigskip

\unitlength=0.90mm \special{em:linewidth 0.5pt}
\linethickness{0.5pt}
\begin{picture}(149.67,30.11)
\put(19.33,30.11){\line(1,0){67.00}}
\put(106.33,30.11){\line(1,0){36.00}}
\put(142.33,13.11){\line(-1,0){35.67}}
\put(86.33,13.11){\line(-1,0){67.00}}
\put(33.00,21.11){\line(1,0){50.00}}
\put(110.00,20.78){\line(1,0){19.33}}
\put(30.00,8.78){\vector(1,1){10.33}}
\put(52.33,8.44){\vector(0,1){10.00}}
\put(76.33,8.11){\vector(-1,1){10.33}}
\put(105.66,7.78){\vector(1,2){5.33}}
\put(132.66,8.11){\vector(-1,1){10.00}}
\put(16.00,21.11){\makebox(0,0)[cc]{$e$}}
\put(29.66,25.44){\makebox(0,0)[cc]{$\pi_1$}}
\put(60.33,25.44){\makebox(0,0)[cc]{$\pi_2$}}
\put(133.00,25.78){\makebox(0,0)[cc]{$\pi_n$}}
\put(145.33,21.44){\makebox(0,0)[cc]{$f$}}
\put(77.00,25.11){\makebox(0,0)[cc]{$e_2$}}
\put(122.00,25.44){\makebox(0,0)[cc]{$e_{n-1}$}}
\put(100.33,25.44){\makebox(0,0)[cc]{$S$}}
\put(96.33,30.11){\makebox(0,0)[cc]{$\dots$}}
\put(96.33,13.11){\makebox(0,0)[cc]{$\dots$}}
\put(26.66,4.78){\makebox(0,0)[cc]{$l(\pi_1,e_1)$}}
\put(52.66,4.11){\makebox(0,0)[cc]{$l(\pi_2,e_1)$}}
\put(78.33,4.11){\makebox(0,0)[cc]{$l(\pi_2,e_2)$}}
\put(104.66,3.78){\makebox(0,0)[cc]{$l(\pi_{n-1},e_{n-1})$}}
\put(134.00,3.11){\makebox(0,0)[cc]{$l(\pi_n,e_{n-1})$}}
\put(49.33,25.11){\makebox(0,0)[cc]{$e_1$}}
\put(33.33,21.11){\circle*{0.94}}
\put(46.66,21.11){\circle*{1.33}}
\put(60.33,21.11){\circle*{0.94}}
\put(73.66,21.11){\circle*{1.33}}
\put(95.00,20.56){\makebox(0,0)[cc]{...}}
\put(19.33,21.00){\line(1,0){14.00}}
\put(129.67,20.67){\line(1,0){12.67}}
\put(19.00,13.00){\vector(0,1){17.00}}
\put(46.67,13.00){\vector(0,1){17.00}}
\put(73.67,13.00){\vector(0,1){17.00}}
\put(117.00,13.00){\vector(0,1){17.00}}
\put(117.00,20.67){\circle*{1.33}}
\put(129.67,20.67){\circle*{0.67}}
\put(142.33,13.00){\vector(0,1){17.00}}
\put(142.33,20.67){\circle*{1.33}}
\put(19.00,21.00){\circle*{1.33}}
\put(14.00,8.33){\vector(1,1){11.00}}
\put(10.33,5.00){\makebox(0,0)[cc]{$l(\pi_1,e)$}}
\put(148.00,8.00){\vector(-1,1){10.33}}
\put(154.67,3.11){\makebox(0,0)[cc]{$l(\pi_n,f)$}}
\end{picture}

\begin{center}
\nopagebreak[4] Fig. \theppp.

\end{center}
\addtocounter{ppp}{1}

A band is said to be a \label{ktabands}$k$-{\em band}
($\theta$-{\em band}, $a$-{\em band}) if $S$ consists of
$k$-letters ($\theta$-letters, $a$-letters). Similarly, a band is
an $a_i(z)$-band (a $z_j$-band) if $S$ consists of one $a$-letter
$a_i(z)$ (of $z_j(r,i)$- and $z_j(r,i)$-letters) for $z\in \tkk$.

Notice that an $S$-band may contain no cells. In this case the
band is called \label{emptyband}{\em empty}.

We say that two bands {\em cross} if their medians cross. We say
that a band is an \label{annulus}{\em annulus} if its median is a
closed curve. In this case the first and the last cells of the
band coincide. \setcounter{annulus}{\value{ppp}}(see Fig.
\theannulus a).

\bigskip

\begin{center}
\unitlength=1.5mm
\linethickness{0.5pt}
\begin{picture}(101.44,22.89)
\put(30.78,13.78){\oval(25.33,8.44)[]}
\put(30.78,13.89){\oval(34.67,18.00)[]}
\put(39.67,9.56){\line(0,-1){4.67}}
\put(25.67,9.56){\line(0,-1){4.67}}
\put(18.56,14.89){\line(-1,0){5.11}}
\put(32.56,7.11){\makebox(0,0)[cc]{$\pi_1=\pi_n$}}
\put(21.00,7.33){\makebox(0,0)[cc]{$\pi_2$}}
\put(15.89,11.78){\makebox(0,0)[cc]{$\pi_3$}}
\put(43.44,12.89){\line(1,0){4.67}}
\put(43.00,8.44){\makebox(0,0)[cc]{$\pi_{n-1}$}}
\put(19.22,10.89){\line(-1,-1){4.00}}
\put(25.89,20.44){\circle*{0.00}}
\put(30.78,20.44){\circle*{0.00}}
\put(30.33,1.56){\makebox(0,0)[cc]{a}}
\put(35.44,20.44){\circle*{0.00}}
\put(62.78,4.89){\line(0,1){4.67}}
\put(62.78,9.56){\line(1,0){38.67}}
\put(101.44,9.56){\line(0,-1){4.67}}
\put(101.44,4.89){\line(-1,0){38.67}}
\put(82.11,16.22){\oval(38.67,13.33)[t]}
\put(62.78,15.78){\line(0,-1){7.33}}
\put(101.44,16.44){\line(0,-1){8.22}}
\put(82.00,12.22){\oval(27.78,14.67)[t]}
\put(67.89,12.89){\line(0,-1){8.00}}
\put(95.89,13.11){\line(0,-1){8.22}}
\put(67.89,13.78){\line(-1,1){4.67}}
\put(95.67,14.00){\line(1,1){5.11}}
\put(73.67,9.56){\line(0,-1){4.67}}
\put(89.67,9.56){\line(0,-1){4.67}}
\put(76.56,21.33){\circle*{0.00}}
\put(81.44,21.33){\circle*{0.00}}
\put(86.11,21.33){\circle*{0.00}}
\put(65.22,12.22){\makebox(0,0)[cc]{$\pi_1$}}
\put(65.22,6.89){\makebox(0,0)[cc]{$\pi$}}
\put(71.00,6.89){\makebox(0,0)[cc]{$\gamma_1$}}
\put(92.78,7.11){\makebox(0,0)[cc]{$\gamma_m$}}
\put(98.56,7.11){\makebox(0,0)[cc]{$\pi'$}}
\put(98.33,12.67){\makebox(0,0)[cc]{$\pi_n$}}
\put(94.11,20.22){\makebox(0,0)[cc]{$S$}}
\put(85.89,6.89){\makebox(0,0)[cc]{$T$}}
\put(76.56,6.89){\circle*{0.00}} \put(79.44,6.89){\circle*{0.00}}
\put(82.33,6.89){\circle*{0.00}}
\put(83.44,1.56){\makebox(0,0)[cc]{b}}
\end{picture}
\end{center}
\begin{center}
\nopagebreak[4]

Fig. \theppp.
\end{center}
\addtocounter{ppp}{1}

Let $\cal B$ be an $S$-band with common $S$-edges $e_1,
e_2,\ldots, e_n$ which is not an annulus. Then the first cell has
an $S$-edge $e$ which forms a pair with $e_1$ and the last cell of
$\cal B$ has an edge $f$ which forms a pair with $e_n$. Then we
shall say that $e$ is the {\em start edge} of $\cal B$ and $f$ is
the {\em end edge} of $\cal B$.  If $p$ is a path in $\Delta$ then
we shall say that a band starts (ends) on the path $p$ if $e$
(resp. $f$) belongs to $p$.

Let $S$ and $T$ be two disjoint sets of letters, let ($\pi$,
$\pi_1$, \ldots, $\pi_n$, $\pi'$) be an $S$-band and let ($\pi$,
$\gamma_1$, \ldots, $\gamma_m$, $\pi'$) is a $T$-band. Suppose
that:
\begin{itemize}
\item the parts of the medians of these bands lying between their
intersections, form a simple closed curve,
\item on the boundary  of $\pi$ and on the boundary of $\pi'$ the pairs of $S$-edges separate the pairs of $T$-edges,
\item the start and end edges of these bands are not contained in the
region bounded by the medians of the bands.
\end{itemize}
Then we say that these bands form an {\em $(S,T)$-annulus} and the
simple curve formed by the medians of these bands is the {\em
median} of this annulus (see Fig. \theannulus b).
\setcounter{nonannulus}{\value{ppp}} Fig. \thenonannulus\ shows
that a multiple intersection of an $S$-band and a $T$-band  does
not necessarily produce an $(S,T)$-annulus. Again, for simplicity,
we say that this annulus is a $(\theta,k)$-{\em annulus}, if $S$
is a $\theta$-band and $T$ is a $k$-band. Similarly, the notion of
a $(\theta,a)$-{\em annulus} is defined.
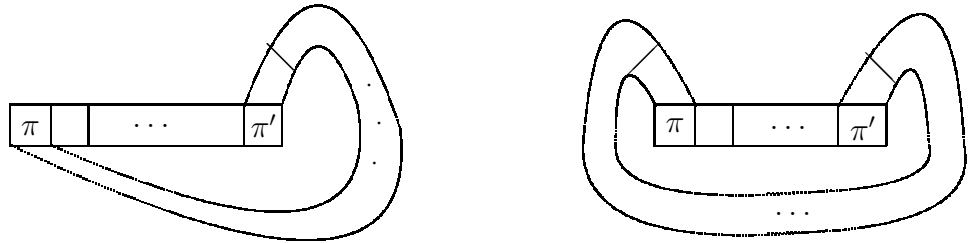
\begin{center}
\unitlength=1.50mm
\linethickness{0.5pt}
\begin{picture}(97.67,32.89)
\put(11.67,10.67){\framebox(24.00,3.56)[cc]{}}
\bezier{88}(35.67,14.22)(38.78,24.67)(42.33,14.22)
\bezier{200}(42.33,14.22)(45.22,-2.67)(15.22,10.67)
\bezier{168}(32.33,14.22)(39.22,32.89)(46.11,12.22)
\bezier{240}(46.11,12.22)(48.33,-6.89)(11.67,10.67)
\put(15.22,10.67){\line(0,1){3.56}}
\put(32.33,14.22){\line(0,-1){3.56}}
\put(13.44,12.22){\makebox(0,0)[cc]{$\pi$}}
\put(34.11,12.45){\makebox(0,0)[cc]{$\pi'$}}
\put(24.11,12.45){\makebox(0,0)[cc]{$\ldots$}}
\put(36.78,17.33){\line(-1,1){2.40}}
\put(43.44,16.00){\circle*{0.00}}
\put(44.33,12.67){\circle*{0.00}} \put(43.89,9.11){\circle*{0.00}}
\put(18.56,14.22){\line(0,-1){3.56}}
\put(68.78,10.67){\framebox(20.44,3.56)[cc]{}}
\put(72.33,10.67){\line(0,1){3.56}}
\put(72.33,14.22){\line(0,0){0.00}}
\put(75.67,14.22){\line(0,-1){3.56}}
\put(85.00,14.22){\line(0,-1){3.56}}
\put(70.56,12.44){\makebox(0,0)[cc]{$\pi$}}
\put(87.22,12.22){\makebox(0,0)[cc]{$\pi'$}}
\bezier{80}(89.22,14.22)(92.33,22.00)(93.00,10.67)
\bezier{80}(93.00,10.67)(94.11,6.44)(78.78,6.44)
\bezier{100}(84.33,6.67)(65.00,4.89)(65.22,10.67)
\bezier{72}(65.22,10.89)(65.22,20.89)(68.78,14.22)
\bezier{168}(85.00,14.22)(95.00,32.00)(96.11,10.67)
\bezier{108}(96.11,10.67)(97.67,2.89)(78.78,2.89)
\bezier{128}(83.89,3.11)(62.11,0.89)(62.33,10.67)
\bezier{156}(62.33,10.67)(63.89,30.89)(72.33,14.22)
\put(90.11,16.22){\line(-1,1){2.67}}
\put(66.33,16.89){\line(1,1){2.89}}
\put(80.56,12.22){\makebox(0,0)[cc]{$\ldots$}}
\put(81.00,4.67){\makebox(0,0)[cc]{$\dots$}}
\end{picture}
\end{center}

\begin{center}
\nopagebreak[4] Fig. \theppp.
\end{center}
\addtocounter{ppp}{1}

If  $\ell$ is the median of an $S$-annulus or an $(S,T)$-annulus
then the maximal subdiagram of $\Delta$ contained in the region
bounded by $\ell$ is called  the \label{innerdiag}{\em inner
diagram} of the annulus.

The union of cells of an $S$-band $\cal B$ forms a subdiagram. The
reduced boundary of this diagram, which we shall call the
\label{boundband}{\em boundary of the band}, has the form $e^{\pm
1}pf^{\pm 1}q\iv$ (recall that we trace boundaries of diagrams
clockwise). Then we say that $p$ is the \label{toppp}{\em top
path} of $\cal B$, denoted by \label{toppq}$\topp({\cal B})$, and
$q$ is the \label{bottt}{\em bottom path} of $\cal B$, denoted by
\label{bottq}$\bott({\cal B})$.


We shall call an $S$-band \label{maximalb}{\em maximal} if it is
not contained in any other $S$-band. If an $S$-band $\ww$ starts
on the contour of a cell $\pi$, does not contain $\pi$  and is not
contained in any other $S$-band with these properties then we call
$\ww$ a {\em maximal $S$-band starting on the contour of $\pi$}.

We now list all types of bands in diagrams over the presentation
$\rrr'$, that will be considered further.

\medskip
1. By definition, any $\Theta$-cell can be included in a
$\Theta$-band ($\Theta$-annulus) of a diagram $\Delta$ over
$\rrr'$. The cells in these bands correspond to the relations
(\ref{mainrel}), (\ref{auxtheta}). Similarly  $\bar\Theta$-bands
and $\theta$-bands are defined.

It is easy to see from (\ref{mainrel}) and (\ref{auxtheta}) that
if a $\theta$-band contains a $\theta(\tau,z)$-cell, for some
$\tau\in \csss$, then all cells in that band correspond to the
same rule $\tau$. Thus we can talk about $\Theta(\tau)$-bands and
$\bar\Theta(\tau)$-bands. Since every $\theta$-cell has two
$\theta$-edges, any maximal $\theta$-band, which is not an
annulus, must start on a boundary edge of the diagram.

2. By definition, an $a_j(z)$-band, where $z\in \tkk $,
$a_j\in\aaa\cup\bar\aaa$, is constructed of
$(a_j(z),\theta)$-cells and $(a_j(z),x)$-cells (see relations
(\ref{auxtheta}), (\ref{auxax})).

These bands can start and end on the boundary of a diagram,
on $(\theta,k)$-cells and on $\gr$-cells. Thus, an $a$-band is a
$a_j(z)$-band for some $j, z$.

3. For every $k\in \tkk$, a $k$-band is constructed of
$(z,\theta)$- and $(z,x)$-cells where $z\in \kk\cup\tkk$ and $z$
projects onto $k$ (see \ref{hardware}). These bands can start and
end on the boundary of the diagram or on hubs.

Notice that the reducibility of the top and the bottom labels of a
$\theta$-band can be easily achieved without changing the boundary
or the type of the diagram (see \cite{SBR}).

Indeed, if the label of the path $\bott(\bb)$ is not reduced then
it contains a subpath $e_1e_2$ where $e_1$, $e_2$ are edges
labelled by $a$ and $a\iv$ respectively, where $a$ is a letter.
\setcounter{pdtwo}{\value{ppp}} These two edges must belong to two
different cells $\pi_1$ and $\pi_2$ of $\bb$ as in the diagram on
the left in the Figure \thepdtwo.

Then we fold $e_1$ and $e_2$, producing a new edge $e$ labelled by
$a$ and we introduce a new edge $f$ with label $a$ which has a
common end vertex with $e$ so that cells $\pi_1$ and $\pi_2$ have
a common edge $e$. The other cells which were attached to $e_1$
and $e_2$ will be attached to the edge $f$. The bottom path of the
$\theta$-band $\bb$ in the new diagram is shorter (by two edges)
than the bottom path of $\bb$ in $\Delta$. The bottom paths of the
other $\theta$-bands are not affected by this operation. Thus
after a finite number of such operations we get a diagram with the
same boundary label as $\Delta$, in which all bands $\bb$ have top
and bottom paths with reduced labels.

\unitlength=0.70mm \special{em:linewidth 0.5pt}
\linethickness{0.5pt}
\begin{picture}(205.66,92.33)
\bezier{684}(0.33,20.67)(47.00,92.33)(92.66,20.67)
\bezier{528}(0.33,20.67)(80.00,-20.00)(92.66,20.67)
\put(37.33,23.33){\line(1,0){12.00}}
\put(49.33,23.33){\line(2,3){5.55}}
\bezier{160}(55.00,31.33)(38.00,42.67)(37.50,23.33)
\put(49.33,23.33){\line(1,-1){5.67}}
\bezier{160}(37.50,23.33)(44.33,2.00)(54.66,18.33)
\put(44.66,28.67){\makebox(0,0)[cc]{$\pi_1$}}
\put(44.66,18.00){\makebox(0,0)[cc]{$\pi_2$}}
\put(54.33,26.67){\makebox(0,0)[cc]{$a$}}
\put(54.33,21.33){\makebox(0,0)[cc]{$a$}}
\put(49.33,23.33){\vector(2,3){5.50}}
\put(49.33,23.55){\vector(1,-1){5.50}}
\put(97.33,20.67){\vector(1,0){10.00}}
\bezier{684}(113.33,20.34)(160.00,92.00)(205.66,20.34)
\bezier{528}(113.33,20.34)(193.00,-20.33)(205.66,20.34)
\put(150.33,23.00){\line(1,0){12.00}}
\bezier{160}(168.00,31.00)(151.00,42.34)(150.66,23.00)
\bezier{160}(150.66,23.00)(157.33,1.67)(167.66,16.00)
\put(157.66,28.34){\makebox(0,0)[cc]{$\pi_1$}}
\put(157.66,17.67){\makebox(0,0)[cc]{$\pi_2$}}
\put(162.00,23.00){\vector(1,0){11.67}}
\bezier{44}(174.67,23.00)(172.33,29.33)(168.00,31.00)
\bezier{40}(174.33,23.00)(167.33,16.67)(167.33,16.00)
\put(186.33,23.00){\vector(-1,0){12.67}}
\put(167.67,25.67){\makebox(0,0)[cc]{$a$}}
\put(179.33,25.67){\makebox(0,0)[cc]{$a$}}
\put(163.67,23.00){\circle*{1.33}}
\end{picture}

\begin{center}
\nopagebreak[4] Figure \theppp.
\end{center}
\addtocounter{ppp}{1}

{\bf Thus we will change the definition of reduced diagrams by
demanding in addition that the top and bottom labels of any
$\theta$-band in a reduced diagram are reduced words.}

\subsection{Forbidden annuli}

Similarly to \cite{talk}, \cite{OlSaBurns} the absence of annuli
of various kinds in reduced diagrams without hubs is important in
our present considerations. The following lemma is similar to
Lemma 3.1 of \cite{OlSaBurns}.

\begin{lemma} \label{NoAnnul}
Let $\Delta$ be a reduced (in particular, a minimal) van Kampen
diagram over the group $\hhh_1$. Then $\Delta$ has no

(1) $\theta$-annuli;

(2) $k$-annuli;

(3) $a$-annuli;

(4) $(\theta, k)$-annuli;

(5)$(\theta, a)$-annuli;
\end{lemma}

\proof We prove the lemma by contradiction. Let us choose a
minimal annulus $\ttt$ in $\Delta$ among all counterexamples to
any of the assertions (1) - (5). Here ``minimal" means that the
inner subdiagram of the annulus is of the smallest possible type.

Let $\ttt$ be a $\theta$-annulus. Suppose it contains a
$(\theta,k)$-cell. Notice that in every $(\theta,k)$-cell the pair
of $\theta$-edges separates cyclically the pair of $k$-edges.
Hence both top and bottom sides of $\ttt$ contain $k$-edges. Thus
the inside subdiagram of $\ttt$ contains a $k$-band starting on
the top or bottom of $\ttt$. Since $\Delta$ does not have hubs,
there exists a maximal  $k$-band $\LL$ starting and ending on the
contour of the inside diagram. This gives rise to a
$(\theta,k)$-annulus with a smaller inside subdiagram than $\ttt$
which contradicts the choice of $\ttt$ (as the smallest
counterexample) and property (4) of the lemma. Hence all cells in
$\ttt$ are $(\theta,a)$-cells corresponding to the auxiliary
relations (\ref{auxtheta}).

Similarly, the inside subdiagram $\Delta'$ of $\ttt$ contains no
$k$-edges. Hence the only possible cells in $\Delta'$ where an
$a$-band can end would be $\gr$-cells (whose contours consist of
$a_i(P_1)$-edges). Hence if $\Delta'$ contains an $a_i(z)$-edge,
and $z\ne P_1$, the maximal $a_i(z)$-band of $\Delta'$ containing
that cell must start and end on the boundary of $\Delta'$. This
gives rise to a $(\theta,a)$-annulus ruled out by property (5) of
the lemma. Thus all $a$-edges in $\Delta'$ are $\aaa(P_1)$-edges.
Hence all cells in $\ttt$ are $(\Theta, \aaa(P_1))$-cells
corresponding to relations $(\ref{auxtheta})$. Notice that all
these relations are commutativity relations because
$\alpha_{\tau\iv}(a_i(P_1))=\alpha_\tau(a_i(P_1))=a_i(P_1)$. Hence
the labels of the top and the bottom paths of $\ttt$ coincide
which contradicts the assumption that $\Delta$ is reduced (see
part (iv) of the definition of non-reduced diagrams).

  Now, let ${\cal T}$ be a $(\theta,k)$-annulus, ${\cal T}_{\theta}$ its
$\theta$-band, and ${\cal T}_k$ its $k$-band. We denote by $\pi$
and $\pi'$ their common cells, i.e. the {\em corner cells} of the
annulus ${\cal T}$. The inner diagram $\Gamma$ of $T$ is bounded
by a path $p_{\theta}p_k$ where $p_{\theta}$ is a subpath of
$\partial({\cal T}_{\theta})$ and $p_k$ is a subpath of $\partial
({\cal T}_k)$.

Property (4), which holds for smaller annuli, ensures that there
are no $\theta$-edges in $p_k$ and there are no $k$-edges in
$p_{\theta}$. Hence every non-corner cell of ${\cal T}_k$ is a
cell corresponding to relations (\ref{auxkx}). There exist no such
cells if ${\cal T}_k$ is a $P_j$- or $R_j$-band for any $j$. Then
$\pi$ and $\pi'$ must have a common $k$-edge. Moreover, $\pi$ is a
mirror copy of $\pi'$, since these cells belong to the same
$\theta$-band ${\cal T}_{\theta}$ as well. But this contradicts
the assumption that the diagram $\Delta$ is reduced. Hence we may
assume that ${\cal T}_k$ is neither $P$- nor $R$-band.

Since every non-corner cell of ${\cal T}_k$  corresponds to
(\ref{auxkx}), we have $\phi(p_k)\equiv U$ where $U$ is a word in
$\xxx$ (see relations (\ref{auxkx}). If $p_k$ contains $a$-edges
of the corner cells, then we just make the path $p_k$ shorter and
make $p_{\theta}$ longer, respectively, keeping the boundary path
of $\Gamma$ equal to $p_{\theta}p_k$.

Again, by (4), for smaller diagrams, every non-corner cell of
${\cal T}_{\theta}$ is a $(\theta,a)$-cell, and the forms of
relations (\ref{auxtheta}), (\ref{mainrel}) show that
$\phi(p_{\theta})\equiv \alpha_{\tau^{\pm 1}}(V)$ for some $V$ in
$a$-letters. The word $V$ has no letters from ${\cal A}(P_1)$
because $\pi$ and $\pi'$ are neither $P_1$- nor $R_1$-cells.

The equality $U\alpha_{\tau^{\pm 1}}(V)=1$ given by $\Gamma$ holds
over the group \label{hhh2}${\cal H}_2$ given by all the $x$-,
$k$-, and $a$-generators, and  $x$-relations (\ref{auxax}) and
(\ref{auxkx}) (but not (\ref{auxtheta})). Indeed, by statements
(1) - (5) for smaller annuli, $\Gamma$ may contain only cells
corresponding to relations of ${\cal H}_2$ and $\gr$-cells. If we
collapse all $\aaa(P_1)$-edges in $\Gamma$, we obtain a diagram
over $\hhh_2$ with the same boundary label (since
$\partial(\Gamma)$ does not contain $\aaa(P_1)$-edges). Thus we
can assume that $\Gamma$ does not contain $\gr$-cells.

The group ${\cal H}_2$ is a (multiple) HNN-extension of the free
group generated by $x$-letters where non-$x$-letters are stable.
Hence $V$ and $\alpha_{\tau^{\pm 1}}(V)$ are freely trivial words.
Then so is the word $U$.

Since there are no reducible pairs of $(k,x)$-cells in ${\cal
T}_k$, and $U=1$ in the free group, we see that $\pi$ and $\pi'$
must have a common $k$-edge, and we come to a contradiction as
above.

Similar arguments work in other cases (see lemma 3.1 in
\cite{OlSaBurns} and lemma 6.1 in \cite{Ol97} for
details).\endproof

\subsection{Surgeries involving $\gr$-cells}

As in our previous papers, we now describe certain surgeries which
can be applied to a diagram to reduce its type.

\begin{lemma}\label{aa-cells}
(1) ({\bf Moving a $\gr$-cell along an $a$-band.}) Assume $\Delta$
is a \vk diagram formed by a $\gr$-cell $\Pi$ and an
$\aaa(P_1)$-band $\ttt$ attached to $\Pi$ along an
$\aaa(P_1)$-edge $e$. Let $e^{-1}p$ and $eqe'q'$ be clockwise
contours of $\Pi$ and $\ttt$, respectively, where $e$ and $e'$ are
$a$-edges. Then there is a reduced diagram $\Delta'$ with contour
$\tilde p\tilde q\tilde e'\tilde q'$, where the labels of $\tilde
p, \tilde q, \tilde e',\tilde q'$ coincide with the labels of $p,
q, e',q'$, respectively, such that $\Delta'$ is formed by a copy
$\Pi'$ of the $\gr$-cell $\Pi$, attached to $\bar e',$ and by
several $(\Theta, \aaa(P_1))$-cells.

 (2) ({\bf Merging two $\gr$-cells connected by an $a$-band.})
Let $\Delta$ be a \vk diagram consisting of two distinct
$\gr$-cells $\Pi_1$, $\Pi_2$ and an $a_i(P_1)$-band $\ttt$
connecting them. Then there is a reduced diagram $\Delta'$ with
the same boundary label, which has at most one $\aaa$-cell and
several $(\Theta,\aaa(P_1))$-cells; and therefore the diagram
$\Delta$ is not reduced.

(3) There does not exist a reduced \vk diagram over $\hhh_1$ which
contains a $\gr$-cell $\Pi$ and an $\aaa$-band $\ttt$ that starts
and ends on $\Pi$.
\end{lemma}

\proof (1) Consider the cell $\pi_1$ of $\ttt$ containing $e$ in
the boundary. By relations (\ref{auxtheta}), the $\Theta$-letter
$\theta$ labelling two edges of $\pi_1$, commutes with all
$\aaa(P_1)$-letters in $\hhh_1$. Hence $\Pi$ can be surrounded by
a $\theta$-annulus (i.e. there exists a $\Theta$-annulus
containing $\pi$ such that the boundary label of the internal
subdiagram of this annulus coincides with the boundary label of
$\Pi$, so we can glue the annulus around $\Pi$). Next we can glue
the resulting diagram inside a $\Theta$-annuli starting containing
the second cell $\pi_2$ of $\ttt$, and so on. Notice that the
boundary labels of all diagrams constructed this way coincide with
the boundary label of $\Pi$ because all cells in the surrounding
annuli are commutativity cells.

Thus, we obtain a diagram $\Delta_0$ with the same boundary label
as $\Pi$, such that $\Delta_0$ contains attached to $e'$. This
subdiagram can be replaced by a copy $\Pi'$ of $\Pi$. Finally, to
restore the boundary label of $\Delta$, we can attach a subdiagram
$\Gamma$ to $\Pi'$, where $\Gamma$ is the mirror copy of
$\Delta_0\backslash\Pi$. The resulting diagram is what we need.
Figure \theppp\ illustrates the procedure used in this proof.

\bigskip

\unitlength=1.00mm \special{em:linewidth 0.5pt}
\linethickness{0.5pt}
\begin{picture}(143.82,59.33)
\put(23.50,49.00){\oval(25.00,20.67)[]}
\put(19.67,38.67){\line(0,-1){28.00}}
\put(19.67,10.67){\line(1,0){8.33}}
\put(28.00,10.67){\line(0,1){28.00}}
\put(19.67,32.67){\line(1,0){8.33}}
\put(19.67,17.00){\line(1,0){8.33}}
\put(26.00,10.67){\vector(-1,0){5.00}}
\put(19.67,19.00){\vector(0,1){10.33}}
\put(11.00,46.00){\vector(0,1){6.00}}
\put(23.00,48.33){\makebox(0,0)[cc]{$\Pi$}}
\put(6.67,48.33){\makebox(0,0)[cc]{$p$}}
\put(15.33,24.33){\makebox(0,0)[cc]{$q'$}}
\put(31.00,24.33){\makebox(0,0)[cc]{$q$}}
\put(23.67,6.33){\makebox(0,0)[cc]{$e'$}}
\put(37.33,34.33){\makebox(0,0)[cc]{$\Delta$}}
\put(28.00,29.00){\vector(0,-1){9.33}}
\put(24.00,27.33){\circle*{0.67}}
\put(24.00,24.33){\circle*{0.67}}
\put(24.00,21.67){\circle*{0.67}}
\put(23.67,36.33){\makebox(0,0)[cc]{$e$}}
\put(52.00,28.00){\makebox(0,0)[cc]{$\longrightarrow$}}
\put(70.00,10.67){\line(1,0){71.67}}
\bezier{312}(70.00,10.67)(69.00,56.67)(101.00,58.33)
\bezier{348}(101.00,58.33)(140.67,58.33)(141.33,10.67)
\bezier{276}(75.00,10.67)(74.00,51.00)(102.33,54.33)
\bezier{304}(102.67,54.33)(136.67,52.67)(136.00,10.67)
\bezier{136}(95.00,10.67)(90.33,29.67)(104.67,31.33)
\bezier{140}(105.00,31.33)(120.67,30.00)(116.67,10.67)
\bezier{108}(99.33,10.67)(94.00,25.00)(105.33,27.00)
\bezier{108}(105.00,27.00)(116.67,26.00)(113.00,10.67)
\put(105.33,17.67){\makebox(0,0)[cc]{$\Pi'$}}
\put(98.00,14.33){\line(-6,1){28.00}}
\put(97.67,18.00){\line(-3,1){27.00}}
\put(98.15,21.62){\line(-5,3){25.03}}
\put(99.68,24.49){\line(-1,1){21.78}}
\put(113.44,12.83){\line(6,1){27.71}}
\put(113.82,17.23){\line(4,1){26.37}}
\put(113.44,20.86){\line(3,1){25.60}}
\put(139.04,29.46){\line(0,0){0.00}}
\put(111.72,24.11){\line(2,1){25.03}}
\put(100.64,37.48){\circle*{0.85}}
\put(104.84,37.48){\circle*{0.76}}
\put(108.66,37.48){\circle*{0.76}}
\put(134.46,10.73){\vector(-1,0){11.46}}
\put(111.72,10.73){\vector(-1,0){11.27}}
\put(90.32,10.73){\vector(-1,0){12.04}}
\put(128.73,7.10){\makebox(0,0)[cc]{$q$}}
\put(105.80,7.10){\makebox(0,0)[cc]{$e'$}}
\put(83.82,7.48){\makebox(0,0)[cc]{$\tilde q'$}}
\put(143.82,38.06){\makebox(0,0)[cc]{$\Delta'$}}
\put(70.64,42.07){\makebox(0,0)[cc]{$\tilde p$}}
\put(98.00,58.20){\vector(1,0){4.33}}
\end{picture}

\begin{center}
\nopagebreak[4] Figure \theppp.
\end{center}
\addtocounter{ppp}{1}

(2) Indeed, the label of the top path of the connecting $a$-band
commutes with all letters in $\aaa(P_1)$. Hence we can apply Lemma
\ref{Minimal} (see also part (iii) of the definition of
non-reduced diagrams).

(3) Proving by contradiction, consider the subdiagram $\Gamma$
bounded by $\partial(Pi)$ and $\ttt$. It follows from Lemma
\ref{NoAnnul}, that $\Gamma$ has no cells except for $\gr$-cells,
and that $|\ttt|=0.$ If $\Gamma$ contains a $\gr$-cell, it can be
merged with $\Pi$. If $\Gamma$ does not contain cells, we obtain a
contradiction with the assumption that the boundary label of $\Pi$
is cyclically reduced.
\endproof

\subsection{Shifting indexes}

The following lemma utilizes the fact that our $\sss$-machines
work the same way in all parts of words
$\Sigma_{r,i}(w_1,w_2,w_3,w_4)$ between two consecutive
$K_j$-letters.

Recall that every letter in the list of generators of $\hhh$ has
the form $z(r,i)$ or $a_i(z)$ or $\theta(\tau,z)$, or
$x(a_i(z),\tau)$ or is a $\bar{ }$-brother of such a letter, where
$z\in\{\lel_j, L_j ,P_j, R_j\}$, $j=1,...,N$, $r\in \bee,
i=1,2,3,4,5$. That $j$ will be called the {\em index} of the
letter. We say that a word $W$ has index $j$ if all letters in $W$
have index $j$. If $W$ has index $j$ and $j'=1,...,N$, then
$\varepsilon_{j'}(W)$ denotes the word $W$ where the index $j$ is
replaced by $j'$ in all letters. The word
$\bar\varepsilon_{j'}(W)$ coincides with $\varepsilon_{j'}(W)$ if
$j'\ne 1$ and is obtained from $\varepsilon_{j'}(W)$ by removing
all $\aaa$-letters if $j'=1$.

\begin{lemma}\label{j-j'}
Suppose an equality $W=1$ holds in the group $\hhh_1$ where $W$
has index $j\ne 1$ and has no $\bar\Theta$-letters, $K_l$-letters
and $\bar K_l$-letters ($l=1,...,N$). Then $\varepsilon_{j'}(W)=1$
in $\hhh_1$ for any $j'=1,...,N$.
\end{lemma}

\proof Let $\Delta$ be a minimal \vk diagram with boundary label
$W$. By Lemma \ref{NoAnnul}, $\Delta$ has no $\bar\Theta$-edges
since its boundary has no $\bar\Theta$-edges.

Similarly, the diagram $\Delta$ does not contain $K_l$-edges and
$\bar K_l$-edges for any $l=1,...,N$, and if $e$ is a $z$-edge in
$\Delta$, $z\in\{L_{j_0},P_{j_0},R_{j_0}\}$, then $j_0=j$. Hence
if a maximal $\aaa(z)$-band in $\Delta$ terminates on a
$\Theta$-cell inside $\Delta$, then $z\in\{L_{j},P_{j},R_{j}\}$.

Since all $a$-edges of the boundary of $\Delta$ belong to
$\aaa(z)\cup\bar\aaa(z)$, $z\in\{L_{j},P_{j},R_{j}\}$, all
$a$-edges inside $\Delta$ must have labels from
$\aaa(z)\cup\bar\aaa(z)\cup\aaa(P_1)\cup\bar\aaa(P_1)$. If
$\Delta$ contains $\aaa(P_1)\cup\bar\aaa(P_1)$-edges, it must
contain $\gr$-cells. By Lemma \ref{aa-cells} (3) any $a$-band
starting on any $\gr$-cell must end on another $\gr$-cell which by
Lemma \ref{aa-cells} (2), (3) contradicts minimality of $\Delta$.
Hence $\Delta$ does not contain $\gr$-cells, and all $a$-edges of
$\Delta$ are $\aaa(z)\cup\bar\aaa(z)$-edges, $z\in\{L_j, P_j,
R_j\}$.

Since every $x$-cell contains a $\aaa$-letter or a $z$-letter,
$z\in \kkk$ and the indexes of all letters in an $a$-cell (resp. a
$z$-cell, $z\ne K_l$), are the same, all labels of $x$-edges in
$\Delta$ have index $j$.

Thus the indexes of the labels of all edges in $\Delta$ are equal
to $j$. When we replace the index $j$ by $j'$ in all letters of
$\Delta$, we obtain a new diagram $\Delta'$ with boundary label
$\varepsilon_{j'}(W)$. Since the $\varepsilon_{j'}$ takes
relations of $\hhh_1$ without $K_l$-letters, $\bar K_l$-letters
and $\bar\Theta$-letters to relations of $\hhh_1$, the diagram
$\Delta'$ is a \vk diagram over $\hhh_1$, so the equality
$\varepsilon_{j'}(W)=1$ holds in $\hhh_1$ as required.
\endproof

\begin{lemma} \label{j-j'(bar)}
Suppose an equality $W=1$ holds in the group $\hhh_1$, where the
word $W$ has index $j$ for some $j\ne 1$ and does not contain
$\Theta$-letters, $K_l$-letters,  and $\bar K_l$-letters for any
$l=1,...,N$. Assume that a reduced \vk diagram over $\hhh_1$ for
the equality $W=1$ contains no $(\aaa,\xxx)$-cells. Then
$\bar\varepsilon_{j'}(W)=1$ in $\hhh_1$ for any $j'=1,...,N$.
\end{lemma}

\proof  By Lemma \ref{NoAnnul}, $\Delta$ has no $\bar K_l$-edges
and $K_l$-edges, since there are no such edges on the boundary of
$\Delta$.Therefore the indexes of all labels of edges on the
contour of every cell in $\Delta$ are the same. Similarly, it does
not contain $\Theta$-edges, and the labels of all
$\bar\Theta$-edges in $\Delta$ have index $j$.

Since every $a$-band in $\Delta$ ends either on the boundary of
$\Delta$ or on a $\bar\Theta$-cell, all $a$-edges have index $j$
and they are $\bar\aaa$-edges. Similarly all $x$-edges in $\Delta$
(if any) have index $j$ as well.

Now if $j'\ne 1$, when we replace the index of every label in
$\Delta$ by $j'$, we get a \vk diagram over $\hhh_1$ with boundary
label $\bar\varepsilon_{j'}(W)$ as required.

Let $j'=1$. It is easy to see that when we apply
$\bar\varepsilon_1$ to every defining relation of $\hhh_1$
occurring in $\Delta$, we get a relation of $\hhh_1$ again. Hence
if we apply $\bar\varepsilon_1$ to all labels of $\Delta$, we
obtain a \vk diagram over $\hhh_1$ with boundary label
$\bar\varepsilon_1(W)$, so $\bar\varepsilon_1(W)=1$ in $\hhh_1$ as
required. \endproof

\section{The group $\hhh_2$}

\label{xrel}

Recall that $\hhh_2$ is the auxiliary group given by all the $x$-,
$k$-, and $a$-generators, and by $x$-relations (\ref{auxax}) and
(\ref{auxkx}) (but not (\ref{auxtheta})).

\begin{lemma} \label{K4}
 The word problem is solvable for the group
$\hhh_2$.
\end{lemma}

\proof The group $\hhh_2$ is an HNN-extension of the free group
generated by $x$-letters where non-$x$-letters are stable. The
associated subgroups are either subgroups generated by $\xxx(z)$,
$z\in \tkk$, or subgroups generated by fourth powers of elements
from $\xxx(z)$, $z\in \tkk$. Hence the membership problem is
decidable for all associated subgroups. So the claim of the lemma
follows from Britton's lemma \cite{LS}.\endproof

We say that a non-empty cyclically reduced word $w$ is
\label{uniformw}{\it uniform} if it is written in $\xxx(z)$ for
some $z\in \tkk$, $z\ne P_j$ . Two uniform words $w_1$ and $w_2$
are said to be \label{relatedu}{\it related} if one of them can be
obtained from another one, by a sequence of substitutions of the
following two forms, or their inverses:

(1) Substitution $x\rightarrow x^4$ applied to every letter of the
word;

(2) Every letter in the word is replaced by the corresponding
letter in $\xxx(z_-)$ if $z, z_-\ne P_j$.

It is obvious that the problem of whether two words in $\xxx$ are
related is decidable.

\begin{lemma}\label{uniform}
Let $\bb$ be an $a$- or a $k$-band in a diagram over $\hhh_2$.
Then the label of the top (bottom) path of $\bb$ is a uniform
word. The labels of the top and the bottom paths are related.
\end{lemma}

\proof This immediately follows from relations (\ref{auxax}) and
(\ref{auxkx}).
\endproof

\begin{lemma} \label{X-words} Two non-empty cyclically reduced words $w_1$ and $w_2$ in $\xxx$
are conjugate in $\hhh_2$ if and only if one of them is a cyclic
permutation of another one, or a cyclic shift of $w_1$ and a
cyclic shift of $w_2$ are relative uniform words.
\end{lemma}

\proof The condition is sufficient for the conjugacy of $w_1$ and
$w_2$ this follows from the conjugacy relations (\ref{auxax}) and
(\ref{auxkx}) in $\hhh_2$.

Now suppose that $w_1$ and $w_2$ are conjugate in $\hhh_2$.
Consider a minimal annular diagram $\Delta$ for the conjugacy of
$w_1$ and $w_2$, $w_1$ (resp. $w_2$) is the label of the inner
(outer) boundary of $\Delta$. We may assume that $w_2$ is not a
cyclic shift of $w_1$, and so the diagram $\Delta$ contains cells.
These cells belong to concentric $a$- and $k$-annuli $\bb_1,\dots,
\bb_d$, counted where the inner boundary to the outer boundary of
$\Delta$.

Therefore by Lemma \ref{uniform} each of the words $w_1, w_2$ must
be uniform, being a side label of the annuli $\bb_1$ or of
$\bb_d$. Also by Lemma \ref{uniform}, a cyclic shift of $w_1$ and
$w_2$ are related.
\endproof

\begin{lemma} \label{Conjugacy}
 The conjugacy problem is solvable in the
group $\hhh_2$.
\end{lemma}

\proof It suffices to find an effective upper bound for the number
of cells in an annular diagram $\Delta$ over $\hhh_2$, in terms of
the sum of lengths of the inner and outer boundaries $p_1, p_2$ of
$\Delta$.

By Lemma \ref{K4}, we may assume that the boundary labels
$w_1=\phi(p_1)$ and $w_2=\phi(p_2)$ of $\Delta$ are non-empty
words which have no pinches over the HNN-extension $\hhh_2$. This
implies that there is no $a$- or $k$-band in $\Delta$, with both
ends lying on $p_1$ (resp. on $p_2$). Indeed, otherwise there
would exist a band $\bb$ with ends $e_1$, $e_2$ lying, say, on
$p_1$, such that, the bottom or the top of $\bb$ is a subpath of
$p_1$ (since each cell of $\Delta$ belongs to a maximal $a$- or
$k$-band); then label of this side of $\bb$ would be a pinch in
$w_1$.

Obviously, Lemma \ref{X-words} solves the conjugacy problem for
words in $\xxx$. So we can assume that one of the words $w_1$ or
$w_2$ contains non-$\xxx$-letter. Hence every maximal $a$- and
$k$-band of $\Delta$ connects $p_1$ with $p_2$.

Let us enumerate these ``radial" bands cyclically
(counter-clockwise): $\bb_1,\dots,\bb_d$. For each $t=1,...,d$,
the boundary of $\bb_t$ has the form $e_t\topp(\bb_t)f_t\iv
\bott(\bb_t)\iv$, where the start edge $e_{t}$ belongs to $p_1$,
the end edge $f_t$ belongs to $p_2$.

Notice that the labels of the paths $\topp(\bb_t)$ and
$\bott(\bb_t)$ are reduced uniform words in $\xxx$ by Lemma
\ref{uniform}. Notice also that $\phi(\topp(\bb_t))$ coincides
with $\phi(\bott(\bb_{t+1}))$ (indices are taken modulo $d$)
except for some prefixes and suffixes whose total length is
bounded from above by a constant $c_0 =  |w_1|+|w_2|$.   Thus, in
particular, $||\bott(\bb_{t+1}| - |\topp(\bb_t)||\le c_0$.

Relations (\ref{auxax}) and (\ref{auxkx}) show that for every
$t=1,...,d$ we have $$|\topp(\bb_t)|=c_t|\bott(\bb_t)|$$ where
$c_t\in \{1/4, 1, 4\}$. If $c_t=1/4$, we call the band $\bb_t$
\label{contracting}{\em contracting}, if $c_t=4$, we call it {\em
expanding}.

Therefore, $|\topp(\bb_{t+1})|=c_{t+1}|\topp(\bb_t)|+c'_{t+1}$ for
some $c'_{t+1}\le c_0$. Hence
$$|\topp(\bb_{t+2})|=c_{t+2}c_{t+1}|\topp(\bb_t)|+c_{t+2}c'_{t+1}+c'_{t+2}.$$
Continuing this way, we finally obtain
$$|\topp(\bb_1)|\le 4^\sigma |\topp(\bb_1)| + c$$ where $\sigma$ is
the number of expanding bands minus the number of contracting
bands among $\bb_1,...,\bb_d$, and $c$ is some constant which is
bounded by $c_0(1+4+...+4^{d-1})<c_04^d/2$.

If $\sigma\ne 0$, then $|\topp(\bb_1)| \le c_04^d/2.$ Hence
$\topp(\bb_1)$ has a recursively bounded length (in terms of
$|w_1|+|w_2|$). Cutting $\Delta$ along $\topp(\bb_1)$, we obtain
an (ordinary, simply connected) \vk diagram $\Gamma$ with
perimeter recursively bounded in terms of $|w_1|+|w_2|$. Clearly,
$\Gamma$ is a minimal diagram over the presentation of $\hhh_2$
(otherwise there would be another diagram $\Gamma'$ with the same
boundary label and a smaller type; by gluing together two copies
of $\topp(\bb_1)$ on the boundary of $\Gamma'$, we obtain a
annular diagrams $\Delta'$ with boundary labels $w_1$, $w_2$,
whose type is smaller than that of $\Delta$, a contradiction). By
Lemma \ref{K4} the number of cells in the diagram $\Gamma$ is
recursively bounded in terms of $|w_1|+|w_2|$. Since $\Gamma$ and
$\Delta$ have the same number of cells, the number of cells in
$\Delta$ is recursively bounded.

Thus we may assume that $\sigma=0$.

Let us choose an arbitrary vertex $o=o_1$ on the path
$\topp(\bb_1)$ far enough from the initial and terminal vertices
of $\bb_1$. Here ``enough" means ``farther than some constant
$C>c_0$ whose value will become clear later".

Denote by $l_1>C$ the distance from $o_1$ to the initial vertex of
$\topp(\bb_1)$. Since $l_1>c_0$, the vertex $o_1$ can be connected
with a vertex $o_2$ of the path $\topp(\bb_2)$ by a path of the
length at most 4 (going at most three edges along $\topp(\bb_1)$
toward its initial vertex and then across the band $\bb_2$ along a
connecting edge. Let $l_2$ be the distance from $o_2$ to the
initial vertex of $\topp(\bb_2)$. It is clear that
$l_2=c_2l_1+c''_2$ for some recursively bounded constant $c''_2$,
$c_2\in \{4, 1/4, 1\}$. We can assume that $C$ is big enough so
that $l_2
>c_0$, so we can continue and find a vertex $o_3$ on
$\topp(\bb_3)$ connected with $o_2$ by a path of length at most 4.
Thus we can construct a sequence of points $o_i\in \topp(\bb_i)$,
$i=1,...,d+1$ (where $\bb_{d+1}=\bb_1$) and a path $p$ connecting
these points such that $|p|\le 4d$, and the distance
$|l_{d+1}-l_1|$ from $o_{d+1}$ to $o_1$ on $\topp(\bb_1)$ is equal
$4^\sigma l_1+C'$ for some recursively bounded constant $C'$ (see
Figure \theppp). Let $p'$ be the loop obtained by concatenation of
$p$ and a portion of $\topp(\bb_1)$ between $o_1$ and $o_{d+1}$.
Since $\sigma=0$, the length of $p'$ is bounded by a constant
$C''$.

\begin{center}
\unitlength=1.2 mm \special{em:linewidth 0.6pt}
\linethickness{0.5pt}
\begin{picture}(76.33,65.67)
\put(39.50,35.00){\oval(73.67,61.33)[]}
\put(22.33,35.50){\oval(26.67,49.67)[]}
\put(35.67,17.00){\line(1,0){40.67}}
\put(35.67,22.33){\line(1,0){40.67}}
\put(35.67,28.67){\line(1,0){40.67}}
\put(35.67,34.67){\line(1,0){40.67}}
\put(35.67,41.33){\line(1,0){40.67}}
\put(35.67,47.33){\line(1,0){40.67}}
\put(58.33,28.67){\rule{-3.33\unitlength}{0.00\unitlength}}
\put(56.67,29.00){\line(-1,0){4.67}}
\put(52.00,29.00){\line(0,1){6.00}}
\put(52.00,35.00){\line(-1,0){3.67}}
\put(48.33,35.00){\line(0,1){6.67}}
\put(48.33,41.67){\line(1,0){10.67}}
\put(59.00,41.67){\line(0,1){6.00}}
\put(59.00,47.67){\line(1,0){4.00}}
\put(63.00,29.00){\line(-1,0){6.67}}
\put(63.00,29.00){\line(0,-1){6.33}}
\put(63.00,22.67){\line(-1,0){6.00}}
\put(57.00,22.67){\line(0,-1){5.33}}
\put(57.00,17.33){\line(-1,0){9.33}}
\put(47.67,17.33){\line(0,-1){7.67}}
\bezier{216}(63.00,47.33)(64.33,65.00)(28.00,63.33)
\bezier{160}(28.00,63.33)(3.33,65.33)(4.33,50.33)
\bezier{228}(4.33,50.33)(2.67,5.67)(15.33,6.33)
\bezier{164}(47.67,10.33)(49.33,3.33)(15.33,6.33)
\put(57.67,29.00){\circle*{0.67}}
\put(57.67,32.00){\makebox(0,0)[cc]{$O_1$}}
\put(32.67,32.00){\makebox(0,0)[cc]{$\bb_1$}}
\put(53.33,36.67){\makebox(0,0)[cc]{$O_2$}}
\put(32.67,36.67){\makebox(0,0)[cc]{$\bb_2$}}
\put(47.00,44.00){\makebox(0,0)[cc]{$O_3$}}
\put(32.67,44.00){\makebox(0,0)[cc]{$\bb_3$}}
\put(11.67,32.00){\makebox(0,0)[cc]{$q$}}
\put(0.67,32.00){\makebox(0,0)[cc]{$p$}}
\put(58.00,49.00){\makebox(0,0)[cc]{$O_4$}}
\put(64.00,31.67){\makebox(0,0)[cc]{$O_{d+1}$}}
\put(56.00,24.00){\makebox(0,0)[cc]{$O_d$}}
\put(47.00,19.00){\makebox(0,0)[cc]{$O_{d-1}$}}
\put(39.67,30.00){\vector(-1,0){3.33}}
\put(53.00,30.00){\vector(1,0){4.67}}
\put(45.67,30.33){\makebox(0,0)[cc]{$l_1$}}
\put(39.33,27.00){\vector(-1,0){3.00}}
\put(57.67,27.00){\vector(1,0){4.67}}
\put(48.00,26.33){\makebox(0,0)[cc]{$l_{d+1}$}}
\end{picture}
\end{center}
\begin{center}
\nopagebreak[4] Figure \theppp.
\end{center}
\addtocounter{ppp}{1}

If the length of $\topp(\bb_1)$ is very large, then by choosing
different starting points $o_1$ on $\topp(\bb_1)$, we can built
very many loops $p'$ whose lengths are bounded by $C''$. Notice
that since the process of building the loop $p'$  starting at
$o_1$ is deterministic, if two of these loops share a common
vertex, they will share all the vertices coming after that one.
Hence if we take initial vertices of these loops sufficiently far
apart (say, farther than $2C''$), these loops won't intersect. The
number of different labels of such loops is recursively bounded
since their lengths are bounded. Hence if $|\topp(\bb_1)|$ were
bigger than that bound, two of the constructed paths would have
the same labels. But this would contradict the minimality of
$\Delta$ (in fact $\Delta$ would not be even reduced). Therefore
$|\topp(\bb_1)|$ is recursively bounded, and so we again found a
cut of $\Delta$ which has a recursively bounded length and turns
$\Delta$ into a simply connected \vk diagram with recursively
bounded perimeter.\endproof

\begin{lemma} \label{slender} Let $\Delta$ be a reduced diagram over $\hhh_1$ with
boundary $ep_1e'q_1(fp_2f')\iv q_2\iv,$ and the following
conditions hold:

\begin{enumerate}

\item[(1)] $e, e', f, f'$ are $\Theta$-edges,

\item[(2)] $\phi(p_1), \phi(p_2)$ are reduced words in $\xxx$,

\item[(3)] $q_1=\topp(\ttt_1)$, $q_2=\bott(\ttt_2)$ for some
$\Theta$-bands $\ttt_1, \ttt_2$ without $k$- and $a(P_j)$-cells,

\item[(4)] The bands $\ttt_1, \ttt_2$ are decomposed as
$\ttt_1=\ttt_1'\ttt_1''$, $\ttt_2=\ttt_2'\ttt_2''$ such that the
words  $\Lab(\bott(\ttt_1'))_a$, $\Lab(\bott(\ttt_1''))_a\iv$,
$\Lab(\topp(\ttt_2'))_a$ and $\Lab(\topp(\ttt_2''))_a\iv$ are
positive,

\item[(5)] The diagram $\Delta\backslash (\ttt_1\cup\ttt_2)$
contains no $\theta$-edges.
\end{enumerate}

Then
\begin{enumerate}

\item[(i)] the length of every $a$-band connecting
$\bott(\ttt_1')$ and $\topp(\ttt_2')$ is between
$\frac{|p_1|}{6^{|\ttt_1'|}}-2$ and $|p_1|$,

\item[(ii)] the length of an $a$-band, connecting
$\bott(\ttt_1'')$ and $\topp(\ttt_2'')$ is between
$\frac{|p_2|}{6^{|\ttt_1''|}}-2$ and $|p_2|$.

\item[(iii)] $|\ttt_1'| \le 2|p_1|+2$.
\end{enumerate}
\end{lemma}

\proof Parts (i) and (ii) of the lemma are symmetric so it is
enough to prove part (i). Denote by $q'_1$ (by $q'_2$) the bottom
(the top) of $\ttt'_1$ (of $\ttt''_2$). Let $\CC_1$, $\CC_2$ be
two neighbor maximal $a$-bands in $\Delta$ corresponding to two
consecutive $a$-letters in the words $\phi(q'_1)_a$ and
$\phi(q'_2)_a$. Denote by $U_1$ and $U'_1$ (by $U_2$ and $U'_2$)
the top and bottom labels of $\CC_1$ (of $\CC_2$), respectively.
By the positiveness condition (4) and the definition of
$\alpha_\tau$, we have that $U_1$ is freely equal to $x^{\pm
1}U'_2(x')^{\pm 1}$ for some $x$-letters $x,x'$, and $U'_2$ is
obtained from $U_2$ by replacing all letters by their fourth
powers (see relations (\ref{auxax})). Hence the length of $\CC_1$
is strictly greater than the length of $\CC_2$ if $|U_2|\ne 0$. If
$U_2$ is empty, then $x\ne x'$ because otherwise the
$\theta$-cells of $\ttt_1$ and $\ttt_2$, connected by $\CC_2$,
correspond to the same rule $\tau$ (defined by $x$) and form a
reducible pair of cells. Hence, in any case $|\CC_1|>|\CC_2|$.
Since the length of the first maximal $a$-band in $\Delta$
(counting from $p_1$) is at most $|p_1|$, we obtain the upper
estimate in part (i). Also we have from above considerations, that
$|U_1|\le 6|U_2|$ if $|U_2|\ne 0$, and this immediately implies
the lower estimate in parts (i) of the lemma.

Part (iii) is obtained as follows. The number of maximal $a$-bands
starting on $q_1'$ cannot exceed $|p_1|+1$ because the length of
the next $a$-band is smaller than the length of the previous one.
Also the length of $q_1'$ does not exceed twice the number of
$a$-edges on $q_1'$. Therefore $|q_1'|\le 2|p_1|+2$.
\endproof

\begin{lemma}\label{qqq} Let $\Delta$ be a diagram over $\hhh_2$
with reduced boundary $pq_1q_2q_3$, and for some $\tau\in\sss^+$,
$\phi(q_1)=\alpha_{\tau^{\pm 1}}(u)$,
$\phi(q_3)=\alpha_{\tau'^{\pm 1}}(v)$ for some reduced words $u,
v$ in $\aaa$, and $\phi(q_2)$ is a reduced word in $\xxx$. Suppose
further that every $a$-band that starts on $q_1$ ends on $p$. Then
the number of cells in $\Delta$ and the length of $|q_1q_2q_3|$
are recursively bounded in terms of $|p|$.
\end{lemma}

\proof We can assume that there are no $a$- or $k$-bands starting
and ending on $p$. Otherwise we could cut off one such band whose
side is a subpath of $p$, producing a diagram of smaller type
satisfying the conditions of the lemma (with side $p$ replaced by
a path $p'$ whose length is at most $4|p|$).

Thus by Lemma \ref{NoAnnul} we can assume that every cell in
$\Delta$ belongs to an $a$-band starting on $q_1\cup q_3$ and
ending on $p$. The number $l$ of such bands is at most $|p|$. Let
us number these bands: $\CC_1, \CC_2,...,\CC_{l_1},..., \CC_l$ so
that $\CC_i$ starts on the $i$-th $a$-edge of $q_1\cup q_3$, $l_1$
is the number of $a$-edges on $q_1$. Let $U_i$ be the label of the
top side of $\CC_i$ and $V_i$ be the label of the bottom side of
$\CC_i$ (see Figure \theppp). Without loss of generality we can
assume that $l_1>1$.

Since $\phi(q_1)=\alpha_{\tau^{\pm 1}}(u)$,
$\phi(q_3)=\alpha_{\tau'^{\pm 1}}(v)$, $|U_1|$ is bounded by
$|p|+1$, $|V_1|$ is at most $4|U_1|$ (see relations
(\ref{auxax})), $U_2$ is bounded by $2+|V_1|+|p|$, $|V_2|\le
4|U_2|$, etc. Hence the lengths of all bands $\CC_1,...,\CC_{l_1}$
are recursively bounded in terms of $|p|$. Similarly the lengths
of all bands $\CC_l,\CC_{l-1},...,\CC_{l_1+1}$, which start on
$q_3$, are recursively bounded in terms of $|p|$. Hence indeed,
the total number of cells in $\Delta$ is recursively bounded in
terms of $|p|$.
\bigskip

\begin{center}
\unitlength=1mm \special{em:linewidth 0.5pt} \linethickness{0.5pt}
\begin{picture}(102.67,65.00)
\put(12.00,6.33){\framebox(87.33,55.00)[cc]{}}
\bezier{300}(42.67,61.33)(53.67,34.67)(99.33,35.00)
\bezier{268}(48.00,61.33)(58.67,37.67)(99.67,39.00)
\bezier{148}(74.00,61.33)(77.67,46.67)(99.33,45.33)
\bezier{108}(83.00,61.33)(83.00,51.00)(99.33,49.67)
\bezier{60}(90.67,61.33)(92.33,54.00)(99.33,53.67)
\put(53.00,61.33){\vector(1,0){13.00}}
\put(71.00,6.33){\vector(-1,0){19.00}}
\put(12.00,26.33){\vector(0,1){19.67}}
\put(99.33,30.00){\vector(0,-1){12.00}}
\bezier{312}(99.33,33.00)(45.67,30.00)(40.00,6.33)
\bezier{268}(99.33,28.33)(57.67,28.67)(45.00,6.33)
\bezier{152}(99.33,20.33)(73.33,17.67)(71.67,6.33)
\bezier{76}(99.33,14.00)(88.00,13.67)(87.00,6.33)
\put(63.00,49.67){\makebox(0,0)[cc]{$U_i$}}
\put(50.33,43.67){\makebox(0,0)[cc]{$V_i$}}
\put(96.00,58.33){\makebox(0,0)[cc]{$\CC_1$}}
\put(89.00,56.00){\makebox(0,0)[cc]{$\CC_2$}}
\put(90.67,37.00){\makebox(0,0)[cc]{$\CC_i$}}
\put(74.00,46.33){\makebox(0,0)[cc]{${\bf \dots}$}}
\put(37.00,13.67){\makebox(0,0)[cc]{$\CC_{l_1+1}$}}
\put(81.00,11.00){\makebox(0,0)[cc]{$\CC_{l-1}$}}
\put(93.67,9.00){\makebox(0,0)[cc]{$\CC_l$}}
\put(7.33,32.67){\makebox(0,0)[cc]{$q_2$}}
\put(54.67,65.00){\makebox(0,0)[cc]{$q_3$}}
\put(102.67,33.33){\makebox(0,0)[cc]{$p$}}
\put(58.00,2.67){\makebox(0,0)[cc]{$q_1$}}
\end{picture}
\end{center}
\begin{center}
\nopagebreak[4] Figure \theppp.
\end{center}
\addtocounter{ppp}{1}

Since every $a$-band in $\Delta$ ends in $p$, the number of double
edges of the boundary (i.e. edges which appear twice in the
contour of the diagram and do not belong to the boundaries of
cells of the diagram) cannot exceed a constant times $p$. Hence
the length of the contour of $\Delta$ does not exceed a constant
times the number of cells in $\Delta$ plus a constant times $|p|$.
Thus $|q_1q_2q_3|$ is recursively bounded in terms of $|p|$.
\endproof

\section{The word problem in $\hhh_1$}

Recall that $\hhh_1$ denotes the group given by all generators of
$\hhh$, by all its defining relations, except for the hub
$\Sigma$, and by all $\gr$-relations.

\begin{lemma} \label{K0}
The words problem in $\hhh_1$ is decidable.
\end{lemma}

\proof It suffices to find an effective upper bound for the number
of cells and perimeters of all $\gr$-cells in a reduced van Kampen
diagram $\Delta$ over $\hhh_1$, depending on the perimeter
$|\partial(\Delta)|$ of $\Delta$. Indeed, since $\gr$ is embedded
into $\hhh$ by Lemma \ref{Emb}, $\gr$ is embedded into $\hhh_1$
(which satisfies fewer relations than $\hhh$ but which does
satisfy all $\gr$-relations. If we know the bound of perimeters of
all $\gr$-cells, we can list them all since the word problem in
$\gr$ is decidable. So in order to check if a word $w$ is equal to
1 in $\hhh_1$, we would have to consider finitely many \vk
diagrams with perimeter $|w|$.

First, arguing as in \cite{SBR} and \cite{talk}, we bound the
number of non-$\hhh_2$-cells in $\Delta$ as follows.

The number of maximal $\theta$-bands of $\Delta$ is at most
$|\partial(\Delta)|$ since these bands must start and end on the
boundary of $\Delta$ by Lemma \ref{NoAnnul} (1). Similar upper
bounds valid for the number of maximal $k$-bands, because the
diagrams over $\hhh_1$ have no hubs. Since every $(\theta,k)$-cell
is the intersection of a $\theta$-band and a $k$-band, and two
such bands can have at most one intersection by Lemma
\ref{NoAnnul} (4), we have a quadratic upper bound (in terms of
$|\partial(\Delta)|$) for the number of $(\theta,k)$-cells of
$\Delta$.

Since a maximal $a$-band in $\Delta$ cannot be an annulus by Lemma
\ref{NoAnnul}(3), and cannot connect two $\gr$-cells by Lemma
\ref{aa-cells}, each maximal $a$-band in $\Delta$ starts or ends
either on the boundary of $\Delta$ or on the boundary of a
$(\theta,k)$-cell. Hence the number of maximal $a$-bands in
$\Delta$ is quadratically bounded in terms of $|\partial\Delta|$.
Since every $(a,\theta)$-cell is an intersection of a
$\theta$-band and an $a$-band, and two such bands can intersect
only once (Lemma \ref{NoAnnul} (5)), the number of
$(a,\theta)$-cells in $\Delta$ is cubically bounded.

The cubic upper bound is also true for the sum of perimeters of
$\gr$-cells of $\Delta$. Indeed, each edge on the boundary of a
$\gr$-cell is the start edge of a $a$-band which ends either on
the boundary of $\Delta$ or on the boundary of a
$(\theta,k)$-cell.

Notice further, that if $\Gamma$ is a subdiagram of $\Delta$
bounded by a simple loop without $\theta$- and $\aaa(P_1)$-edges,
then $\Gamma$ is a diagram over $\hhh_2$. Indeed, $\Gamma$ has no
$\theta$-cells by Lemma \ref{NoAnnul}(1). Also $\Gamma$ has no
$\gr$-cells, because otherwise the boundary of $\Gamma$ would
contain an $\aaa(P_1)$-edge (the end edge of a $a$-band starting
on the $\gr$-cell), contrary to the assumption.

This implies that every maximal connected subdiagram $\Gamma$ of
$\Delta$ without $\theta$- and $\aaa(P_1)$-edges on the boundary
is simply connected. The maximality of $\Gamma$ implies that every
edge of the contour of $\Gamma$ belongs either to the contour of
$\Delta$ or to a contour of a non-$\hhh_2$-cell.

Hence the sum of perimeters of all such $\Gamma$ does not exceed
the sum of perimeters of all non-$\hhh_2$-cells of $\Delta$ plus
the perimeter of $\Delta$, i.e. it is effectively bounded. Now we
can apply Lemma \ref{K4}, which provides an effective upper bound
for the number of $\hhh_2$-cells in $\Delta$.
\endproof

Our main goal is to obtain similar result for the conjugacy
problem in $\hhh_1$. In this section, we use Lemma \ref{K0} for
getting some preliminary results.

\begin{lemma} \label{xdist} Let $\Delta$ be a reduced diagram over
the presentation consisting of relations (\ref{auxtheta}),
(\ref{auxax}) and $\gr$-relations with contour $pq$ where
$\Lab(p)$ is reduced a word in $\xxx$. Then $|p|$ is at most
$2|q|2^{\mu_a(q)}$ where $\mu_a(q)$ is the number of $a$-edges in
$q$, and the number and perimeters of cells in $\Delta$ is
recursively bounded as function of $|q|$.
\end{lemma}

\proof  By Lemmas \ref{NoAnnul} and \ref{aa-cells}, part (2),
every maximal $a$-band in $\Delta$ starts and ends on $q$ or
starts on $q$ and ends on a $\gr$-cell. This implies that the
number and the total perimeter of all $\gr$-cells in $\Delta$ are
bounded by $|q|$, since an $a$-band cannot start and end on
$\gr$-cells by part (iv) of the definition of reduced diagram and
Lemma \ref{aa-cells}(3). Using Lemma \ref{aa-cells}, part (1), we
can move the $\gr$-cells along the $a$-bands connecting these
cells and $q$, to $q$ and then remove them from the diagram
increasing the length of $q$ by at most a factor of 2. So we can
assume that $\Delta$ does not contain $\gr$-cells, and $q$ is
replaced by $q_0$ with $|q_0|\le 2|q|$ and the number $n(q_0)$ of
$a$-edges in $q_0$ for $a\notin\aaa(P_1)$ does not exceed the
number of such edges in $q$. Thus all maximal $a$-bands in
$\Delta$ start and end on $q_0$.

Let us use induction on the number $n_a(q_0)$ to prove that
$|p|\le |q_0|2^{n(q_0)}$ . If $n(q_0)=0$, then $\Delta$ contains
no cells and the statement is obvious. Since every cell in
$\Delta$ belongs to an $a$-band, there exists a maximal $a$-band
whose side is contained in $q_0$. Cutting this $a$-band off
$\Delta$, we get a diagram $\Delta'$ with boundary $pq'$ where
$(q')=n(q_0)-2$, $|q'|< 4|q_0|$. By the induction hypothesis,
$$|p|\le |q'|2^{n(q')}< 4|q_0|2^{n(q_0)-2}=|q_0|2^{n(q_0)}.$$
Since we have recursively bounded the number of $a$-bands in
$\Delta$ and the length of each $a$-band, the number of cells in
$\Delta$ is recursively bounded.
\endproof

\begin{lemma} \label{tdist} Let $\Delta$ be a reduced diagram over
the presentation consisting of relations (\ref{auxtheta}),
(\ref{auxax}) and $\gr$-relations. Suppose $\partial(\Delta)=pq$
where $\Lab(p)$ is the label of a side of a reduced $k$-band (the
$k$-band is not contained in $\Delta$, of course). Suppose also
that all $\theta$-bands starting on $p$ end on $q$. Then the
number and the perimeters of cells in $\Delta$ are recursively
bounded in terms of $|q|$.
\end{lemma}

\proof Let $\ttt_1,...,\ttt_n$ be the maximal $\theta$-bands
starting on $p$ (counted from the beginning of $p$ to the end of
$p$). Then $n\le |q|$. Let $\Delta_0$ be the subdiagram bounded by
$\partial(\Delta)$ and $\bott(\ttt_1)$ (the one that does not
contain $\ttt_1$), for every $i=1,...,n-1$ let $\Delta_i$ be the
subdiagram bounded by $\partial(\Delta)$, $\topp(\ttt_i)$ and
$\bott(\ttt_{i+1})$. Finally let $\Delta_n$ be the subdiagram
bounded by $\topp(\ttt_n)$ and $\partial(\Delta)$, not containing
$\ttt_n$. Then for every $i=0,...,n$, the $a$-edges on
$\partial(\Delta_i)\cap p$ belong to the initial or terminal
segments of $\partial(\Delta_i)\cap p$ of constant length (not
exceeding the length of a relator of the form (\ref{mainrel}) ).
Hence by Lemma \ref{NoAnnul}, all but a constant number of maximal
$a$-bands in $\Delta_i$ ($i=0,...,n$) starting on the boundary of
$\ttt_{i+1}$ ($i=0,...,n-1$), end either on $q$ or on a side of
$\ttt_{i}$ or on the boundary of a $\gr$-cell.

Suppose that one of the $a$-bands in $\Delta_i$ starting on a side
of $\ttt_{i+1}$ ($i=0,1,...$) ends on the boundary of a
$\gr$-cell. Then the label of an $a$-edge of $\ttt_{i+1}$ belongs
to $\aaa(P_1)$. This implies that the $\Theta$-edges of
$\ttt_{i+1}$ belong to $\Theta(P_1)$. This means that $p$ is a
side of a $P_1$- or $R_1$-band (see relations (\ref{mainrel}).
Since there are no $P_1$- or $R_1$-relations of the form
(\ref{auxkx}), the length of $p$ is at most a constant times the
number of $\theta$-edges in $p$, i.e. does not exceed a constant
times $|q|$. Thus the perimeter of $\Delta$ is linearly bounded in
terms of $|q|$, so the statement of the lemma follows from Lemma
\ref{K0}.

Thus we can assume that none of the $a$-bands in $\Delta_i$
starting on a side of $\ttt_{i+1}$ ends on the boundary of a
$\gr$-cell. Therefore all but a constant number of the $a$-bands
in $\Delta_i$ starting on the boundary of $\ttt_{i+1}$ end on
$\ttt_i$ or on $q$. Therefore the length of $\ttt_1$ is
recursively bounded in terms of $|q|$, and for each of each
$i=2,...,n$, $|\ttt_i|$ is recursively bounded in terms of $|q|$
and $|\ttt_{i-1}|$. Hence the lengths of all $\ttt_i$ are
recursively bounded in terms of $|q|$.

Thus for every $i=0,...,n$, the diagram $\Delta_i$ consists of
cells corresponding to relations (\ref{auxax}), and
$\gr$-relations the word $\Lab(\partial(\Delta_i)\cap p)$ has the
form $uvu'$ where $u$, $u'$ have lengths bounded by a constant,
and $v$ is a word in $\xxx$, and the complement of
$\partial(\Delta_i)\cap p$ to $\partial(\Delta_i)$ has length
recursively bounded in terms of $|q|$. Therefore by Lemma
\ref{xdist}, the length of $\partial(\Delta_i)\cap p$ is also
recursively bounded in terms of $|q|$, and the number and
perimeters of cells in $\Delta_i$ are also recursively bounded.
This implies that $|p|$ and the number of cells in $\Delta$ are
recursively bounded in terms of $|q|$.
\endproof

Let us introduce the following desirable properties of a boundary
component $p$ (inner or outer contour) of an annular diagram
$\Delta$ over $\hhh_1$:

\label{RRRR}
\begin{enumerate}
\item[(R1)] $p$ does not contain the boundary
of a $\gr$-cell as a subpath;

\item[(R2)] no $a$-band in $\Delta$ starts and ends on $p$;
\item[(R3)] no $\theta$-band in $\Delta$ starts and ends on $p$;

\item[(R4)] no $k$-band starts and ends on $p$.
\end{enumerate}

The next lemma shows that in some cases we can achieve these
properties without making too many changes in the diagram.

For every reduced path $p$ in a diagram over $\hhh_1$, denote by
$\mu_k(p)$ (resp. $\mu_\theta(p)$, $\mu_a(p)$, $\mu_x(p)$) the
number of all $k$-edges (resp. $\theta$-edges, $a$-edges,
$x$-edges) in $p$. The vector $(\mu_k(p), \mu_\theta(p), \mu_a(p),
\mu_x(p))$ will be called the \label{typep}{\em type} of $p$. We
order types lexicographically.

\begin{lemma}\label{rollsprep1} Let $\Delta$ be a reduced annular diagram over $\hhh_1$ with
boundary components $p$ and $q$. Then by removing some cells from
$\Delta$, one can construct an annular diagram $\Delta_1$
satisfying condition (R4) with boundary components $p_1$ and $q$
where $p_1$ satisfies condition (R4), such that $|p_1|$, the
number and perimeters of cells in the diagram
$\Delta\backslash\Delta_1$ are recursively bounded in terms of
$|p|$.
\end{lemma}

\proof Indeed, suppose that $\Delta$ contains a $k$-band starting
and ending, say, on $p$. Consider such a $k$-band $\bb$ which is
the closest to $p$, so that the subdiagram $\Gamma$ of $\Delta$
bounded by $\topp(\bb)$ and $p$ contains no $k$-bands. By Lemma
\ref{NoAnnul}, no $\theta$-band can start and end on $\bb$.
Therefore, by Lemma \ref{tdist}, the number and perimeters of
cells in $\Gamma$ are recursively bounded in terms of $|p|$.
Removing $\Gamma$ and $\bb$ we obtain a diagram with smaller
number of $k$-edges on the boundary. Thus we can continue that
operation at most $|p|$ times and obtain the desired diagram
$\Delta_1$.\endproof

\begin{lemma}\label{rollsprep} Let $\Delta$ be a reduced annular diagram over $\hhh_1$ with a
boundary components $p$ and $q$, and $p$ satisfy the condition
(R4). Then there exists another (annular) diagram $\Delta'$ with
contours $p'$ and $q$ satisfying (R1), (R2), (R3), (R4) and such
that:

(R0) The path $p'$  satisfies the following properties:
\begin{itemize}
\item The type of $p'$  is not greater than the type
 of $p$, and $|p'|$ is recursively bounded in
 terms of $|p|$.

\item The word
$\Lab(p')_\theta$  is freely conjugate of the word
$\Lab(p)_\theta$ ;

\item If $p$ do not contain $k$-edges then the word
$\phi(p')_a$ is conjugate of the word $\phi(p)_a$ modulo the
$\gr$-relations from $\Delta$;

\item The words $\phi(p)$ and  $\phi(p')$ are conjugate modulo relations
used in $\Delta$. The number and lengths of relations used to
deduce these conjugacies are recursively bounded in terms of
$|p|$;
\end{itemize}
\end{lemma}

\proof Let us describe three operations which can be applied to
$\Delta$.

1. If a $p$ contains a subpath which is the boundary of a
$\gr$-cell, we can remove that cell from $\Delta$ lowering the
length of the boundary. This operation does not increase the
number of $k$-edges, $\theta$-edges or $a$-edges on this boundary
component of the diagram. The $\theta$-projections of the labels
of the boundary components of the diagram do not change. Notice
that the resulting diagram satisfies (R0).

2. Suppose that an $a$-bands $\CC$ in $\Delta$ starts and ends on
 $p$. Let $\Delta_1$ be the subdiagram of $\Delta$
bounded by $\CC$ and $p$. It has no $k$-edges by property (R4) and
Lemma \ref{NoAnnul}, because $\CC$ cannot possess a $k$-cell.
Since the top/bottom of $\CC$ has no $a$-edges, we can apply lemma
2.12, as in the proof of Lemma \ref{xdist}, and remove all
$\gr$-cells of $\Delta_1$ increasing the length of the $p$-part of
its boundary at most twice, but not changing $\CC$.

Every maximal $a$-band in $\Delta_1$ starts and ends on $p$ since
$\Delta_1$ has neither $k$- nor $\gr$-cells. Therefore we could
find an $a$-band $\CC$ that starts and ends on, say, $p$, and such
that there are no cells between the top path of the band and $p$.
Then the length of this band would be bounded by $|p|$. Notice
that the number of $\theta$-edges on the top path of an $a$-band
is equal to the number of $\theta$-edges on the bottom path of the
band, and the corresponding numbers of $\xxx$-edges can differ at
most 4 times. Therefore we could delete that band and reduce the
number of $a$-edges in $\partial(\Delta)$ but not increase the
number of $\theta$-edges there and preserve property (R0) (in
particular, since $\CC$ connects two consecutive $a$-edges on $p$,
$\phi(p')_a$ is obtained from $\phi(p)_a$ by removing a subword of
the form $aa\iv$, so $\phi(p)_a=\phi(p')_a$ in the free group).

3. Suppose that a $\theta$-band $\ttt$ in $\Delta$ starts and ends
on $p$, and there are no $\theta$-cells in the simply connected
subdiagram $\Gamma$ of $\Delta$ between the top side of that band
an $p$. If $\topp(\ttt)$ is a part of $p$, then we can remove
$\ttt$ as in part 2, reducing the number of $\theta$-edges on the
boundary of the diagram, cancelling two consecutive mutually
inverse $\theta$-letters in the cyclic word $\Lab(p)_\theta$, and
not increasing the number of $k$-edges there, thus preserving (R0)
(notice, in particular, that if there are no $k$-edges on $p$,
then $\Lab(\topp(\ttt))_a\equiv \Lab(\bott(\ttt))_a$, whence
$\phi(p')_a\equiv\phi(p)_a$).

Assume now that there is a $\gr$-cell $\Pi$ in $\Gamma$, connected
with $\ttt$ by an $a$-band. This $a$-band has no cells, since
$\Gamma$ contains no $\theta$-cells. Hence the contour $p_\Pi
q_\Pi$ of $\Pi$ has a common part $q_\Pi$ of positive length with
$\ttt$. Then , by relations (\ref{auxtheta}), $\phi(p_\Gamma)$ and
$\phi(q_\Gamma)$ commute with the $\theta(P_1)$-letters labelling
the $\theta$-edges of the subband $\ttt_0$ of $\ttt$ with top side
$q_\Pi$.
 Therefore one can construct an auxiliary $\theta$-band $\ttt'$,
whose top label is obtained from that of $\ttt$ by replacing the
subwords $\Lab(q_\Pi)$ by $\Lab(p_\Pi)$. Then one can paste
$\ttt'$ to $\Gamma\backslash (\Pi\cup\ttt_0)$ reducing the number
of $\gr$-cells in $\Gamma$, but preserving property (R0). Notice
that this transformation does not affect the boundary labels of
$\Delta$.

Hence one may assume further that no $\gr$-cell of the reduced
subdiagram $\Gamma$ is connected with $\ttt$ by an $a$-band.
Therefore the sum of perimeters of all $\gr$-cells in $\Gamma$ is
not greater than the length of the $p$-part $p_\Gamma$ of its
boundary. Then, as in part 2, one can remove all $\gr$-cells from
$\Gamma$ increasing the number of $a$-edges in $p_\Gamma$ at most
twice and preserving (R0). Thus we may assume that $\Gamma$ has no
$\gr$-cells.

Let $\bb_1,...,\bb_r$ be all the maximal $a$- and $k$-bands in
$\Gamma$ starting on $p$, of non-zero length.
 Each of the bands $\bb_1,\dots \bb_r$
ends either on $\ttt$ or on $p$ (an $a$-band $\bb_i$ cannot end on
a $k$-cell not from $\ttt$ because every $k$-cell containing an
$a$-edge also contains a $\theta$-edge). Therefore there is one of
them (if $r>0$), say $\bb_1$, such that the top or bottom side
$p^1$ of $\bb_1$ has no edges in common with the sides of
$\bb_2,\dots,\bb_r$ and has at most 4 common $x$-edges with
$\ttt.$  ($\ttt$ can have at most 2 consecutive $x$-letter in its
boundary label, as it follows from the form of
$\theta$-relations.) Since all the edges of $p^1$ (with at most 4
exceptions) belong to $p$,  the length of $\bb_1$ is linearly
bounded. Removing this band from $\Gamma$ does not violate
property (R0). Considering the diagrams $\Gamma\backslash\bb_1$,
$\Gamma\backslash(\bb_1\cup\bb_2) \dots$ we can recursively bound
the length of every $\bb_2,\dots \bb_r$, because the number $r$
was linearly bounded above, and remove them all from $\Gamma$. But
the case where $\Gamma$ has no cells was considered in the
beginning.

Now repeating operations 1, 2, 3 recursively bounded (in terms of
$|p|$) number of times we will get the desired annular diagram.
\endproof

We call a word $W$ \label{minimalw}$\theta $-{\em minimal} if the
number of $\theta$-letters  in $W$ are not greater than the
similar numbers for any word $W'$ which is a conjugate of $W$ in
$\hhh_1$. Similarly we define $k$-{\em minimal} words. A word $W$
is called $\theta k$-{\em minimal} if it is both $\theta$- and
$k$-minimal. A word $W$ is said to be $\theta a$-{\em minimal}, if
the number of $\theta$-letters and the number of non-$\aaa(P_1)$
$a$-letters in $W$ are not greater than the corresponding numbers
for any word $W'$ which is a conjugate of $W$ modulo the relations
of $\hhh_1$ except for $(\theta,k)$-relations. A boundary
component of an (annular) diagram $\Delta$ over $\hhh_1$ is called
\label{minimalb}minimal if its label is $\theta k$-minimal and, if
$\Delta$ has no $(\theta,k)$-cells, it is also $\theta a$-minimal.
Let us consider reduced annular diagrams with boundary labels $W$
and $W'$ for various words $W'$. The immediate analysis of the
transformations used in lemmas \ref{rollsprep1} and
\ref{rollsprep} proves

\begin{lemma} \label{rollsprep2} For any word $W$ representing an element of $\hhh_1$, there is an
annular diagram $\Delta_0$ over $\hhh_1$ with boundary labels $W$
and $W_0$, where $W_0$ is $\theta k$-minimal ($\theta a$-minimal,
and $\Delta_0$ has no $(\theta,k)$-cells), such that the length of
$W_0$, the number, and the perimeters of cells in $\Delta_0$ are
recursively bounded as functions of $|W|$.
\end{lemma}

\section{Some special diagrams}
\label{SomeSpecial}

\subsection{$\theta$-bands and trapezia}

\label{bandsand}

 Let $W$ be a reduced word in the generators of $\hhh$. The {\em
base} of $W$ is the projection of $W$ onto $\tkk$.

Let $\tau\in\csss$, and $W\equiv U_1z_1U_2\dots z_nU_{n+1}$ a word
where $z_1,\dots,z_n$ are basic letters, $U_1,\dots,U_{n+1}$ are
words in $\aaa\cup\bar\aaa\cup\xxx$. Denote $z_0=(z_1)_-$,
$z_{n+1}=(z_n)_+$. For every word $V$ let $V_{\backslash x}$ be
the word $V$ with all $x$-letters deleted. We say that the word
$W$ is \label{taur}$\tau$-{\it regular} if

\begin{itemize}

\item for every $i=1,\dots,n+1$ the word $z_{i-1}(U_i)_{\backslash x}z_i$ is
admissible for $\csss$ and in the domain of $\tau$;

 \item if $\tau\in\sss$, then $W\equiv \alpha_\tau (W_{\backslash x})$,
 and if $\tau\in \bar sss$, then $W$ does not contain $x$-letters.

\end{itemize}

Thus for every admissible word $W$ of $\sss$ to which
$\tau\in\sss$ is applicable the word $\alpha_\tau(W)$ is
$\tau$-regular, and if $\tau\in\bsss$ is applicable to an
admissible word $W$ of $\bsss$ then $W$ itself is $\tau$-regular.

The following lemma immediately follows from the form of
$\theta$-relations from $\rr'$.

\begin{lemma} \label{Band}
Let $\bb$ be a reduced $\Theta(\tau)$-band for some
$\tau\in\csss$. Then the reduced labels of the top and the bottom
paths of $\bb$ have the same bases which will be called
\label{baseb}{\em the base} of $\bb$, and are $\tau$- or
$\tau\iv$-regular words.
\end{lemma}

For future references it is convenient to list all possible
2-letter bases of $\tau$-regular words (i.e. $2$-letter subwords
of $\tilde\Sigma^{\pm 1}$).

\begin{lemma} \label{Base}
Let $yz$ be a 2-letter base of a reduced $\Theta(\tau)$-band,
$y,z\in \tkk\cup\tkk\iv$. Then $yz$ or $(yz)\iv$ has one of the
following forms: $\lel_jL_j$, $L_jP_j$, $P_jR_j$, $R_j\rer_j$,
$z_jz_j\iv$, $z_j\iv z_j$ where $z\in \{K,L,P,R\}$, $j=1,...,N$.
In addition if $\tau\in\bsss$, $yz$ cannot have the form
$z_1z_1\iv$ or $z_1\iv z_1$ where $z\in\{L,R,P\}$, and it cannot
be equal to either $K_1K_1^{-1}$ or to $K_2^{-1}K_2$ . Moreover if
$\tau$ locks $zz_+$-sectors then $yz$ cannot have the form
$zz\iv$.
\end{lemma}

\proof Only the ``moreover" statement needs explanation. If $\tau$
locks $zz_+$-sectors, then there are no corresponding relations of
the form (\ref{auxtheta}).  Hence if the base has the form
$zz\iv$, the band has two consecutive mirror image
$(\theta,z)$-cells that cancel. This contradicts the assumption
that $\bb$ is reduced.
\endproof

The following lemma immediately follows from relations
(\ref{mainrel}), (\ref{auxtheta}).

\begin{lemma} \label{VV'}
Let $yz$ be the base of a $\Theta(\tau)$-band $\bb$. Let $V$ be
the label of the top path of $\bb$, $V'$ be the label of the
bottom path of $\bb$. Let $V=\alpha_{\tau^{\pm
1}}(W_1y(r,i)W_2z(r,i)W_3)$, $V'=\alpha_{\tau^{\mp 1}
}(W_1'y(r',i')W_2'z(r',i')W_3')$ provided $\tau\in\sss$ and let
$V=W_1y(r,i)W_2z(r,i)W_3$, $V'=W_1'y(r',i')W_2'z(r',i')W_3'$
provided $\tau\in\bsss$. Let $u(y_-), v(y), u(y), v(z)$ be the
words associated with $\tau$ and defined in Section
\ref{notation}. Then, in the free group, $W_1'=W_1v(y_-)$,
$W_2'=u(y)W_2v(y)$, $W_3'=u(z)W_3$.
\end{lemma}

A \label{trapeziumt}{\it trapezium} (see Figure \theppp) is a
reduced van Kampen diagram $\Delta$ over the group $\hhh_1$ whose
boundary path is factorized as $p_1p_2p_3p_4$, where

(1) the paths $p_2$ and $p_4$ have no $\theta$-edges;

(2) $p_1$ and $p_3\iv$ are sides of $k$-bands $\bb$ and $\bb'$
starting on $p_4$ and ending on $p_2$ (i.e. $p_1=\topp(\bb)$,
$p_3\iv = \bott(\bb')$).

(3) every maximal $k$-band in $\Delta$ contain at least one
$\theta$-cell and connects $p_2$ and $p_4$.
\medskip

\begin{center}
\unitlength=1mm \special{em:linewidth 0.5pt} \linethickness{0.5pt}
\begin{picture}(88.00,77.00)
\put(25.33,21.33){\rule{51.67\unitlength}{3.00\unitlength}}
\put(18.00,33.33){\rule{66.33\unitlength}{3.00\unitlength}}
\put(22.00,46.33){\rule{53.33\unitlength}{3.00\unitlength}}
\put(34.33,59.33){\rule{30.33\unitlength}{3.67\unitlength}}
\put(23.33,8.67){\line(1,0){57.33}}
\put(29.33,73.67){\line(1,0){41.00}}
\bezier{60}(29.33,73.67)(30.00,63.00)(34.33,63.00)
\bezier{76}(34.33,59.33)(22.00,55.00)(22.00,49.33)
\bezier{48}(22.00,46.33)(18.00,44.00)(18.00,36.33)
\bezier{56}(18.00,33.33)(18.00,26.67)(25.33,24.33)
\bezier{52}(25.33,21.33)(23.00,17.00)(23.33,8.67)
\bezier{64}(80.33,8.67)(80.33,21.00)(77.00,21.33)
\bezier{68}(77.00,24.33)(84.00,23.67)(84.33,33.33)
\bezier{72}(84.33,36.33)(83.67,46.33)(75.33,46.33)
\bezier{64}(75.33,49.33)(71.67,57.00)(64.67,59.33)
\bezier{64}(64.67,63.00)(70.33,63.33)(70.33,73.67)
\put(18.00,39.00){\vector(1,3){1.50}}
\put(45.00,73.67){\vector(1,0){12.67}}
\put(69.00,57.50){\vector(1,-1){3.00}}
\put(61.67,8.67){\vector(-1,0){11.33}}
\put(14.67,38.33){\makebox(0,0)[cc]{$p_1$}}
\put(50.33,77.00){\makebox(0,0)[cc]{$p_2$}}
\put(83.67,47.00){\makebox(0,0)[cc]{$p_3$}}
\put(54.67,5.00){\makebox(0,0)[cc]{$p_4$}}
\bezier{52}(27.00,8.67)(26.33,15.67)(28.67,21.00)
\bezier{60}(29.00,24.33)(19.67,29.67)(21.67,33.00)
\bezier{52}(22.00,36.67)(22.67,45.00)(26.67,46.33)
\bezier{80}(26.67,49.33)(24.67,52.67)(39.33,59.33)
\bezier{56}(38.67,63.00)(33.67,66.00)(34.00,73.67)
\bezier{60}(66.33,73.67)(67.00,65.67)(60.67,63.00)
\bezier{64}(60.33,59.33)(70.67,53.67)(70.67,49.33)
\bezier{72}(71.67,46.33)(83.00,42.67)(81.33,36.67)
\bezier{60}(81.67,33.33)(81.33,26.67)(73.33,24.67)
\bezier{64}(73.33,21.33)(79.33,16.00)(77.67,8.67)
\put(25.00,5.33){\makebox(0,0)[cc]{$\bb$}}
\put(79.33,5.33){\makebox(0,0)[cc]{$\bb'$}}
\put(84.67,22.67){\makebox(0,0)[cc]{$\ttt_1$}}
\put(88.00,34.67){\makebox(0,0)[cc]{$\ttt_2$}}
\put(69.67,61.00){\makebox(0,0)[cc]{$\ttt_d$}}
\put(46.67,53.67){\makebox(0,0)[cc]{$\dots$}}
\end{picture}
\end{center}

\begin{center}
\nopagebreak[4] Figure \theppp.
\end{center}
\addtocounter{ppp}{1}

By definition, $p_2$ and $p_4$ are the \label{topt}{\it top} and
the \label{bott}{\it bottom} of the trapezium, respectively. As
follows from Lemma \ref{NoAnnul}, every maximal $\theta$-band of
the trapezium connects $p_1$ and $p_3$. We shall usually enumerate
maximal $\theta$-bands starting on $p_1$ from the bottom to the
top: $\ttt_1,...,\ttt_d$. The number of them is called the
\label{heightt}{\em height} of the trapezium. The paths
$p_1,p_2,p_3,p_4$, the bands $\bb$, $\bb'$ and $\ttt_1,...,\ttt_d$
form the \label{datat}{\em data} associated with the trapezium. By
Lemma \ref{NoAnnul}, each of these $k$-bands intersects each of
the $\theta$-bands exactly once. Therefore by Lemma \ref{Band} the
bases of the $\theta$-bands $\ttt_1,...,\ttt_d$ are the same.
Hence the base of the label of the bottom $p_4$ will be called the
\label{baset}{\it base} of the trapezium (it is equal to the base
of each of $\ttt_i$). Notice that every $k$-band is a trapezium
with a 1-letter base.

A reduced annular diagram $\Delta$ is called a \label{ringgg}{\it
ring} with boundary components $p_2$ and $p_4$ if:

(1) the labels of $p_2$ and $p_4$ contain no $\theta$-edges;

(2) the boundaries $p_2$ and $p_4$ are minimal;

(3) $\Delta$ has at least one $(k,\theta)$-cell.

It follows from the definition, that a ring $\Delta$ contains at
least one $k$-band having a $(\theta,k)$-cell, every maximal
$k$-band of $\Delta$ connects $p_2$ and $p_4$, and a ring is
obtained from a trapezium by identifying two bands $\bb$ and
$\bb'$ (in this case $\bb$ must be a copy of $\bb'$ in the
definition of trapezium). The $\theta$-bands in a ring are annuli
surrounding the hole of the ring. The top and bottom paths of the
trapezium turn into the {\em outer} and {\em inner} paths of the
ring. The base and the height of a ring are defined as for
trapezia, but the base is considered as a cyclic word.

If the bands $\bb$ and $\bb'$ from the definition of trapezia are
allowed to be $a$-bands, $a\notin\aaa(P_1)$, one obtains the
definition of \label{quasitrapezium}{\it quasitrapezium}.

A \label{quasiringq}{\it quasiring} is either a ring or a reduced
annular diagram with contours $p$ and $q$ such that (1) $\Delta$
contains no $(\theta,k)$-cells but contains a $(\theta,a)$-cell
with $a\ne \aaa(P_1)$; (2) the labels of the contour $p$ and $q$
have neither $\theta$ nor $k$-letters; (3) the boundaries $p$ and
$q$ are minimal.

 The definitions of the base and the history of a quasitrapezium
(quasiring) are similar to those for trapezia (ring), but the base
of a quasitrapezium can be empty. Of course if the base of a
quasiring is not empty then it is a ring.

\begin{lemma} \label{ring123} For any reduced
annular diagram $\Delta$ over $\hhh_1$ whose contours $p$ and $q$
do not contain $\theta$-edges, there exists a ring or a quasiring,
or a diagram with minimal boundaries and having no cells,
$\Delta'$ whose contours $p'$ and $q'$ satisfy properties
(R0)-(R4) of Lemma \ref{rollsprep}, in which the lengths $p'$ and
$q'$ are recursively bounded in terms $|p|+|q|$, and there exists
a conjugacy diagram for $\phi(p)$ (for $\phi(q)$) and $\phi(p')$
(and $\phi(q')$) with recursively bounded number and perimeters of
cells.
\end{lemma}

\proof By Lemma \ref{rollsprep2} one can assume that $p$ is a
minimal boundary. We may also assume that $q$ enjoys the same
property.

If we have a $(\theta,k)$-cells in $\Delta$, then the diagram is a
roll. So we assume that there are no $(\theta,k)$-cells in it.

If there is a $k$-edge in $\partial\Delta$, then a maximal
$k$-band, containing no $\theta$-cells, connects $p$ and $q$.
Therefore $\Delta$ has no $\theta$-edges by Lemma \ref{NoAnnul}.
Hence it is a diagram over the free product of the group $\hhh_2$
and $\gr$, and the statement follows from Lemma \ref{Conjugacy}.

A reduced diagram cannot possess a $\theta$-annulus consisting of
$(a(P_1),\theta)$-cells (with equal boundary labels). Hence, if
$p$ and $q$ have neither $\theta$- nor $k$-edges, then $\Delta$ is
a $quasiring$ as desired, or it has no $\theta$-cells. If $\Delta$
has no $\theta$-cells, then we can repeat the argument of the
previous paragraph.
\endproof

 Let $\ttt_i$, $i=1,...,d$ be a
$\Theta(\tau_i)$-band where $\tau_i\in\csss$. Then the word
$h=\tau_1\tau_2...\tau_d$ is called the \label{historyt}{\em
history} of the trapezium (ring). Notice that the length of the
history is equal to the height of the trapezium (ring). By Lemma
\ref{NoAnnul}, the history of any $k$-band in a trapezium (ring)
is equal to the history of the trapezium (ring).


\begin{lemma} \label{No a-cells}
(1) Let $y_1\dots y_s$ be the base of a trapezium $\Delta$ whose
top and bottom labels have no letters from $\aaa(P_1)$. Assume the
words $y_1\dots y_{s-1}$, $(y_2\dots y_s)^{-1}$ contain no
occurrence of the positive letter $P_1$, and the words $y_2\dots
y_s$, $(y_1\dots y_{s-1})^{-1}$ contain no occurrence of the
positive letter $R_1$. Then $\Delta$ has no $\gr$-cells.

(2) If the base of a quasiring does not contain $P_1^{\pm 1}$ and
$R_1^{\pm 1}$ and the inner and outer paths of the ring do not
contain $\aaa(P_1)$-edges then the ring does not have $\gr$-cells.
\end{lemma}

\proof (1) Indeed, by Lemma \ref{aa-cells}, we can assume that
every $\aaa(P_1)$-band starting on a $\gr$-cell must end either on
the boundary of the diagram or on the boundary of a $P_1$-cell or
on the boundary of a $R_1$-cell. The first case would imply that
the bottom or the top of the trapezium contains $\aaa(P_1)$-edges.
The second case would imply that one of the letters
$y_1,...,y_{s-1}, y_2\iv,...,y_s\iv$ is $P_1$. The third case
would imply that one of the letters $y_2,...,y_s,
y_1\iv,...,y_{s-1}\iv$ is $R_1$.

Statement (2) is proved similarly.
\endproof

\subsection{Trapezia with 2-letter bases}
\label{trapeziawith}

Let $\Delta$ be a trapezium with a 2-letter base $yz$. We use the
notations from Section \ref{bandsand} and from the definition of a
trapezium for the data of $\Delta$. Let $V_0\equiv
\phi(p_4^{-1})$, $U_{d+1}\equiv \phi(p_2)$, and for all
$i=1,...,d$ let $U_i\equiv\phi(\bott(\ttt_i))$,
$V_i\equiv\phi(\topp(\ttt_i))$.

For every word $W$, we denote the projection
$W_{\aaa\cup\bar\aaa}$ of $W$ onto the alphabet $\aaa\cup\bar\aaa$
by $W_a$. Similarly we denote $W_{\Theta\cup\bar\Theta}$ by
$W_\theta$.

\begin{lemma}\label{V2V3}
(1) Suppose that $\Delta$ does not have $\gr$-cells. Then $(U_1)_a
=(V_0)_a$, $(V_d)_a = (U_{d+1})_a$ in the free group, and
$(V_i)_a\equiv (U_{i+1})_a$, $i=1,...,d-1$.

(2) If $\Delta$ contains $\gr$-cells then
$(U_1)_a=(V_0)_a$,$(V_d)_a = (U_{d+1})_a$, $(V_i)_a =
(U_{i+1})_a$, $i=1,...,d-1$, modulo the $\gr$-relations.
\end{lemma}

\proof Let $\Gamma_0$ be the subdiagram of $\Delta$ bounded by
$p_1, p_3, p_4, \bott(\ttt_1)$, $\Gamma_d$ be the subdiagram
bounded by $p_1,p_3, p_2, \topp(\ttt_d)$, and for every
$i=1,...,d-1$ let $\Gamma_i$ be the subdiagram bounded by
$p_1,p_3$, $\topp(\ttt_i)$ and $\bott(\ttt_{i+1})$.

(1) Suppose that $\Delta$ does not contain $\gr$-cells. A maximal
$a$-band in $\Gamma_i$ which starts on a side of a $\theta$-band
cannot end up on the same side because otherwise the top/bottom
labels of the $\theta$-band would not be reduced. This implies
that $(V_1)_a\equiv (U_2)_a$, ..., $(V_{d-1})_a\equiv (U_d)_a$.

Also every $a$-band which starts on $\bott(\ttt_1)$ ends on $p_4$,
every $a$-band which starts on $\topp(\ttt_d)$ ends on $p_2$. If
an $a$-band starts and ends on $p_4$ (resp. $p_2$) then the label
of the $a$-path between the start and end edges of that band must
be freely equal to 1. Hence $(V_0)_a=(U_1)_a$ and
$(V_d)_a=(U_{d+1})_a$ in the free group.

(2) Cells in $\Gamma_i$ do not contain $\theta$-edges. Hence they
correspond to either relations of the form (\ref{auxax}) or
$\gr$-relations. If we remove all $\xxx$-letters in relations
(\ref{auxax}), the relations become trivial. Thus if we collapse
all $\xxx$-edges in $\Gamma_i$, $i=0,...,d$, that diagram becomes
a diagram with $\gr$-cells only. (More precisely, to the boundary
label of $\Gamma_i$, we apply the homomorphism which preserves
letters $a_i(P_1)$, $1\le i\le m$, and sends other $a$-letters and
$\xxx$-letters to 1.) Since the portions of $p_1$ and $p_3$ on the
contour of $\Gamma_i$ consist of $\xxx$-edges, these portions
collapse to vertices. Now \vk lemma implies that
$(U_1)_a=(V_0)_a$,$(V_d)_a = (U_{d+1})_a$, $(V_i)_a =
(U_{i+1})_a$, $i=1,...,d-1$, modulo the $\gr$-relations, as
required.
\endproof

The following lemma immediately follows from the fact that there
are no $(P_j,\xxx)$- and $(R_j,\xxx)$-relations of the form
(\ref{auxkx}).

\begin{lemma} \label{PR} One of the sides of every $P_j$-band ($R_j$-band) in a \vk diagram
over $\hhh_1$ contains $a$- and $\theta$-edges only.
\end{lemma}

\begin{lemma} \label{History}
(1) If the base of a trapezium  $\Delta$ contains a $P_j$- or a
$R_j$-letter, then the history is a reduced word.

(2) The history of an arbitrary ring is a reduced word.
\end{lemma}

\proof (1) Let $\tau\tau\iv$ be a 2-letter subword of the history
of a trapezium $\Delta$. Consider the corresponding pair $\ttt_i$,
$\ttt_{i+1}$  of consecutive $\theta$-bands in $\Delta$. The
intersections of these bands with a $P_j$- or $R_j$-band $\ccc$ of
$\Delta$ are two cells $\pi$, $\pi'$ which are consecutive cells
in $\ccc$ by Lemma \ref{PR}. Since $\pi$ and $\pi'$ correspond to
$\tau$ and $\tau^{-1}$, they form a reducible pair of cells which
contradicts the assumption that $\Delta$ is reduced.

(2) Let $\Delta$ be a ring. Without loss of generality we can
assume that the history of $\Delta$ is $\tau\tau^{-1}$, so
$\Delta$ is of height 2. By (1), we can assume that the base of
$\Delta$ has neither $P_j$- nor $R_j$-letters. We can also suppose
that $\Delta$ is constructed of two $\theta$-bands $\ttt_1$,
$\ttt_2$ and a of the subdiagram $\Gamma$ bounded by
$\topp(\ttt_1)$ and $\bott(\ttt_2)$, which have no $\theta$-cells.


By Lemma \ref{No a-cells}, applied to all subtrapezia of $\Delta$
with 2-letter bases, $\Delta$ has no $\gr$-cells. Therefore, by
Lemma \ref{V2V3} (1) we have that the words obtained from the
labels of $\topp(\ttt_1)$ and $\bott(\ttt_2)$ after deletion of
all $\xxx$-letters, are freely equal. The same is true for the
labels of $\bott(\ttt_1)$ and $\topp(\ttt_2)$ by Lemma \ref{VV'},
since bands $\ttt_1$ and $\ttt_2$ correspond to mutually inverse
rules. Hence $\phi(\bott(\ttt_1))\equiv\phi(\topp(\ttt_2))$ by
Lemma \ref{Band}. But this contradicts the assumption that
$\Delta$ is reduced (see part (iv) of the definition of
non-reduced diagram).

\endproof

\subsection{Trapezia simulate the work of $S$-machines}
\label{trapeziesimulate}

The following lemma shows that every trapezium simulates the work
of $\csss$.

\begin{lemma}\label{simulatecsss} Let $\Delta$ be a
trapezium of height $d$ with $\theta$-bands $\ttt_1,...,\ttt_d$,
$W$ (resp. $W'$)  be the projection of the label of the bottom
(resp. top)  path of $\Delta$ onto
$\aaa\cup\bar\aaa\cup\kkk\cup\bar\kkk$, $U_i$ (resp. $V_i$),
$i=1,...,d$, be the projection of  $\phi(\bott(\ttt_i))$ (resp.
$\phi(\topp(\ttt_i))$ on the same set. Let $h=\tau_1\cdot...\cdot
\tau_d$ be the history of $\Delta$. Then  $U_i, V_i$, $i=1,...,d$,
are admissible words for $\csss$. In addition

\begin{equation}\label{scsss}
W=U_1 (\mod \gr),  V_1=U_1\circ\tau_1, U_2=V_1 (\mod \gr),...,
V_d=U_d\circ\tau_d, W'=V_d (\mod \gr),
\end{equation}

\noindent (here \mod $\gr$ means that the equality is true modulo
$\gr$-relations, other equalities are true in the free group). If
$\Delta$ does not contain $\gr$-cells then one can remove (\mod
$\gr$) from the previous statements. We also have
$||U_i|-|V_i||\le Bc$ where $B$ is the length of the base of
$\Delta$ and $c$ is the maximum of lengths of words in $\bar{\cal
E}$.

\end{lemma}

\proof Suppose first that $\Delta$ contains just one $\theta$-band
$\ttt$. Let $U$ and $V$ be the projections of $\phi(\bott(\ttt))$
and $\phi(\topp(\ttt))$ respectively on
$\aaa\cup\bar\aaa\cup\kkk\cup\bar\kkk$. Lemmas \ref{Band} and
\ref{VV'} imply that $U$ and $V$ are admissible words for $\csss$
and $V=U\circ \tau$. By Lemma \ref{V2V3} $W = U$ and $W'=V$ in the
free group if $\Delta$ does not have $\gr$-cells, otherwise these
equalities hold modulo $\gr$-relations. The fact that
$||V_i|-|U_i||$ follows from Lemma \ref{VV'}. Now the proof can be
finished by a simple induction on the height of $\Delta$.
\endproof

We shall call (\ref{scsss}) the \label{compt}{\em computation
associated with trapezium} $\Delta$. The {\em length} of the
computation is the height of the trapezium.

\begin{lemma} \label{x-power}
(1) Let $x_1, y_1,x_2,y_2$ be letters in $\xxx$, $x_1\ne y_1^{\pm
1}$, $x_2\ne y_2^{\pm 1}$, and let $U_1, U_2$ be non-empty reduced
words which are products of fourth powers of letters from $\xxx$.
Then $U_1x_1y_1\iv U_2y_2x_2\iv\ne 1$ in the free group.

(2) Let $x_1, x_2\in \xxx^{\pm 1}$, $U$ and $x_1Ux_2$ are products
of fourth powers of letters from $\xxx$ then $x_2\equiv x_1\iv$
and $U$ is a power of $x_1$.

(3) Let $x_1, x_2\in \xxx^{\pm 1}$, $U$ is a non-empty power of a
letter $x\in \xxx$, $x_1Ux_2$ is a product of fourth powers of
letters of $\xxx$. Then $x_1,x_2\in \{x, x\iv\}$.
\end{lemma}

\proof (1) Suppose $U_1x_2y_2\iv U_2y_1x_1\iv=1$ in the free
group. Then $U_1x_2y_2\iv=x_1y_1\iv U_2\iv $. Since $U_1$ is a
nonempty product of fourth powers of letters, the reduced form of
$U_1x_2y_2\iv$ ends with $x_2^{\pm 1}y_2\iv$. But since $U_2$ is a
non-empty product of fourth powers of letters, the word $x_1y_1\iv
U_2$ ends with a third power of a letter, a contradiction.

The other statements of the lemma are proved similarly. \endproof

A (quasi-)trapezium or a ring $\Delta$ is of the
\label{fsmtype}{\it first (second, mixed) type} if its history is
a word in $\sss$ (resp. $\bar\sss$, $\csss$ but not $\sss$ or
$\bsss$).

\begin{lemma} \label{-+}
 Let $\Delta$ be a quasitrapezium of the first or mixed type with
history of length 2. Let $\ttt_1$, $\ttt_2$ be the two maximal
$\theta$-bands of $\Delta$ counting from the bottom up. Assume
that $\ttt_1$ is a $\Theta$-band. Let $V$ be the label of
$\topp(\ttt_1)$. Then $V_a$ cannot contain subwords $a_i(z)\iv
a_{i'}(z)$ if $z\in \{K_j, L_j\}$ and it cannot contain subwords
$a_i(z)a_{i'}(z)\iv$ if $z=R_j$.
\end{lemma}

\proof We shall consider only the case $z=L_j$ because other cases
are similar. Suppose that $V_a$ contains a subword $a_i(L_j)\iv
a_{i'}(L_j)$. Let $\bb_1$, $\bb_2$ be the two neighbor $a$-bands
starting on $\topp(\ttt_1)$ and ending on $\bott(\ttt_2)$
corresponding to this subword, i.e. $\bb_1$ is a
$a_i(L_j)\iv$-band, $\bb_2$ is a $a_{i'}(L_j)$-band.

Let $\Gamma$ be the subdiagram of $\Delta$ which is situated
between $\bb_1$ and $\bb_2$, i.e. it is bounded by $\bott(\bb_1)$
and $\topp(\bb_2)$, a portion of $\topp(\ttt_1)$, and a portion of
$\bott(\ttt_2)$. Let $\partial(\Gamma)=u_1qu_2 p\iv$ be the
decomposition of the boundary of $\Gamma$ where
$u_1=\bott(\bb_1)$, $u_2=\topp(\bb_2)\iv$.

Notice that the start (end) edges of $\bb_1$ and $\bb_2$ on
$\topp(\ttt_1)$ (resp. $\bott(\ttt_2)$) belong to two different
cells $\pi_1$, $\pi_1'$ (resp. $\pi_2,\pi_2'$) because in every
$\theta$-cells a edges labelled by $a$-letters in opposite
exponents are always separated by a $\theta$-edge: this is obvious
for relations (\ref{auxtheta}) and relations (\ref{mainrel})
corresponding to rules not from $\sss(34)\cup\bsss(34)$; for the
remaining relations it follows from the assumption that all words
in $\bee$ are positive.

Notice also that by definition, the word $\alpha_{\tau^{\pm 1}
}(a_i(z))$ is completely determined by its $x$-letter. Suppose
that $|p|\le 1$. This means that the two $x$-letters in the
relations corresponding to the cells $\pi_1, \pi_1'$ cancel. That
implies, by Lemma \ref{Band}, that the labels of the start edges
of $\bb_1, \bb_2$ must be mutually inverse which contradicts the
assumption that $\phi(\topp(\ttt_1))$ is a reduced word. Hence
$|p|=2$, $\phi(p)=x_1y_1\iv$, $x,y\in \xxx^{\pm 1}$, $x_1\ne
y_1^{\pm 1}$. Similarly if $\ttt_2$ corresponds to a rule from
$\sss$, $|q|=2$, $\phi(q)=x_2y_2\iv$, $x_2,y_2\in \xxx^{\pm 1}$,
$x_2\ne y_2^{\pm 1}$.

Let $U_1\equiv \phi(u_1)$, $U_2\equiv \phi(u_2)$. Since $\Gamma$
contains no cells (by Lemma \ref{NoAnnul}), the equality
$U_1\phi(q)U_2\phi(p)\iv=1$ must be true in the free group.

\medskip

{\bf Case 1.} Suppose that $\ttt_2$ corresponds to a rule from
$\bsss$. Then  $U_1U_2y_1x_1\iv=1$ in the free group where $U_1,
U_2$ are reduced products of fourth powers of letters from
$\xxx^{\pm 1}$. Considering the homomorphism of the free group
onto $\mathbb{Z}$ which kills all letters except $x_1$, we
immediately get a contradiction.

\medskip

{\bf Case 2.} Suppose now that $\ttt_2$ corresponds to a rule from
$\sss$. Then the equality $U_1x_2y_2\iv U_2y_1x_1\iv=1$ is true in
the free group. Taking the projections on $\la x_1\ra$, $\la
x_2\ra$ as in Case 1, we deduce that $x_1\equiv x_2, y_1\equiv
y_2$.

\medskip

{\bf Case 2.1.} Suppose that $U_1, U_2$ are not empty. Then Lemma
\ref{x-power} (1) immediately gives a contradiction.

\medskip

{\bf Case 2.2.} Suppose that one of the words $U_1$ or $U_2$ is
empty. Then the other word must be empty as well (these words are
freely conjugate and reduced). Notice that the paths $u_1$, $u_2$
cannot contain consecutive edges with mutually inverse edges
because otherwise the cells containing these edges cancel (they
correspond to relations of the form (\ref{auxax}) with the same
$a$-letter). Hence the paths $u_1$ and $u_2$ are empty. Therefore
the cells $\pi_1$ and $\pi_2$ have a common $x$-edge and a common
$a$-edge. Since these two cells do not cancel ($\Delta$ is
reduced), one of them corresponds to a relation of the form
\ref{mainrel} and the other one corresponds to a relation of the
form \ref{auxtheta}. (If these cells are both $k$-cells, then it
is clear from the structure of the $a$-bands between $\pi_1$ and
$\pi_2$, that $\bb_1$ determines equal occurrences of $a$-letters
in the boundary labels of $\pi_1$ and $\pi_2$.) Without loss of
generality we can assume that $\pi_1$ corresponds to a relation of
the form (\ref{mainrel}) and $\pi_2$ corresponds to a relation of
the form (\ref{auxtheta}). The $\Theta$-edges on these cells have
the same labels (up to the direction) because the label of the
$\Theta$-edge is encoded in the $x$-edges of a relation
(\ref{mainrel}) or (\ref{auxtheta}). Hence $\ttt_1$ and $\ttt_2$
correspond to mutually inverse rules of $\sss$.

\begin{center}
\unitlength=1mm \special{em:linewidth 0.5pt} \linethickness{0.5pt}
\begin{picture}(109.33,54.67)
\put(105.00,54.67){\line(-1,0){100.00}}
\put(5.00,54.67){\line(0,-1){6.67}}
\put(5.00,48.00){\line(1,0){99.67}}
\put(86.33,48.00){\line(0,1){6.67}}
\put(51.00,48.00){\line(0,1){6.67}}
\put(105.00,11.00){\line(-1,0){100.00}}
\put(5.00,11.00){\line(0,-1){6.67}}
\put(5.00,4.33){\line(1,0){99.67}}
\put(86.33,4.33){\line(0,1){6.67}}
\put(51.00,4.33){\line(0,1){6.67}}
\put(38.00,50.33){\makebox(0,0)[cc]{$\pi_1'$}}
\put(64.67,50.67){\makebox(0,0)[cc]{$\pi_2'$}}
\put(38.00,6.33){\makebox(0,0)[cc]{$\pi_1$}}
\put(64.67,6.67){\makebox(0,0)[cc]{$\pi_2$}}
\put(109.33,51.33){\makebox(0,0)[cc]{$\ttt_2$}}
\put(109.33,7.33){\makebox(0,0)[cc]{$\ttt_1$}}
\put(39.67,48.00){\vector(1,0){10.67}}
\put(62.33,48.00){\vector(-1,0){10.67}}
\put(39.67,11.00){\vector(1,0){10.67}}
\put(62.33,11.00){\vector(-1,0){10.67}}
\put(45.67,50.00){\makebox(0,0)[cc]{$x_2$}}
\put(56.33,50.00){\makebox(0,0)[cc]{$y_2$}}
\put(56.33,8.33){\makebox(0,0)[cc]{$y_1$}}
\put(45.67,8.33){\makebox(0,0)[cc]{$x_1$}}
\bezier{36}(42.00,48.00)(44.00,45.00)(49.67,45.00)
\put(49.67,44.67){\line(0,-1){30.00}}
\bezier{40}(49.67,15.00)(43.33,14.33)(42.00,11.00)
\bezier{36}(52.00,45.00)(58.67,45.33)(59.00,48.00)
\put(52.00,45.00){\line(0,-1){30.00}}
\bezier{40}(52.00,15.33)(59.00,14.33)(59.33,11.00)
\put(49.67,23.67){\vector(0,1){11.67}}
\put(52.00,35.33){\vector(0,-1){12.33}}
\put(46.00,29.00){\makebox(0,0)[cc]{$U_1$}}
\put(54.67,29.00){\makebox(0,0)[cc]{$U_2$}}
\put(7.00,48.00){\line(0,-1){37.00}}
\put(12.33,11.00){\line(0,1){37.00}}
\put(23.00,11.00){\line(0,1){37.00}}
\put(34.33,48.00){\line(0,-1){37.00}}
\put(68.67,48.00){\line(0,-1){37.00}}
\put(28.33,27.33){\makebox(0,0)[cc]{$\bb_0$}}
\put(39.67,37.00){\makebox(0,0)[cc]{$\bb_1$}}
\put(61.00,37.00){\makebox(0,0)[cc]{$\bb_2$}}
\put(9.00,22.67){\makebox(0,0)[cc]{$\bb$}}
\end{picture}

\end{center}

\begin{center}
\nopagebreak[4] Figure \theppp.
\end{center}
\addtocounter{ppp}{1}

Let $\bb$ be the $k$-band starting on the $k$-edge of the cell
$\pi_1$ (that belongs to $\topp(\ttt_1)$) and ends on a $k$-edge
of $\bott(\ttt_2)$. Assume that there is an $a$-band $\bb_0$
connecting $\ttt_1$ and $\ttt_2$ between $\bb$ and $\bb_1$. Assume
that it is the closest one to $\bb_1$. Notice that it starts with
an edge of $\topp(\ttt_1)$ labelled by a negative $a$-letter,
since the all $a$-labels of $a$-edges of $\pi_1$, situated between
a $k$-edge and a $\theta$-edge, are positive or negative
simultaneously.

Consider the subdiagram $\Gamma^0$ between $\bb_0$ and $\bb_1$.
Its contour has decomposition $u_1^0q^0u_2^0(p^0)^{-1}$simlar to
the decomposition $u_1qu_2p^{-1}$ of the contour of $\Gamma$. Now
we obtain the equality $U_1^0xU_2^0x^{-1}=1$ in the free group for
the boundary label of $\Gamma^0$. Here $U_2^0$ must be empty since
the word $U_1$ was empty. Hence $U_1^0$ is empty being the product
of forth powers of $x$-letters. There are fewer maximal $a$-bands
between $\bb$ and $\bb_0$ than between $\bb$ and $\bb_1$.
Therefore, arguing in this way, we finally conclude that
$\phi(\bott(\bb))$ is empty too. But then $\pi_1$ and $\pi_2$ must
have a common $k$-edge, and they cancel because $\ttt_1$ and
$\ttt_2$ correspond to mutual inverse rules of $\sss$. This
contradicts the fact that $\Delta$ is reducible.
\endproof

If the base of a trapezium $\Delta$ is a 2-letter word, and its
height is equal to 2, then $\Delta$ is called a \label{smallt}{\it
small} trapezium.

\begin{lemma} \label{Sm Tr}
Let $\Delta$ be a small trapezium of the first or mixed type with
reduced history. Then the projection of the words
$V_1=\phi(\topp(\ttt_1))$, $V_2=\phi(\bott(\ttt_2))$ onto
$\aaa\cup\bar\aaa\cup\kkk$ are admissible words for $\sss$.
\end{lemma}

\proof By Lemma \ref{simulatecsss}, we need only check the
positivity conditions of the definition of admissible words. Also
it is clear that it is enough to prove the statement for $V_1$.
Indeed, by Lemma \ref{NoAnnul} every $a$- and $k$-band starting
on $\topp(\ttt_1)$ ends on $\bott(\ttt_2)$ and vice versa; so the
projections of $V_1$ and $V_2$ onto $\aaa\cup\bar\aaa\cup\kkk$
are the same in the free group (in particular they do not contain
$\bar\aaa$-letters).

Hence it is enough to consider the cases when $\ttt_1$
corresponds to a rule from $\sss$, $V_1$ is a $zL_j$-, $zP_j$- or
$R_jz$-sector. All three cases are similar, so we consider only
one of them when $V_1$ is a $zL_j$-sector. We need to show that
the projection of $V_1$ onto $\aaa$ is positive.

Suppose that $V_1$ contains a negative letter $a\iv$, $a\in
\aaa(\lel_j)$. First let $a\iv $ be the last $\aaa$-letter of
$V_1$. Since $V_1$ and $V_2$ are reduced, every $a$-band starting
on $\topp(\ttt_1)$ ends on $\bott(\ttt_2)$ and vise versa. In
particular the $a$-band $\CC$ starting on the last $\aaa$-edge of
$\topp(\ttt_1)$ ends on the last $\aaa$-edge of $\bott(\ttt_2)$,
hence the last $\aaa$-letter of $V_2$ is also $a\iv$.

Relations (\ref{mainrel}), (\ref{auxtheta}) show that then $V_1$
ends with $a\iv x$ where $x=x(b,\tau)^{\pm 1}$, and $V_2$ ends
with $a\iv x'$ where $x'=x(c,\tau')^{\pm 1}$, for some $b,c\in
\aaa$. Let $\Gamma$ be the diagram bounded by $\topp(\CC)$,
$\bott(\bb')$, a portion of $\topp(\ttt_1)$ and a portion of
$\bott(\ttt_2)$. By Lemma \ref{NoAnnul} and \ref{No a-cells},
$\Gamma$ contains no cells. Both words $\phi(\topp(\CC))$ and
$\phi(\bott(\bb'))$ are products of fourth powers of letters from
$\xxx$, and $x$ and $x'$ are not mutually inverse because
$\tau'\ne\tau\iv$ by an assumption of the lemma. This immediately
leads to a contradiction by Lemma \ref{x-power} (2).

Thus the last $a$-letter of $V_1$ (resp. $V_2$) is positive. Since
$V_1$ contains a negative $a$-letter $a\iv$, $(V_1)_a$ must
contain a subword of the form $a\iv b$ for some $b\in\aaa$. But
this contradicts Lemma \ref{-+}.
\endproof

The next lemma shows how to get rid of $\gr$-cells in some
trapezia.

\begin{lemma} \label{refined}
Let  $\Delta$ be a quasitrapezium of the first type. Suppose that
the history of $\Delta$ does not contain
$\tau\in\sss(34)\cup\sss(4)\cup\sss(45)$ or the base of $\Delta$
is empty. Then there exists a quasitrapezium $\Delta'$ of the
first type such that

\begin{itemize}

\item $\Delta'$ has the same labels of the bottom, left and right sides
as $\Delta$,

\item $\phi(\topp(\Delta))=\phi(\topp(\Delta'))$ modulo
$\gr$-relations

\item $\Delta'$ has no $\gr$-cells and the type of $\Delta'$ is
not higher than that of $\Delta$.
\end{itemize}
\end{lemma}

\proof We only consider the case when $\Delta$ is a trapezium with
non-empty base. The other cases are similar.

If $\Delta$ does not contain $\gr$-cells then we can take
$\Delta'=\Delta$, so suppose that $\Delta$ contain $\gr$-cells.

Suppose there is a small subtrapezium $\Gamma$ in $\Delta$ with
base $yz$ where $yz$ is one of the words $(P_1R_1)^{\pm 1}$,
$P_1P_1^{-1}$, $R_1^{-1}R_1$. Assume that $\Gamma$ is crossed by
maximal $\theta$-bands $\ttt_t$ and $\ttt_{t+1}$. Denote
$q_1=\bott(\ttt_t)$,...,$q_4=\topp(\ttt_{t+1})$. Let
$V_i=\phi(q_i)$, $i=1,...,4$. Notice that $V_2,V_3, V_4$ are words
over  $\aaa(P_1)$ by Lemmas \ref{Band}, \ref{VV'} and $V_2=V_3$
modulo $\gr$-relations by Lemma \ref{V2V3}.

Denote by $\Gamma_0$ the subdiagram of $\Gamma$ with the boundary
label $V\equiv V_2V_3^{-1}$. Assume $\Gamma_0$ contains at least
one $\gr$-cell. Then by Lemma \ref{Emb} it has exactly one
$\gr$-cell and no other cells.

Let $\tau$ be the history of $\ttt_{t+1}$. We cut $\ttt_{t+1}$
along a $\theta(\tau, P_1)$-edge and the boundary of $\Gamma_0$ to
obtain a closed subpath $p$ labelled by the word
 $\theta(\tau,P_1)^{-1}V\theta(\tau,P_1)$. It is freely equal to
$\theta(\tau,P_1)^{-1}V\theta(\tau,P_1)V^{-1}V$. Removing the
subdiagram bounded by this loop, we get a hole in the diagram
$\Delta$.

Recall that $\theta(\tau,P_1)$ commutes with $\aaa(P_1)$-letters
by relations (\ref{mainrel}), since $\tau$ does not belong to
$\sss(34)\cup\sss(4)\cup\sss(45)$ by the assumptions of the lemma.

Hence there exists a $\theta$-band with boundary label
$\theta(\tau,P_1)^{-1}V\theta(\tau,P_1)V^{-1}$. If we connect a
$\gr$-cell by a vertex to the end of the first $\theta$-cell of
this band, we get a diagram $\Gamma_1$ with boundary label freely
equal to $\phi(p)$. Now we can fill the hole bounded by $p$ by the
diagram $\Gamma_1$. As a result the $\theta$-bands of $\Delta$,
except $\ttt_{t+1}$ do not change, the band $\ttt_{t+1}$ gets
longer, and the $\gr$-cell moves outside $\Gamma$ toward the top
of $\Delta$. Let $\tilde\Delta$ be the resulting diagram, and
$\Delta_1$ be the diagram obtained by reducing $\tilde\Delta$.

Two $k$-cells of the same $k$- or $\theta$-band in $\tilde\Delta$
cannot form a reducible pair of cells (see (ii) in the definition
of non-reduced diagram) because otherwise diagram $\Delta$ were
not reduced. Hence $k$-cells are not removed in the process of
reducing $\tilde\Delta$.

This implies that $\Delta_1$ satisfies the definition of a
trapezium. The height of $\Delta_1$ is the same as the height of
$\Delta$, the diagram between $\ttt_1$ and $\ttt_{t+1}$ in
$\Delta_1$ contains fewer $\gr$-cells than the corresponding
subdiagram of $\Delta$. Hence after a number of such
transformations, we get a trapezium $\Delta_s$ with the same
boundary label as $\Delta$, and all $\gr$-cells in $\Delta_s$ are
between $\topp(\ttt_d)$ and the top path $p_2$. Notice that
$\Delta_s$ has the same number of $k$- and $\gr$-cells as $\Delta$
but possibly bigger number of auxiliary $(\theta,a)$-cells.

For every $\gr$-cell $\pi$ in $\Delta_s$, there exists a path in
$\Delta_s$ which does not cross any $\theta$-bands and connects
$\pi$ with $p_2$. Hence all $\gr$-cells in $\Delta_s$  can be cut
off from $\Delta_s$ without changing the bottom, left and right
sides of $\Delta_s$. The label of the top side does not change
modulo $\gr$-relations. After we remove all $\gr$-cells from
$\Delta_s$, we get a trapezium $\Delta'$ satisfying the first two
of the required conditions. Since $\Delta$ contains $\gr$-cells,
the type of $\Delta'$ is smaller than the type of $\Delta$: to
obtain $\Delta'$ from $\Delta$, we add auxiliary $\theta$-cells
and remove $\gr$-cells.\endproof

\begin{lemma}\label{quasiring} For every quasiring $\Delta$ with empty base and contours $p$ and
$q$ there exists a quasiring with the same boundary labels and
recursively bounded number of cells.
\end{lemma}

\proof
  By Lemma \ref{rollsprep2} we can assume that
$\Delta$ satisfies condition (R1).

 By the definition, there is an $a(z)$-band, connecting $p$
and $q$, in $\Delta$. Since the base of $\Delta$ is empty, every
$\theta$-cell of $\Delta$ must be an $\aaa(z)$-cell, and $z\ne
P_1$ by the definition of a quasiring.

Denote by $\Delta'$ the minimal annular subdiagram of $\Delta$
which contains all the $\theta$-annuli of $\Delta$. It contains no
$\gr$-cells by Lemma \ref{No a-cells}(2). Hence every $a$-band
starting on $p$ ends on $q$ and vice versa. Therefore the length
of arbitrary $\theta$-annulus is linearly bounded in terms of
$|p|, |q|$, and we have the bounded number of possibilities for
their labels. Since $\Delta$ is reduced, we conclude that the
number of its $\theta$-annuli is also recursively bounded. The
statement is obtained after the application of Lemma
\ref{Conjugacy} to the remaining annular subdiagrams of $\Delta$
which are diagrams over $\hhh_1$.

\endproof

The following lemma shows that trapezia of the first type
simulates the work of $\sss$. Notice that this lemma does not
follow from Lemma \ref{simulatecsss} because admissible words for
$\sss$ differ from admissible word for $\csss$.

\begin{lemma} \label{simulatesss}
Let $\Delta$ be a trapezium of the first type with base $B$ and
the reduced history $h=\tau_1\cdot...\cdot \tau_d$, $\tau_i\in
\sss$, $d\ge 2$, with $\theta$-bands $\ttt_1,...,\ttt_d$. Let $W$
(resp. $W'$) be the projection of the label of the bottom (resp.
top) path of $\Delta$ onto $\aaa\cup\bar\aaa\cup\kkk\cup\bar\kkk$,
$U_i$ (resp. $V_i$), $i=1,...,d$, be the projection of
$\phi(\bott(\ttt_i))$ (resp. $\phi(\topp(\ttt_i))$ on the same
set. Then $U_2,...,U_d, V_1,...,V_{d-1}$ are admissible words for
$\sss$. In addition

$$\begin{array}{l} W=U_1 (\mod \gr),  ||V_1|-|U_1||\le c|B|,
U_2=V_1 (\mod \gr), V_2=U_2\circ\tau_2, \\ ||V_2|-|U_2||\le c|B|,
..., V_{d-1}=U_{d-1}\circ\tau_{d-1}, ||V_{d-1}|-|U_{d-1}|\le
c|B|,\\ U_d = V_{d-1} (\mod \gr), ||V_d|-|U_d||\le c|B|,  W'=V_d
(\mod \gr)
\end{array}$$ where $c$ is the maximum of lengths of words in
$\bee$. If $\Delta$ does not contain $\gr$-cells then one can
remove ``(\mod $\gr$)" from the previous statement (i.e. all
equalities will be true in the free group).
\end{lemma}

\proof The proof is similar to the proof of Lemma
\ref{simulatecsss}. The fact that $U_2,...,U_d, V_1,...,V_{d-1}$
are admissible words for $\sss$ follows from Lemma \ref{Sm Tr}.
\endproof

The definition of admissible words for $\bar\sss$ does not have
the positivity conditions. This makes analysis of computations of
$\bsss$ different from the analysis of computations of $\sss$. On
the one hand the analysis of computations of $\bsss$ is more
difficult because there are more possible bases of admissible
words for $\bsss$. On the other hand the following lemma which is
an immediate corollary of Lemma \ref{W1}, part 2, allows us to
replace any computation of the form $U\bullet h = V$ of $\bsss$ by
any other computation with the same initial and terminal words. In
the case of computations over $\sss$ the possibilities of
modifying the computation were very limited: we could only replace
$h$ by a reduced form of $h$ or we could remove certain subwords
of $h$ (which fix certain admissible words).

We call a computation $U\bullet h=V$ \label{reducedc}{\em reduced}
if the word $h$ is freely reduced.

\begin{lemma} \label{Delta(1)}
For every reduced computation $U\bullet h=V$ of $\bsss$ there
exists a trapezium with top label $V$, bottom label $U$, and
history $h$. If $U=U'z$, $V=V'z_1$ where $z$ (resp. $z_1$)
coincides with the first letter of $U'$ (resp. $V'$) then for
every computation $U\bullet h=V$ of $\bsss$ there exists a ring of
type 2 with boundary labels $U'$ and $V'$.
\end{lemma}

\proof Assume that $|h|=1$, and $h=\tau\in \sss^+$. We consider a
path $p$ on the plane labelled by $U$. To every its $a$-edge we
attach a cell corresponding to an appropriate $\bar\Theta$-analog
of relation (\ref{auxtheta}) which in turn corresponds to the rule
$\tau$. Similarly, to every $k$-edge of $p$ we attach a cell
corresponding to (\ref{mainrel}). Then the $\bar\Theta$-edges of
the neighbor cells should be identified. We get a
$\bar\Theta$-band $\ttt$. (Of course we identify some other
labels, if necessarily, to obtain a reduce label of the
$\topp(\ttt)$.) The form of these relations shows that we obtain a
band $\ttt$ with top label graphically equal to $V$.

Two $k$-cells of $\ttt$ cannot form reducible pair of $k$-cells
since the admissible words $U$ and $V$ have no subwords $z^{\pm
1}z^{\mp 1}$ for a $k$-letter $z$. Hence $\ttt$ is a reduced
diagram.

Similarly, one can construct $\ttt$ if $\tau\in \sss^-$.

If $|h|=d>1$, then we should construct the a diagram $\Delta$
gluing together $d$ bands $\ttt_1,\dots,\ttt_n$, corresponding to
the letters of the history $h$. Two cells from distinct bands
cannot form a reducible pair since the history $h$ is reduced.
Hence $\Delta$ is a reduced diagram, and it is the desired
trapezium.

To prove the second statement, we notice that in this case the
resulting trapezium is bounded by two $k$-bands which are copies
of each other. Hence we can identify them to get the desired ring
(or
 if the resulting annular diagram is compressible, then there is a
 ring of smaller type) with the same boundary labels. However we should verify
 that the words $U'$, $V'$ are not conjugates in $\hhh_1$ of words having fewer
 occurrences of $k$-letters . Assume the is a reduced annular diagram for such a
 conjugation with contours $p$ and $q$ where $\Lab(p)\equiv U$.
Then there must be a $k$-band $\bb$ which starts and ends on $p$
such that there are no $k$-bands in the simply connected
subdiagram $\Gamma$ bounded by $\bb$ and $p$. The subdiagram
$\Gamma$ contains no $\theta$-edges by Lemma \ref{NoAnnul} since
$\Lab(p)$ has no $\theta$-letters. Therefore the word $u$ written
between corresponding mutual inverse $k$-letters of $U$ is equal
to a word having no $a$-letters. Hence if $|u|>0$, then must be an
$a$-band in $\Gamma$, which starts and ends on $p$, and gives a
cancellation in $u$. If $|u|=0$ then the $k$-letters cancel. In
both cases we get a contradiction to the fact that the admissible
word $U$ must be reduced.
\endproof

\section{Computations of $\csss$}
\label{Rings}

Notice that if we cut a ring along the top side of a $k$-band
$\bb$, we get a rectangular \vk diagram. If we then attach a copy
of $\bb$ along its top side to one of the sides of the rectangle,
we get a trapezium $\Delta$ with boundary $p_1p_2p_3p_4$ which
satisfies one additional property that $p_1$, $p_3\iv$ are sides
of two copies of the same $k$-band (the band $\bb$). The top and
bottom sides of $\Delta$ are the outer and inner boundary
components of the initial ring.

As we know (see Lemma \ref{simulatecsss}), with every trapezium
$\Delta$ with bottom label $U$ and top label $V$, we can associate
a computation $U_{\backslash x}=W_0, W_1=W_0\circ \sigma_1,...,
W_n=W_{n-1}\circ\sigma_n=V_{\backslash x}$ of the machine $\csss$
(modulo $\gr$ if the trapezium contains $\gr$-cells). Here
$U_{\backslash x}$ and $V_{\backslash x}$ are obtained from $U$
and $V$ by deleting their $x$-letters. The words $W_i$ are the
$\aaa\cup\bar\aaa\cup\kkk\cup\bar\kkk$-projections of the labels
of sides of $\theta$-bands in $\Delta$.

Computations modulo $\gr$ will be called \label{gcomp}$\gr$-{\em
computations}. Previously considered computations will be called
\label{freec}{\em free}. The concept of admissible word $W$ is
defined naturally for $\gr$-computations: we add to the
reducibility of such a word that $W$ has  subwords neither of the
form $P_1uP_1^{-1}$ nor of the form $R_1^{-1}uR_1$ for a word $u$
in $\aaa(P_1)$-letters which is equal to 1 modulo $\gr$-relations.
Locked sectors are defined naturally for $\gr$-computations as
well (if a rule locks a $zz_+$-sector then all admissible words in
the domain of the rule must have $zz_+$-sectors with empty modulo
$\gr$-relations inner parts). Notice that $\gr$-computations can
be viewed as computations of $S$-machines obtained from $\csss$ by
replacing the free product of $\gr$ and the free group over
$\{a_{m+1}(P_1),\dots,a_{\bar m}(P_1)\}$ for the free subgroup
generated by $\aaa(P_1)$ in the hardware.

Notice that if $U^{\pm 1}$ does not contain $P_1z$- or
$zR_1$-sectors then free computations and computations modulo
$\gr$ are the same. Also any computation of $\bsss$ is free.

Given a ring $\Delta$,  by Lemma \ref{simulatesss}, we get a
$\gr$-computation $U_1\bullet \sigma_1...\sigma_d=V_d$. Suppose
that we have recursively bounded the number $d$ in terms of the
lengths of the contours of the diagram. If the base of the ring
contains $P_j$ or $R_j$ for some $j$ then we can cut the ring
along a side of the corresponding $k$-band which does not contain
$x$-edges. This side is of bounded length, the resulting diagram
will be simply connected. Now using Lemma \ref{K0}, we can
recursively bound the number of cells in a ring.

If the base does not contain $P_j$ and $R_j$ for any $j$ then, to
bound the number of cells in $\Delta$ from above, we can assume by
lemmas \ref{ring123} and \ref{No a-cells} that the ring does not
contain $\gr$-cells. Therefore the words $U_i$, $V_i$ are of
recursively bounded lengths. By Lemma \ref{Conjugacy}, the number
of cells between any two consecutive $\theta$-annuli in $\Delta$
is bounded.

In both cases the history of the the resulting ring coincides with
the history of $\Delta$ because it cannot be shorten by condition
(iv) in the definition of a reduced diagram $\Delta$. Thus we
obtain

\begin{lemma}\label{ring321}
If the length of the computation corresponding to a ring $\Delta$
is recursively bounded, then there exists a ring  with the same
boundary labels and the same history as $\Delta$ whose number of
cells is recursively bounded in terms of lengths of boundaries of
$\Delta$.
\end{lemma}

If a computation (a $\gr$-computation) corresponds to a ring we
call it a \label{ringcom}{\it free ring computation} (a {\it ring
computation}). A ring computation is a reduced computation by
Lemma \ref{History}.


For every word $h$ and every natural $t\le |h|$ we denote the
prefix of $h$ of length $t$ by $h[t]$.

\begin{lemma}\label{ringc} If $U\bullet h=V$ is a ring
computation then $U\bullet h[t]$ is also a ring computation for
every $t=1,...,|h|$).
\end{lemma}

\proof The computation $U\bullet h[t]$ corresponds to a subring of
the ring corresponding to $U\bullet h=V$ containing the first $t$
$\theta$-bands of that ring.
\endproof

By definition the \label{basec}{\em base of a computation} is the
base of the starting word $U$ (which is the same as the base of
$V$), and the \label{hiscomp}{\em history of computation}
$U\bullet h$ is $h$. This agrees with the definitions of the base
and the history of a trapezium.

\begin{rk} \label{remsim} {\rm When we cut a ring to make
a trapezium, we can choose any of the $k$-bands of the ring. Thus
any letter of the base of the ring can be chosen to be the first
letter (and the last letter) of the base of the computation. Also
by taking a mirror image of a trapezium, we get a trapezium with
inverse base. For simplicity let us always assume that if the base
of a ring computation contains a letter $K_1^{\pm 1}$, then the
base starts and ends with $K_1$.}
\end{rk}

Here is a translation of Lemma \ref{refined} into the language of
computations. We leave it as an exercise to the reader to
translate the proof of Lemma \ref{refined} as well.

\begin{lemma}\label{refined1}
Let $U\bullet h=V$ be a ring computation of $\sss$. Suppose that
$h$ does not contain letters from
$\sss(34)\cup\sss(4)\cup\sss(45)$. Then there exists a free ring
computation $U\bullet h=V'$ where $V=V'$ modulo $\gr$.
\end{lemma}

\subsection{Brief history}

\label{bh} \setcounter{equation}{0}

\begin{lemma} \label{Monoton}
Let $U\bullet h=V$ be a reduced free computation of $\sss(g)$ or
of $\bar\sss(g)$, $g=1,...,5$, with $|h|=2$, $h=\sigma_1\sigma_2$.
Suppose that $U$ is a $zz_+$-sector. Then inequalities
$|U|<|U\circ\sigma_1|>|V|$ cannot happen simultaneously.
\end{lemma}

\proof Indeed, the definition of $\csss$ shows that if
$|U|<|U\circ\sigma_1|$, $\sigma_1\in \csss$,  then, depending on
$g$, $(U\circ\sigma_1)_a=U_ab$, or $(U\circ\sigma)_a=bU_a$,
$b\in\aaa\cup\bar\aaa$, where $U_ab$ (resp. $bU_a$) is a freely
reduced word, and $\sigma_1$ is completely determined (by Lemma
\ref{csss} among the rules from $\sss(g)$ by the letter $b$. If in
addition $|V|<|U\circ\sigma_1|$, then $V_a=U_abb\iv$ (resp.
$V_a=b\iv bU_a$) since $V=(U\circ\sigma_1)\circ\sigma_2$.
Therefore $\sigma_2=\sigma_1\iv$, a contradiction with the
assumption that the computation $U\bullet h=V$ is reduced.
\endproof

Lemma \ref{Monoton} immediately implies the following

\begin{lemma} \label{height}
Let $U\bullet h=V$ be a reduced free computation of $\sss(g)$ or
of $\bar\sss(g)$, $g=1,\dots,5$, and $U$ is a $zz_+$-sector.
Suppose that rules of $\sss(g)$ are active with respect to
$zz_+$-sectors. Then the length $|h|$ is at most $|U|+|V|$.
\end{lemma}

Let $U\bullet h=V$ be a computation. A maximal nonempty subword of
$h$ consisting of rules from $\sss(\omega)$ for a given
$\omega\in\Omega=\{1,2,3,4,5\}$ will be called an \label{age}{\em
age} of $(\omega)$ in the computation. The word $h$ has a unique
decomposition $g_1h_1g_2h_2...$ where $g_i$ are transition rules
($g_1$ may be missing) and $h_i$ are ages ($i=1,2,...$). If
$g_i\in \sss(\omega_i')$, $h_i$ is an age of $\omega_i$,
$i=1,2,...$, then the word
$(\omega_1')(\omega_1)(\omega_2')(\omega_2)...$ will be called the
\label{brh}{\em brief history} of the computation $U\bullet h=V$
and is denoted by \label{brh1}$br(h)$.


\begin{lemma} \label{(12)(2)(12)}
Let $U\bullet h=V$ be a free reduced computation of $\sss$ or
$\bsss$. Suppose that $U$ is a $zz_+$-sector, and for some
$\omega, \omega'\in \Omega$, rules from $\sss(\omega)$ lock
$zz_+$-sectors and rules from $\sss(\omega')$ are active with
respect to $zz_+$-sectors. Then the brief history $\br(h)$ does
not contain subwords $(\omega)(\omega')(\omega)$.
\end{lemma}

\proof Suppose that
$\br(h)=\beta_1(\omega)(\omega')(\omega)\beta_2$ for some prefix
$\beta_1$ and some suffix $\beta_2$. Let $h=h_1h_2h_3h_4h_5$ be
the corresponding decomposition of $h$. Since rules from
$\sss(\omega)$ lock $zz_+$-sectors, we have that the words
$(U\circ h_1h_2)_a$ and $(U\circ h_1h_2h_3h_4)_a$ are empty. Since
 rules from $h_3$ are active with respect to $zz_+$-sectors, and
are from the same $\sss(\omega')$, there exists a word $u$ in
$\aaa\cup\bar\aaa$ such that $(U\circ h_1h_2\sigma_3)_a$ is freely
equal to either $(U\circ h_1h_2)_au$ or $u(U\circ h_1h_2)_a$.
Hence $u$ is freely equal to 1. Since by Lemma \ref{csss} each
letter in $u$ and $\omega'$ uniquely determines a rule in $\sss$
(or $\bsss$), we conclude that $h_3$ is not freely reduced or
empty, a contradiction.
 \endproof

\medskip

Consider a brief history of the form $h\equiv
(12)(2)(23)(3)(34)(4)(45)(5)(51)$. Any word of the form $h^{\pm
1}$ will be called a \label{histper}{\it historical period}. (We
set here $(1)^{-1} = (1), (12)^{-1}=(12),\dots$, though one can
prefer formula $(12)^{-1}=(21)$.)

\begin{lemma} \label{Br Hist}
Suppose that the base of a free reduced computation $U\bullet h=V$
of $\sss$ or $\bsss$ contains a subword $\lel_jL_jP_jR_j\rer_j$
for $j\ne 1$. Then $\br(h)$ is a subword of a word of the form
$$f_0h_1f_1\dots h_sf_s$$ where $h_1,...,h_s$ are historical
periods and each of (possibly empty) words $f_0,\dots,f_s$ are
ages of $(1)$.
\end{lemma}

\proof Indeed, notice that for every $g=2,3,4,5$ rules in
$\sss(g)$ are active with respect to $zz_+$-sectors for at least
two letters $z$ occurring in $\lel_jL_jP_jR_j\rer_j$, and each
rule from
$\sss(gg')\cup\bar\sss(gg')\cup\sss(g'g)\cup\bar\sss(g'g)$ locks
one of these two sectors. Now using Lemma \ref{(12)(2)(12)}, we
can conclude that $\br(h)$ does not contain subwords of the form
$(g'g)(g)(g'g)$ for any $g$. This immediately implies that
$\br(h)$ has the desired form.
\endproof

\begin{lemma} \label{Ring Br Hist}
Let $U\bullet h=V$ be a ring computation of $\sss$ (and the first
and the last letters of the base $B$ of $U$  the same). Then
$\br(h)$ does not contain subwords $(12)(2)(12)$, $(23)(2)(23)$,
$(23)(3)(23)$, $(34)(4)(34)$, $(45)(4)(45)$, and $(51)(5)(51)$.
\end{lemma}

\proof Notice that since the first and the last letters in the
base $B$ are the same, for every letter $z$ occurring in $B$, $U$
contains a $zz'$-sector (for some $z'$).

Let $\br(h)$ contain a subword $(12)(2)(12)$. First assume that
the base $B$ contains $L_j$ or $\lel_j$ for some $j$. Since
$\br(h)$ contains $(12)$ and rules from $\sss(12)$ lock
$zL_j$-sectors, we can conclude that $B$ contains a subword
$\lel_jL_j$ (here we use the assumption that the first and the
last letters in $B$ are the same, otherwise we could have $L_j$ as
the first letter in $B$, and $B$ could contain no copies of
$\lel_jL_j$). Since rules of $\sss(2)$ are active with respect to
$\lel_jL_j$-sectors, we get a contradiction with Lemma
\ref{(12)(2)(12)}. Therefore we can assume that $B$ contains
neither $L_j$ nor $\lel_j$. By passing to a subcomputation, we can
assume that $\br(h)=(12)(2)(12)$.

Notice that the rules from $\sss(12)\cup\sss(2)$ do not change
modulo $\gr$ the $a$-projections of admissible words without
$zL_j$- and $L_jz$-sectors. Hence the $a$-projections of
corresponding sectors in $U$ and $V$ are the same modulo $\gr$.
The $\bee$ and $\Omega$-coordinates of $U$ and $V$ are
$(1,\emptyset)$. Hence $U$ and $V$ are equal modulo $\gr$, and the
labels of the corresponding $\theta$-annuli of the ring are equal
modulo $\gr$-relations too by Lemma \ref{Band}. Hence the ring is
not reduced (see condition (iv) in the definition of reduced
diagram), a contradiction.

The other five statements can be proved by contradiction quite
analogously.
\endproof

\subsection{Standard computations} \label{roft}
\setcounter{equation}{0}

A computation $U\bullet h=V$ of $\sss$ (of $\bar\sss$) is said to
be \label{standardc}{\em standard} or $j$-{\em standard}
computation if its base has the form $(\lel_jL_jP_jR_j\rer_j)^{\pm
1}$ and the computation is reduced. If $U$ is an admissible word
with this base then for each $z\in \{\lel_j, L_j, P_j, R_j\}$ let
$z(U)$ be the $a$-projection of the $zz_+$-sector in $U$.

\begin{lemma}\label{Monoton1} Let $U\bullet h=V$ be a reduced free computation of $\sss(g)$,
$g\in \{2,3,4,5\}$. Suppose that the base of the computation is
$\lel_jL_jP_j$,

(1)  The equality $$|\lel_j(U)|+|L_j(U)|= |\lel_j(U\circ
h[t])|+|L_j(U\circ h[t])|$$ holds for every $t=1,...,|h|$.

(2) Let $z\in\{\lel_j, L_j\}$. Then the sequence of lengths
$|z(U\circ h[t])|$, $t=0,1,... $ is monotone.
\end{lemma}

\proof Indeed, since all words $U\circ h[t]$ are admissible for
$\sss$, the words $\lel(U\circ h[t])$ and $L_j(U\circ h[t])$ are
positive for $t=1,...,|h|$. If $g\ne 2,4$ then the rules
$\tau\in\sss(g)$ are not active with respect to $\lel_jL_j$- and
$L_jP_j$-sectors. Hence in that case the words $U\circ h[t]$ are
all the same ($t=1,...,|h|$) which proves (1) and (2).

If $g\in \{2,4\}$ then by Lemma \ref{csss} for every rule
$\tau\in\sss(g)$ there exist two $a$-letters
$b_1=b_1(\tau),b_2=b_2(\tau)$ such that for every admissible word
$W$, with base $\lel_jL_jP_j$, $\lel_j(W\circ\tau)=\lel_j(W)b_1$,
$L_j(W\circ\tau)=b_2L_j(W)$ and the rule $\tau$ is completely
determined by each of the letters $b_1, b_2$. Moreover $b_1$ is a
positive letter if and only if $b_2$ is a negative letter. Since
the words $\lel_j(U\circ h[t])$, $L_j(U\circ h[t])$ are always
positive and $h$ is reduced, we see that either for all $\tau$ in
$h$ the letter $b_1(\tau)$ is positive and the letter $b_2(\tau)$
is negative, or for all letters $\tau$ of $h$, $b_1(\tau)$ is
negative and $b_2(\tau)$ is positive. Therefore either for every
$t=1,...,|h|-1$, $|\lel_j(U\circ h[t+1])|=|\lel(U\circ h[t])|+1$
and $|L_j(U\circ h[t+1])|=|L_j(U\circ h[t])|-1$ or for every
$t=1,...,|h|-1$, $|\lel_j(U\circ h[t+1])|=|\lel(U\circ h[t])|-1$
and $|L_j(U\circ h[t+1])|=|L_j(U\circ h[t])|+1$. This implies
parts (1) and (2) again.
\endproof

For every word $W$, we denote by \label{diff}$\diff(W)$ the
difference between the number of positive, and the number of
negative occurrences of $a$-letters in a word $W.$ Also let, as
before, $c$ be the maximal length of relations in $\bee$.

Let $h$ be a word over $\sss$. For every $t\le |h|$ let $s(t)$ be
the number of occurrences of rules from $\sss(34)$ in $h[t]$.

\begin{lemma} \label{+sc}
 Let $U\bullet h=V$ be a $j$-standard computation of $\sss$,
$j\ne 1$. Suppose that for some $z\in\{\lel_j, L_j, R_j\}$, and
some $t\le |h|$
\begin{equation}|z(U\circ h[t])|
> s(t)c+|U|.\label{sc}\end{equation}

Then, for every $d\ge t$,

(a) the $d$-th rule of $h$ does not belong to
$\sss(12)\cup\sss(1)\cup\sss(34)\cup\sss(51)$ if $z\in
\{\lel_j,R_j\}$ and the $d$-th rule of $h$ does not belong to
$\sss(34)$ if $z=L_j$,

(b) if $z\in \{\lel_j, R_j\}$ then the word $z(U\circ h[d])$ is
not empty; if $z=L_j$ then at least one of the words
$\lel_j(U\circ h[d])$ or $L_j(U\circ h[d])$ is not empty,

(c) $|U\circ h[d||\ge d-t$.
\end{lemma}

\proof We may assume that the base of $U$ has the form
$\lel_jL_jP_jR_j\rer_j$ because otherwise we can replace $U$ by
$U\iv$.

{\bf Case 1.} Suppose that $z=\lel_j$. Without loss of generality
we can assume that $t$ is the minimal number such that (\ref{sc})
holds. From (\ref{sc}), we have that $t>0$.

Then the $t$-th rule of $h$ belongs to
$\sss(2)\cup\sss(34)\cup\sss(4)$ because otherwise $\lel_j(U\circ
h[t-1])=\lel_j(U\circ h[t])$ which contradicts the minimality of
$t$. Similarly this rule cannot belong to $\sss(34)$ because in
that case $s(t-1)=s(t)-1$ and
$$|\lel_j(U\circ h[t-1])|\ge |\lel_j(U\circ h[t])|-c>
s(t)c-c+|U|\ge s(t-1)c+|U|$$ contrary to the minimality of $t$.

Furthermore, the minimality of $t$ and the form of the rules from
$\sss(2)\cup\sss(4)$ imply that $|\lel_j(U\circ
h[t])|=|\lel_j(U\circ h[t-1])|+1$.

{\bf Case 1.1.} Suppose first that the $t$-th rule in $h$ belongs
to $\sss(2)$. If for every $t'\ge t$ the $t'$-rule in $h$ is from
$\sss(2)$ then

$$|\lel_j(U\circ h[t])|<|\lel_j(U\circ h[t+1])|<...
<|\lel_j(U\circ h)|$$ by Lemma \ref{Monoton}. (The series
increases because of the minimality of $t$.) This implies parts
(a), (b), (c) of the conclusion of the lemma since $\sss(12)$
locks $\lel_j$-sector.

Suppose that for some $t'>t$ the $t'$-th rule in $h$ is not from
$\sss(2)$. Let $t'$ be the first such number. Then by Lemma
\ref{Br Hist} the $t'$-th rule in $h$ belongs to
$\sss(12)\cup\sss(23)$. Since rules from $\sss(12)$ lock
$\lel_jL_j$-sectors, the $t'$-th rule in $h$ cannot belong to
$\sss(12)$ by Lemma \ref{Monoton}. So it belongs to $\sss(23)$.

Clearly the claims of the lemma are true if $t'=|h|$.

Suppose that $t'<|h|$. Then by Lemma \ref{Br Hist} $h\equiv
h[t']h'h''$ where $h'$ is a maximal subword of $h$ consisting of
rules from $\sss(3)$. We will show that $h''$ is empty.

Indeed, suppose that $h''$ is not empty. Then by Lemma \ref{Br
Hist}, $h''$ starts with a rule $\tau$ from $\sss(34)$. Let
$U'=U\circ h[t']h'$. Since $\tau$ locks $P_jR_j$-sectors, and
since all words $U\circ h[u]$ are admissible for $\sss$
($u=0,1,...,|h|$), we have that $P_j(U')=\emptyset$ and the words
$\lel_j(U')$, $L_j(U')$ and $R_j(U')$ are positive.

Therefore, by Lemma \ref{Monoton1},
\begin{equation}\label{6.2}\diff(U')=|U'_a|\ge
|\lel_j(U')|=|\lel_j(U\circ h[t'])|\ge |\lel_j(U\circ
h[t])|+(t'-t)>s(t)c+|U_a|+(t'-t).\end{equation} But it is easy to
see that for every rule $\tau$  in $\sss\backslash\sss(34)$ and
for every admissible word $W$ in the domain of $\tau$,
$\diff(W)=\diff(W\circ\tau)$, and for every $\tau\in\sss(34)$,
$\diff(W\circ\tau)\le \diff(W)+c$. Since there are at most $s(t)$
rules from $\sss(34)$ in $h[t']h'$, we can conclude that

\begin{equation}\label{6.2.2}
\diff(U')\le \diff(U)+s(t)c\le|U_a|+s(t)c,\end{equation} a
contradiction with (\ref{6.2}).

Thus $h''$ is empty, and part (a) of the lemma is established.
Since rules from $\sss(3)$ are not active with respect to
$\lel_jL_j$-sectors, and $\lel_j(U\circ h[t])$ is not empty, part
(b) follows as well.

Recall that $h'$ is not empty because $t'\ne |h|$. Notice that the
$R_j\rer_j$-sector $W$ of $U\circ h[t']$ is in the domain of $h'$.
By Lemma \ref{height} applied to the computation $W\bullet h'$, we
have
\begin{equation}\label{6.3.1}
 |h'|=|h|-t'\le |R_j(W)|+|R_j(W\circ h')|= |R_j(U\circ
h[t'])|+|R_j(U\circ h)|.
\end{equation}

 Notice that rules in $\sss(23)$ lock
$R_j\rer_j$-sectors, so $|R_j(U\circ h[t'])|=0$ (since the last
rule in $h[t]$ is from $\sss(23)$), and therefore $R_j(U\circ
h)|\ge d-t'$ by (\ref{6.3.1}). Notice also that $\lel_j(U\circ
h[d])|\ge |\lel_j(U\circ h[t])+t'-t$ by (\ref{6.2}) and the
equality $\lel_j(U\circ h[d])=\lel_j(U\circ h[t'])$. Hence for
every $d\ge t'$ $$\begin{array}{l}|U\circ h[d]|> |\lel_j(U\circ
h[d])|+|R_j(U\circ h[d])|\ge\\ (|\lel_j(U\circ
h[t])|+t'-t)+(d-t')\ge s(t)c+d-t+|U|\ge d-t,\end{array}$$ which
proves part (c).

{\bf Case 1.2.} The case when the $t$-th rule in $h$ belongs to
$\sss(4)$ is similar. Let $h=h[t]h'$. Then the brief history of
$h'$ cannot start with $(4)(45)$ by Lemma \ref{Monoton} since
rules from $\sss(45)$ lock $\lel_jL_j$-sectors.

The brief history also cannot start with $(4)(34)$. Indeed,
suppose that $h'=h''\tau h'''$ where $h''$ is a word over
$\sss(4)$, $\tau\in\sss(34)$. Then as in Case 1.1, we denote
$U'=U\circ h[t]h''$. Then we can use the fact that $\tau$ locks
$P_jR_j$-sectors, and that words $\lel_j(U')$ and $R_j(U')$ are
positive, we deduce inequalities
\begin{equation}\label{6.2.1}\diff(U')=|U'_a|\ge
|\lel_j(U')|=|\lel(U\circ h[t]h'')|\ge |\lel(U\circ
h[t])|+|h''|>s(t)c+|U_a|+|h''|.\end{equation} On the other hand we
have the inequality (\ref{6.2.2}) as in Case 1.1, a contradiction
with (\ref{6.2.1}).

{\bf Case 2}. The case when $z=R_j$ is completely similar to Case
1.

{\bf Case 3}. Let $z=L_j$. As in Case 1, we assume $t$ being
minimal. Then the $t$-th rule in $h$ belongs to
$\sss(2)\cup\sss(4)\cup\sss(1)$.

{\bf Case 3.1}. Suppose  the $t$-th rule in $h$ belongs to
$\sss(2)$. Let $h=h[t]h_1$. If all rules in $h_1$ belong to
$\sss(2)$, we prove the statement as in Case 1.

If $h_1=h_2\tau h_3$ where $h_2$ is a word over $\sss(2)$,
$\tau\in\sss(23)$ then by Lemma \ref{Monoton} $|L_j(U\circ
h[t-1])|<|L_j(U\circ h[t]h_2)|$ which is impossible since
$|L_j(U\circ h[t]h_2)|=0$ (because rules from $\sss(23)$ lock
$L_jP_j$-sectors).

Hence we can assume that $\tau\in\sss(12)$.

{\bf Case 3.1.1}. Assume that $h[t]h_2$ is not a word in
$\sss(2)$. Then it contains a rule from $\sss(12)$ or from
$\sss(23)$. The first option is impossible by Lemma \ref{Br Hist}.
Hence $h[t]= h[q]\tau'h'$ where $\tau'\in\sss(23)$, $h'$ is a word
in $\sss(2)$.

By Lemma \ref{Monoton}, the value $|\lel_j(U\circ h[p])|$ must
decrease when $p$ runs from $q$ to $t+|h_2|$ because
$|\lel_j(U\circ h[t]h_2)|=0$ (since $\tau$ locks
$\lel_jL_j$-sectors). Hence the sum $|\lel_j(U\circ
h[p])|+|L_j(U\circ h[p])|$ cannot increase. Therefore
$|\lel_j(U\circ h[q])|\ge |L_j(U\circ h[t])|>s(t)c+|U|$ since
$|L_j(U\circ h[t]h_2)|=0$ ($\tau'$ locks $L_jP_j$-sectors).

Since $s(q)\le s(t)$, the number $q$ satisfies the hypothesis of
the lemma for $z=\lel_j$ (Case 1 of the proof). But this
contradicts part (a) of the lemma which we already proved in Case
1 (by this part, in Case 1, the brief history of the suffix of
$h$ that starts after $h[q]$ cannot contain (12)).

{\bf Case 3.1.2}. Now assume that $h[t]h_2$ is a word in
$\sss(2)$. Hence $s(t)=0$ (since $h[t]$ does not contain rules
from $\sss[34]$). By the assumption of the lemma
$$|\lel_j(U\circ h[t])|+|L_j(U\circ h[t])|\ge
|L_j(U\circ h[t])|>|U|\ge |\lel_j(U)|+|L_j(U)|.$$

The strict inequality means that $|\lel(U\circ h[p])|$ cannot
monotonically decrease when $p$ runs from 0 to $t+|h_2|$ (an age
of $(2)$). Hence $|\lel_j(U\circ h_th_2|>0$ by Lemma
\ref{Monoton}. But this contradicts the fact that $\tau$ (as every
rule in $\sss(12)$) locks $\lel_jL_j$-sectors. This contradiction
completes Case 3.1.

{\bf Case 3.2}. The case when the $t$-th rule in $h$ belongs to
$\sss(4)$ is similar, but we need to substitute $\sss(4),
\sss(45), \sss(34)$ for $\sss(2), \sss(12), \sss(23)$.

{\bf Case 3.3}. Suppose that the $t$-th rule in $h$ belongs to
$\sss(1)$. Let $h=h[t]h'$. If all rules in $h'$ belong to
$\sss(1)$ then we can argue as in Case 1, apply Lemma
\ref{Monoton}, and deduce parts (a), (b) of the conclusion of the
lemma, and the inequality

$$|L_j(U\circ h[d])|\ge |L_j(U\circ h[t])|+(d-t)$$
for every $d\ge t$ which implies part (c).

{\bf Case 3.3.1}. Assume that the letter $(1)$ in the brief
history of $h'$ is followed by $(12)$. Let $h'=h'[t']h_1$ where
$h'[t']$ is the maximal prefix of $h'$ which is a word over
$\sss(1)$. Since $|h_1|>0$, $d\ge t+t'+1$. As above, we obtain by
Lemma \ref{Monoton},
\begin{equation}
|L_j(U\circ h[t+t'+1])|=|L_j(U\circ h[t+t'])|= |L_j(U\circ
h[t])|+t'\ge |U|+s(t)c+t'+1 \ge t'+1>0, \label{t't}
\end{equation}
and all three parts of the statement of the lemma are true if
$d=t+t'+1$. So suppose that $t+t'+1<d$.

Then
$$\diff(U\circ h[t+t'+1]) >
|U_a|+s(t)c+t'-|P_j(U\circ h[t+t'+1])|,$$ but $\diff(U\circ
h[t+t'+1])\le |U_a|+s(t)c$ by the same argument as in Case 1.
Therefore \begin{equation}|P_j(U\circ h[t+t'+1])|\ge
t'+1.\label{ppp}\end{equation}

Then, by Lemma \ref{Br Hist}, the brief history of $h_1$ starts
with $(12)(2)$ (it cannot start with $(12)(23)$ since $|L_j(U\circ
h[t+t'+1])|>0$ by (\ref{t't})). Let $h_1=h_1[t_1]h_2$ where
$h_1[t_1]$ is the maximal prefix of $h_1$ with brief history
$(12)(2)$.

{\bf Case 3.3.1.1}. Suppose that $h_2$ is empty. Then part (a) of
the lemma is obviously true. Part (b) is also true because
$|\lel_j(h[q])|+|L_j(h[q]|$ stays the same for all $q\ge t+t'+1$
(Lemma \ref{Monoton1} (1)), and $|L_j(U\circ h[t+t'+1])|>0$  by
(\ref{t't}).

In order to prove part (c), notice that by Lemma \ref{Monoton1}
(2) the sequences $|L_j(U\circ h[p])|$, $p=t+t'+1,t+t'+2...$, is
monotonically decreasing. Therefore $d-t-t'\le |L_j(U\circ
h[t+t'+1])|$. Therefore by Lemma \ref{Monoton1} (1) and by
(\ref{ppp})

\begin{equation}\label{3311}\begin{array}{ll} |U\circ h[d]|= \\ &
|U\circ h[t+t'+1]|> |L_j(U\circ
h[t+t'+1])|+|P_j(U\circ h[t+t'+1])| \ge \\ & (d-t-t')+ (t'+1)\ge\\
d-t+1\end{array}\end{equation} which proves part (c).

{\bf Case 3.3.1.2}. Now suppose that $h_2$ is not empty, and $d>
t+t'+t_1$. If $|h_2|=1$ then $h_2\in \sss(23)$ and parts (a), (b)
follow from Case 3.3.1.1, part (c) follows from (\ref{3311}). So
suppose that $|h_2|>1$, $h_2=h_2[t_2]h_3$ where $h_2[t_2]$ is the
maximal prefix of $h_2$ with brief history $(23)(3)$.

Since the brief history of $h_2$ starts with (23) and rules of
$\sss(23)$ lock $L_jP_j$-sectors, $|L_j(U\circ h[t+t'+1+t_1]|=0$.
As in Case 3.3.1.1 we have that $|\lel_j(U\circ h[t+t'+1+t_1])|\ne
0$. Since rules from $\sss(3)$ are not active for
$\lel_j$-sectors, the length $|\lel_j(U\circ h[q])|$ stays the
same (and $>0$) while $q$ runs from $t+t'+1+t_1$ to
$t+t'+1+t_1+t_2$.

Suppose that $h_3$ is not empty, then it starts with a rule in
$\sss(34)$. Since rules in $\sss(34)$ lock $P_jR_j$-sectors,
$P_j(U\circ h[t+t'+1+t_1+t_2])$ is empty.

Therefore we obtain by Lemma \ref{Monoton1}
$$
\begin{array}{l} s(t)c+|U| <  \\
\qquad |\lel_j(U\circ h[t])|+|L_j(U\circ h[t])| \le |\lel_j(U\circ
h[t+t'+1])|+|L_j(U\circ h[t+t'+1])|= \\ \qquad |\lel_j(U\circ
h[t+t'+1+t_1+t_2])| + |L_j(U\circ h[t+t'+t_1+t_2])|\le \\
\diff(U\circ h[t+t'+1+t_1+t_2]).\end{array}$$

But, as in Case 1, $\diff(U\circ h[t+t'+1+t_1+t_2])\le
\diff(U)+s(t)c$, a contradiction.

Thus $h_3$ is empty. Parts (a) and (b) of the lemma follow as in
Case 3.3.1.1. Also notice that since rules in $\sss(3)$ lock
$L_jP_j$-sectors, we have

\begin{equation}\label{10}
|U\circ h[d]|\ge |\lel_j(U\circ h[d])|+|P_j(U\circ h[d]
)|.\end{equation}

Since rules in $\sss(3)$ are not active for $\lel_jL_j$-sectors,
$$|\lel_j(V)|=|\lel_j(U\circ h[t+t'+1+t_1])|.$$ By Lemma
\ref{Monoton1}, since rules in $\sss(12)$ lock $\lel_jL_j$-sectors
and rules in $\sss(23)$ lock $L_jP_j$-sectors, we have
$$|\lel_j(U\circ h[t+t'+1+t_1])|=|L_j(U\circ h[t+t'+1])|\ge
|L_j(U\circ h[t])|+t'>|U|+s(t)c+t'.$$ Hence

\begin{equation}\label{11}
|\lel_j(U\circ h[d])|>|U|+s(t)c+t'.
\end{equation}

Also notice that by Lemma \ref{Monoton},
\begin{equation}\label{12}|\lel_j(U\circ h[d])|=t_1-1.\end{equation}

By Lemma \ref{Monoton}, since rules of $\sss(2)\cup\sss(23)$ lock
$R_j\rer_j$-sectors, we have that

\begin{equation}\label{13}
R_j(U\circ h[d])=|h_2|-1=d-t-t'-t_1-1.\end{equation}

Also notice that $$|\lel_j(U\circ h[d])|+|R_j(U\circ h[d]
)|-|P_j(U\circ h[d])|\le \diff(U\circ h[d])=\diff(U\circ h[t]))\le
\diff(U)+s(t)c.$$ Hence, by (\ref{11}) and (\ref{13}),

\begin{equation}\label{14}\begin{array}{l}
|P_j(U\circ h[d])|\ge \\ \qquad (|\lel_j(U\circ
h[d])|-s(t)c-\diff(U))+
|R_j(U\circ h[d])|\ge\\
\qquad
t'+d-t-t'-t_1-1+|U|-\diff(U)\ge \\
d-t-t_1+1.\end{array}
\end{equation}

Now using (\ref{12}) and (\ref{14}), we get

$$|U\circ h[d]|\ge |\lel_j(U\circ h[d])|+|P_j(U\circ h[d])|\ge (t_1-1)+d-t-t_1+1=d-t$$ as required.

{\bf Case 3.3.2}. The remaining case is when the letter $(1)$ in
the brief history of $h$ is followed by $(51)$. This case can be
analyzed as above, but the analysis is shorter, because the
longest possible suffix of the brief history under consideration
now is $(1)(51)(5)$ (instead of $(1)(12)(2)(23)(3)$ in Case
3.3.1). This subword cannot be followed by $(45)$ by the same
argument as in Case 3.3.1.2 (see the argument showing that the
subword $(1)(12)(2)(23)(3)$ cannot precede the letter $(34)$)
because rules from $\sss(45)$ lock $P_jR_j$-sectors. The lemma is
proved.
\endproof

\begin{lemma} \label{Soon or Never}
Let $U\bullet h=V$ be a $j$-standard computation of $\sss$, $j\ne
1$. Then

(1) if $h$ contains a letter from $\sss(34)$ then such a letter
occurs in the prefix $h[30(|U|+c)]$.

(2) if $|h|\ge |V|+30(|U|+c)$ then $h[30(|U|+c)]$ contains a
subword starting with a rule from $\sss(34)$ and ending with a
rule from $\sss(12)\cup\sss(51)$.
\end{lemma}

\proof (1) Suppose that $h$ contains rules from $\sss(34)$. Let
$l=30(|U|+c)$. Suppose, by contradiction, that $h[l]$ does not
contain a rule from $\sss(34)$.

If there exists $t\le l$ that satisfies the conditions of Lemma
\ref{+sc} (for some $z\in \{\lel_j, L_j, R_j\}$) then $s(t)=0$,
and by Lemma \ref{+sc} $h$ does not contain letters from
$\sss(34)$, a contradiction.

Therefore for every $t=1,...,l$, $$|z(U\circ h[t])|\le |U|$$ for
every $z\in \{\lel_j, L_j, R_j\}$. Then by Lemma \ref{height}, the
length of an arbitrary subword of $h[l]$ over $\sss(g)$, $g\in
\{1,2,3,4,5\}$ is at most $2|U|$. Since $(34)$ is not in the brief
history of $h[l]$, the brief history of $h[l]$ contains at most 5
occurrences of $g\in \{1,2,3,4,5\}$ by Lemma \ref{Br Hist}. Thus
the length $l$ of $h[l]$ does not exceed $10|U|+4<30(|U|+c)$
contrary to the definition of $l$.

(2) If there are two letters from $\sss(34)$ in $h[l]$ then there
is a letter from $\sss(12)$ and a letter from $\sss(51)$ between
them by Lemma \ref{Br Hist}, and the claim of the lemma is true.

Therefore, we will suppose that $h[l]$ contains exactly one letter
from $\sss(34)$. By Lemma \ref{Br Hist}, the length of the brief
history of $h[l]$ is at most $19$, and the brief history contains
at most $10$ maximal subwords over $\sss(1)\cup...\cup\sss(5)$.

Suppose that some number $t\le |h|$ satisfies an assumption of
Lemma \ref{+sc} for some $z\le\{\lel_j, L_j, R_j\}$, and $s(t)=1$.
Take the minimal $t$ with this property. Then $|z(U\circ h[p])|\le
c+|U|$ for $p<t$, for any  $z\in \{\lel_j, L_j, R_j\}$.

Then, as in the proof of part (1),  $t-1\le
10(2|U|+2c)+9<30(|U|+c)$ because the number of ages of  $(1),
(2),\dots (5)$ is at most $10$ in the history of $h[t-1]$. This
inequality together with the inequality $|U\circ h|\ge |h|-t$
given by Lemma \ref{+sc} (part (c)) implies that $|h|\le |U\circ
h|+30(|U|+c)$ contrary to the assumption of part (2).

Hence there is no number $t$ between 1 and $|h|$ such that
$|z(U\circ h[t])|\ge c+|U|$ for some $z\in \{\lel_j, L_j, R_j\}$.
Then, as in (1), Lemmas \ref{height} and \ref{Br Hist} provide us
with the inequality $|h|\le 10(2|U|+2c)+9<30(|U|+c)$ contrary to
the assumption of the lemma. The lemma is proved.
\endproof

\begin{lemma} \label{Determ}
Let $U\bullet h=V$ be a $j$-standard computation of $\sss$ and
$j\ne 1$ or a $j$-standard computation of $\bar\sss$ (for any
$j$).

Assume that $h$ contains an occurrence of a rule $\tau_1$ from
$\sss(34)\cup\bar\sss(34)$ preceding an occurrence of a rule
$\tau_2$ from
$\sss(12)\cup\sss(51)\cup\bar\sss(12)\cup\bar\sss(51)$:
$h=h_1\tau_1h_2\tau_2h_3$. Then the words $z(U\circ
h_1\tau_1h_2)$, $z\in \{L_j, P_j\}$, and $z(U\circ h_1)$, $z\in
\{\lel_j,R_j\}$ are completely determined by the word $\tau_1h_2$.
\end{lemma}

\proof If $j=1$ and $h$ is a word over $\bar\sss$, then all the
words listed in the formulation of the lemma are empty, so the
statement is trivially true.

Suppose that $j\ne 1$. We may suppose that $h_2$ does not contain
rules from $\sss(34)\cup\sss(12)\cup\sss(51)\cup\bar\sss(34)\cup
\bar\sss(12)\cup\bar\sss(51)$. By Lemma \ref{Br Hist}, the brief
history of $\tau_1h_2\tau_2$ is either $(34)(3) (23)(2)(12)$ or
$(34)(4)(45)(5)(51)$. These cases are similar, so let the brief
history be equal to $(34)(3) (23)(2)(12)$.

Let $h_2=h'\tau h''$ where $\tau\in\sss(23)$ (or $\tau\in
\bar\sss(23)$). Since rules from $\sss(23)$ lock $L_jP_j$-sectors,
and each rule $\tau(2,r,i)^{\pm 1} $ of $\sss(2)$ (of
$\bar\sss(2)$) multiplies inner parts of the $L_jP_j$-sectors of
admissible words by the letter $a_i(L_j)^{\mp 1}$ (resp. by $\bar
a_i(L_j)^{\mp 1}$) on the left (Lemma \ref{csss}), the word
$L_j(U\circ h'\tau h'')$ is a copy of $h''$ written from right to
left. Therefore the word $L_j(U\circ h_1\tau h_2)$ is completely
determined by $h_2$. Similarly, the word $P_j(U\circ h_1\tau h_2)$
is completely determined by $h'$ (and, in turn, by $h_2$) because
rules from $\sss(34)$ lock $P_jR_j$-sectors, rules from $\sss(3)$
are active with respect to $P_jR_j$-sectors, and rules from
$\sss(23)\cup\sss(2)$ are not active with respect to these
sectors.

Similarly, since $\tau_2$ (as all rules from $\sss(12)$) locks
$\lel_j$-sectors and $R_j$-sectors, rules from $\sss(2)$ are
active with respect to these sectors, and rules from
$\sss(3)\cup\sss(23)$ are not active with respect to these
sectors, the word $\lel_j(U\circ h_1)$ (resp. $R_j(U\circ h_1)$)
is completely determined by $\tau_1h_2$.
 \endproof

\begin{lemma} \label{Posit Band}
Let $U\bullet h=V$ be a $j$-standard computation of $\sss$ and
$j\ne 1$. Suppose that $|h| \ge |V|+30(|U|+c)$. Then
$h[30(|U|+c)]$ can be decomposed as
$h[30(|U|+c)]=h_1\tau_1h_2\tau_2h_3$ where $\tau_1\in \sss(34)$,
$\tau_2\in \sss(12)\cup\sss(51)$, $h_2$ contains no letters from
$\sss(34)\cup\sss(12)\cup\sss(51)$. Moreover the word $P_j(U\circ
h_1\tau_1h_2\tau_2)$ is positive for any such decomposition of
$h[30(|U|+c)]$.
\end{lemma}

\proof By Lemma \ref{Soon or Never} there exists a decomposition
$h[30(|U|+c)]=h_1\tau_1h_2\tau_2h_3$ where $\tau_1\in \sss(34)$,
$\tau_2\in \sss(12)\cup\sss(51)$. Taking $h_2$ of minimal length,
we can assume that $h_2$ contains no letters from
$\sss(34)\cup\sss(12)\cup\sss(51)\cup\sss(1)$. So we only need to
show that $P_j(U\circ h_1\tau_1h_2\tau_2)$ is positive. Consider
the homomorphism \label{gamma}$\gamma $ from the free group
generated by $\aaa\cup\kk$ into the free group generated by
$\{a_1,...,a_{\bar m}\}$ which maps $a_i(z)$ into $a_i$ for every
$z$, and maps all other letters to 1. It is clear that for every
rule $\tau\in\sss\backslash\sss(34)$, and every admissible word
$W$ in the domain of $\tau$, $\gamma(W\circ\tau)=\gamma(W)$. Let
$U'$ be the subword of $U$ stating with a $P_j$-letter and ending
with a $\rer_j$-letter. Then $U'\circ h=V'$ where $V'$ is the
corresponding  subword in $V$. Moreover for every $t$ between $1$
and $|h|$, and $z\in\{P_j,R_j\}$, the $zz_+$-sector in $U'\circ
h[t]$ coincides with the $zz_+$-sector in $U\circ h[t]$. Since
rules in $\sss(34)$ lock $P_jR_j$-sectors, $$\gamma(U'\circ
h_1\tau_1)=\gamma(R_j(U'\circ h_1\tau_1)).$$ The word in the right
hand side of this equality is a copy of $R_j(U'\circ h_1\tau_1)$
and is positive by definition of admissible words. Hence
$\gamma(U'\circ h_1\tau_1h_2)=\gamma(R_j(U'\circ h_1\tau_1)$ (we
use here that $h$ does not contain rules from $\sss(1)$) is a
positive word. But the rules in $\sss(12)\cup\sss(51)$ lock
$R_j\rer_j$-sectors, so $\gamma(U'\circ h_1\tau_1h_2)$ is a copy
of $P_j(U'\circ h_1\tau_1h_2)=P_j(U\circ h_1\tau_1h_2)$. Hence
$P_j(U\circ h_1\tau_1h_2)$ is a positive word as required.
\endproof

To formulate the next lemma, we introduce a natural homomorphism
$\lambda$ from the free group $F_0$ with basis $\aaa\cup\bar\aaa$
onto the factor group $H_0$ of $F_0$ over the normal closure of
all $\gr$-relations. Thus $H_0$ is the free product of a copy of
$\gr$ generated by $a_1(P_1),...,a_m(P_1)$ and a free group
generated by all other $a$-letters.

\begin{lemma} \label{j=1}
Let $U\bullet h=V$ be a $1$-standard $\gr$-computation of $\sss$.
Suppose that $h=h_1\tau_1h_2\tau_2h_3$ where $\tau_1\in\sss(34)$,
$\tau_2\in\sss(12)\cup\sss(51)$. Assume that $h$ satisfy the
conclusion of Lemma \ref{Br Hist}. Then the words $z(U\circ
h_1\tau_1)$, $z\in\{\lel_1,R_1\}$, $L_1(U\circ h_1\tau_1h_2)$ and
the $\lambda$-image of $P_1(U\circ h_1\tau_1h_2)$ are determined
by the word $\tau_1h_2\tau_2$.
\end{lemma}

\proof The proof of the assertion about $z(U\circ h_1\tau_1)$,
$z\in\{\lel_1,R_1\}$ and $L_1(U\circ h_1\tau_1h_2)$ coincides with
the corresponding proof in Lemma \ref{Determ}. Let $U'$ be the
$P_1R_1$-sector in $U$. Then $P_1(U\circ h[t])=P_1(U'\circ h[t])$
for every $t$. Notice also that $P_1(U'\circ h_1)$ is equal to 1
modulo $\gr$ because rules from $\sss(34)$ lock the
$P_1R_1$-sectors. Hence, modulo $\gr$, $U'\circ h_1$ is a 2-letter
word over $\kkk$. Therefore $U'\circ h_1\tau_1h_2=(U'\circ
h_1\tau_1)\circ h_2$ is determined by $h_2$ modulo $\gr$ as
desired. Hence the $\lambda$-image of $P_1(U\circ h_1\tau_1)$ is
determined by the word $\tau_1h_2\tau_2$.
\endproof

\subsection{Tame and wild computations}
\label{twc}

Recall that for every positive word $w$ in $\{a_1,...,a_{\bar
m}\}$ we have defined a word $\Sigma(w)$ in (\ref{Sigma(w)}) which
has the base $\tsigma$ defined in (\ref{basehub}).

Let $U$ be an admissible word for $\sss$ with a base
$\tsigma^sK_1$. Then the word $U$ and any ring computation
$U\bullet h=V$ over $\csss$  will be called \label{tamec}{\em
tame}. All non-tame admissible words $U$ whose bases start and end
with the same letter will be called \label{wild}{\em wild}. If $U$
is a wild admissible word then any ring computation $U\bullet h=V$
will be called {\em wild} as well.

Recall that every admissible word for $\csss$ has an
$\bee$-coordinate and an $\Omega$-coordinate. All $k$-letters in
an admissible word $U$ have the same coordinates $r$, $i$, so we
can write $U=U'(r,i)$ where $U'$ is the word $U$ with $(\bee,
\Omega)$-coordinates removed from all $k$-letters.

\begin{lemma} \label{Tame Ring}
Let $U\bullet h=V$ be a tame computation of $\sss$, with the base
$U=\tsigma^sK_1$. Suppose that $|h|\ge |V|+30(|U|+c)$. Then the
prefix $h[30(|U|+c)]$ has the form $h_1\tau h_2$ where
$\tau\in\sss(12)\cup\sss(51)$ and $U\circ h_1$ or $U\circ h_1\tau$
has the form $\Sigma(u,v)^sK_1$ where $u$ and $v$ are positive
words over $\{a_1,...,a_{\bar m}\}$, and $\Sigma(u,v)$ is an
admissible word of the following form
\begin{equation}\label{6.3}
\begin{array}{l}(\prod_{j=1}^{sN/2}
K_{2j-1}L_{2j-1}v_{2j-1}P_{2j-1}u_{2j-1}R_{2j-1} K_{2j}\iv
R_{2j}^{-1}u_{2j}^{-1}P_{2j}^{-1}v_{2j}^{-1}L_{2j}^{-1})(\emptyset,1)
\end{array}\end{equation} where $v_k$, $u_{k'}$
are copies of $v$ and $u$ under the substitutions $a_i\mapsto
a_{i}(L_k)$, $a_i\mapsto a_i(P_{k'})$ for every $k$ and for $k'\ne
1 (mod N)$; $u_{1+Nl}$ is equal to a copy of $u$ modulo $\gr$.
\end{lemma}

\proof Let $U'$ be any admissible subword of $U$ with base
$\lel_jL_jP_jR_j\rer_j$. Let $V'$ be the corresponding subword in
$V$. Then $U'\bullet h=V'$ is a $j$-standard computation, and we
can apply Lemma \ref{Posit Band}. This gives us the desired
occurrence of $\tau\in\sss(12)\cup\sss(51)$ in the prefix
$h[30(|U|+c)]$. Then, depending on whether $\tau\in
\sss^-(12)\cup\sss^+(51)$ or $\tau\in\sss^+(12)\cup\sss^-(51)$,
either $U\circ h_1\tau$ or $U\circ h_1$ has the form (\ref{6.3}).
The fact that all words $u_i$ are copies of each other and all
words $v_i$ are copies of each other follows from Lemmas
\ref{Determ} and \ref{j=1} (in case $j=1$) since the histories of
the $j$-standard computations are all equal to $h$ (notice that
Lemma \ref{j=1} can be applied because $h$ satisfies the
conclusion of Lemma \ref{Br Hist} since there exists a
$j$-standard computation $U'\bullet h=V'$ for $j\ne 1$). The fact
that $u_i$ and $v_i$ are positive words follows from Lemma
\ref{Posit Band} and from the fact that words $U\circ h[t]$ are
admissible for $\sss$.
\endproof

\begin{lemma} \label{Est' sector}
Let $|\sss|$ be he number of rules in $\sss$. Let $U\bullet h=V$
be a wild ring computation of $\sss$. Suppose that $|h|\ge
|\sss|+1$, $h$ contains a rule from $\sss(34)$. Then $U^{\pm 1}$
contains an admissible subword with base $\lel_jL_jP_jR_j\rer_j$
for some $j\ge 1$. If in addition the brief history of $h$
contains one of the letters $(1), (12), (2), (23), (45), (5),
(51)$, then $U$ contains such an admissible subword with base
$\lel_jL_jP_jR_j\rer_j$ for some $j>1$.
\end{lemma}

\proof If the base of $U$ does not contain $L$-, $P$- and
$R$-letters, then the words $U\circ h[t]$, $t=1,...,|h|$ differ
only by the $(\bee,\Omega)$-coordinates (since rules of $\sss$ do
not change $K_jK_j\iv$- and $K_j\iv K_j$-sectors of admissible
words). Since the computation $U\bullet h=V$ is a ring
computation, the labels of the annuli of the ring, which
correspond to words $U\circ h[t]$, $U\circ h[t']$ ($t\ne t'$) can
not be equal modulo $\gr$-relations. Then lemma \ref{Band} shows
that the number of the annuli  cannot be greater than the number
of rules of $\sss$. Hence $|h|$ does not exceed $|\sss|+1$.

Assume now that the base contains $L_j$-, $P_j$- or $R_j$-letter
for some $j$. Then since rules from $\sss(34)$ lock $L_jP_j$ and
$P_jR_j$-sectors, the base of $U^{\pm 1}$ must contain a subword
$L_jP_jR_j$. By Lemma \ref{noPivP} this subword must be inside a
subword of the form $\lel_jL_jP_jR_j\rer_j$ and the first claim
of the lemma true.

Suppose that $j=1$. Notice that $\rer_1=\rer_2^{-1}$. Now if the
brief history of $h$ contains letter from $\{(1), (12), (2), (23),
(51)\}$ then the corresponding rules in $h$ lock
$R_2\rer_2$-sectors. Therefore the base of $U^{\pm 1}$ contains
$R_2\iv$ which implies, as before, that it contains a subword of
the form $\lel_2L_2P_2R_2\rer_2$, as required.

Similarly notice that $\lel_1=\lel_N$, and rules from
$\sss(45)\cup\sss(5)$ lock $\lel_jL_j$-sectors for every $j$. So
if the brief history of $h$ contains $(45)$ or $(5)$ then the base
of $U^{\pm 1}$ contains the subword $\lel_NL_NP_NR_N\rer_N$.
\endproof

\begin{lemma} \label{(34) and (12)} Let $U\bullet h=V$ be a wild or tame ring computation of $\sss$,
$|h|$ exceeds $|\sss|+1$. Suppose that $h$ contains a rule from
$\sss(34)$ and a rule from $\sss(1)\cup\sss(12)\cup\sss(51)$. Then
the computation is tame.
\end{lemma}

\proof The presence of a rule locking $\lel_jL_j$-sector (locking
$R_j\rer_j$-sector), implies that the letter $\lel_j$ (letter
$\rer_j$) can occur in the base of $U$ only in subwords
$(\lel_jL_j)^{\pm 1}$ (resp. in subwords $(R_j\rer_j)^{\pm 1}$).
Together with the argument from the proof of Lemma \ref{Est'
sector}, this implies that arbitrary letter of the base of $U$ has
a uniquely determined right neighbor, which coincides with the
right neighbor of this letter in the base of a tame ring.
\endproof

\begin{lemma} \label{(34)}
Let $U\bullet h=V$ be a wild ring computation of $\sss$, $h$
contains a rule from $\sss(34)$. Then $|h|$ is bounded by a linear
function in $|U|+|V|$.
\end{lemma}

\proof By Lemma \ref{(34) and (12)}, we may assume that the brief
history $\br(h)$ of $h$ does not contain the letters $(12)$ and
$(51)$. Suppose in the beginning, that $\br(h)$ has at least one
of the letters $(1),(2), (23), (45), (5)$. Then, applying Lemmas
\ref{Est' sector} and \ref{Soon or Never}, we conclude that $|h|$
is bounded by a linear function in $|U|+|V|$.

So, we may assume that $\br(h)$ contains no letters from the list
$$\{(1), (12), (2), (23), (45), (5), (51)\},$$ i.e. it is a word
in letters $(3), (34), (4)$. By Lemma \ref{Ring Br Hist} it has no
subwords $(34)(4)(34)$. Thus, we may suppose that $\br(h)$ is a
subword of $(4)(34)(3)(34)(4)$. Hence
$$h=h_{(4)}h_{(34)}h_{(3)}h'_{(34)}h'_{(4)}$$
where $h_\alpha$ is a word in $\sss(\alpha)$ (the first and the
last few subwords in this decomposition of $h$ may be empty),
$|h_{(34)}|,|h'_{(34)}|\le 1$.

By Lemma \ref{Est' sector}, we may assume that $U^{\pm 1}$
contains an admissible subword $U'$ with the base
$\lel_jL_jP_jR_j\rer_j$ for some $j$. Consider the computation
$U'\bullet h=V'$. By Lemma \ref{height} applied to the
$L_jP_j$-sector $U''$ of $U'$, $|h_{(4)}|$ is bounded by $|U|+|V|$
since $L_j(U''\circ h_{(4)})$ is empty.

The same is true for $|h'_{(4)}|$. Hence $|U\circ h_{(4)}|$ and
$|U\circ h_{(4)}h_{(34)}h_{(3)}|$ are bounded by a linear function
in $|U|+|V|$ (the application of every rule can increase the
length of an admissible word with a given base only by a
constant). Now again by Lemma \ref{height} $|h_{(3)}|$ is bounded
by a linear function in $|U|+|V|$.\endproof

\begin{lemma} \label{-(34)}  Let $U\bullet h=V$ be a wild computation of $\sss$. Suppose
that $h$ does not contain rules from $\sss(34)$.  Then $|h|$ is
bounded by a linear function in $|U|+|V|$.
\end{lemma}

\proof Lemma \ref{Ring Br Hist} (and the absence of (34)-letters
in the brief history $B$ of $h$) shows that $\br(h)$ does not have
any of the six prohibited words mentioned in Lemma \ref{Ring Br
Hist} and also $(34)(3)(34)$. Therefore one can apply the argument
of Lemma \ref{Br Hist}, and conclude that $\br(h)$ must be
obtained by removing some letters in one of the following four
words:

\begin{enumerate}
\item[(i)] $(4)(45)(5)(45)(4)$,
\item[(ii)] $(3)(23)(2)(12)(1)(12)(2)(23)(3)$,
\item[(iii)] $((3)(23)(2)(12)(1)(51)(5)(45)(4))^{\pm 1}$,
\item[(iv)] $(4)(45)(5)(51)(1)(51)(5)(45)(4)$
\end{enumerate}

{\bf Case (i)}. In this case, as above, we have the following
decomposition: $$h=h_{(4)}h_{(45)}h_{(5)}h'_{(45)}h'_{(4)}$$ where
some of the factors can be empty.

Suppose that the base of $U$ does not contain $L$-letters. Since
rules in $\sss(4)$ are not active with respect to $zz'$-sectors
where $z,z'\ne L_j$, every rule in $\sss(4)$ fixes $U$. Therefore
the words $h_{(4)}$ and $h'_{(4)}$ are empty (otherwise we get a
contradiction with the assumption that the computation $U\bullet
h=V$ is a ring computation (rings are incompressible) , so all
words $U\circ h[t]$, $t=0,1,..,|h|$ are supposed to be different).

If the base of $U$ contains $L_j$ then it contains $\lel_jL_j$ by
Lemma \ref{noPivP}. By Lemma \ref{height}, this gives a linear
upper bound for the lengths of $h_{(4)}$ and $h'_{(4)}$ in terms
of the lengths of $U$ and $V$  by the locking argument for
$\sss(34)$.

If the base of $U$ contains $R_j$ then it contains a subword
$R_j\rer_j$ and the length of $h_{(5)}$ is also linearly bounded
from above by Lemma \ref{height}. Otherwise $|h_{(5)}|=0$ because
the computation is absolutely reduced. Thus all five subwords in
the decomposition of $h$ have linearly bounded length.

{\bf Case (ii)}. By Lemma \ref{refined1}, we can turn $U\bullet
h=V (\mod \gr)$ into a free computation $U\bullet h=V'$ where
$V'=V (\mod \gr)$. All sectors in $V'$ but the $P_1R_1$-sector are
freely equal to the corresponding sectors in $V$. Using Lemma
\ref{height} to the computations $z(U)\bullet h=z(V')$,
$z\in\{\lel_j,L_j,R_j\}$, we linearly bound the lengths of maximal
subwords of $h$ over $\sss(g)$ ($g=1,2,3$) in terms of
$|z(U)|+|z(V')|\le |U|+|V|$. For example, if $g=1$, $z=L_j$, etc.

{\bf Cases (iii),(iv)} are similar to cases (1) and (ii).
\endproof

\begin{lemma} \label{Non-tame-1}
Let $U\bullet h=V$ be a wild ring computation of $\sss$. Then
$|h|$ is bounded by a linear function in $|U|+|V|$.
\end{lemma}

\proof The statement follows from lemmas \ref{(34)} and
\ref{-(34)}. \endproof

We call an admissible word $U$ for $\bsss$ (resp. $\sss$)
\label{acceptedw}{\em accepted} by $\bsss$ (by $\sss$) if there
exists a ring computation of the form $U\circ
h=\Sigma^s(w)K_1(\emptyset,1)$, $s\ne 0$, of $\bsss$  (resp. of
$\sss$ mod $\gr$) for some (for some positive) word $w$.

\begin{lemma} \label{To Sigma1}
(1)If a tame admissible word $U$ whose base starts and ends with
the same letter is not accepted by $\sss$ then the length of the
brief history of any ring computation $U\bullet h=V$ of $\sss$ is
at most $|\sss|+20$.

(2)
  If $U$ is accepted by $\sss$, then the length $d$ of any ring
  computation
 $$U\circ h=\Sigma^s(w)K_1(\emptyset,1),$$ such that
  $U\circ h[d']\ne\Sigma^s(w')K_1(\emptyset,1)$ for a word $w'$
  and $d'<d$, is recursively bounded in terms of $|U|$.

\end{lemma}

\proof (1) Suppose that the word $U$ is not accepted but the
length of the brief history of some ring computation $U\bullet h$
is $>|\sss|+20$. Then by Lemma \ref{Br Hist} $\br(h)$ contains an
occurrence of $(34)$ and one of the letters $(1), (12), (2), (23),
(45), (5), (51)$. Since the computation is tame, the base of $U$
is $\tsigma^sK_1$ for some $s\ne 0$. Then by Lemma \ref{Br Hist},
the $\br(h)$ contains two occurrences of $(34)$, and an occurrence
of $(12)$ or $(51)$ between them: $h=h_1\tau_1h_2\tau_2h_3$ where
$\tau_1\in\sss(34)$, $\tau_2\in\sss(12)\cup\sss(51)$. Then by
Lemmas \ref{Determ} and Lemma \ref{j=1}, we conclude that all
sectors of $U\circ h_1\tau_1h_2\tau_2$ are determined by
$\tau_1h_2$, so they are copies of each other (or, in the case of
$P_1z$-sectors, copies modulo $\gr$). Therefore, depending on
whether $\tau_2\in \sss^+(12)$ or $\tau_2\in \sss^-(12)$, the word
$U\circ h_1\tau_1h_2$ or $U\circ h_1\tau_1h_2\tau_2$ has the form
$\Sigma_{\emptyset,1}(\emptyset,u,v,\emptyset)K_1(\emptyset,1)$
(see Definition \ref{sigmari}) where $u$ and $v$ are positive
words by Lemma \ref{Posit Band}. Let us denote this word by $U'$.

By Lemma \ref{sigma}, part  a) (for $i=1$), there exists a ring
computation $$U'\bullet h' =\Sigma_{\emptyset,1}(\emptyset,
\emptyset, uv,
\emptyset)^sK_1(\emptyset,1)=\Sigma(uv)^sK_1(\emptyset,1)$$ of
$\sss(1)$ and, by lemmas \ref{W} and \ref{Band}, the corresponding
ring can be united with the ring which corresponds to $U\bullet
h_1\tau_1h_2$ or $U\bullet h_1\tau_1h_2\tau_2$. Now by Lemma
\ref{ringc}, the computation $U\bullet h_1\tau_1h_2h'$ or
$U\bullet h_1\tau_1h_2\tau_2h'$ (or another computation if the
corresponding annular diagram is compressible) is a ring
computation accepting $U$, a contradiction.

(2) First we need to check if $U$ is tame. If it is not tame, it
cannot be accepted. So assume that $U$ is tame. By the proof of
part (1) (and the definitions of compressible and reduced
diagrams), in order to check if a tame word $U$ is accepted, one
needs to check only computations $U\bullet h$ whose brief history
does not contain two occurrences of $(34)$ and an occurrence of
$(12)$ or $(51)$ between them. By Lemma \ref{Br Hist},
$|\br(h)|\le 19$. If for some $z\in \{\lel_j, L_j, R_j\}$, $j\ne
1$, and some $t\le |h|$, the length $|z(U\circ h[t])|$ exceeds
$c+|U|$ then by Lemma \ref{+sc}, part (b), one of the words
$\lel_j(U\circ h[d])$, $L_j(U\circ h[d])$ or $R_j(U\circ h[d])$ is
not empty for every $d\ge t$. Since the words
$\lel_j(\Sigma^s(w)K(\emptyset,1))$,
$L_j(\Sigma^s(w)K(\emptyset,1))$ or
$R_j(\Sigma^s(w)K(\emptyset,1))$ are empty, we can consider only
computations where $|z(U\circ h[t])|\le c+|U|$. By Lemma
\ref{height}, the length of every subword of $h$ over $\sss(g)$,
$g\in \{1,2,3,4,5\}$ is bounded by $2|U|+c$. This gives a bound on
the length of $h$.
\endproof

\begin{lemma}\label{1sttype} There exists an algorithm to check, given two admissible words $U,
V$ for $\sss$, if there exists a ring computation $U\bullet h=V$
of $\sss$. The length of arbitrary such ring computation is
recursively bounded in terms of $|U|+|V|$.
\end{lemma}

\proof By lemmas \ref{Non-tame-1} and \ref{ring321}, we have a
linear bound for any wild ring computation of $\sss$ which can
connect $U$ and $V$.

Let us assume that there is a tame ring computation connecting $U$
and $V$. By Lemmas \ref{Tame Ring} and \ref{ring321}, we can
assume that $U$ and $V$ have the form (\ref{6.3}) for some $s\ne
0$. Let us use the notation from Lemma \ref{Tame Ring}. Since the
word problem in $\gr$ is decidable, $u_{1+Nl}$, $l=1,...,s$, is a
copy of $u$ modulo $\gr$, we can assume that $u_{1+Nl}$ is a copy
of $u$ in the free group. (Clearly if $U'=U (\mod \gr)$, then a
computation $U\bullet h=V$ exists if and only if a computation
$U'\bullet h =V$ exists, and the corresponding rings differ only
in a $\gr$-cell attached to the boundary of one of them.)
Furthermore we may suppose that the words $U$ and $V$ are accepted
by $\sss$ (see the proof of Lemma \ref{To Sigma1}), so we can
assume by Lemma \ref{To Sigma1}, part (2), that $U\equiv
\Sigma(w)^sK_1(\emptyset,1)$, $V\equiv
\Sigma(w')^sK_1(\emptyset,1)$.

Notice that for every rule $\tau\in\sss$ and every admissible word
$U'$ with the base $W_j=\lel_jL_jP_jR_j\rer_j$ which is in the
domain of $\tau$, $\delta(U')=\delta(U'\circ \tau)$ for the
homomorphism $\delta$ introduced before Lemma \ref{alpha}. By
Lemma \ref{Emb}, we have similar equality for a ring computation
(which compute modulo $\gr$) of arbitrary length. Therefore if
$\delta(w')\ne \delta(w)$, there could be no computation of $\sss$
connecting $U$ and $V$.

Assume that $\delta(w')=\delta(w)$. Then a ring computation
$U\bullet h=V$ of $\sss$ does exist by lemmas \ref{Positive} and
\ref{Special}. Since the equality $\delta(w)=\delta(w')$ is
effectively verifiable (the word problem in $\bgr$ is decidable),
we can effectively verify (using Lemma \ref{ring321}) if there
exists a ring computation of $\sss$ connecting $U$ and $V$. By
Lemma \ref{Band}, all rings of the first type corresponding to
computations $U\bullet h'=V$ for some $h'$, have boundary labels
$U$ and $V$. (These boundary labels have no $x$-letters by Lemma
\ref{Band} since  $U\equiv \Sigma(w)^sK_1(\emptyset,1)$, $V\equiv
\Sigma(w')^sK_1(\emptyset,1)$.) It follows from the definition of
compressible and reduced diagrams, that the heights of all of
these rings must be equal. This height is recursively bounded.
Indeed, the above construction recursively bounds the number of
${\cal S}(34)$-rules in any ring computation connected $U$ and
$V$, because there are no $(\theta,k)$-cells which are higher than
those corresponding to rules from $\sss(34)$. Then arguing as in
the proof of Lemma \ref{To Sigma1}, we get a recursive upper bound
for the length of the history. The lemma is proved.
\endproof

\subsection{Computations of $\bar\sss$} \label{rost}
\setcounter{equation}{0}

\begin{lemma} \label{Equations}
For any $s \ge 1$, and any two admissible words $U, V$ of $\bsss$
there is an algorithm to decide whether there exists a computation
$U\bullet h =V$ of $\bsss$, and $\br(h)$ is of length $\le s$. The
lengths of all ring computations of this form are recursively
bounded by a function of $|U|+|V|$.
\end{lemma}

\proof If $U$ and $V$ have different bases then there is no
computation connecting $U$ and $V$. So we can assume that $U$ and
$V$ have the same bases.

We can assume that $h$ is reduced and write it as $h\equiv
h_1...h_s$ where $h_i$ are words over $\bsss(\omega_i)$,
$\omega_i\in \{1, 12,\dots 51\}$, $i=1,...,s$. There are only
finitely many choices for the sequence $\omega_1,...,\omega_s$. So
we can fix one of them. We can also assume that if
$\sss(\omega_i)$ consists of transition rules then $h_i$ is of
length $1$, so there are finitely many choices for such $h_i$. Fix
one of these choices.

By Lemma \ref{bsss}, for each of the other $h_i$, $i=1,...,s$,
there exist two sequences of numbers $\{\varepsilon_{i,z},
\delta_{i,z}\in \{-1,0,1\}\ |\ z\in \tkk\}$, such that for every
admissible word $W$ of the form $V(g,r)$, and any group word $w_i$
in $\{a_1,...,a_{\bar m}\}$ there exists a word $h_i$ in
$\bsss(g)$ satisfying the following properties:

\begin{enumerate}
\item $h_i$ is a copy of $w_i$;
\item $W$ is in the domain of $h_i$;
\item $W\circ h_i$ is obtained from $W$ by all replacements
of the form $$\bar z(r,g)\to w_i(z_-)^{\varepsilon_{i,z}} \bar
z(r,g)w_i(z)^{\delta_{i,z}},$$ where $z$ runs over all letters in
the base of $W$, $\varepsilon_{i,z}\in\{-1,1\}$.
\end{enumerate}

Thus, by Lemmas \ref{csss} and \ref{bsss}, applying $h$ to $U$
results in multiplying the interior of each sector of $U$ on the
left and the right by products of copies of unknown words
$w_i(z)$, their inverses, and constant words corresponding to the
transition rules. The result of this multiplication should be the
interior of the corresponding sector in $V$. Hence the existence
of $h$ such that $U\circ h=V$ is equivalent to the existence of a
solution of one of a finite systems of equations in the free group
(each system corresponds to the choice of the brief history of $h$
and the choice of transition rules in $h$; the equations
correspond to the ages in the brief history and the unknowns are
the words $w_i(z)$).

It remains to apply a theorem of Makanin \cite{Makanin} saying
that there exists an algorithm checking if a system of equations
over a free group has a solution. The heights of all ring
computations of this form must be equal by the definition of
compressible and reduced diagrams. This height is recursively
bounded by Lemma \ref{Delta(1)}, because the above construction
gives the list of  computations (if any exists) with arbitrary
distribution of at most $s$ transition rules; and this allows us
to effectively select those whose rings are incompressible. The
lemma is proved.
\endproof

\begin{lemma} \label{To Sigma}
(1) If a tame word $U$ is not accepted by $\bsss$ then the length
of the brief history of any ring computation $U\bullet h=V$ of
$\bsss$ is at less than 19 with  at most one occurrence of (34).

(2) There is an algorithm which checks whether a given admissible
word is accepted by $\bsss$, this algorithm also finds an
accepting computation if it exists.
\end{lemma}

\proof (1) Let the base of $U$ be $\tsigma^lK_1$ for some $l\ne
0$. Then using Lemma \ref{Br Hist} as in the proof of Lemma
\ref{To Sigma1}, if the length of the brief history of a
computation $U\bullet h$ is bigger than $19$ or if the number of
occurrences of (34) in the brief history is at least 2, there
exists $t$ between $1$ and $|h|$ such that $$U\circ
h[t]=\bar\Sigma_{\emptyset,1}(\emptyset, u, v, \emptyset)^l\bar
K_1(\emptyset,1)=\Sigma_{\emptyset,1}(\emptyset, u, v,
\emptyset)^lK_1(\emptyset,1).$$ By Lemma \ref{Sigma}, part b),
there exists a computation $U\bullet
h[t]h'=\Sigma(uv)^lK_1(\emptyset,1).$ By Lemma \ref{Delta(1)},
this computation (or another computation if the corresponding
annular diagram is not reduced) is a ring computation, so $U$ is
accepted. This contradiction proves the first part of the lemma.

(2) By Part 1 of the lemma and Lemma \ref{Br Hist}, we can
restrict ourselves to computations $U\bullet h$ with $\br(h)$
containing at most one occurrence of $(34)$ and $|\br(h)|\le 19$.

Suppose that $U$ is accepted, $U\circ h=\Sigma(w)^lK_1$ be the
corresponding computation. As in Lemma \ref{Posit Band}, introduce
a homomorphism $\gamma$ of the free group with basis
$\bar\aaa\cup\bar\kkk$ onto free group over $\{a_1,\dots,a_{\bar
m}\}$, by the rule $\bar a_i(z)\mapsto a_i$, $z\to 1$, $z\in
\bar\kkk$. Then for every $t=0,...,|h-1|$, $\gamma(U\circ
h[t])=\gamma(U\circ h[t+1])$ if $h[t+1]$ does not end with a rule
from $\sss(34)$. Otherwise
$$||\gamma(U\circ h[t])|-|\gamma(U\circ h[t+1])||\le cl.$$
Hence the length of $w^l=\gamma(\Sigma(w)^lK_1(\emptyset,1))$ is
bounded by $|U|+cl$ (since $h$ contains only one rule from
$\bar\sss(34)$). Hence $|w|\le |U|+cl$. Therefore
$|\Sigma(w)^lK_1(\emptyset,1)|$ does not exceed $l(|U|+lc)+4N$.
Now Lemma \ref{Equations} gives us an algorithm to check if an
accepting ring computation exists.
\endproof

\begin{lemma} \label{Non-tame} Assume $U\bullet h$ be a wild ring computation of $\bar\sss$.
Then $|\br(h)|\le 13$.
\end{lemma}

\proof We examine several cases taking into account Lemma
\ref{Base}. Let $B$ be the base of $U$.

{\bf Case 1.} Suppose that $B$ contains letter $P_j^{\pm 1}$ for
some $j\ne 1$.

{\bf Case 1.1}. Suppose that $B$ contains no letters except
$P_j^{\pm 1}$. Since rules from $\bar\sss(\omega)$, $\omega\in \{
(23), (3), (34), (4), (45)\}$ lock $P_jR_j$-sectors, these rule
cannot occur in $h$. Rules from  $\bsss(2)\cup\bsss(5)$ are not
active for $P_jP_j\iv$- and $P_j\iv P_j$-sectors, so existence of
these rules in $h$ would contradict the assumption that we have a
ring computation. Hence $\br(h)$ does not contain $(2)$ and $(5)$.
This implies that $\br(h)$ cannot be longer than the words
$(51)(1)(12)$, $(51)(1)(51)$ and $(12)(1)(12)$, $|\br(h)|\le 3$.

{\bf Case 1.2.}. Suppose that $B$ contains a subword $yz$ where
$y$ (or $z$) is equal to $P_j^{\pm 1}$ but $z$ (resp. $y$) is not
$P_j^{\pm 1}$. We can assume that the exponent of $P_j^{\pm 1}$ in
$yz$ is 1 (otherwise we can switch from $U$ to $U\iv$). Therefore
by Lemma \ref{Base}, $yz$ is equal to either $L_jP_j$ or $P_jR_j$.

{\bf Case 1.2.1}. Let $yz\equiv L_jP_j$, $j\ne 1$. Then by Lemma
\ref{Base}, since $B$ starts and ends with the same letter, there
are two possibilities: either $B$ has a subword $L_jP_jP_j\iv$ or
it contains $L_jP_jR_j$.

{\bf Case 1.2.1.1}. Suppose there is a subword $L_jP_jP_j\iv$ in
$B$ for $j\ne 1$. Then there are no letters from $\{(34), (4),
(45)\}$ in $\br(h)$ since the corresponding rules lock
$P_jR_j$-sectors. The word $B$ also has no subwords $(23)(2)(23)$
by Lemma \ref{(12)(2)(12)}. If $B$ contains a subword
$L_j^{-1}L_jP_jP_j^{-1}$, then $\br(h)$ cannot have letters $(12)$
or $(51)$ since the corresponding rules lock $\lel_jL_j$-sectors.
If $B$ has a subword $\lel_jL_jP_jP_j^{-1}$, then $\br(h)$ cannot
have subwords $(12)(2)(12)$, by Lemma \ref{(12)(2)(12)}. Thus, in
any case, subwords $(12)(2)(12)$ do not occur in $\br(h)$.
Consider two subcases.

{\bf Case 1.2.1.1.1}. Suppose that $B$ does not contain $R_p$ for
$p\ne 1$. Then $\br(h)$ contains neither $(3)$ nor $(5)$ because
$U\bullet h$ is a ring computation. This implies that $\br(h)$
cannot be longer than the words $(23)(2)(12)(1)(12)(2)(23)$ or
$(51)(1)(12)(2)(23)$, $|\br(h)|\le 7$.

{\bf Case 1.2.1.1.2}. Suppose that $B$ contains $R_p$, $p\ne 1$.
If $B$ has no subwords $R_p\rer_p$, then letters
$(1),(12),(2),(23),(51)$ do not occur in $\br(h)$ (the locking
argument again). Recall that by Case 1.2.1.1, (45) does not occur
in $\br(h)$.

If $B$ has a subword $R_p\rer_p$, $p\ne 1$, then $\br(h)$ cannot
contain a subwords $(51)(5)(51)$ or $(23)(3)(23)$ by Lemma
\ref{(12)(2)(12)}.  In view of the previous restrictions, the
longest possible brief history cannot be longer than the word
$(3)(23)(2)(12)(1)(12)(2)(23)(3)$, $|\br(h)|\le 9$.

{\bf Case 1.2.1.2}. Suppose that there exists a subword
$L_jP_jR_j$ in $B$, $j \ne 1$. Then $\br(h)$ has no subwords
$(23)(2)(23)$, $(34)(3)(34)$, $(34)(4)(34)$ and $(45)(5)(45)$ by
Lemma \ref{(12)(2)(12)}. The subwords $(12)(2)(12)$ and
$(45)(4)(45)$ can be excluded as in Case 1.2.1.1. The subwords
$(23)(3)(23)$ and $(51)(5)(51)$ are impossible by a similar reason
(but one has to consider the alternative whether $R_j\rer_j$
occurs in the base or not).

The restrictions we have obtained show that there is one of the
letters $(12), (51)$ in $\br(h)$ if the length of the brief
history is at least 8. Therefore $B$ has a subword
$\lel_jL_jP_jR_j$ by the locking argument. Similarly, if the
length of $\br(h)$ is at least 6, it has one of the letters $(12),
(23),(51)$, and so $B$ contain subword $\lel_jL_jP_jR_j\rer_j$.
Continuing the extension of the subword of the base using the
locking argument and Lemma \ref{Base}, we see that the computation
is tame if the length of $\br(h)$ is at least 8, contrary the
assumption of the lemma.

{\bf Case 1.2.2}. Let $yz\equiv P_jR_j$. Then by Lemma \ref{Base},
one can assume that $B$ has a subword $P_j^{-1}P_jR_j$ (otherwise
we have Case 1.2.1.2). Hence $\br(h)$ has neither of the letters
$(23), (3), (34)$ (the locking argument). This brief history
cannot contain subwords $(51)(5)(51)$ by Lemma \ref{(12)(2)(12)}.

If there are no letters $L_p$ in $B$ for $p\ne 1$, then there are
no letters $(2)$, $(4)$ in $\br(h)$ because rules from
$\bsss(2)\cup\bsss(4)$ are not active with respect to
$zz'$-sectors if $z,z'\ne L_p$. Otherwise there are no subwords
$(12)(2)(12)$ or $(45)(4)(45)$ in $\br(h)$ by Lemma
\ref{(12)(2)(12)}. Therefore such subwords are absent in any case.
Hence $\br(h)$ cannot be longer than words
$(2)(12)(1)(51)(5)(45)(4)$, $(2)(12)(1)(12)(2)$,
$(4)(45)(5)(51)(1)(51)(5)(45)(4)$, $|\br(h)|\le 9$.

{\bf Case 2}. Suppose now that $B$ does not contain $P_j$ for
$j\ne 1$. Since the computation is a ring computation, $\br(h)$
does not contain $(1)$.

If there are nether $L_j$- nor $R_j$-letters in $B$ for $j\ne 1$,
then  $\br(h)$ contains no letters from $\{(1),...,(5)\}$.
Therefore $\br(h)$ must be a subword of a product of several
factors of the form $(12)(23)(34)(45)(51)$ or
$(51)(45)(34)(23)(12)$ because the history of the computation is
reduced. However $\br(h)$ cannot contain two occurrences of $(12)$
or $(51)$ otherwise the computation is not a ring computation.
Indeed the transition rules not from $\sss(34)$ change only
$(\bee,\Omega)$-coordinates of the words, rules from $\sss(34)$
also change only $(\bee,\Omega)$-coordinates in admissible words
without $L_j$ in their bases. Hence if $\br(h)$ contains two
occurrences of $(12)$ or $(51)$, the corresponding ring is not
reduced. Thus, the length of $\br(h)$ is at most 13 in this case.
Hence, one may assume further that there is either $L_j$- or
$R_j$-letter in $B$ for some $j\ne 1$.

{\bf Case 2.1}. Suppose that $B$ contains $L_j$, $j\ne 1$. Then
the locking argument shows that $\br(h)$ contains no letters from
$\{(23), (3), (34)\}$. Then the standard application of Lemma
\ref{(12)(2)(12)} shows that $\br(h)$ cannot contain either
$(12)(2)(12)$ or $(45)(4)(45)$.

If $B$ has no letter $R_j$ with $j\ne 1$ then $\br(h)$ has no
letter (5) (otherwise the computation is a ring computation since
the corresponding diagram is compressible). Hence by Lemma
\ref{(12)(2)(12)}, $\br(h)$ has no subwords $(51)(5)(51)$ whether
$B$ contains $R_j$ or not. Hence $\br(h)$ cannot be longer than
the words $$(2)(12)(1)(51)(5)(45)(4), (2)(12)(1)(12)(2),
(4)(45)(5)(51)(1)(51)(5)(45)(4),$$ $|\br(h)|\le 9$.

{\bf Case 2.2}. Suppose that $B$ contains $R_j$, $j\ne 1$, but
$L_p$ does not occur in $B$ if $p\ne 1$.. Then $\br(h)$ cannot
contain letters $(34),(4), (45)$ since $P_jR_j$ is not a subword
of $B$. Besides $\br(h)$ has no letters $(2), (4)$ since otherwise
the computation is not a ring computation. The subwords
$(23)(3)(23)$, $(51)(5)(51)$ are impossible by Lemma
\ref{(12)(2)(12)}. Hence $\br(h)$ cannot be longer than the words
$(3)(23)(12)(1)(51)(5)$, $(3)(23)(12)(1)(12)(23)(3)$,
$(5)(51)(1)(51)(5)$, $|\br(h)|\le 7$.

The lemma is proved.\endproof

\begin{lemma} \label{2d type} There exists an algorithm to check, given two admissible words $U,
V$ for $\bsss$, if there exist a ring computation $U\bullet h=V$
of $\bsss$ . The length of $|h|$ is recursively bounded in terms
of $|U|+|V|$.
\end{lemma}

\proof We can assume that the bases of $U$ and $V$ are the same,
the first and the last letters of that base coincide.

First we check, using Lemma \ref{To Sigma}, if both words are
acceptable by $\bsss$.

{\bf Case 1.} Suppose that the answer is ``yes". Let $U\bullet
h_1=\Sigma(w_1)^sK_1(\emptyset,1)$ and $V\bullet
h_2=\Sigma(w_2)^sK_1(\emptyset,1)$ be accepting computations.

Let $\delta$ be the homomorphism $\hhh\to \bar\gr$ from Lemma
\ref{alpha}. By Lemma \ref{Special(bar)} and Lemma
\ref{simulatecsss}, if $\delta(w_1)$ and $\delta(w_2)$ are equal
in $\bar\gr$ then there exists a computation
$\Sigma(w_1)^sK_1(\emptyset,1)\bullet
h_3=\Sigma(w_2)^sK_1(\emptyset,1)$. Then $U\bullet h_1h_3h_2\iv
=V$ (or another computation if the corresponding annular diagram
is not reduced) is a ring computation (by Lemma \ref{Delta(1)}).

Suppose that $\delta(w_1)\ne \delta(w_2)$ in $\bar\gr$. Then by
Lemma \ref{alpha1} there exist no computations connecting
$\Sigma(w_1)^sK_1(\emptyset,1)$ and
$\Sigma(w_2)^sK_1(\emptyset,1)$. Hence there is no computation
connecting $U$ and $V$. Since the word problem in $\bar\gr$ is
decidable, we deduce that, in Case 1, we can check if $U$ and $V$
are connected by a computation of $\bsss$.

{\bf Case 2.} Suppose that $U$ or $V$ are not accepted by $\bsss$.
Then we can assume that both are not accepted (otherwise $U$ and
$V$ obvious  cannot be connected by a computation). Then by Lemmas
\ref{Non-tame} and \ref{To Sigma}, one has to check only
computations whose brief histories have lengths at most $13$ (if
$U$ is wild) or $19$ if $U$ is tame. This can be done by Lemma
\ref{Equations}.

In both cases we obtain a computation (if it exists) of
recursively bounded length. Thus the number of transition rules is
bounded for every ring computation which connect $U$ and $V$. Then
Lemma \ref{Equations} helps us to effectively select a desired
ring computation.
\endproof

\subsection{All ring computations} \label{allrings} \setcounter{equation}{0}

\begin{lemma} \label{I-II-I}
Let $U\bullet h=V$ be a ring computation of $\csss$. Then

(1) $h=h_1h_2h_3$ where $h_1, h_3$ are words in $\sss$, $h_2$ is a
word in $\bsss$.

(2) If $h_2$ is not empty, then  $|\br(h_1)|,|\br(h_3)|$ are
recursively bounded in terms of $|U|, |V|$.
\end{lemma}

\proof (1) Suppose $h=h'h_0h''$ where non-empty word $h_0$ is over
$\sss$, $h'$ ends and $h''$ starts with a rule from $\bsss$. Let
$U'=U\bullet h'$, and consider the computation $U'\bullet h_0$.
Denote by $\Delta$ the ring corresponding to the computation
$U\bullet h$. Notice that a rule of $\sss$ as well as a rule from
$\bar\sss$ can be applied to $U'$ and to $U'\circ h_0$. Hence
there are no $a$-letters in $U'$ except for $a(P_j)$-letters.
Therefore the subring $\Delta_0$ of $\Delta$, which correspond to
$U'\bullet h_0$ has no $\xxx$-letters by Lemma \ref{Band}.

Notice that by definition of $\bsss$ if in every rule $\tau$ of
$\sss$, we remove $a$-letters in subrules $z\to uz'v$ where $z\in
\{K_1, L_1, P_1, R_1\}$ and put $\bar{ }$ over every letter of the
rule, we get a rule $\bar\tau$ from $\bsss$. Moreover, if for any
admissible word $W$ in the domain of $\tau$, we remove $a$-letters
from $zz'$-sectors where $z,z'\in \{K_1, L_1, P_1, R_1, K_2\}$,
and add $\bar{ }$ to every letter, we get an admissible word $\bar
W$ for $\bsss$, and $$\bar W\circ \bar\tau=\overline{W\circ
\tau}.$$

Since the words $U'$ and $U'\circ h_0$ are in domains of rules
from both $\bsss$ and $\sss$, $\lel_jz$, $L_jz$-, $R_jz$-sectors
of these words do not contain $a$-letters, for any $j$, the
$\lel_1$-, $L_1z$-, $P_1z$-, $R_1z$-sectors do not contain
$a$-letters, the $a$-letters in $P_jz$-sectors are from
$\{a_1(P_1),...,a_m(P_1)\}=\{\bar a_1(P_1),...,\bar a_m(P_1)\}$,
and the $(\bee, \Omega)$-coordinates of these words are
$(\emptyset,1)$. Hence $\bar U'=U$ and $\overline{U'\circ
h_0}=U'\circ h_0$. Therefore we can replace $h_0$ by a word $\bar
h_0$ over $\bsss$ so that $U'\circ \bar h_0=U'\circ h_0$. By
Lemmas \ref{Delta(1)} and \ref{ringc} and the lack of the
$\xxx$-edges in the boundary of $\Delta_0$ , the computation
$U\bullet h'\bar h_0h''=V$ or another computation if the
corresponding ring diagram is not reduced, is a ring computation,
but some $\Theta$-letters of $h$ are replaced by
$\bar\Theta$-letters in $h'\bar h_0h''$. This contradicts the
definition of reduced and compressible diagrams, because
$(\Theta,k)$-cells are higher than $(\bar\Theta,k)$-cells.

 (2) Assume $h_2$ is non-empty and consider, for
instance, $\br(h_1)$. If $U$ is wild and there are no letters
$(34)$ in $\br (h_1)$, then the forms (i) - (iv) from the proof of
Lemma \ref{-(34)} show that $|\br(h_1)|\le 9$. If $U$ is wild, and
$\br(h_1)$ contains $(34)$ and one of the letter $(1), (12), (2),
(23), (45), (5), (51)$, then, by Lemma \ref {Est' sector}, we can
assume that $U^{\pm 1}$ contains an admissible subword with base
$\lel_jL_jP_jR_j\rer_j$ for some $j\ne 1$. If $U$ is wild, and
every letter of $\br (h_1)$ belongs to $\{(3),(34),(4)\}$, then
$|\br (h_1)|\le 5$ as this was shown in the proof of Lemma
\ref{(34)}. Thus, we may assume further that the base of $U^{\pm
1}$ contains an admissible subword with base
$\lel_jL_jP_jR_j\rer_j$ for some $j\ne 1$.

By Lemma \ref{Br Hist}, we may assume that  $\br(h_1)$ contains
each of the letters $(12),(23),(34)$, $(45),(51).$ Using the
locking argument as in Lemma \ref{(34) and (12)}, we conclude that
$U$ is tame. Then we may assume by Lemma \ref{To Sigma}(1), that
$U$ is accepted by $\sss$: $U\circ g_1 = \Sigma^s(w)K_1(0, 1)$ for
some ring computation $U\bullet g_1$ of $\sss$.

Let $U_1$ be an admissible subword of the tame word $U^{\pm 1}$
with base $K_1L_1P_1R_1K_2^{-1}$. Then the standard application of
the homomorphism $\delta$ (introduced before Lemma \ref{alpha}) to
$U_1\bullet g_1$ shows that $\delta (U_1)=\delta (w)$ in the group
$\bar\gr$. The same homomorphism applied to $U_1\bullet h_1$ shows
that $\delta(U_1)=1$ in $\bar\gr$ since the $\bar\sss$-rules of
$h_2$ lock $K_1L_1$-, $L_1P_1$-, $P_1R_1$- and
$R_1K_2^{-1}$-sectors. Hence the positive word $w$ represents the
identity of $\bar\gr$.

Now, by lemmas \ref{Positive} and \ref{Special}, there exists a
ring computation $(U\circ g_1)\bullet g_2=\tsigma$ of $\sss$. Set
$g=g_1g_2$, then the length of the ring computation $U\bullet
g=\tsigma$ (or of some other ring computation if the annular
diagram, corresponding to $U\bullet g$ is not reduced) is
recursively bounded as function of $|U|$. (Since we already know
that such a computation exists, and the number of its ${\cal
S}(34)$-rules is bounded, we get an upper bound for its length as
in the proof of Lemma \ref{To Sigma1}(2).)

The computation $U\bullet g g^{-1}h_1$ corresponds to a diagram
consisting of three reduced rings. It contains the subcomputation
$\tsigma\bullet(g^{-1}h_1)=U\circ h_1$. Since
$\tsigma\equiv\bar\tsigma$ and $U\circ h_1$ is in the domain of
$h_2$, we have, as in part (1) of the proof, a ring computation of
$\bar\sss$: $\tsigma\bullet f =U\circ h_1$. Hence the ring
computation $U\bullet h'$ can be replaced by the  computation
$U\bullet g(fh_2)h_3$ with recursively bounded length of $\br(g)$
where computation $(U\circ g)\bullet fh_2$ is  a computation of
$\bar\sss$). Therefore the number of rules from $\sss(34)$ is also
recursively bounded in $h_1$ since otherwise $\Delta$ were
compressible. Then, by Lemma \ref{Br Hist}, we have a recursive
upper bound for the length of $br(h_1)$. (Again $U\circ)\bullet
fh_2$ can be replaced by a ring computation if the corresponding
annular diagram is not reduced.) \endproof

According to Lemma \ref{I-II-I}, we will consider  ring
computations of the form $U\bullet h_1h_2h_3=V$ where $h_1,h_3$
are words over $\sss$, $h_2$ is a word over $\bsss$, and if $h_2$
is not empty, then $|\br(h_1)|, |\br(h_3)|$ are recursively
bounded in terms of $|U|$ and $|V|$. We fix this notation for the
next two lemmas.

\begin{lemma} \label{Short}
Assume the word $h_2$ is non-empty. Then the words $h_1$ and $h_3$
have lengths bounded by a recursive function of $|U|$ and $|V|$
respectively.
\end{lemma}

\proof Consider the computation $U\bullet h_1$ (the computation
$V\bullet h_3\iv$ is similar).

{\bf Case 1.} First suppose that the base of the computation
contains the $j$-standard base $\lel_jL_jP_jR_j\rer_j$, $j\ne 1$.
Then we can consider the restriction of the computation to a
subword of $U$ with the $j$-standard base. We know that the length
of $\br(h_1)$ is recursively bounded. Therefore the number of
occurrences of $(34)$ in $h$ is also bounded. Notice also that
since $U\circ h_1$ is in the domain of a rule of $\bsss$, the
$\lel_jL_j$-, $L_jP_j$, $R_j\rer_j$-sectors of $U\circ h_1$ do not
contain $a$-letters. Therefore, by Lemma \ref{+sc}(b),
$|\lel_j(U\circ h_1[t])|+|L_j(U\circ h_1[t])|+|R_j(U\circ
h_1[t])|$, $t=1,...,|h_1|$ is bounded by a recursive function in
$|U|$. Then Lemma \ref{height} gives a recursive upper bound for
lengths of arbitrary maximal subword of $h_1$ over $\sss(g)$,
$g=1,...,5$. Since the number of such subwords is recursively
bounded, we get a recursive bound for $|h_1|$ in terms of $|U|$
and the claim is proved.

{\bf Case 2.} Assume now that the base of $U^{\pm 1}$ does not
contain $W_j=\lel_jL_jP_jR_j\rer_j$ as a subword for any $j\ne 1$.

Notice, however, that since $U\circ h_1$ is in the domain of a
rule from $\sss$ and a rule from $\bsss$, the word $U\circ h_1$
does not contain letters from
$\aba{K_j}\cup\aba{L_j}\cup\aba{R_j}$ since such letters does not
belong to $\aaa\cap\bar\aaa$ . Therefore, by Lemma \ref{noPivP},
every $\lel_j$-, $L_j$-, $P_j$-, $R_j$- and $\rer_j$-letter in the
base of $U$ must occur in subwords of the form $(\lel_jL_j)^{\pm
1}$, $(R_j\rer_j)^{\pm 1}$, $(L_jP_j)^{\pm 1}$.  Hence, by Lemma
\ref{Base} the base of the corresponding ring $\Delta$ must
contain a letter $L_j$ or a letter $R_j$ for $j\ne 1$.

Let us consider only the first case, because one can apply the
same argument to the second one. In this case, subword
$(P_jR_j)^{\pm 1}$ cannot be a neighbor in the base to the
occurrence of $L_j^{\pm 1}$, since otherwise we get a subword
$W_j$ in the base of $U$. Consequently, $\br(h_1)$ does not
contain the letter $(34)$ (because rules from
$\sss(34)\cup\bsss(34)$ lock $L_jP_j$- and $P_jR_j$-sectors) and
the base of $U^{\pm 1}$ contains a subword
$W_j'=\lel_jL_jP_jP_j\iv L_j\iv(\lel_j)\iv$. Let $U'$ be a subword
of $U$ with base $W_j'$.

Consider the computation $U'\bullet h_1$. The word $U'\circ h_1$
has the form
$$(\lel_jL_jP_jvP_j^{-1}L_j^{-1}\lel_j^{-1})(r,i)$$
($(r,i)$ are the $(\bee,\Omega)$-coordinates), where $v$ is a word
in $\{a_1(P_j),...,a_m(P_j)\}$. (Recall that $a_i(P_j)$ is
identified with a $\bar a_i(P_j)$ only if $i\le m$.)

Let $\gamma$ be the homomorphism defined in the proof of Lemma
\ref{Posit Band}. Then $|\gamma(v)|\le |U'|$.  This means that we
have a linear bound of $|U'\circ h_1|$ in terms of $|U'|$. Since
this argument works for every subword of $U^{\pm 1}$ with base
$W_j'$, $j\ne 1$, and $U\circ h_1$ contains no $a$-letters except
letters from $\aaa(P_j)\cup\bar\aaa(P_j)$, $j\ne 1$, we obtain a
linear bound for $|U\circ h_1|$ in terms of $|U|$. It remains to
apply Lemma \ref{Non-tame-1}.
\endproof

\begin{lemma} \label{Main}There is an algorithm to check whether, given two admissible words
$U,V$ for $\csss$,  there exists a ring computation $U\bullet
h=V$. The length of the history $h$ is recursively bounded in
terms $|U|, |V|$ for arbitrary such ring computation.
\end{lemma}

\proof By Lemma \ref{1sttype}, we can assume that there are no
computations of $\sss$ connecting $U$ and $V$. By Lemma \ref{2d
type}, we can assume that there is no ring computation of $\bsss$
connecting $U$ and $V$. (Recall that all ring computations
$U\bullet f =V$ of $\csss$ must have equal length $|f|$ depending
on $U$ and $V$, by the definition of a reduced diagram.)

Hence, by Lemma \ref{I-II-I}, we have to analyze only the ``mixed"
situation when $h_2$ is not empty and $h_1$ or $h_3$ is not empty
(see notations introduced before Lemma \ref{Short}).

The sum $|h_1|+|h_3|$ is recursively bounded by Lemma \ref{Short}.
Therefore the lengths $|U\circ h_1|$ and $|U\circ h_1h_2|$ are
recursively bounded modulo $\gr$-relations (i.e., if we substitute
the subwords in ${\cal A}(P_1)$ by equal modulo $\gr$ shorter
subwords). Therefore by Lemma \ref{2d type}, we can recursively
bound $|h_2|$. Hence we can recursively bound the sum
$|h|=|h_1|+|h_2|+|h_3|$. \endproof

\begin{lemma} \label{Main1} Let $\Delta$ be any annular diagram over $\Delta$ whose contours
$p, q$ do not contain $\theta$-edges. Then there exists another
annular diagram $\Delta'$ with the same boundary labels and
recursively bounded number of cells (in terms of $|p|+|q|$).
Moreover,  $\Delta'$ and $\Delta$ have equal histories if $\Delta$
is a ring.
\end{lemma}

\proof By Lemma \ref{ring123}, we can assume that $\Delta$ is a
ring or a quasiring. Lemma \ref{quasiring} treats the quasiring
case, so we can assume that $\Delta$ is a ring.  By Lemma
\ref{Main}, we can bound the size of the corresponding ring
computation. Finally Lemma \ref{ring321} bounds the number of
cells in $\Delta$.\endproof

\section{Spirals}
\setcounter{equation}{0}
 \label{spirals}

We need the following properties of words in a free group $F$.

\begin{lemma} \label{powers}
For arbitrary elements $u,v,w$ of $F$ and any integer $d\ge 0$,
the length of arbitrary product $u^{t}wv^{t}$ in $F$ is not
greater than $|u|+|v|+|w|+|u^{d}wv^{d}|$ provided $0\le t\le d$.
Furthermore, either equality $u^twv^t=w$ is true in $F$ for every
$t$ or $|u^twv^t|\ge t -(|u|+|v|+|w|)$.
\end{lemma}

\proof Since, for any word $U$, a power $U^s$ is freely equal to
$YU_0^sY^{-1}$ where $U_0$ is the cyclically reduced form of $U$
and $|Y|\le |U|$, it suffices to assume that the words $u$ and $v$
are cyclically reduced. Moreover one of them can be assumed to be
non-empty, and words $u$ and $v$ can be assumed to be not equal to
proper powers of some non-trivial elements of $F$.

Of course, there may be cancellations in the words $u^{t}wv^{t}$.

Suppose first, that the cancellations are not substantial in the
following sense: the length of the suffix (resp. prefix)  of
$u^tw$ (resp $wv^t$) that is cancelled in the product $u^twv^t$
does not exceed $|w|+|u|+|v|$ for every $t$. Then the sequence of
lengths $|u^twv^t|$ increases linearly with the slope $\ge
|u|+|v|\ge 1$ for $t\ge |w|+|u|+|v|$. Hence $|u^{t}wv^{t}|\le
|u|+|v|+|w|+|u^dwv^d|$, as desired.

Now suppose that there exists a substantial cancellation in the
above mentioned sense. Then both words $u^t$ and $v^{-t}$ contain
a common subword $T$ with $|T|\ge |u|+|v|$. This implies (see,
e.g. \cite{FW}) that $u$ and $v^{-1}$ coincide up to a cyclic
shift: $v^{-1}\equiv u_2u_1$ and $u\equiv u_1u_2$ for some words
$u_1, u_2$. Furthermore the word $wu_1^{-1}$ must be freely equal
to to a power $u^l$ for some $l$. Therefore $u^twv^{t}=u^tu^l
u^{1-t}u_2^{-1}=u^{l+1}u_2^{-1}= wu_1^{-1} u_1u_2u_2^{-1}=w$ in
$F$ for arbitrary $t$. The lemma is proved.
\endproof

By Britton's lemma, every $x$-letter has infinite order in the
group $\hhh_2$. It is also convenient to introduce auxiliary
letters $x^{1/4}, x^{1/16},\dots$ for every letter $x\in \xxx$. In
other words, we embed our group $\hhh_2$ into the multiple
amalgamated product of $\hhh_2$ and infinite cyclic groups
$C_{x,i}=\la x^{\frac{1}{4^i}}$, $x\in \xxx\ra, i=0, 1, 2,...$
where $C_{x,0}=\la x\ra$, $x\in \xxx$, and for every $i\ge 1$,
$C_{x,i}$ is identified with the subgroup of $C_{x,i+1}$ generated
by $(x^{\frac{1}{4^{i+1}}})^4$. Letters $x^{\frac{1}{4^n}}$ will
be called \label{fractional}{\em fractional letters}, words in
these letters will be called {\em fractional words}.

\begin{lemma} \label{2 bands}
Let $\Delta$ be a quasitrapezium  with the standard decomposition
of the boundary $p_1p_2p_3p_4$. Assume that the height of $\Delta$
is 2, $\Delta$ has no $\gr$-cells and $w_1$, $w_3$ are maximal
subwords in the alphabet $\xxx$ of $\phi(p_1)$, $\phi(p_3^{-1})$,
respectively. (Recall that $p_1=\topp(\bb)$,
$p_3^{-1}=\bott(\bb')$ as in the definition of trapezium or
quasitrapezium.) Then $w_3$ is freely equal to $uw'_1v$ where
$w'_1$ is obtained from $w_1$ by an injective  substitution of the
form $x\mapsto (x')^{4^c}$ ($x,x'\in \xxx$), the integer $c$ (same
for all $x\in \xxx$) and the fractional words $u$, $v$ depend
recursively only on the boundary labels of the maximal
$\theta$-bands $\ttt_1$, $\ttt_2$ of $\Delta$.
\end{lemma}

\proof Consider the subdiagram $\Gamma$ of $\Delta$ bounded by
$\topp(\ttt_1)$ and $\bott(\ttt_2)$. Let $\bar p_1\bar p_2 \bar
p_3 \bar p_4$ be its  boundary, where $\phi(\bar p_1)\equiv w_1$
and $\phi(\bar p_3)^{-1}\equiv w_3$. Consider the maximal system
of maximal $z$-bands $\bb_1,\dots, \bb_n$ which start on $\bar
p_4$ and end on $\bar p_2$, where $z$ is one of $k$- or
$a$-letters . In particular, $\bb_1$ and $\bb_n$ are subbands of
$\bb$ and $\bb'$ respectively. Let $p_1^ip_2^ip_3^i p_4^i$ be the
standard decomposition of the boundary $\partial(\bb_i)$, where
$p_2^i$, $p_4^i$ are edges from the paths $\bar p_2$ and $\bar
p_4$, respectively, and the labels $w_1^{(i)}, w_3^{(i)}$ of
$p_1^i$, $(p_3^{i})^{-1}$ are some words in $\xxx$. In particular,
$w_1^{(1)}\equiv w_1$ and $w_3^{(n)} \equiv w_3$. The words
$w_1^{(i)}, w_3^{(i)}$ are empty if $\bb_i$ is a $P_j$- or
$R_j$-band, since there are no $(P_j,x)$- or $(R_j,x)$-relations
of the form (\ref{auxkx}).

\begin{center}
\unitlength=1mm \special{em:linewidth 0.4pt} \linethickness{0.4pt}
\begin{picture}(101.67,50.67)
\put(5.33,6.33){\framebox(92.67,40.00)[cc]{}}
\put(5.33,40.00){\line(1,0){92.67}}
\put(5.33,13.00){\line(1,0){92.67}}
\put(47.00,46.33){\vector(1,0){16.00}}
\put(41.00,6.33){\vector(-1,0){11.00}}
\put(5.33,22.00){\vector(0,1){9.33}}
\put(98.00,17.67){\vector(0,1){11.33}}
\put(1.33,26.00){\makebox(0,0)[cc]{$w_1$}}
\put(101.67,26.00){\makebox(0,0)[cc]{$w_3$}}
\put(54.67,48.67){\makebox(0,0)[cc]{$p_2$}}
\put(35.33,3.33){\makebox(0,0)[cc]{$p_4$}}
\put(74.00,9.67){\makebox(0,0)[cc]{$\ttt_1$}}
\put(74.00,43.33){\makebox(0,0)[cc]{$\ttt_2$}}
\put(24.67,13.00){\line(0,1){27.00}}
\put(73.00,13.00){\line(0,1){27.00}}
\put(24.67,21.67){\vector(0,1){12.33}}
\put(73.00,18.33){\vector(0,1){15.67}}
\put(21.00,28.67){\makebox(0,0)[cc]{$w_1^{(i)}$}}
\put(79.33,28.67){\makebox(0,0)[cc]{$w_3^{(i+1)}$}}
\put(45.33,13.00){\line(0,1){27.00}}
\put(45.33,19.33){\vector(0,1){14.67}}
\put(41.00,27.00){\makebox(0,0)[cc]{$w_3^{(i)}$}}
\put(32.67,18.00){\makebox(0,0)[cc]{$\bb_i$}}
\put(63.33,18.00){\makebox(0,0)[cc]{$\bb_{i+1}$}}
\bezier{28}(52.67,40.00)(51.00,38.00)(47.00,37.67)
\put(47.00,37.33){\line(0,-1){21.00}}
\bezier{24}(47.00,16.33)(51.00,14.67)(52.33,13.00)
\put(47.00,21.67){\vector(0,1){12.00}}
\put(53.00,27.00){\makebox(0,0)[cc]{$w_1^{(i+1)}$}}
\end{picture}
\end{center}

\begin{center}
\nopagebreak[4] Fig. \theppp.

\end{center}
\addtocounter{ppp}{1}

Notice, that by $x$-relations (\ref{auxax}) and (\ref{auxkx}), the
word $w_3^{(i)}$ is obtained from the word $w_1^{(i)}$ by an
injective substitution of the form $x\mapsto (x')^l$ on the set of
$x$-letters of the word $w_1^{(i)}$ where $l\in \{1,4, 1/4\}$, $l$
does not depend on an $x$-letter $x$.

There are no cells of $\Gamma$ in the maximal subdiagram of
$\Delta$ whose boundary consists of $x$-edges, and which is
situated between $\bb_i$ and $\bb_{i+1}$. Indeed,  by Lemma
\ref{NoAnnul}, otherwise any non-$\gr$-cell would be included in
an annulus, and there are at most two edges of $\bar p_2$
($x$-edges) between $p_2^i$ and $p_2^{i+1}$ and at most 2 edges of
$\bar p_4$ between $p_4^i$ and $p_4^{i+1}$ (it follows from Lemma
\ref{Band}). Hence the word $w_1^{(i+1)}$ is freely equal to
$u_iw_3^{(i)}v_i$ for $|u_i|, |v_i|\le 2$, and $w_1^{(i+1)}$ is
obtained from from $w_1^{(i)}$ by the above mentioned substitution
$x\mapsto (x')^l$ and multiplication by the words $u_i$, $v_i$
which depend only on the bands $\ttt_1, \ttt_2$. The obvious
induction on $i\le n$, shows that $w_3$ is freely equal to
$uw_1'v$ where $w'$ is obtained from $w_1$ by a substitution of
the form $x\mapsto (x')^{4^c}$ with $|c|\le n$ ($c$ depends on the
bands $\ttt_1$, $\ttt_2$), and $u$, $v$ are fractional words which
can be effectively calculated given $\ttt_1, \ttt_2$.
\endproof

Let $\Delta$ be a trapezium with at least 2-letter base. Suppose
that the side bands $\bb$ and $\bb'$ of $\Delta$  are partitioned
into subbands $\bb_1$, $\bb_2$ and $\bb'_1$, $\bb'_2$
respectively, such that the band $\bb_1$ is the copy of $\bb'_2$
(or $\bb_2$ is the copy of $\bb'_1$), and each of the subbands
$\bb_1,\bb_2, \bb_1',\bb'_2$ has a positive number of
$(\theta,k)$-cells.

We can identify $\bb_1$ with $\bb'_2$ (or $\bb_2$ with $\bb'_1)$.
If the resulting annular diagram $\Gamma$ is reduced and has
minimal boundaries, it is called a \label{spiralsp}{\it spiral}
(see Figure \theppp).

Thus, a spiral $\Gamma$ contains a $k$-band $\bb''$, which
connects its inner contour $q$ with the outer contour $p$; this
band $\bb''$ is obtained from $\bb$ and $\bb'$ by the partial
identification. The $k$-band $\bb''$ will be called the
\label{mainsp}{\em main} $k$-band of the spiral. Notice that the
length of the identified portion of the bands $\bb$ and $\bb'$
does not exceed the lengths of the contours of the spiral. Our
goal will be to estimate the number of cells in a spiral with the
given boundary labels in terms of the lengths of the contours.

\begin{center}
\unitlength=1mm \special{em:linewidth 0.4pt} \linethickness{0.4pt}
\begin{picture}(145.00,60.66)
\put(5.33,10.33){\framebox(64.00,37.67)[cc]{}}
\put(12.33,48.00){\line(0,-1){38.00}}
\put(62.00,48.00){\line(0,-1){38.00}}
\put(5.33,35.00){\line(1,0){7.00}}
\put(5.33,41.67){\line(1,0){64.00}}
\put(11.67,35.00){\line(1,0){57.67}}
\put(5.33,28.00){\line(1,0){64.00}}
\put(5.33,21.33){\line(1,0){64.00}}
\put(5.33,16.00){\line(1,0){64.00}}
\put(5.33,34.33){\rule{7.00\unitlength}{0.67\unitlength}}
\put(12.00,10.33){\rule{0.33\unitlength}{37.67\unitlength}}
\put(5.33,10.33){\rule{0.67\unitlength}{37.67\unitlength}}
\put(62.00,10.33){\rule{0.67\unitlength}{37.67\unitlength}}
\put(69.00,10.33){\rule{0.33\unitlength}{37.67\unitlength}}
\put(8.67,40.33){\makebox(0,0)[cc]{$\bb_2$}}
\put(9.00,22.67){\makebox(0,0)[cc]{$\bb_1$}}
\put(65.67,14.33){\makebox(0,0)[cc]{$\bb_1'$}}
\put(65.67,38.67){\makebox(0,0)[cc]{$\bb_2'$}}
\put(32.33,31.33){\makebox(0,0)[cc]{$\Delta$}}
\put(62.33,20.67){\rule{6.67\unitlength}{0.67\unitlength}}
\bezier{32}(111.33,28.00)(107.33,26.33)(109.00,23.00)
\bezier{72}(109.33,22.67)(117.33,19.00)(118.67,28.00)
\bezier{48}(118.67,28.00)(118.00,34.00)(112.00,33.00)
\bezier{60}(111.33,33.00)(103.33,32.33)(103.00,25.00)
\bezier{76}(103.00,24.33)(103.33,16.00)(114.33,16.00)
\bezier{92}(114.00,16.00)(126.00,15.67)(126.00,26.67)
\bezier{108}(126.00,26.33)(125.33,40.67)(113.00,40.00)
\bezier{136}(113.00,40.00)(96.00,40.67)(96.00,24.00)
\bezier{136}(96.33,27.67)(96.00,10.67)(113.33,10.00)
\bezier{148}(113.33,10.00)(133.00,9.67)(132.67,26.67)
\bezier{160}(132.67,26.67)(133.33,45.33)(112.33,47.00)
\bezier{180}(112.33,47.00)(89.33,47.33)(89.67,25.33)
\bezier{164}(89.67,25.33)(90.67,5.33)(111.67,3.67)
\bezier{200}(111.67,3.67)(138.33,2.33)(140.33,26.00)
\bezier{220}(140.33,26.33)(141.00,52.33)(111.67,55.00)
\bezier{32}(111.67,26.00)(109.00,23.67)(113.67,23.67)
\bezier{64}(113.33,23.67)(120.33,28.00)(113.00,30.67)
\put(111.67,28.00){\line(-1,-6){0.33}}
\bezier{52}(113.00,30.67)(106.00,31.33)(105.33,25.00)
\bezier{56}(105.33,25.00)(106.33,18.67)(113.67,18.00)
\bezier{76}(113.33,18.00)(123.00,18.00)(123.00,27.00)
\bezier{96}(123.00,27.00)(124.00,38.00)(110.67,37.00)
\bezier{96}(110.33,37.00)(98.67,36.33)(98.67,24.33)
\bezier{116}(99.00,26.67)(97.00,14.33)(113.00,12.33)
\bezier{128}(113.00,12.33)(128.33,10.33)(130.67,26.67)
\bezier{148}(130.67,26.67)(131.33,43.33)(111.33,44.33)
\bezier{168}(111.00,44.33)(90.67,44.33)(92.00,22.67)
\bezier{136}(92.00,22.67)(95.67,6.00)(112.67,6.00)
\bezier{176}(112.33,6.00)(136.33,5.67)(137.33,26.00)
\bezier{208}(137.33,26.00)(139.67,49.67)(111.67,52.67)
\put(111.67,52.33){\line(0,1){2.67}}
\put(111.33,28.00){\rule{0.67\unitlength}{27.33\unitlength}}
\put(114.33,27.67){\rule{1.00\unitlength}{27.00\unitlength}}
\bezier{16}(111.67,26.00)(112.00,26.67)(115.33,27.33)
\put(110.67,26.67){\circle*{2.11}}
\put(112.33,26.67){\circle*{0.67}}
\put(112.33,27.00){\circle*{1.89}}
\put(113.33,26.50){\circle*{4.89}}
\put(115.33,26.50){\circle*{4.89}}
\put(112.33,25.50){\circle*{4.89}}
\put(111.67,26.67){\circle*{1.49}}
\put(113.33,27.00){\circle*{1.89}}
\put(109.33,24.67){\circle*{2.11}}
\put(109.67,25.67){\circle*{2.00}}
\put(110.67,23.00){\circle*{2.40}}
\put(113.67,22.67){\circle*{2.75}}
\put(112.00,22.67){\circle*{1.33}}
\put(111.67,22.67){\circle*{2.00}}
\put(115.33,23.33){\circle*{2.67}}
\put(116.33,24.33){\circle*{2.11}}
\put(117.00,25.33){\circle*{1.49}}
\put(117.33,26.00){\circle*{1.33}}
\put(117.33,26.67){\circle*{1.49}}
\put(117.67,27.33){\circle*{1.33}}
\put(117.67,28.00){\circle*{1.49}}
\put(117.33,28.67){\circle*{2.00}}
\put(117.33,29.67){\circle*{2.00}}
\put(117.33,30.67){\circle*{2.11}}
\put(116.33,31.00){\circle*{1.49}}
\bezier{28}(72.33,25.67)(75.67,28.00)(78.00,26.00)
\bezier{24}(78.00,26.00)(80.67,24.00)(83.00,26.00)
\put(82.67,26.00){\vector(1,0){2.00}}
\put(145.00,29.33){\makebox(0,0)[cc]{$\Gamma$}}
\put(109.18,23.32){\circle*{0.97}}
\put(115.61,29.29){\circle*{1.53}}
\put(115.76,30.05){\circle*{0.97}}
\put(116.07,25.15){\circle*{1.53}}
\put(116.37,26.07){\circle*{0.92}}
\put(116.68,27.14){\circle*{0.68}}
\put(117.75,25.92){\circle*{0.87}}
\put(118.06,26.68){\circle*{1.10}}
\put(118.36,27.60){\circle*{0.87}}
\put(118.36,28.52){\circle*{0.97}}
\put(118.36,29.29){\circle*{0.87}}
\put(118.36,30.05){\circle*{0.87}}
\put(112.24,21.79){\circle*{1.22}}
\put(111.32,21.79){\circle*{0.87}}
\put(108.87,23.62){\circle*{0.97}}
\put(110.40,25.15){\circle*{0.68}}
\put(110.86,25.61){\circle*{0.92}}
\put(110.25,24.39){\circle*{0.97}}
\put(115.76,31.43){\circle*{1.94}}
\put(115.61,31.89){\circle*{1.86}}
\put(116.53,31.58){\circle*{1.37}}
\put(113.31,60.66){\makebox(0,0)[cc]{$\bb''$}}
\end{picture}
\end{center}

\begin{center}
\nopagebreak[4] Figure \theppp.
\end{center}
\addtocounter{ppp}{1}

It is easy to see that every maximal $\theta$-band of $\Gamma$ is
a union of $\theta$-bands of $\Delta$, it starts on the inner
contour $q$ of $\Gamma$ and ends on the outer contour $p$. It may
intersect the $k$-band of $\Gamma$ many times, and looks like a
spiral. (However it can intersect every $k$-band of $\Gamma$ only
in one direction: this follows from Lemma \ref{NoAnnul}.) Hence
the number of maximal $\theta$-bands of the spiral is at most
$\min (|q|, |p|)$.

Let $h_1, h_2$ be the histories of the bands $\bb_1$, $\bb_2$.
Then $\bb$ and $\bb'$ have the history $h_1h_2$. On the other hand
$h_2$ is a prefix of the history of $\bb'$. Hence $h_1h_2=h_2h_3$
for some $h_3$. If $|h_1|>|h_2|$ , then this immediately implies
that $h_1h_2$ is a periodic word with period $h_2 $(i.e. $h_1h_2$
is a subword of $h_2^l$ for some $l>0$). Notice that $|h_2|$ does
not exceed the length of the inner contour of the spiral. So
$|h_1h_2|<2|h_2|\le |q|+|p|$ if $|h_1|\le |h_2|$, because every
$\theta$-band in $\Gamma$ connects $q$ with $p$. Hence i any case,
the history of $\Delta$ is a periodic word whose period does not
exceed $|q|+|p|$.

\label{qspirals}{\it Quasispirals} are defined like spirals when
one takes a quasitrapezium with empty base as the original diagram
instead of a trapezium.

\begin{lemma}\label{sp1234} Let $\Delta$ be a reduced annular diagram over $\hhh_1$ with
minimal boundaries, whose labels contain at least one and
 a $k$- or $a$-edge, $a\notin\aaa(P_1)$. Let the corresponding maximal
 $k$- or $a$-band $\CC$ have recursively
bounded (in terms of $|p|+|q|$) number of $\theta$-cells. Then
either (1) there exists a reduced annular diagram $\Delta'$ over
$\hhh_1$ with the same boundary labels as $\Delta$ and with
recursively bounded (in terms of $|p|+|q|$) number and perimeters
of cells, or (2) there is a reduced annular diagram $\Delta'$
without $\gr$-cells and with minimal boundaries $p'$ and $q'$
which are at most twice longer than $p$ and $q$ respectively, and
$\Lab(p)$ and $\Lab(p')$ (resp. $\Lab(q), \Lab(q')$) are conjugate
and the number and lengths of relations needed to deduce these
conjugacies are recursively bounded.
\end{lemma}

\proof By Lemma \ref{rollsprep2}, we can assume that $\Delta$ has
minimal boundaries. If there exists an $a$-band starting on a
$\gr$-cell in $\Delta$ and ending on a $k$-band of $\Delta$ then
this must be a $P_1$- or an $R_1$-band, and so we would be able to
turn $\Delta$ into a simply connected \vk diagram by cutting along
a side of this $k$-band. The length of this cut is recursively
bounded as this follows from the lemma assumption and the lack of
$x$-cells in $P_1$- and $R_1$-bands. After that we could apply
Lemma \ref{K0}. Thus we can assume that every $a$-band starting on
a $\gr$-cell of $\Delta$ ends on $p$ or $q$. Therefore by Lemma
\ref{aa-cells}, parts (1), (2), we can remove all the $\gr$-cells
from $\Delta$ increasing the length of the boundary at most by a
factor of 2. The total perimeter of all $\gr$-cells that we remove
is bounded by $|p|+|q|$.
\endproof

\begin{lemma}\label{sp123} Let $\Delta$ be a reduced annular diagram over $\hhh_1$ with
minimal boundaries $p$ and $q$ having a $\theta$-band $\ttt$ and
at least one $k$- or $a$-edge for an $a$-letter $a\in
(\aaa\cup\bar\aaa)\backslash\aaa(P_1)$. Let $\CC$ be corresponding
$k$ or $a$-band connecting the contours $q$ and $p$. Then there
exists a diagram $\Delta'$ with the same boundary label as
$\Delta$, and such that $\Delta'=\Delta_1\cup\Delta_2\cup\Delta_3$
where the number and perimeters of cells in $\Delta_1, \Delta_3$
are recursively bounded and $\Delta_2$ is either empty or a
spiral, or a quasispiral.
\end{lemma}

\proof By Lemmas \ref{rollsprep1}, \ref{rollsprep}, we can obtain
a minimal diagram $\Delta'$ satisfying (R1). Then all maximal
$\theta$-bands in $\Delta$ connect $p$ and $q$, so their number
does not exceed $|p|+|q|$. Also all maximal $k$-bands in $\Delta'$
connect $p$ and $q$. If there are no $k$-bands, then every
$a$-band, for $a\notin \aaa(P_1)$, connects $p$ and $q$.


 All maximal $a$-bands starting on $p$ or $q$ end either
on $q$ and $p$ respectively, or on a $k$-band or on a $\gr$-cell
(recall that different $a$-bands cannot cross and different
$k$-bands cannot cross).

Suppose first that every $\theta$-band in $\Delta$ intersects
every $k$-band and every $a$-band, $a\notin \aaa(P_1)$, in
$\Delta$ at most twice. Since the number of maximal $\theta$-bands
in $\Delta$ is bounded by $|p|+|q|$, the number of $\theta$-cells
in $\CC$ is at most $2(|p|+|q|)$. By Lemma \ref{sp1234} we can
assume that $\Delta$ does not contain $\gr$-cells. Since the
number of maximal $k$-bands in $\Delta$ is bounded by $|p|+|q|$,
the number of $\theta$-cells in $k$-bands of $\Delta$ is bounded
by $2(|p|+|q|)^2$. Therefore the number of maximal $a$-bands in
$\Delta$ is bounded by $2(|p|+|q|)^2+(|p|+|q|)<3(|p|+|q|)^2$.
Since each of the $a$-bands intersects each of the $\theta$-bands
at most twice, the number of $\theta$-cells in each of them is
bounded by $2(|p|+|q|)$. Therefore the total number of
$\theta$-cells in $\Delta$ is bounded by
$6(|p|+|q|)^3+2(|p|+|q|)^2$. Therefore the $\theta$-band $\ttt$
has a recursively bounded length. Thus we can cut $\Delta$ along a
side of $\ttt$, obtain a simply connected \vk diagram with
perimeter recursively bounded in terms of $|p|+|q|$, and use Lemma
\ref{K0} to bound the number of cells in $\Delta$.

Notice that if a $\theta$-band intersects an $a$- or $k$-band
$\CC'$ three times then by Lemma \ref{NoAnnul} it must go around
the hole of the diagram twice and so by Jordan's lemma it
intersects the $a$- or $k$-band $\CC$ connecting $p$ and $q$ at
least twice.

Thus we can assume that the $\theta$-band $\ttt$ in $\Delta$
intersects $\CC$ twice. Moreover we can conclude that there is a
subband $\CC_0$ of $\CC$ with linearly (in terms of $|p|+|q|$)
bounded number of $\theta$-cells in $\CC\backslash \CC_0$, such
that for every $\theta$-cell $\pi$ of $\CC_0$ there is a
$\theta$-band $\ttt_{\pi}$ which starts with $\pi$ and whose last
cell $\pi'$ (and no other cell) also belongs to $\CC$, $\pi'\ne
\pi$. These bands $\ttt_{\pi}$, for $\pi\in \CC_0$, (and the
auxiliary cells between them) form a subdiagram $\Delta_2$ which
is a spiral or a quasispiral with main band $\CC_0$ by the
construction and the assumptions on $p$ and $q$. If there exists a
$\gr$-cell in $\Delta$ outside $\Delta_2$ and an $a$-band
connecting it with the boundary of $\Delta_2$, we can use Lemma
\ref{aa-cells}, move the $\gr$-cell to the boundary of $\Delta_2$,
and then include that cell in $\Delta_2$. Thus we can assume that
there are no $\gr$-cells in $\Delta$ outside $\Delta_2$ which are
connected with the boundary of $\Delta_2$ by an $a$-band.


\begin{center}
\unitlength=1mm \special{em:linewidth 0.4pt} \linethickness{0.4pt}
\begin{picture}(144.00,113.00)
\put(70.33,54.33){\circle*{5.20}}
\put(69.00,56.67){\line(0,1){54.00}}
\put(72.33,56.00){\line(0,1){54.33}}
\bezier{356}(72.33,99.33)(28.67,99.33)(26.67,54.00)
\bezier{324}(26.67,54.00)(29.67,15.33)(72.33,17.00)
\bezier{288}(72.33,17.00)(104.67,17.67)(108.00,57.67)
\bezier{224}(108.00,58.00)(107.00,76.33)(69.00,73.33)
\bezier{328}(72.33,92.33)(34.33,97.00)(31.33,53.00)
\bezier{288}(31.33,53.00)(32.67,20.67)(72.33,21.67)
\bezier{244}(72.33,21.67)(100.67,23.33)(101.33,55.67)
\bezier{200}(101.33,55.67)(105.67,68.67)(69.00,66.67)
\put(71.08,92.38){\line(1,0){1.32}}
\put(70.20,99.22){\line(1,0){2.21}}
\put(69.09,73.41){\line(1,0){3.31}}
\put(69.09,66.79){\line(1,0){3.31}}
\bezier{80}(69.09,63.48)(58.28,63.04)(57.40,53.99)
\bezier{96}(57.40,53.77)(58.06,42.74)(70.64,44.06)
\bezier{76}(70.64,44.06)(79.90,46.27)(80.35,55.76)
\bezier{64}(80.35,55.76)(80.35,60.61)(69.09,59.73)
\put(70.64,63.26){\makebox(0,0)[cc]{$\CC_0$}}
\put(71.00,70.00){\makebox(0,0)[cc]{$\pi'$}}
\put(71.00,96.00){\makebox(0,0)[cc]{$\pi$}}
\put(71.00,105.00){\makebox(0,0)[cc]{$\CC$}}
\put(108.00,68.00){\makebox(0,0)[cc]{$\ttt_\pi$}}
\bezier{392}(69.00,111.00)(15.00,113.00)(6.00,70.00)
\bezier{164}(6.00,70.00)(3.00,43.00)(12.00,33.00)
\bezier{112}(12.00,33.00)(29.00,27.00)(30.00,17.00)
\bezier{196}(30.00,17.00)(34.00,1.00)(66.00,4.00)
\bezier{168}(66.00,4.00)(101.00,13.00)(102.00,7.00)
\bezier{208}(102.00,7.00)(111.00,0.00)(136.00,32.00)
\bezier{240}(136.00,32.00)(144.00,45.00)(130.00,88.00)
\bezier{292}(130.00,88.00)(122.00,106.00)(69.00,111.00)
\put(72.33,101.33){\line(0,1){9.33}}
\put(135.67,68.33){\vector(1,-4){0.67}}
\put(138.33,69.67){\makebox(0,0)[cc]{$p$}}
\put(13.33,55.33){\makebox(0,0)[cc]{$\Delta_1$}}
\put(42.67,53.00){\makebox(0,0)[cc]{$\Delta_2$}}
\put(62.00,52.00){\makebox(0,0)[cc]{$\Delta_3$}}
\put(75.00,52.00){\makebox(0,0)[cc]{$q$}}
\bezier{52}(69.00,73.33)(60.00,72.67)(57.00,70.33)
\bezier{44}(69.00,66.67)(65.33,67.33)(58.33,66.33)
\bezier{56}(71.67,92.33)(85.33,91.67)(84.67,91.33)
\bezier{52}(72.33,99.33)(81.67,99.00)(84.67,98.00)
\end{picture}
\end{center}

\begin{center}
\nopagebreak[4] Fig. \theppp.

\end{center}
\addtocounter{ppp}{1}

The diagram $\Delta$ is a union of $\Delta_2$ and two annular
diagrams $\Delta_1$, $\Delta_3$ which are situated between
$\Delta_2$ and $p$ and between $\Delta_2$ and $q$ respectively.

Our purpose is to obtain a recursive upper bound for the number of
cells in $\Delta_1$ and $\Delta_3$ as a function in $|p|, |q|$. We
will consider only the diagram $\Delta_1$ ($\Delta_3$ is similar).

By the choice of $\CC_0$, each of the maximal $\theta$-bands in
$\Delta$ has linearly bounded (in terms of $|p|+|q|$) number of
intersections outside $\Delta_2$ with every maximal $a$- or
$k$-band of $\Delta$. Notice also that an $a$-band of $\Delta_1$
which has at least one $\theta$-cell, cannot start and end on the
boundary of $\Delta_2$ by Lemma \ref{NoAnnul}. Hence every maximal
$a$-band of $\Delta_1$ must start and end on $p$ or on a
$(k,\theta)$-cell, and therefore the number of $\theta$-cells in
$\Delta_1$ is recursively bounded in terms of $|p|+|q|$. By our
assumption, every $a$-band starting on a $\gr$-cell of $\Delta_1$
ends on $p$ or on a $\theta$-cell of $\Delta_1$. Therefore the
total perimeter of all $\gr$-cells in $\Delta_1$ is recursively
bounded.

Since $\theta$-cells of $\Delta_1$ are arranged in $\theta$-bands
starting on $p$, we may delete all (recursively bounded number) of
them from $\Delta_1$ (the result will be again an annular diagram)
and assume that every cell of $\Delta_1$ is an $(a,x)$-cell or a
$(k,x)$-cell corresponding to the relations (\ref{auxax}),
(\ref{auxkx}), or a $\gr$-cell. Then the $\gr$-cells can also be
deleted from $\Delta_1$ by Lemma \ref{aa-cells} (and such
transformations do not change the boundary of $\Delta_2$).

Since there are no $a$-, $k$-, $(a,\theta)$-, $(k,\theta)$-annuli
in $\Delta$ by Lemma \ref{NoAnnul}, all cells of $\Delta_1$ belong
to maximal $a$- and $k$-bands $\CC_1, \CC_2, \dots$ starting on
$p$ (and ending on the other boundary component of $\Delta_1$).
Hence it suffices to bound the length of each of these $a$- and
$k$-bands in $\Delta_1$.


Notice that the common boundary component, say $u$, of $\Delta_1$
and $\Delta_2$ is the boundary component of a spiral (or
quasispiral). By the above construction of $\Delta_2$, we have
$u=u_1u_2$ where $u_1$ belongs to the boundary of $\CC_0$, the
last edge of $u_1$ is a $\theta$-edge,
 and the label of $u_2$ is the label of a $\theta$-band (modulo
 $\gr$).

Suppose that there is a band $\CC_i$ among $\CC_1, \CC_2,\dots$,
which ends on $u_2$. We may assume that $i$ is chosen so that the
end of $\CC_i$ is the nearest to the beginning of $u_2$ among all
$\CC_j$ ending on $u_2$. Then $\bott(\CC_i)$ cannot have common
edges with any $\CC_j$ since the last edge of $u_1$ is a
$\theta$-edge, which also belong to $p$. (Recall that now we have
no $\theta$-cells in $\Delta_1$.) The bottom of $\CC_i$  has at
most 2 common edges with $u_2$ because $u_2$ can have at most 2
consecutive $x$-letter in its boundary label by Lemma \ref{Band}.
Hence the length of $\CC_i$ is recursively bounded by $|p|$.

Then suppose that there is a band $\CC_i$ ending on $u_1$ (and so
$\CC$ is a $k$-band). It is clear from the form of
$(k,\theta)$-relations that there is a subpath $evf$ of $u_1^{\pm
1} $ such that $e$ is the end of $C_i$, $f$ is a $\theta$-edge,
and $v$ contains no two consecutive $x$-letters in its label. We
may assume that no $C_j$ ends on $v$. Then, as in the first case,
the length of $C_i$ is recursively bounded by $p$.

Thus, in any case, there is a band, say $\CC_1$, whose perimeter
and whose number of cells are linearly bounded in terms of $|p|$.
We delete $\CC_1$ from $\Delta_1$ and then choose a band in the
remaining subdiagram of $\Delta_1$, whose perimeter and whose
number of cells are recursively bounded as well. Then we delete
$\CC_2$, etc. Since the number of $a$-, $k$-bands $\CC_1,\dots$ is
recursively bounded, we conclude that there is a recursive bound
for the number of cells in $\Delta_1$, and the lemma is proved.
\endproof

\begin{rk} \label{rk123} {\rm The proof of Lemma \ref{sp123} shows that for every
maximal $k$-band $\CC$ in a spiral $\Delta$, one can obtain
another spiral $\Delta'$ with contours $p'$ and $q'$, such that
(a) the words $\phi(p)$ and  $\phi(p')$ ($\phi(q)$ and $\phi(q')$)
are conjugate modulo relations used in $\Delta$; (b) the number
and lengths of relations used to deduce these conjugacies are
recursively bounded in terms of $|p|+|q|$; (c) the set of
$\theta$-cells of the main band of $\Delta'$ is a subset of the
set of cells of the band $\CC$.
 Thus if the base of a spiral has letters $P_j$ or $R_j$ we shall
always assume that the main $k$-band of the spiral is a $P_j$- or
a $R_j$-band.}
\end{rk}

\begin{lemma} \label{historyspiral}
Let $h$ be a periodic word in $\csss$ with period $h_0$ and let
$W\bullet h=W'$ be a free computation of $\csss$. Then the number
of different words among $W\circ h[t]$, $t=1,...,|h|$ is
recursively bounded in terms of $|h_0|+|W|+|W'|$.
\end{lemma}

\proof By lemma \ref{csss}, for every $zz'$-sector of $W$ there
exist two words $u$ and $v$ depending on $h_0$ whose size is
linearly bounded in terms of $|h_0|$, such that the inner part of
the corresponding $zz'$-sector in each $W\circ h_0^t$ is equal in
the free group to $u^tW_{zz'}v^t$ where $W_{zz'}$ is the inner
part of the $zz'$-sector of $W$, $t=1,...,s$. By Lemma
\ref{powers}, we can assume that $s$ is big enough so that
$s-|u|-|v|>|W'|$, so $u^tW_{zz'}v^t=W_{zz'}$ for every
$t=1,...,s$. Hence $W=W\circ h_0=W\circ h_0^2=...$. Replacing $W$
by $W\circ h'$ for any prefix $h'$ of $h_0$, we get that $W\circ
h'=W\circ hh'=W\circ h^2h'=...$. Therefore the number of different
words in the computation associated with $\Delta$ is recursively
bounded in terms of $|W|+|W'|$ and $|h_0|$.\endproof

\begin{lemma} \label{spiral} Let a minimal diagram $\Gamma$  be a spiral obtained from a
trapezium $\Delta$. Then the minimum of lengths of all the maximal
$k$-bands of $\Gamma$ is bounded by a recursive function of
$|q|+|p|,$ where $q, p$ are inner and outer contours of $\Gamma$.
 \end{lemma}

\proof Let $W$ be the projection of the label of the bottom side
of the trapezium $\Delta$ corresponding to the spiral $\Gamma$
onto $\aaa\cup\bar\aaa\cup\kk\cup\bar\kk$, and let $h$ be the
history of $\Delta$ corresponding to $\Gamma$. Notice that
$h=h_0^sh'$ for some $s$ and some prefix $h'$ of $h_0$, and
$|W|,|h_0|\le |q|+|p|$. We need to bound $s$ in terms of
$|q|+|p|$. Indeed, suppose we have bounded $s$. Suppose the spiral
contains a $\gr$-cells. If there exists $a$-band starting on a
$\gr$-cell there ends on a $P_1$- or $R_1$-band, and the spiral
has a cut of a recursively bounded length along a side of that
band since it has no $\xxx$-cells. Thus we can assume, by Lemma
\ref{aa-cells}, that all $a$-bands starting on $\gr$-cells end on
the boundary of the spiral, in which case the total perimeter of
these cells is recursively bounded, and we can move these cells
out of the spiral using Lemma \ref{aa-cells}. Therefore we can
assume that the spiral does not contain $\gr$-cells. Since the
number of $k$-cells is recursively bounded in terms of $|h|$ and
$|q|+|p|$, we have  that the number of maximal $a$-bands is
bounded too. Hence the length of $\theta$-bands in $\Delta$ are
bounded as well, because every such a band has recursively bounded
number of intersections with every $a$-band. Therefore there is a
cut of $\Delta$ with recursively bounded length, which makes
possible to apply Lemma \ref{K0}.

 Let us assume
that $s>|q|+|p|$.

Let us enumerate $k$-bands in $\Gamma$ clockwise:
$\bb_1,...,\bb_c$. Let $\ttt$ be a maximal $\theta$-band of
$\Gamma$. We can decompose $\ttt$ into several subbands $\ttt(1),
\ttt(2),\dots, \ttt(r_c)$ such that every $\ttt(l)$ (with possible
exception $\ttt(1)$) starts with a cell of the band $\bb_1$ and
contains at most one common cell with each of the $k$-bands
$\bb_1,\dots \bb_c.$ The bands $\ttt(2),\ttt(3),\dots$ will be
called the {\it turns} of $\ttt$. For every $l=2,...$ there is a
path $t(l)$ going along $\bb_1$, which connect the initial vertex
of the top $p(l)$ of the turn $\ttt(l)$ with that of the top of of
$\ttt(l+1)$. Notice that $p(l+1)_- = p(l)_+$. Denote by $q(l)$ the
loop $t(l)p(l)^{-1}$ associated with $\ttt(l)$. Notice also that
for every $l=1,2..$ words
$\phi(p(l))_{\aaa\cup\bar\aaa\cup\kkk\cup\bar\kkk}$ coincides for
some $i$, with the word $W\circ h[i]$ with omitted the last
$k$-letter $z$. The label of $t(l)$ consists of at most $|h_0|$
$\theta$-letters and a number of $x$-letters.

Consider several cases.

{\bf Case 1.} There is a $P_j$- or $R_j$-letter in the base of
$\Delta$.

Then we assume that $\bb_1$ is a $P_j$- or a $R_j$-band by Remark
\ref{rk123}. We have the following alternative.

{\bf Case 1.1}. Suppose that the history of $\Delta$ contains a
rule from $\bsss\cup\sss(34)\cup\sss(4)\cup\sss(45)$.

Then we can assume that $\ttt$ corresponds to one of these rules
$\tau$. Let $W_{zz'}$ be a sector of $W$ which is not a $P_1R_1$-,
$P_1P_1\iv$- or $R_1\iv R_1$-sector. Then the computation
$W_{zz'}\bullet h$ is a free computation and, by Lemma
\ref{historyspiral}, the number of different words among
$W_{zz'}\circ h[i]$, $i=1,2...$ is recursively bounded in terms of
$|q|+|p|$.  Notice that the rule $\tau$ locks $P_1R_1$-sectors.
Hence the $P_1z$- and $zR_1$-sectors of
$\Lab(p(l))_{\aaa\cup\bar\aaa\cup\kkk\cup\bar\kkk}z$ are of length
$2$. Therefore by Lemma \ref{Band}, the number of different words
among $\Lab(p(l))$, $l=1,2,...$, is recursively bounded.

Since the sides of $P_j$-, and $R_j$-bands do not contain
$x$-edges, the lengths of the paths $t(l)$, $l=1,2,...$ do not
exceed $|h_0|\le |q|+|p|$.

Therefore the number of different labels of the loops
$t(l)p(l)^{-1}$, $l=1,2,...$ is recursively bounded in terms of
$|q|+|p|$. Notice that for different $l$'s these loops do not
intersect each other. Since $\Gamma$ is minimal, it cannot contain
two non-intersecting loops with the same labels. Hence $|h|$ is
recursively bounded in terms of $|q|+|p|$. Hence the length of
$\bb_1$ is recursively bounded, as desired.

{\bf Case 1.2}. Suppose that no rules from
$\bsss\cup\sss(34)\cup\sss(4)\cup\sss(45)$ occur in $h$.

The history $h$ is reduced by Lemma \ref{History}. By Lemma
\ref{Sm Tr}, then $W\bullet h$ is a computation of $\sss$ and $W$
is an admissible word for $\sss$. Then one can apply Lemma
\ref{refined1} to the computation $W\bullet h=W' (\mod \gr)$
corresponding to the trapezium $\Delta$ and obtain a free
computation $W\bullet h=W''$ where $W''=W' (\mod \gr)$. Notice any
two $P_j$-bands (resp. $R_j$-bands) which have the same histories
consist of the same cells since there are no $(P_j,x)$- and
$(R_j,x)$-relations of the form (\ref{auxkx}). Therefore since $W$
starts and ends with the same $P_j$- or $R_j$-letter, the free
computation $W\bullet h$ corresponds to a spiral whose inner
contour has the same label as the inner contour of $\Gamma$ and
the outer contour has the same label $\mod \gr$ as the outer
contour of $\Gamma$. Therefore the inner parts of non-$P_1z$- and
non $zR_1$-sectors of $W\circ h$ coincide with the corresponding
subwords of $W'$. We would be able to apply the same argument as
in Case 1.1 to the spiral corresponding to the free computation
$W\bullet h$ if we bound the sizes of $P_1z$- and $zR_1$-sectors
of the word $W\circ h$.

Let us prove that these sizes are recursively bounded. We can
certainly assume by Remark \ref{rk123}, that the base $W^{\pm 1}$
starts and ends on $P_1$ or $R_1$. If the base of $W^{\pm 1}$
contains both $P_1$ and $R_1$, we can assume that it starts with
$P_1$.

By Lemma \ref{noPivP} there are the following possibilities for
the base of $W^{\pm 1}$:

\begin{itemize}
\item the base of $W^{\pm 1}$ has a prefix $P_1R_1K_2$ and suffix
$K_1L_1P_1$;
\item the base of $W^{\pm 1}$ has a prefix $P_1P_1\iv$ and suffix
$K_1L_1P_1$;
\item the base of $W^{\pm 1}$ has a prefix $R_1\iv R_1K_2$.
\end{itemize}

These three cases are completely similar, so let us consider only
the first one. Then $W^{\pm 1}$  has the form
$P_1V_1R_1V_2K_2...K_1V_3L_1V_4P_1$ where $V_1$ (resp $V_2$,
$V_3$, $V_4$) is a word in $\aaa(P_1)$ (resp. $\aaa(R_1),
\aaa(K_1)$, $\aaa(L_1)$). By Lemma \ref{csss} there exist words
$u(L_1,h_0)$, $v(L_1,h_0)$, $u(P_1,h_0)$, $v(P_1,h_0)$,
$u(R_1,h_0)$, $v(R_1,h_0)$ such that

\begin{equation}\begin{array}{l}\beta(u(L_1,h_0))\equiv
\beta(v(L_1,h_0))\iv,\\ \beta(u(P_1,h_0))\equiv
\beta(v(P_1,h_0))\iv, \beta(u(R_1,h_0))\equiv
\beta(v(R_1,h_0))\iv\end{array}\label{eq123}\end{equation} and for
every $t=1,...,s$ $$\begin{array}{l}W\circ h_0^t=\\
P_1u(P_1,h_0)^tV_1v(R_1,h_0)^tR_1u(R_1,h_0)^tV_2K_2^{-1}...K_1V_3v(L_1,h_0)^tL_1u(L_1,h_0)^tV_4
v(P_1,h_0)^tP_1.\end{array}$$ Since we assumed that $s>|q|+|p|$,
$|V_3|< |q|$ and $V_3v(L_1,h_0)^t\le |p|$, we must have
$v(L_1,h_0)=1$ in the free group. Hence by (\ref{eq123})
$u(L_1,h_0)=1$ in the free group. Since $|V_4|<|q|$ and
$|V_4v(P_1,h_0)^s|<|p|$, we have that $v(P_1,h_0)=1$ in the free
group. Then $u(P_1,h_0)$ is freely equal to 1. Similarly,
considering the words $V_2$ and $v(R_1,h_0)^sV_2$, we conclude
that $v(R_1,h_0)=1$ in the free group. Hence $u(R_1,h_0)$ is
freely trivial. But then
$|u(P_1,h_0)^sV_1v(R_1,h_0)^s|=|V_1|<|q|$, and we have a (linear)
bound of the size of the $P_1R_1$-sector of the word $W\circ h$ as
desired.

{\bf Case 2}. Suppose now that there are neither $P_j$- nor
$R_j$-letters in the base of $\Delta$.

Then the trapezium $\Delta$ has no $\gr$-cells by Lemma \ref{No
a-cells}. Again, it suffices to give a recursive upper bound for
the length of the loops $q(l)$, because this would give a
recursive estimate for the number of their labels. As in Case 1,
we can bound the sizes of the paths $p(l)$. Since the history of
$\Delta$ is periodic and the lengths of the words $\phi(p(l))$
(whose $(k,a)$-projections are the intermediate words, without the
last letter, of the corresponding computation ) are bounded, the
words $\phi(p(l))$ form a periodic sequence whose period $n$ is
recursively bounded in terms of $|q|+|p|$.

Thus it remains to obtain an upper bound for the lengths of the
paths $t(l)$. Since the number of $\theta$-edges in $t(l)$ is at
most $\min(|q|, |p|)$, it suffices to obtain an upper estimate for
the length of each of the $x$-subpaths $t_{lj}$ of $t(l)$
connecting neighbor $\theta$-bands of $\Delta$.

The $\theta$-bands in $\Gamma$ subdivide this diagram into several
(at most $|q|$) subdiagrams $\Gamma_i$ containing no
$\theta$-cells. The path $t_{lj}$ is a side of one of the
$k$-bands in one of the subdiagrams $\Gamma_i$. Thus it is to show
that $k$- and $a$-bands in $\Gamma_i$ have bounded lengths. Notice
that the boundary of the \vk diagram $\Gamma_i$ has the form
$q^1q^2q^3q^4$ where $\Lab(q^1)$ and $\Lab(q^3)$ are subwords in
$\xxx$ of the boundary labels of $\Delta$ and $\Lab(q^2),
\Lab(q^4)$ are subwords of the top and the bottom labels of
$\theta$-bands in $\Gamma$.

\begin{center}
\unitlength=1mm \special{em:linewidth 0.4pt} \linethickness{0.4pt}
\begin{picture}(126.67,96.67)
\put(55.67,37.00){\framebox(3.67,56.00)[cc]{}}
\put(55.67,83.33){\rule{25.33\unitlength}{1.33\unitlength}}
\put(55.67,74.67){\rule{25.67\unitlength}{1.67\unitlength}}
\put(55.67,83.00){\rule{25.67\unitlength}{1.33\unitlength}}
\bezier{48}(81.33,84.67)(90.00,83.33)(91.33,80.33)
\bezier{48}(81.33,82.67)(88.67,81.67)(91.33,78.33)
\bezier{48}(81.33,76.33)(90.00,75.67)(91.00,72.00)
\bezier{48}(81.33,74.33)(89.00,73.33)(91.00,69.67)
\put(59.33,54.67){\rule{-20.00\unitlength}{-0.67\unitlength}}
\put(38.67,53.00){\rule{20.67\unitlength}{2.00\unitlength}}
\put(38.33,44.33){\rule{21.00\unitlength}{1.67\unitlength}}
\bezier{56}(38.33,55.00)(28.00,53.67)(26.33,50.67)
\bezier{56}(38.33,53.00)(30.00,52.67)(26.33,48.33)
\bezier{60}(38.00,46.00)(27.33,46.00)(26.00,42.33)
\bezier{52}(38.00,44.33)(30.67,44.67)(26.00,40.67)
\bezier{268}(55.67,70.67)(13.67,68.00)(7.67,43.33)
\bezier{376}(7.67,43.33)(4.67,4.00)(59.33,5.00)
\bezier{456}(59.33,4.67)(126.67,7.00)(125.33,54.00)
\bezier{432}(125.33,54.00)(124.33,96.67)(59.33,93.00)
\put(32.33,49.67){\makebox(0,0)[cc]{$q^2$}}
\put(38.67,48.33){\makebox(0,0)[cc]{$q^4$}}
\put(82.67,80.00){\makebox(0,0)[cc]{$q^2$}}
\put(66.33,78.33){\makebox(0,0)[cc]{$q^4$}}
\put(73.33,78.67){\makebox(0,0)[cc]{$\Gamma_i$}}
\put(55.67,77.00){\vector(0,1){4.00}}
\put(52.33,78.00){\makebox(0,0)[cc]{$q^1$}}
\put(59.33,51.33){\vector(0,-1){6.67}}
\put(62.00,48.67){\makebox(0,0)[cc]{$q^3$}}
\put(52.00,48.67){\makebox(0,0)[cc]{$\Gamma_i$}}
\bezier{132}(55.67,37.00)(48.33,28.33)(70.00,28.00)
\bezier{280}(70.00,28.00)(94.67,35.67)(59.33,62.67)
\put(78.67,49.67){\makebox(0,0)[cc]{$q$}}
\put(40.00,61.00){\makebox(0,0)[cc]{$\Delta$}}
\put(119.33,47.33){\makebox(0,0)[cc]{$p$}}
\put(125.33,55.00){\vector(0,-1){3.67}}
\end{picture}

\end{center}

\begin{center}
\nopagebreak[4] Fig. \theppp.

\end{center}
\addtocounter{ppp}{1}

The words $\Lab(q^2)$, $\Lab(q^4)$ are periodic with the length of
the period  bounded in terms of $n$. Hence $\Gamma_i$ can be cut
by its $k$- or $a$-bands into subdiagrams $\Psi_1,\Psi_2,\dots
\Psi_r$ with contours $q_i^1q_i^2q_i^3q_i^4$ where $q_i^1$,
$q_i^3$ are words in $\xxx$, $q_i^2$ (or $q_i^4$) have the same
labels of length $n$ for all $i$, with the possible exception for
$q_r^2, q_r^4$ whose lengths are recursively bounded in terms of
$n$. By Lemma \ref{2 bands}, the subdiagram $\Psi_r$ has bounded
number of cells and can be ignored. Also Lemma \ref{2 bands}
implies that it suffices to bound $|q_i^1|$ for all $i$ in terms
of $|q_1^1|, |q_{r-1}^1|$.

Let $w_i\equiv \Lab(q_i^1)$. By Lemma \ref{2 bands} for $\Psi_i$,
we have that $w_{i+1}$ is freely equal to $uw'_iv$ where the
fractional words $u, v$ depend on  $\Lab(q_i^2), \Lab(q_i^4)$, and
therefore, do not depend on $i$, and so they belong to some finite
effectively determined set of words; $w'_i$ results from $w_i$
after a substitution of the form $x\mapsto (x')^{4^d}$ for a
recursively bounded $d$.

If $d=0$, then we have, by the obvious induction, that $w_i$ is
freely equal to a copy of $u^{i-1}w_1v^{i-1}$. Then, by Lemma
\ref{powers}, $|w_i|$ is effectively bounded by $|w_1|+|w_{z-1}|$,
as desired. Hence we may assume that $d\ge 1$. (We would
interchange the top and the bottom of the spiral if $d\le -1$.)

Since $d\ge 1$, $|w_{i+1}|\ge 4|w_i|-C$ where $C=|u|+|v|$. Then
either $|w_i|\le C$ for every $i$, or, if $|w_i|>C$ for some $i$,
the series $|w_i|, |w_{i+1}|,\dots $ monotonically increases, and
so $|w_i|\le |w_{z-1}|$.

In any case the lengths of the words $w_i$ are recursively bounded
in terms of $|q|+|p|$ and the lemma is proved.
\endproof

\begin{lemma} \label{qspiral}
Let a minimal diagram $\Delta$ be a quasispiral having no
$k$-cells. Then the length of a shortest $a$-band of $\Delta$,
which connects $p$ and $q$ is bounded by a recursive function of
$|q|+|p|,$ where $q, p$ are the inner and the outer contours of
$\Delta$.
 \end{lemma}

\proof One can define turns of the quasispiral as for spirals (in
the beginning of the proof of Lemma \ref{spiral} (but one should
choose a maximal $a$-band instead of $\bb_k$). All the turns have
equal boundary labels since every maximal $a$-band has one common
cell with every turn. (Recall that there are no terminal cells for
$a$-bands in $\Delta$, because the base of $\Delta$ is empty.)
Then one can apply Lemma \ref{2 bands} exactly as in Case 2 of the
proof of Lemma \ref{spiral} to complete the proof.
\endproof

\section{Rolls}
\setcounter{equation}{0} \label{Rolls}

A reduced annular diagram $\Delta$ over the group $\hhh_1$ is said
to be a \label{roll}{\it roll} if its inner contour $q$ and its
outer contour $p$ are minimal boundaries and have no $k$-letters
(see Figure \theppp).

\begin{center}
\unitlength=.95 mm \special{em:linewidth 0.5pt}
\linethickness{0.5pt}
\begin{picture}(86.00,84.33)
\bezier{344}(44.00,84.00)(86.00,84.33)(86.00,40.33)
\bezier{352}(38.33,6.00)(85.67,6.00)(85.67,46.33)
\put(44.67,42.33){\circle{14.00}}
\bezier{252}(14.33,42.67)(14.33,74.67)(45.33,74.67)
\bezier{240}(46.00,74.67)(74.67,74.00)(76.33,42.67)
\bezier{224}(14.33,42.67)(18.00,13.67)(45.00,15.33)
\bezier{240}(45.67,15.33)(76.67,14.00)(76.33,42.67)
\bezier{228}(18.00,41.67)(17.67,71.33)(45.00,71.33)
\bezier{212}(45.67,71.33)(71.33,69.00)(72.33,42.00)
\bezier{188}(18.00,41.67)(21.67,18.33)(45.00,19.33)
\bezier{208}(72.33,42.00)(69.67,15.33)(45.00,19.33)
\bezier{152}(28.00,42.33)(28.33,63.33)(45.67,64.00)
\bezier{144}(45.67,64.00)(60.67,62.67)(62.67,42.00)
\bezier{148}(63.00,42.00)(63.00,23.00)(45.00,23.33)
\bezier{128}(43.33,23.33)(28.00,25.00)(28.33,42.00)
\bezier{116}(30.67,43.00)(33.00,60.33)(44.67,61.00)
\bezier{120}(45.00,61.00)(56.67,60.67)(59.33,42.33)
\bezier{124}(31.00,43.00)(31.33,26.00)(45.00,26.33)
\bezier{116}(45.67,26.00)(58.33,26.00)(59.67,42.33)
\put(44.33,49.33){\line(0,1){34.67}}
\put(46.67,49.00){\line(0,1){35.00}}
\put(50.67,46.00){\line(3,4){22.33}}
\put(51.67,44.00){\line(3,4){22.67}}
\put(51.33,40.00){\line(5,-1){33.00}}
\put(50.33,38.67){\line(5,-1){33.67}}
\put(46.67,35.67){\line(0,-1){29.33}}
\put(44.33,35.33){\line(0,-1){29.33}}
\put(39.00,38.33){\line(-5,-4){26.00}}
\put(38.00,40.67){\line(-5,-4){26.67}}
\put(11.33,19.33){\line(0,0){0.00}}
\put(38.33,45.67){\line(-2,1){29.67}}
\put(39.00,47.00){\line(-2,1){29.67}}
\put(21.00,68.33){\makebox(0,0)[cc]{$\bb_2$}}
\put(30.67,59.67){\makebox(0,0)[cc]{$\bb_1$}}
\put(44.33,46.33){\makebox(0,0)[cc]{$\ttt_1$}}
\put(47.67,42.67){\makebox(0,0)[cc]{$\ttt_2$}}
\bezier{196}(8.33,60.28)(22.78,83.89)(43.61,83.89)
\bezier{124}(12.78,17.78)(21.67,6.11)(38.33,6.11)
\bezier{188}(12.78,17.78)(0.83,37.22)(8.33,60.28)
\end{picture}
\end{center}
\begin{center}
\nopagebreak[4] Figure \theppp.
\end{center}
\addtocounter{ppp}{1}

It follows that every maximal $\theta$-band of a roll $\Delta$
connects its inner boundary $q$ and outer boundary $p$. The cyclic
\label{historyr}{\em history} of a roll is the projection of
$\phi(p)$ on $\csss$ (it is equal to the projection of $\phi(q)$
on $\csss$). Clearly (by Lemma \ref{NoAnnul}) maximal $k$-bands of
a roll $\Delta$ form concentric annuli surrounding the hole of the
diagram.  It follows from the definition that arbitrary annular
subdiagram  of a roll bounded by some sides of $k$- or $a$-annuli
is a roll if it contains a $(\theta,k)$-cell. In any case it is a
roll if it is bounded by sides of two $a$-annuli, since there are
no $a$-edges on these sides.

 We count the $k$-annuli of $\Delta$ from $q$ to
$p$: $\bb_1,...,\bb_r$. The bases of $\bb_1,...,\bb_r$ form the
\label{baser}{\em base} of the roll $\Delta$.

Our goal is to replace every roll by a roll with the same boundary
labels and  recursively bounded number of cells (in terms of
$|p|+|q|$).

\subsection{Rolls without a base}

Let $z\in \tilde{\cal K}$. We say that an $x$-edge ($a$-edge,
$\theta$-edge) $e$ is of {\em type} $z$ if the label of $e^{\pm
1}$ belongs to the set $\xxx(z)$ (to $\aaa(z)\cup\bar\aaa(z)$, to
$\Theta (z)\cup\bar\Theta(z)$).

\begin{lemma}\label{bez1} Let $\Delta$ be a roll with empty
base. Suppose that either

(1) the contours $p$ and $q$ are simple loops without common
vertices or

(2) all edges in $p\cup q$ are of the same type $z$ for some
$z\in\tilde{\cal K}$.

Then the labels of all edges in $\Delta$ are of the same type.

\end{lemma}

\proof It is clear from the list of relations of $\hhh_1$, that
all edges of a non-$k$-cell have the same type. We say that two
cells of $\Delta$ are neighbor if they have a common edge. It
follows that the transitive closure of this relation is an
equivalence, and all edges of the cells of an equivalence class
$\Gamma$ are of the same type.

Let $\Gamma_1,\dots,\Gamma_s$ be all the classes of $\Delta$. If
$s>1$, then the deletion of finitely many  vertices from
$\cup_{i=1}^s \Gamma_i$ makes $\Delta$ disconnected, a
contradiction in case (1). Hence $s=1$ in this case, and the
statement is true because every edge of $p\cup q$ must belong to a
cell of $\Gamma_1$. In case (2) the statement follows from the
observation that any boundary edge of each of the components
$\Gamma_1,\dots,\Gamma_s$ must lie on $p\cup q$.

\endproof

\begin{lemma}\label{bez bazy}
There is a recursive function $f$ with the following property.

Let $\Delta$ be a roll with an empty base. Then there exists a
roll $\Delta'$ with the same boundary labels as $\Delta$, which
has a path of length at most $f(|p|+|q|)$, connecting the
components $p$ and $q$ of the boundary of $\Delta'$.
\end{lemma}

\proof We can assume that $\Delta$ contains $\theta$-edges since
otherwise one can refer to Lemma \ref{Conjugacy}. By Lemma
\ref{rollsprep}, we can assume that $\Delta$ satisfies properties
(R1) - (R4).

{\bf Case 1}: Assume that there is at least one $\aba{P_1}$-edge
in $\Delta$. By lemmas \ref{bez1} and \ref{K0}, we can assume that
all $a$-edges in $\Delta$ belong to $\aba{P_1}$, all
$\theta$-edges belong to $\ata{P_1}$. Since the set $\xxx(P_1)$
does not exist, $\Delta$ does not have $x$-edges.

Thus all cells in $\Delta$ correspond to the commutativity
relations of the form (\ref{auxtheta}) and $\gr$-relations.
Therefore $\Delta$ is a diagram over the direct product of the
free group generated by $\Theta$ and the free product of the copy
of $\gr$ generated by $\aaa(P_1)$ and the free group with basis
$\{a_{m+1},...,a_{\bar m}\}$. Since the conjugacy problem is
decidable in both factors of this direct product, we can replace
$\Delta$ by a diagram with the same boundary label which has a
path, connecting the inner and outer contours, of length bounded
by a recursive function in $|p|+|q|$, as required.

{\bf Case 2.} Suppose that $\Delta$ has no $\aaa(P_1)$-edges.

Again by Lemma \ref{bez1}, for some $z\in \tkk$, $z\neq P_1$, all
edges from $\Delta$ are of type $z$.

{\bf Case 2.1.} First we consider the case where $p$ and $q$ have
no $a$-edges. By Lemma \ref{NoAnnul}, $\Delta$ is a union of
concentric $a$-annuli $\CC^1, \CC^2, \dots$.  If we find a
recursive upper bound for the lengths of the $a$-annuli, we
recursively bound their number because if the sides of two such
annuli have equal labels, we can remove the annuli and all the
cells between them to obtain a smaller roll. Since all the
$a$-annuli have equal numbers of $\theta$-cells (each of them
intersects each of the $\theta$-bands only once), it suffices to
bound the lengths of the subbands of $\CC^i$ between any two
neighbor maximal $\theta$-bands $\ttt$ and $\ttt'$. Let $\Gamma$
be the subdiagram situated between $\ttt$ and $\ttt'$. Clearly all
cells in $\Gamma$ correspond to the relations (\ref{auxax}).

We can assume that $z\ne P_j$ for any $j$ since there are no
auxiliary $(\aaa(P_j),x)$-relations (\ref{auxax}). Hence $z\in
\{K_j, L_j, R_j\}$ for some $j$.

We consider only the first of these cases because all other cases
are similar. Thus assume that $z=K_j$ for some $j$.

The bands $\ttt$ and $\ttt'$ are reduced, and so $\Gamma$ has
boundary label of the form $U_1V_1U_2V_2$ where $V_1$ and $V_2$
are labels of sides of $\ttt$ and $\ttt'$, and $U_1$ and $U_2$ are
words in $\xxx$. Clearly the $a$-projections of $V_1$ and $V_2\iv$
are graphically the same since every $a$-band starting on $\ttt$
ends on $\ttt'$ and vise versa. Let $V=(V_1)_a$.

Lemma \ref{-+} can be applied to a minimal subdiagram $\Gamma'$ of
$\Gamma$, containing two neighbor $a$-bands. It implies that the
word $V$ is a product $V'(V'')^{-1}$ for some positive words $V',
V''$. Now Lemma \ref{VV'} implies that all conditions of Lemma
\ref{slender} hold for the union of $\Gamma'$ and the
corresponding portions of the $\theta$-bands $\ttt$, $\ttt'$.
Lemma \ref{slender} gives the desired recursive bound for the
lengths of $a$-bands in $\Gamma$ in terms $|U_1|, |U_2|$.

\begin{center}
\unitlength=1mm \special{em:linewidth 0.4pt} \linethickness{0.4pt}
\begin{picture}(113.33,79.67)
\bezier{156}(39.33,50.00)(21.00,47.00)(21.00,26.67)
\bezier{128}(21.33,30.33)(23.33,15.00)(39.67,15.33)
\bezier{112}(39.67,15.33)(53.67,15.33)(55.67,29.00)
\bezier{148}(55.67,29.00)(57.33,47.67)(39.67,50.00)
\put(47.00,43.67){\makebox(0,0)[cc]{$q$}}
\bezier{356}(86.67,70.67)(113.33,40.67)(86.67,0.00)
\put(44.00,49.00){\line(2,1){42.67}}
\bezier{60}(86.67,70.33)(81.33,79.00)(76.67,79.67)
\put(50.33,46.00){\line(2,1){40.00}}
\put(55.67,27.67){\line(3,-1){39.33}}
\put(52.67,20.67){\line(3,-1){39.00}}
\put(87.33,12.67){\makebox(0,0)[cc]{$\ttt$}}
\put(83.33,65.00){\makebox(0,0)[cc]{$\ttt'$}}
\put(73.33,53.00){\makebox(0,0)[cc]{$V_2$}}
\put(74.00,25.67){\makebox(0,0)[cc]{$V_1$}}
\put(58.67,36.67){\makebox(0,0)[cc]{$U_2$}}
\put(82.67,38.00){\makebox(0,0)[cc]{$\Gamma$}}
\put(105.67,39.00){\makebox(0,0)[cc]{$U_1$}}
\end{picture}

\end{center}

\begin{center}
\nopagebreak[4] Fig. \theppp.

\end{center}
\addtocounter{ppp}{1}

{\bf Case 2.2} Now let us consider the case when
$\partial(\Delta)$ contains $a$-edges. In this case, there is an
$a$-band $\CC$ starting on $p$ and ending on $q$ since $\Delta$
satisfies (R2). So we can apply Lemma \ref{sp123} and represent
$\Delta$ as a union of a quasispiral, and two subdiagrams with
recursively bounded number of cells. Then the sub-quasispiral of
$\Delta$ will have a recursively bounded perimeter, and we can
apply Lemma \ref{qspiral}.\endproof

\begin{lemma} \label{a-annulus} Let $\Delta$ be a roll with the empty base, satisfying (R3),
which has no $\gr$-cells. Let $p$ and $q$ be the outer and inner
components of $\partial(\Delta)$ and $\Lab(p)_a$ freely equal to
the empty word. Then either there is an $a$-annulus between $p$
and $q$ surrounding the hole of $\Delta$, or the lengths of the
$\theta$-bands of $\Delta$ are recursively bounded in terms of the
number of $a$-letters in $\Lab(p)$ and $\Lab(q)$.
\end{lemma}

\proof By Lemma \ref{rollsprep}, we can assume that no $a$-bands
in $\Delta$ start and end on $p$ (resp. $q$). Therefore we can
assume that all $a$-bands in $\Delta$ connect $p$ and $q$
(otherwise it would contain an $a$-annulus).

Since the $a$-projection of $\phi(p)$ is freely trivial, there is
a subpath $p'$ of $p$ which starts with an edge $e$, ends with
$e'$ and $\phi(e)=\phi(e')\iv\in \aaa\cup\bar\aaa$, and there are
no $a$-edges in $p'$ between $e$ and $e'$. Consider the $a$-bands
$\CC$ and $\CC'$ starting on $e$ and $e'$ respectively. Suppose
that both bands intersect the same $\theta$-band $\ttt$. Then the
intersection cells $\CC\cap\ttt$ and $\CC'\cap\ttt$ must have a
common $\theta$-edge, and their $a$-edges must have mutually
inverse labels. Thus these cells cancel, a contradiction. Hence,
say, $\CC$ does not intersect all maximal $\theta$-bands in
$\Delta$. This and the Jordan lemma implies that none of the
$a$-bands in $\Delta$ intersects any of the maximal $\theta$-bands
there more than ones. Therefore by Lemma \ref{NoAnnul} each
$a$-band intersects each $\theta$-band in $\Delta$ at most once.
Hence the length of each $\theta$-band in $\Delta$ does not
exceed the number of $a$-edges on $p$.\endproof

\subsection{Rolls with $L$-annuli}
\label{rla}

\begin{lemma} \label{M-a} Let a reduced annular diagram $\Delta$ be an $L_j$-annulus with an
outer boundary $p$ and inner boundary $q$, whose cyclic history is
a word of the form $\tau\tau^{-1}\tau\tau^{-1}\dots$ for some
$\tau\in\sss^+(2)\cup\sss^+(4)$. Assume that all $x$-edges of $p$
are labelled by the same letter $x=x(a,\tau)$, $a\in \aaa(L_j)$.
Then there exists an annular diagram $\Delta'$ which contains no
$k$-cells with the outer path $p'$ and inner path $q'$, where
$\phi(p')\equiv \phi(p)$ and $\phi(q')$ can be obtained from
$\phi(q)$ by replacing every $a$-, $x$-, and $\theta$-letter by
its brother in $\aaa(L_j)\cup\xxx(L_j)\cup\Theta(L_j)$. Similarly
there also exists an annular diagram $\Delta''$ with outer contour
$p''$ and inner contour $q''$, without $k$-cells where all $a$-,
$x$-, and $\theta$-labels of the boundary belong to
$\aaa(\lel_j)\cup\xxx(\lel_j)\cup\Theta(\lel_j)$, $\phi(q'')\equiv
\phi(q)$, $\phi(p'')$ is a copy of $\phi(p)$.
\end{lemma}

\proof Assume that $\tau\in \sss^+(2)$ (the other case is
similar). Consider an $(L_j,\theta)$-cell $\Pi$ of $\Delta$ with
boundary label of the form
$$\theta_1^{-1}L_j\theta_2x^{-1}a
L_j^{-1}(a')^{-1}x'$$ (see relations (\ref{mainrel})) where $a\in
\aaa(L_j)$, $a'\in\aaa(\lel_j)$, $x\in\xxx(L_j)$, $x'\in
\xxx(\lel_j)$. We substitute the occurrence $L_j$ by $x'a'$, the
occurrence $L_j^{-1}$ by $(a')^{-1}x'$, the letter $\theta_2$ by
its copy $\theta_1$, the letter $a$ by $a'$ and $x$ by $x'$. Then
we obtain a relation of the form (\ref{auxtheta}). We can
substitute the $k$-cell $\Pi$ by the new cell $\Pi'$ without
$k$-edges.

Then, for every $(L_j,x)$-cell $\pi$ (see relations
(\ref{auxkx})), situated between $(L_j,\theta)$-cells $\Pi_1$ and
$\Pi_2$, we change the letters $L_j$ in its label by $x'a'$ if the
subtrapezium, which includes $\Pi_1$ and $\Pi_2$, has history
$\tau^{-1}\tau$, or by $(x')^{-1}a'$ if the history is
$\tau\tau^{-1}$. The letter $x$ in the boundary label of $\pi$
will be replaced by $x'$. We obtain the boundary label of the form
$(x')^4(x')^{\pm 1}a'(x')^{-1}(a')^{-1}x^{\mp 1}$ which is equal
to 1 modulo relations (\ref{auxax}).

Since the cells of $\Delta$ were glued in the $L_j$-band along
$L_j$-edges, the substitutions we have just described provide us
with the desired diagram $\Delta'$.
\endproof

\begin{lemma} \label{MM-a}
Let a roll $\Delta$ with base $L_jL_j\iv$ or $L_j\iv L_j$ be
bounded by two $L_j$-annuli $\bb_1$ and $\bb_2$. Assume that the
non-empty history of the roll is $\tau\tau^{-1}\tau\tau^{-1}\dots$
for some $\tau\in\sss^+(2)\cup\sss^+(4)$. Suppose that all
$x$-edges of the inner boundary of $\bb_1$ have the same labels
and all $x$-edges of the outer boundary of $\bb_2$ have the same
labels $x(a,\tau)$ (for some $a$). Then there exists a roll
$\Delta'$ which has no $k$-cells and the same boundary labels as
$\Delta$.
\end{lemma}

\proof By Lemma \ref{NoAnnul}, every $\theta$-band of $\Delta$
crosses each of $\bb_1$, $\bb_2$ only once. Therefore one can
apply Lemma \ref{No a-cells} (part 2) to conclude that $\Delta$
contains no $\gr$-cells. Besides, by Lemma \ref{bez1} (assumption
2 of Lemma \ref{bez1} holds), the subroll with trivial base
bounded by $\bb_1$ and $\bb_2$, consists of $\aaa(L_j)$-cells or
it consists of $\aaa(\lel_j)$-cells. Notice that for every
$\aaa(L_j)$-cell, can be constructed its $\aaa(\lel_j)$-copy, and
vice versa (see relations (\ref{auxtheta}), (\ref{auxax}) and the
definition of the mapping $\alpha_{\tau}$). Then one can change
all these cells by their $\aaa(\lel_j)$-copies or, respectively,
by their $\aaa(L_j)$-copies and simultaneously transform the bands
$\bb_1$ and $\bb_2$ according to Lemma \ref{M-a}. Reducing the
resulting diagram, we obtain the desired roll $\Delta'$.
\endproof

\begin{lemma} \label{MivM}
Let $\Delta$ be a roll with a non-empty history $h$ which is a
word in $\sss(2)\cup\sss(4)$. Assume that the base of $\Delta$ is
$L_j\iv L_j$ for some $j$, and $\Delta$ is bounded by two
$L_j$-annuli $\bb_1$ and $\bb_2$. Then $h$ is a word of the form
$\tau\tau\iv\tau\dots$ for some $\tau\in\sss$, all $x$-edges in
the boundaries $p$ or $q$ are labelled by the same letter
$x=x(a,\tau)^{\pm 1}$, and all $x$-edges not in $p\cup q$ are also
labelled by the same letter $x'=x(a',\tau)^{\pm 1}$.
\end{lemma}

\proof Assume that the history $h$, as a cyclic word, has two
neighbor letters $\tau, \tau'$ which are not mutually inverse.
Consider the small subtrapezium $\Gamma$ of $\Delta$ with history
$\tau\tau'$ and a two-letter base. Since the history of $\Gamma$
is reduced, by Lemma \ref{Sm Tr}, the projection of the label of
the bottom of the trapezium $\Gamma$ on $\aaa\cup\kkk$ is an
admissible word for $\sss$. But this contradicts Lemma
\ref{noPivP} (an admissible word for $\sss$ cannot have base
$L_j\iv L_j$). Hence $h$ does not have a reduced subword of length
2, whence it has the form $\tau\tau\iv\tau\dots$ for some $\tau$.

Let $\Delta'$ be the roll with empty base obtained by removing
$\bb_1$, $\bb_2$. Since the history of the roll has the form
$\tau\tau\iv...$, there exists $i$ such that all $a$-edges on the
boundary of the roll have labels $a_i(z)^{\pm 1}$ for some $i$
where $z=\lel_j$.
 By Lemma \ref{bez1}(2), all edges of $\Delta'$ are of type
 $\lel_j$.


Let $\Gamma$ be a small subtrapezium of $\Delta$ bounded by two
maximal $\Theta$-bands $\ttt_1, \ttt_2$ of $\Delta$ starting on
$q$ and ending on $p$. Denote $\pi_1=\ttt_1\cap\bb_1$ and
$\pi_2=\ttt_2\cup\bb_1$.

We assume that $\ttt_2$ follows $\ttt_1$ in the clockwise order.
Let $V\equiv\phi(\bott(\ttt_1)), V'\equiv\phi(\topp(\ttt_2))$.
Since every maximal $a$-band of $\Gamma$ starting on
$\topp(\ttt_1)$ ends on $\bott(\ttt_2)$, the $a$-projections $V_a$
and $V'_a$ are identical and each of these words contains a letter
$b^{\pm 1}=\in\aaa(\lel_j)^{\pm 1}$.

Suppose that the first $a$-letter of the word $V_a$ is negative
and the last one is positive. Then $V_a$ has a subword of the form
$a_1\iv a_2$ which is ruled out by Lemma \ref{-+}.

Therefore either the first letter of $V_a$ is positive or the last
letter is negative. These two cases are similar (one can be
obtained from the other one by interchanging of $p$ and $q$).

Thus we may assume that the first letter of $V_a$ is positive.
Hence by Lemma \ref{-+}, $V=V_1V_2$ where $(V_1)_a$ is a non-empty
positive word and $(V_2)_a$ is a negative word.

Again, we consider the first $a_s(\lel_j)$-band $\CC$ in $\Gamma$,
counting from $q$, connecting $\bott(\ttt_1)$ and $\topp(\ttt_2)$.
Let $\Psi$ be the diagram containing $\CC$ and bounded by a part
of $\bott(\ttt_1)$, $\bott(\CC)$, a part of $\topp(\ttt_2)$ and a
part of the outer contour of $\bb_1$.

The boundary label of $\Psi$ is of the form $W\equiv
U_1a_s(\lel_j)^{-1}x^{\pm 1 }U_2x^{\mp 1}a_s(\lel_j)$, where $U_1$
and $U_2$ are some words in $\xxx$, $U_1$ is the top label of the
$a_s(\lel_j)$-band $\CC$ and $U_2$ is read on the outer contour of
$\bb_1$ between two $L_j$-cells. Hence the word $U_2$ is freely
equal to a product of fourth powers of $x$-letters labelling
$x$-edges of the outer side of $\bb_1$ (by relations
(\ref{auxkx})). This word cannot be empty, because otherwise the
two $L_j$-cells $\pi_1, \pi_2$ would form a reducible pair. The
letter $x$ labels either an $x$-edge of the $L_j$-cell $\pi_1$
(where $\ttt_1$ starts) or an $x$-edge of an auxiliary
$\theta$-cell $\pi$ of $\ttt_1$ (that could happen either if the
cell $\pi_1$ has no $a$-edges on its top side, or if the first
auxiliary cell in $\ttt_1$ shares its $x$-edge with $\pi_1$). We
shall refer to these possibilities as {\bf Cases 1} and {\bf 2}.

Notice that the letter $x$ labels an $x$-edge of a $\theta$-cell.
Therefore it corresponds to the $a$-letter appearing as a label of
an $a$-edge of this cell and the $S$-rule $\tau$. We shall prove
that all $x$-edges in $\Delta'$ are labelled by $x^{\pm 1}$, and
that all $a$-edges on the sides of the $\theta$-bands of $\Delta'$
are labelled by $a_s(\lel_j)^{\pm 1}$. Then we shall prove that
$s=i$. This would imply that all $x$-edges on the boundary of
$\Delta$ have the same label.

The word $x^{\pm 1}U_2x^{\mp 1}$ is freely inverse of the bottom
label $U'_1$ of $\CC$. Therefore, by Lemma \ref{x-power}(2),
$U_2=x^{4d}$ for some $d\ne 0$. In particular, the letter
$a_s(\lel_j)$ is determined by the $x$-letter $x'$ such that $x$
is a copy of $x'$ and $(x')^d$ is read on $p$ between $\pi_1$ and
$\pi_2$. Thus, $U_1=x^{-d}$.

\begin{center}
\unitlength=1mm \special{em:linewidth 0.4pt} \linethickness{0.4pt}
\begin{picture}(107.67,61.67)
\put(19.33,61.67){\line(0,-1){54.00}}
\put(19.33,7.67){\line(1,0){88.33}}
\put(19.33,16.00){\line(1,0){88.00}}
\put(19.33,52.67){\line(1,0){87.67}}
\put(19.33,43.33){\line(1,0){87.33}}
\put(29.33,43.33){\line(0,-1){27.33}}
\put(62.67,43.33){\line(0,-1){27.33}}
\bezier{124}(29.33,43.33)(36.67,31.33)(29.33,16.00)
\put(29.33,52.67){\line(0,1){8.67}}
\put(24.00,56.67){\makebox(0,0)[cc]{$\bb_1$}}
\put(25.00,47.00){\makebox(0,0)[cc]{$\pi_1$}}
\put(98.00,48.00){\makebox(0,0)[cc]{$\ttt_1$}}
\put(98.00,11.00){\makebox(0,0)[cc]{$\ttt_2$}}
\put(25.00,11.00){\makebox(0,0)[cc]{$\pi_2$}}
\put(26.00,32.00){\makebox(0,0)[cc]{$U_2$}}
\put(36.00,32.00){\makebox(0,0)[cc]{$U_1'$}}
\put(66.00,32.00){\makebox(0,0)[cc]{$U_1$}}
\put(51.00,22.00){\makebox(0,0)[cc]{$\CC$}}
\put(15.00,32.00){\makebox(0,0)[cc]{$q$}}
\end{picture}
\end{center}

\begin{center}
\nopagebreak[4] Fig. \theppp.

\end{center}
\addtocounter{ppp}{1}

If $|(V_1)_a|\ge 2$, one can almost repeat the argument (but
applying now part (3) of Lemma \ref{x-power} instead of part (2))
and conclude that the next to the first $a$-letter
$a_{s'}(\lel_j)$ of $(V'')_a$ is determined by $x$, and so $s'=s$.
Hence $(V_1)_a$ is a power of $a_s(\lel_j)$ and all $x$-edges of
the $a$-bands starting on the portion of $\bott(\ttt_1)$ labelled
by $V_1$ have the same label $x$. Similarly, if $(V_2)_a$ is
non-empty, it must be a power of some $a_{s'}(\lel_j)$, and all
$x$-edges on $a$-bands starting on the part of $\bott(\ttt_1)$
labelled by $V_2$ have the same label $x''$. The neighbor
$a_s(\lel_j)$- and $a_{s'}(\lel_j)$-bands, connecting
$\bott(\ttt_1)$ and $\topp(\ttt_2)$, have a common boundary
sections labelled by non-empty $x^{t_1}\equiv (x'')^{t_2}$,
because otherwise a cell of $\ttt_1$ and a cell of $\ttt_2$ would
form a reducible pair. We conclude that $x\equiv x''$ and hence
$s=s'$. Therefore $V_2$ is empty. Similarly
$\Lab(\topp(\ttt_2))_a$ is a power of $a_s(\lel_j)$.

Considering the next (in the clockwise order) small subtrapezium
$\Gamma'$ of $\Delta'$, the small subtrapezium after that, and so
on, we conclude that all $a$-edges on the (reduced) top and bottom
sides of all $\theta$-bands in $\Delta'$ have the same label
$a_s(\lel_j)$ and all $x$-edges in all small subtrapezia of
$\Delta'$ are labelled by $x$.

It remains to show that $s=i$. Indeed, if Case 1 holds, that is if
$\pi_1$ does not share its $x$-edge with the first auxiliary cell
$\pi$ of $\ttt_1$, then the first $a$-edge of $\bott(\ttt_1)$ is
an edge of $\pi_1$, and so $s=i$. Suppose that Case 2 holds. Then
the $a$-edges of $\pi$ are labelled by $a_i(\lel_j)^{\pm 1}$.

The bottom side of $\pi_1$ (considered as a $\theta$-band)
consists of one $L_j$-edge (see relations (\ref{mainrel})), so the
bottom side of $\pi$ does not have common edges with the bottom
side of $\pi_1$. Therefore the $a$-edge of the bottom side of
$\pi$ appears on the (reduced) bottom path of $\ttt_1$. Therefore,
again, $s=i$. This completes the proof of the lemma.
\endproof

\begin{lemma} \label{24} Let $\Delta$ be a roll
with non-empty base whose history is a word over
$\sss(2)\cup\sss(4)\cup\bar\sss$. Then either all rules in the
history belong to $\bar\sss$ or all of these rules are from
$\sss(2)\cup\sss(4)$.
\end{lemma}

\proof Indeed, consider a $k$-annulus $\bb$ of $\Delta$. If the
history of $\Delta$ contains rules both from $\bar\sss$ and from
$\sss(2)\cup\sss(4)$ then there exist two neighbor cells in $\bb$,
one corresponding to a rule in $\bar\sss$ and the other
corresponding to a rule from $\sss(2)\cup\sss(4)$. But this is
impossible since rules from $\sss(2)\cup\sss(4)$ do not change the
$\Omega$-coordinate of the base letters, and letters in $\kkk\cap
\bar\kkk$ have $\Omega$-coordinates 1.
\endproof

\begin{lemma} \label{bezMivM}
Let $\Delta$ be a roll with base $L_j^{-1}L_j$ for some $j$.
Assume that the non-empty history $h$ of $\Delta$ is a word in
$\sss(2)\cup\sss(4)\cup\bar\sss$. Then there exists a roll
$\Delta'$ with an empty base, having the same boundary labels as
$\Delta$.
\end{lemma}

\proof One may assume that $\Delta$ is bounded by two $L_j$-annuli
$\bb_1$ and $\bb_2$. By Lemma \ref{24}, there are two
possibilities: either all rules in the history of $\Delta$ are
from $\bar\sss$ or all of them are from $\sss(2)\cup\sss(4)$.

In the first case the inner side of $\bb_1$ has the same label as
the outer side of $\bb_2$ since there exists no $(\bar
a,x)$-cells. Hence the role $\Delta$ can be replaced by an empty
roll.

In the second case, one can apply Lemma \ref{MivM} and then Lemma
\ref{MM-a}. \endproof

\begin{lemma} \label{M-rolls}
Let $\Delta$ be a roll. Assume that the history $h$ of $\Delta$ is
non-empty and the base of $\Delta$ is either (a) equal to
$L_j^{-1}L_j$ or (b) equal to $L_jL_j^{-1}$ or (c) equal to $L_j$
or for some $j$. Then there exists a roll $\Delta'$ with the same
boundary labels as $\Delta$, whose base is obtained from the base
of $\Delta$ by removing some (or none) letters, and and the number
and perimeters of cells are recursively bounded in terms of
lengths of the boundary components $|p|,|q|$ (as usual, $p$ is the
outer and $q$ is the inner boundary components).
\end{lemma}

\proof 1. By Lemma \ref{rollsprep} we can assume that the roll
$\Delta$ satisfies the conditions (R1) and (R2) of that lemma. We
assume that $\Delta$ has the lowest possible number of
$L_j$-annuli and the lowest type among all rolls satisfying the
conditions of the lemma, plus conditions (R1), (R2), and having
the same boundary labels as $\Delta$.

If we remove the $L_j$-bands from $\Delta$, we get two or three
rolls with empty bases. Two of these rolls, say, $\Delta_1$ and
$\Delta_2$ are bounded by an $L_j$-annulus and $q$ and an
$L_j$-annulus and $p$, respectively, and the third roll $\Delta_3$
(if exists) is bounded by two $L_j$-annuli.

Notice that it is enough to be able to replace $\Delta$ by a roll
with the same boundary labels and recursively bounded lengths of
$L_j$-bands. Indeed in that case we would be able to apply Lemma
\ref{bez bazy} to the subrolls $\Delta_1$, $\Delta_2$, $\Delta_3$.

2. Fix a big enough recursive function $f(x)$. From what follows
it will be clear how big $f(x)$ should be: $f(x)$ has to satisfy
certain (finite) number of inequalities of the form $f(x)>g(x)$
where $g(x)$ is a recursive function.

3. Suppose that one can connect the boundary components of
$\Delta_i$, $i=1,2$, by a path of length $\le f(|p|+|q|)$. Then by
Lemma \ref{tdist}, the number of cells in $\Delta_i$ is
recursively bounded in terms of $|p|, |q|$. Thus we can remove
one of the subrolls $\Delta_1$ or $\Delta_2$ and one of the
$L_j$-annuli from $\Delta$ and reduce the problem to the same
problem about a roll with fewer $L_j$-annuli.

Thus we can assume that there are no paths of length $\le
f(|p|+|q|)$ connecting the boundary components of $\Delta_i$,
$i=1,2$ (if $\Delta_i$ does exist). By Lemma \ref{K0}, we may also
assume that $p$ and $q$ are simple loops.By Lemma \ref{bez1}, we
conclude that for some $z\in \{L_j, \lel_j\}$, all edges of the
diagram $\Delta_i$ are of the same type $z$. In particular,
$\Delta$ does not contain $\gr$-cells.

4. As a consequence of 3, we can deduce that the lengths of
maximal $\theta$-bands in $\Delta_1$ and $\Delta_2$ are bigger
than $1$ (since $f(x)>1$). Hence for every $\tau$ in the history
of $\Delta$ there exist relations of the form (\ref{auxtheta})
containing a letter from $\Theta(\lel_j)\cup\bar\Theta(\lel_j)$.
This implies that $\tau$ does not lock $(zL_j)$-sectors and
$(L_jz)$-sectors. So every rule $\tau$ in the history of $\Delta$
belongs to $\sss(2)\cup\sss(4)\cup\bar\sss(2)\cup\bar\sss(4)$.
Hence if one of the rules in that history belongs to $\bar\sss$,
all of them belong to $\bar\sss$, which implies that $\Delta_i$
($i=1,2$) does not contain $x$-edges at all (there are no letters
in $\xxx$ that correspond to $\bar\aaa$-letters). This, in turn,
would imply that the lengths of the $L_j$-annuli in $\Delta$ equal
the number of maximal $\theta$-bands starting on $p$ and ending on
$q$, hence the lengths of $L_j$-annuli in $\Delta$ are bounded in
terms of $|p|,|q|$, as desired. Hence we can assume that all rules
$\tau$ in the history of $\Delta$ belong to $\sss(2)$ or
$\sss(4)$.

5. Now if the base of $\Delta$ is $L_j\iv L_j$ (as in case (a) of
the lemma), we can apply Lemma \ref{bezMivM} (to obtain a roll
with the same boundary labels as $\Delta$ but with empty base)
and then use Lemma \ref{bez bazy}.

6. Thus suppose that the base of $\Delta$ is $L_j$ or $L_jL_j\iv$.
Then $z_1=K_s^{\pm 1}$ for some $s$.

Consider two arbitrary consecutive (in the clockwise order)
maximal $\theta$-bands $\ttt_1$ and $\ttt_2$ in
$\Delta_1\cup\bb_1$ starting on the inner contour $q$. Let
$\Gamma$ be the subdiagram of $\Delta_1\cup\bb_1$ bounded by
$\bott(\ttt_1)$ and $\topp(\ttt_2)$. Fewer than $|q|$ of the
maximal $a$-bands in $\Gamma$ that start on $\topp(\ttt_1)$ end
on $q$. The total number of cells in these bands is recursively
bounded (we can remove these bands one by one starting with the
band which shares one of its sides with $q$). Thus we can remove
these cells from $\Delta$. Similar operation can be performed on
all other subdiagrams of $\Delta_1$ situated between two
consecutive $\theta$-bands, since the number of such diagrams is
at most $|q|$. Thus we can assume that every maximal $a$-band in
$\Gamma$ starting on $\bott(\ttt_1)$ ends on $\topp(\ttt_2)$. Let
$\CC$ be the first of these $a$-bands (counting from $q$ to
$\bb_1$). Then the bottom side of $\CC$ is a part of $q$.

Note that the rules of $\sss(2)\cup\sss(4)$ have the form $[...,
K_s\to K_s,...]$. Therefore the corresponding relations
\ref{mainrel} are of the form

$$ \theta(\tau,(K_s)_-)\iv K_s(r,\omega)\theta(\tau,K_s) = K_s(r, \omega) $$
where $\omega\in \{2,4\}$, $r\in\bee$. Using cells corresponding
to these relations and relations of the form (\ref{auxkx}), we can
build a $K_s$-band $\bb$ whose top side has the same label as the
bottom side of $\CC$. Then we can attach $\bb$ to the bottom side
of $\CC$ and obtain a small trapezium
$\tilde\Gamma=\Gamma\cup\bb$.

7. Suppose that the history of $\tilde\Gamma$ (= the history of
$\Gamma$) is reduced. Then by Lemma \ref{Sm Tr}, the labels of
the top and bottom paths of $\tilde\Gamma$ are admissible words.
In particular it means that the $a$-projection of the word $V$
(resp. $V'$) written on the bottom (resp. top) path of the
trapezium $\tilde\Gamma$ is positive. But then Lemma
\ref{slender} (iii) implies that the length of $\ttt_1$ does not
exceed $2|\CC|+2$. This contradicts our assumption that it is
greater than $f(|p|+|q|)$.

Hence the history of $\Gamma$ is not reduced. Since $\Gamma$ was
chosen arbitrarily, we deduce that the history of $\Delta$ has the
form $\tau\tau\iv...$ for some $\tau\in \sss(2)\cup\sss(4)$.

8. Let $V=\Lab(\bott(\ttt_1))$. By Lemma \ref{-+}, the word the
$a$-projection $V_a$ is equal to $V'V''$ where $V'$ is positive,
and $V''$ is negative. Accordingly, we can subdivide $\ttt_i$
($i=1,2$) into two subbands $\ttt_i'$ and $\ttt_i''$ such that
$$\Lab(\topp(\ttt_1'))_a=\Lab(\bott(\ttt_2')_a=V',
\Lab(\topp(\ttt_1''))_a=\Lab(\bott(\ttt_2'')_a=V'',$$ and we can
subdivide $\Gamma$ into two subdiagrams $\Gamma'$ and $\Gamma''$
where $\Gamma'$ is formed by the maximal $a$-bands of $\Gamma$
starting on $\topp(\ttt_1')$. Notice that by Lemma \ref{slender},
the length of of each of the $a$-bands in $\Gamma'$ is recursively
bounded in terms of $|q|$, and the number of these $a$-bands is
recursively bounded as well. Hence the number of cells in
$\Gamma'$ is recursively bounded, and we can assume that
$\Gamma=\Gamma''$ (and that the same is true for all other
subdiagrams of $\Delta_1$ situated between two consecutive
$\theta$-bands).

9. As in the proof of Lemma \ref{MivM}, we can deduce that the
$x$-edges in $\Delta_1$ have the same labels and the word $V'$ is
a power of some letter $a=a_i(\lel_j)$. Notice that the letter $a$
does not depend on the choice of $\Gamma$. Thus for every maximal
$\theta$-band $\ttt$ of $\Delta_1$, $\Lab(\bott(\ttt))_a=a^k$ or
some negative $k$,  $|k|>f(|p|+|q|)$ depending on $\ttt$. By Lemma
\ref{a-annulus}, there exists an $a$-annulus $\CC$ in $\Delta_1$
such that the number of cells between $\CC$ and $q$ is recursively
bounded. Without loss of generality, we can assume that $q$ is the
inner contour of $\CC$. Since $k>>1$, there exist one more
$a$-annulus $\CC_1$ such that the inner contour of $\CC_1$ is the
outer contour of $\CC$.

10. Now suppose that the condition (b) of the lemma holds, that is
the base of $\Delta$ is $L_jL_j\iv$. Then we can repeat the above
argument for $\Delta_2$ and conclude that all $x$-edges of each
side of the second $L_j$-annulus in $\Delta$ have the same labels.
This allows us to apply Lemma \ref{MM-a}, eliminate the
$L_j$-annuli from $\Delta$, and then use Lemma \ref{bez bazy}.

11. Finally suppose that the base of $\Delta$ is $L_j$ as in case
(c) of the lemma.

Notice that the label of the inner boundary of $\bb_1$ has the
form $$\theta (x\iv a) x^{4k_1} (x\iv a)\iv \theta\iv
x^{4k_2}\theta (x\iv a)...x^{4k_l}=\theta x\iv a x^{4k_1} a\iv
x\theta\iv x^{4k_2}\theta x\iv a...x^{4k_l}$$ for some integers
$k_1,...,k_l$ where $a\in \aaa(z_1), x=x(a,\tau)$. The label of
the outer boundary of $\bb_1$ is freely equal to $$\theta' (x')\iv
a'(x')^{k_1} (a')\iv x'\theta'^{-1} (x')^{k_2}\theta' (x')\iv
a'...(x')^{k_l}$$ where $\theta'$ is the ``brother" of $\theta$,
$a'$ is a ``brother" of $a$, $x'=x(a',\tau)$. Notice that for
every choice of parameters $k_1,...,k_l$ there exists an
$L_j$-annulus with boundary labels as above.

Modulo relations (\ref{auxax}) any word of the form
$ax(a,\tau)^ta\iv$ is equivalent to $x(a,\tau)^{4t}$. Therefore
for every integers $s_1,...,s_l$ divisible by $16$ there exists a
roll with base $L_j\iv$ and without $a$-annuli with outer boundary
label
$$\theta x^{s_1}\theta\iv x^{s_2} \theta x^{s_3} ... x^{s_l}$$
and inner boundary label
$$\theta' (x')^{s_1/4}\theta' (x')^{s_2/4}\theta'
x^{s_3/4}...x^{s_l/4}.$$ Let us denote this roll by
$\Psi(s_1,...,s_l)$. The ``inverse" roll, that is the roll
obtained from $\Psi(s_1,...,s_l)$ by switching the inner and outer
contours, will be denoted by $\Psi\iv(s_1,...,s_l)$. The base of
the inverse role is $L_j$.

\begin{center}
\unitlength=1mm \special{em:linewidth 0.4pt} \linethickness{0.4pt}
\begin{picture}(136.67,73.33)
\put(2.33,48.00){\line(1,0){20.33}}
\put(2.33,27.67){\line(1,0){20.00}}
\put(22.00,27.67){\line(0,1){20.33}}
\put(22.00,48.00){\line(2,3){17.00}}
\put(39.00,73.33){\line(0,-1){70.00}}
\put(39.00,3.33){\line(-2,3){16.33}}
\put(39.00,48.00){\line(1,0){25.33}}
\put(39.00,27.67){\line(1,0){25.00}}
\put(64.00,48.00){\line(0,1){25.33}}
\put(64.00,73.33){\line(4,-5){20.33}}
\put(84.33,48.00){\line(0,-1){20.33}}
\put(84.33,27.67){\line(-5,-6){20.33}}
\put(64.00,3.33){\line(0,1){24.33}}
\put(84.33,48.00){\line(1,0){19.00}}
\put(84.33,27.67){\line(1,0){18.67}}
\put(39.00,59.67){\circle*{1.33}}
\put(39.00,16.00){\circle*{1.33}}
\put(39.00,64.00){\vector(0,1){5.00}}
\put(39.00,56.33){\vector(0,-1){7.00}}
\put(39.00,18.33){\vector(0,1){6.33}}
\put(39.00,13.00){\vector(0,-1){7.00}}
\put(25.04,52.51){\vector(2,3){5.77}}
\put(27.18,20.65){\vector(2,-3){5.77}}
\put(25.67,60.67){\makebox(0,0)[cc]{$\theta'$}}
\put(25.67,16.00){\makebox(0,0)[cc]{$\theta$}}
\put(42.00,11.67){\makebox(0,0)[cc]{$x$}}
\put(42.00,20.67){\makebox(0,0)[cc]{$a$}}
\put(42.00,54.00){\makebox(0,0)[cc]{$a'$}}
\put(42.00,64.67){\makebox(0,0)[cc]{$x'$}}
\put(103.00,27.67){\line(0,1){20.33}}
\put(103.00,48.00){\line(2,3){17.00}}
\put(120.00,73.33){\line(0,-1){70.00}}
\put(120.00,3.33){\line(-2,3){16.33}}
\put(120.00,59.67){\circle*{1.33}}
\put(120.00,16.00){\circle*{1.33}}
\put(120.00,64.00){\vector(0,1){5.00}}
\put(120.00,56.33){\vector(0,-1){7.00}}
\put(120.00,18.33){\vector(0,1){6.33}}
\put(120.00,13.00){\vector(0,-1){7.00}}
\put(106.04,52.51){\vector(2,3){5.77}}
\put(108.18,20.65){\vector(2,-3){5.77}}
\put(106.67,60.67){\makebox(0,0)[cc]{$\theta'$}}
\put(106.67,16.00){\makebox(0,0)[cc]{$\theta$}}
\put(123.00,11.67){\makebox(0,0)[cc]{$x$}}
\put(123.00,20.67){\makebox(0,0)[cc]{$a$}}
\put(123.00,54.00){\makebox(0,0)[cc]{$a'$}}
\put(123.00,64.67){\makebox(0,0)[cc]{$x'$}}
\put(64.00,59.67){\circle*{1.33}}
\put(64.00,16.00){\circle*{1.33}}
\put(64.00,64.00){\vector(0,1){5.00}}
\put(64.00,56.33){\vector(0,-1){7.00}}
\put(64.00,18.33){\vector(0,1){6.33}}
\put(64.00,13.00){\vector(0,-1){7.00}}
\put(61.00,11.67){\makebox(0,0)[cc]{$x$}}
\put(61.00,20.67){\makebox(0,0)[cc]{$a$}}
\put(61.00,54.00){\makebox(0,0)[cc]{$a'$}}
\put(61.00,64.67){\makebox(0,0)[cc]{$x'$}}
\put(79.00,54.67){\vector(-3,4){4.67}}
\put(77.67,20.00){\vector(-3,-4){5.00}}
\put(80.33,58.33){\makebox(0,0)[cc]{$\theta'$}}
\put(80.33,17.00){\makebox(0,0)[cc]{$\theta$}}
\put(50.67,51.33){\makebox(0,0)[cc]{$(x')^{k_1}$}}
\put(50.67,23.67){\makebox(0,0)[cc]{$x^{4k_1}$}}
\put(93.33,23.67){\makebox(0,0)[cc]{$x^{4k_2}$}}
\put(93.33,51.67){\makebox(0,0)[cc]{$(x')^{k_2}$}}
\put(120.00,48.00){\line(1,0){16.67}}
\put(120.00,27.67){\line(1,0){16.33}}
\put(3.00,37.67){\makebox(0,0)[cc]{$\bb_1$}}
\put(64.00,48.00){\line(0,-1){20.33}}
\end{picture}

\end{center}

\begin{center}
\nopagebreak[4] Fig. \theppp.

\end{center}
\addtocounter{ppp}{1}

Notice that for some integers $s_1,...,s_l$ divisible by $16$ the
label of the inner contour of $\Psi\iv(s_1,...,s_l)$ coincides
with the label of the outer contour of $\CC_1$. Therefore we can
do the following surgery: cut the roll $\Delta$ along the outer
contour of $\CC_1$, insert in the hole a copy of
$\Psi(s_1,...,s_l)$ and $\Psi\iv(s_1,...,s_l)$ (the resulting
non-reduced ``roll" will have the base $L_jL_j\iv L_j$). Now
reduce the smallest subroll $\Delta'$ of the resulting diagram
with the base $L_j\iv L_j$.

There are two possibilities. Either after the reduction, the base
becomes $L_j$ and only the $L_j$-annulus of
$\Psi\iv(s_1,...,s_l)$ will remain, or the base of the new (still
non-reduced diagram will be $L_jL_j\iv L_j$ but the smallest
subroll with base $L_j\iv L_j$ will be reduced. Applying Lemma
\ref{bezMivM} to that subroll, we can get back to the first
possibility.

Thus we can construct another roll $\Delta''$ with base $L_j$
with the same boundary labels as $\Delta$, such that there are at
most two $a$-annuli between the $L_j$-annulus of $\Delta''$ and
the inner contour of $\Delta''$. By Lemma \ref{a-annulus} this
means that the number of cells between the $L_j$-annulus of
$\Delta''$ and the inner contour is recursively bounded. This
implies that the length of the $L_j$-annulus is recursively
bounded in terms of $|q|$, which completes the proof. \endproof

\subsection{Rolls with base $L_{j-1}\iv K_jL_j$}
\label{rwb}

A cyclic freely trivial word will be called a Dyck word. If $W$ is
a Dyck word in the alphabet $\{z_1^{\pm 1}, z_2^{\pm 1},...\}$
then for every scheme of cancellation of $W$, we can associate a
pairing of letters of $W$ (i.e. a selection of some pairs of
letters in $W$, each letter occurring in one pair), and a
placement of parentheses in $W$ which show in which order we
should cancel the letters. This pairing will be called
\label{cancpar}{\em cancellation pairing} in $W$.  For example, a
pairing in the cyclic word $b\iv abb\iv a\iv b$ is represented as
$b\iv)(a(bb\iv)a\iv)(b$. Each open parenthesis corresponds to
unique closed parenthesis. A pair of open and closed parentheses
corresponding to each other will be called \label{connectedp}{\em
connected}. It is be convenient to draw the word on the boundary
of a disc on the plane (we always read cyclic words in the
clockwise direction): One can connect the letters of each pair
drawn on the disk boundary, by arcs situated on the plane, say,
outside the disk, without intersections of these arcs.

\begin{center}
\unitlength=.8 mm \special{em:linewidth 0.4pt}
\linethickness{0.4pt}
\begin{picture}(136.00,80.84)
\bezier{248}(72.00,64.00)(39.00,61.00)(38.00,32.00)
\bezier{224}(38.00,32.00)(43.00,5.00)(72.00,5.00)
\bezier{256}(67.00,5.00)(104.00,6.00)(105.00,33.00)
\bezier{240}(105.00,33.00)(106.00,58.00)(72.00,64.00)
\put(69.67,63.67){\circle*{1.33}}
\put(101.00,48.33){\circle*{1.33}}
\put(93.67,12.00){\circle*{1.33}} \put(58.33,6.67){\circle*{1.33}}
\put(39.33,27.00){\circle*{1.33}}
\put(44.33,52.00){\circle*{1.33}}
\put(53.00,58.67){\vector(2,1){2.33}}
\put(84.00,60.88){\vector(3,-1){5.29}}
\put(104.77,30.54){\vector(0,1){3.46}}
\put(74.43,5.50){\vector(4,1){3.87}}
\put(48.57,12.01){\vector(-1,1){3.26}}
\put(40.01,43.58){\vector(-1,-4){1.22}}
\put(55.08,55.59){\makebox(0,0)[cc]{$a$}}
\put(77.48,58.64){\makebox(0,0)[cc]{$b$}}
\put(101.51,31.97){\makebox(0,0)[cc]{$b$}}
\put(75.85,9.16){\makebox(0,0)[cc]{$a$}}
\put(47.75,17.10){\makebox(0,0)[cc]{$b$}}
\put(41.85,38.49){\makebox(0,0)[cc]{$b$}}
\bezier{260}(87.46,59.66)(119.02,62.51)(104.97,32.17)
\bezier{296}(55.08,59.87)(92.55,80.84)(115.36,59.46)
\bezier{440}(115.33,59.33)(136.00,7.00)(82.67,7.00)
\bezier{320}(45.33,15.33)(6.33,17.00)(39.33,41.67)
\end{picture}
\end{center}

\begin{center}
\nopagebreak[4] Fig. \theppp.

\end{center}
\addtocounter{ppp}{1}

Then it is clear when a connected pair of parentheses or  letters
is {\em inside} another pair of parentheses (letters) or, vice
versa, is {\em outside}, i.e. {\em contains} it. From the
definition of a cancellation pairing, it is clear that no two open
(resp. closed) parentheses can stay next to each other because no
letter can cancel twice during the process of cancellation. Thus
there is a letter next in the clockwise (counterclockwise)
direction of every open (closed) parenthesis. If two parentheses
are connected then we call the corresponding letters
\label{connectedl}{\em connected}.

Let us call a pair of connected letters $(z_i^{\pm 1}, z_i^{\mp
1})$ \label{normalp}{\em normal}, if every pair of connected
letters containing this pair consists of different occurrences of
the same letter $z_i^{\pm 1}$.

For every cancellation pairing in a cyclic Dyck word $W$, a pair
of connected letters of the form $(z_i\iv, z_i)$ will be called
\label{minuspair}{\em minus pair}, a pair of the form $(z_i,
z_i\iv)$ will be called a \label{pluspair}{\em plus pair}. We call
a Dyck word $W$ a \label{minusw} {\em minus word} if for some
pairing (called \label{minusp}{\em minus pairing}) every pair of
connected letters  is a minus pair. Similarly one can define a
plus pairing and a plus word.

For example the cyclic word $abb^{-1}a^{-1}abb^{-1}a^{-1}$ is both
a plus word: $$(a(bb^{-1})a^{-1})(a(bb^{-1})a^{-1})$$ and a minus
word: $$a)b)(b^{-1}(a^{-1}a)b)(b^{-1}(a^{-1},$$ but it does not
have a pairing where every connected pair of letters is normal.

\begin{lemma}\label{A} If a minus word
$W$ contains a subword of the form  $z_i^{-1}z_j$, then $i=j$, and
the two letters of this subword are connected in any minus-pairing
of $W$.
\end{lemma}

\proof If these two letters are not connected in a given minus
pairing then they would be separated by a parenthesis. If is clear
that in any minus pairing a negative (positive) letter stays next
in the clockwise (resp. counterclockwise) direction of any open
(resp. closed) parenthesis. Therefore the parenthesis separating
$z_i\iv z_j$ must be both open and closed, a contradiction.
\endproof

Let $A$ and $B$ be alphabets and let $\phi: A\to B\cup\{1\}$ be a
map. Let $W$ be any cyclic group word in the alphabet $A$. Then
any pairing $\CC$ in $\phi(W)$ induces a pairing of some letters
in $W$: we simply pull the arrangement of parentheses in $\phi(W)$
to $W$. Notice that some letters can be not pared in $W$, and the
paired letters can be not mutual inverse, and, in general, the
induced pairing is not a cancellation pairing.

Let $h$ be a word over $\csss$. The projection $h_0$ of $h$ on
$\sss(2)\cup\sss(4)$ will be called the \label{24projection}{\em
$(2,4)$-projection} of $h$. If we further identify in $h_0$ rules
corresponding to the same letter $a_i$, $i=1,...,\bar m$ (such
rules will be called \label{asimilar}{\em $a$-similar}), then we
get a word $h_1$. If $h_1$ is a Dyck word then any cancellation
pairing in $h_1$ induces a pairing of letters in $h$. That pairing
will be called a \label{24pairing}$(2,4)$-{\em pairing} in $h$.
The $(2,4)$-pairing in $h$ will be called {\em $z$-good} for some
$z\in\tkk$ if no rule locking $zz_+$-sectors occurs inside a pair
of connected parentheses. A pairing is called \label{zbest}{\em
$z$-best} if no rule locking $zz_+$-sectors occurs in $h$.
Further, consider any word $W=X_1\theta_1X_2\theta_2...X_s$ where
$\theta_i\in \Theta(z)$, $z\in \{\lel_j, L_j\}$ for some $j$, the
words $X_1,...,X_s$ do not contain $\theta$-letters. Let $h$ be a
history of this word that is the projection of
$\theta_1...\theta_{s-1}$ onto $\csss$. Then any $(2,4)$-pairing
in $h$ induces a pairing in $W$. This pairing in $W$ will be
called \label{zgood}{\em $z$-good} ({\em $z$-best}) if it is
induced by a $z$-good ($z$-best) pairing in $h$.

\begin{lemma} \label{B}
Let $z\in \{\lel_j, L_j\}$ and

\begin{equation}\label{*}
w \equiv X_0 \theta_1 X_1 \dots \theta_s X_s
\end{equation}
for some $s$, where $\theta_i^{\pm 1}$ are letters of the form
$\theta(\tau,z)$ for some $z$, where $\tau$ does not lock
$zz_+$-sectors, $X_j$ are words in $\xxx(z)$. Let $a=a_i(z)$ for
some $i\in \{1,...,\bar m\}$. Then:

(1) Modulo relations (\ref{auxtheta}) and (\ref{auxax}), we have
$awa^{-1} = w'$ where

$$w' \equiv X'_0 \theta_1 X'_1...\theta_s X'_s ,$$
$X'_0,\dots,X'_s$ are words in $\xxx(z)$. The number of
applications of relations (\ref{auxtheta}) in the derivation of
this equality does not exceed $s$. If in addition
$\theta_t=\theta_{t+1}^{-1}=\theta(\tau,z)^{-1}$ for some index
$t$, and $X_t$ is a power of a letter $x=x(a,\tau)$ then $X'_t$
is again a power of $x$.

(2) If $\theta_1$ is a negative letter and $\theta_s$ is a
positive letter then, modulo relations (\ref{auxtheta}) and
(\ref{auxax}) the word $w$ is equal to the word

$$w_1=Y_0 (a^{-1}x\theta_1)Y_1\theta_2 Y_2 \dots
\theta_{s-1}Y_{s-1}(\theta_s (x')^{-1}a) Y_s,$$ where $Y_0,\dots,
Y_s$ are words in $\xxx(z)$, $x=x(a,\tau')$ where
$\theta_1^{-1}=\theta(\tau',z)$ and $x'=x(a,\tau'')$ where
$\theta_s=\theta(\tau'',z)$. The number of applications of
relations (\ref{auxtheta}) in the derivation of this equality is
at most $s$. If $\theta_t\equiv \theta_{t+1}^{-1}
\equiv\theta_l(\tau,z)^{-1},$ and $X_t$ is a power of
$x=x(a,\tau)$, then $Y_t$ is again a power of $x$.
\end{lemma}

\proof (1)   By (\ref{auxax})  $ax = x^4a$ for $x\in \xxx(z)$.
Furthermore

$$a\theta(\tau,z)^{\pm 1}a^{-1}   =   x^{\mp 1}\theta(\tau,z)^{\pm
1}x^{\mp 1}$$ by relations (\ref{auxtheta}). These relations allow
us to move the letter $a$ in $awa\iv$ to the right until it
cancels with $a\iv$. The additional condition for $X'_t$ can be
obtained just as immediately.

     (2) First make two insertions which do not change the value of
$w$ in the free group:

 $$w=X_0 a^{-1} (a \theta_1 X_1 \dots \theta_s a^{-1}) a X_s .$$
Now apply part (1) of the lemma to the word $v$ in parentheses,
utilizing the fact that $a\theta(\tau,z))^{-1} =
x\theta(\tau,z)^{-1}xa$ by (\ref{auxtheta}). As a result, we
obtain a word in the desired form. \endproof

     \begin{lemma} \label{C} Let $\Delta$ be a roll
with base $K_jL_j$ or $K_j^{-1}L_{j-1}$ whose history is a word
over $\sss$. Suppose that the inner contour $q$ of $\Delta$
coincides with the inner side $q$ of the $K_j$-annulus $\bb'$ and
the outer contour of $\Delta$ coincides with the outer contour $p$
of the $L_j$-annulus $\bb$. Then:

     (1) Each maximal $a$-band in $\Delta$ starting on $\bb$ ends on
$\bb$.

     (2) If $U$ is the label of the inner contour $p_0$ of $\bb$
then $U_a$ is a Dyck word and the maximal $a$-bands in $\Delta$
determine a  cancellation pairing of $U_a$.

     (3) If an $a$-band $\CC$ determines a minus pair of letters
$(a\iv,a)$ in $U_a$ then $(a\iv,a)$ are consecutive letters in
$U_a$, and a subband of $\CC$ connects two consecutive
$\theta$-cells $\pi_1$, $\pi_2$ in $\bb$. The corresponding
subword in the history of the roll is $\tau\iv\tau$ where
$\tau\in\sss^+(2)\cup\sss^+(4)$.

     (4) Let $\Gamma$ be a small subtrapezium
in $\Delta$ with history $\tau\iv\tau$ from part (3). Let $\ttt_1$
and $\ttt_2$ be two $\theta$-bands containing the cells $\pi_1,
\pi_2$. Then the maximal subword in $\xxx$ written on $p_0$ (see
part (2)) between $\ttt_1$ and $\ttt_2$ is a power of a letter
$x(a,\tau)$.
\end{lemma}

\proof We consider only the case when the base is $K_jL_j$ (the
other case is similar).

      (1) The statement is obvious since the sides of the $K_j$-band
do not have $a$-edges.

     (2) If an $a$-band starts on an edge $e$ of the path $p_0$,
$\Lab(e)=a^{\pm 1}$ then it ends on an edge of $p_0$ with label
$a^{\mp 1}$. Since the maximal $a$-bands do not intersect, they
determine a pairing in $U_a$.

     (3) Let $\Delta_0$ be the \vk subdiagram of $\Delta$ bounded
by $\CC$ and a part of $p_0$ (and including $\CC$). Suppose that
$\Delta_0$ contains another maximal $a$-band $\CC'$ corresponding
a plus pair. Then we should be able to find a small subtrapezium
$\Gamma$, whose $\theta$-bands $\ttt_1$, $\ttt_2$ are connected by
a subband of $\CC$ and by a subband of $\CC'$. Then consider the
label $V$ of a part of a side of $\ttt_1$ from $\bb'$ to $\bb$.
When we read the word $V$ from left to right (from $\bb'$ to
$\bb$) then we first read a negative $a$-letter (since $\CC$
corresponds to a minus pair) and then a positive $a$-letter since
$\CC'$ corresponds to a plus pair. But this contradicts Lemma
\ref{-+}. Hence the $a$-band $\CC'$ does not exist.

It is clear from the form of relations (\ref{auxtheta}) that the
$a$-band $\CC$ connects two $\theta$-cells $\pi_1$ and $\pi_2$
corresponding to rules $\tau_1\iv$ and $\tau_2$ where the rules
$\tau_1, \tau_2\in \sss^+(2)\cup\sss^+(4)$ are similar. If in the
corresponding subword $h=\tau_1\iv...\tau_2$ of the history of
$\Delta$, there exists a reduced subword of length $2$ then there
exists a small subtrapezium $\Gamma$ with reduced history, such
that $\CC$ crosses $\Gamma$. Therefore the word $V_0$ defined for
$\Gamma$ like the word  $V$ in the previous paragraph, contains a
negative $a$-letter. But this contradicts Lemma \ref{Sm Tr}.

     Therefore $h\equiv \tau^{-1}\tau\tau^{-1}\tau \dots$. As we
have proved before, all pair of edges on $p_0$ defined by
$a$-bands in $\Delta_0$ are minus pairs. Obviously it is possible
only if $h\equiv \tau\iv\tau$.

\begin{center}
\unitlength=1.00mm \special{em:linewidth 0.4pt}
\linethickness{0.4pt}
\begin{picture}(125.00,96.34)
\put(14.67,16.67){\line(0,1){79.67}}
\put(22.00,17.00){\line(0,1){79.00}}
\put(116.33,17.33){\line(0,1){78.67}}
\put(125.00,17.00){\line(0,1){78.67}}
\put(14.67,66.00){\line(1,0){110.33}}
\put(14.67,59.33){\line(1,0){110.33}}
\put(14.67,38.00){\line(1,0){110.33}}
\put(14.67,31.00){\line(1,0){110.33}}
\put(18.00,47.00){\makebox(0,0)[cc]{$\bb'$}}
\put(120.33,47.00){\makebox(0,0)[cc]{$\bb$}}
\put(26.67,62.67){\makebox(0,0)[cc]{$\ttt_1$}}
\put(26.67,34.33){\makebox(0,0)[cc]{$\ttt_2$}}
\put(14.67,19.00){\line(0,-1){15.00}}
\put(22.00,19.00){\line(0,-1){15.00}}
\put(116.33,19.00){\line(0,-1){15.00}}
\put(125.00,19.00){\line(0,-1){15.00}}
\bezier{992}(116.33,92.00)(0.00,53.33)(116.33,5.67)
\bezier{904}(116.33,86.33)(10.00,52.00)(116.33,11.00)
\bezier{620}(116.33,78.67)(44.67,51.00)(116.33,19.33)
\bezier{476}(116.00,72.33)(61.33,50.00)(116.33,26.33)
\put(116.33,86.67){\vector(0,1){5.00}}
\put(116.33,78.67){\vector(0,-1){6.00}}
\put(116.33,19.33){\vector(0,1){6.33}}
\put(116.33,11.00){\vector(0,-1){6.00}}
\put(66.00,59.33){\vector(-1,0){4.67}}
\put(85.33,59.33){\vector(1,0){6.67}}
\put(119.00,89.33){\makebox(0,0)[cc]{$a$}}
\put(119.00,75.33){\makebox(0,0)[cc]{$a'$}}
\put(119.00,22.00){\makebox(0,0)[cc]{$a'$}}
\put(119.00,8.00){\makebox(0,0)[cc]{$a$}}
\put(64.67,61.33){\makebox(0,0)[cc]{$a$}}
\put(90.00,61.33){\makebox(0,0)[cc]{$a'$}}
\put(75.00,61.33){\makebox(0,0)[cc]{$V$}}
\put(60.33,48.67){\makebox(0,0)[cc]{$\CC$}}
\put(84.00,48.67){\makebox(0,0)[cc]{$\CC'$}}
\end{picture}

\end{center}

\begin{center}
\nopagebreak[4] Fig. \theppp.

\end{center}
\addtocounter{ppp}{1}

     (4) By part (3), a subband $\CC_0$ of $\CC$ connects two
$a$-edges which are the closest $a$-edges to the $L_j$-edges of
sides of $\ttt_1$ and $\ttt_2$. Since $\CC$ defines a minus pair
in $U_a$, the word $W$ over $\xxx$ written on the side of $\CC_0$
which is closer to $\bb$ is a product of fourth powers (see
relations (\ref{auxax})). The word $W'$ written on $p_0$ between
$\ttt_1$ and $\ttt_2$ is also a product of fourth powers. Then by
Lemma \ref{Band}, we have the equality $W'=xWx\iv$ in the free
group. By Lemma \ref{x-power} (part (4)), it follows that both $W$
and $W'$ are powers of $x$.
\medskip

\begin{center}
\unitlength=1.00mm \special{em:linewidth 0.4pt}
\linethickness{0.4pt}
\begin{picture}(91.67,55.67)
\put(69.00,55.67){\line(1,0){15.33}}
\put(84.33,55.67){\line(5,-6){7.33}}
\put(91.67,47.00){\line(-5,-6){7.33}}
\put(84.33,38.33){\line(-1,0){15.33}}
\put(69.00,38.33){\line(-3,4){6.33}}
\put(62.67,47.00){\line(3,4){6.67}}
\put(65.67,51.33){\circle*{1.33}}
\put(66.00,51.33){\vector(3,4){2.67}}
\put(65.00,50.33){\vector(-2,-3){1.67}}
\put(69.00,51.33){\makebox(0,0)[cc]{$a$}}
\put(62.67,50.33){\makebox(0,0)[cc]{$x$}}
\put(77.33,47.00){\makebox(0,0)[cc]{$\pi_1$}}
\put(69.00,38.33){\line(0,-1){10.33}}
\put(69.00,28.00){\line(1,0){15.67}}
\put(84.67,28.00){\line(0,1){10.33}}
\put(69.00,35.33){\vector(0,-1){5.67}}
\put(71.00,32.67){\makebox(0,0)[cc]{$W'$}}
\put(68.67,28.00){\line(1,0){15.33}}
\put(84.00,28.00){\line(5,-6){7.33}}
\put(91.33,19.33){\line(-5,-6){7.33}}
\put(84.00,10.67){\line(-1,0){15.33}}
\put(68.67,10.67){\line(-3,4){6.33}}
\put(62.33,19.33){\line(3,4){6.67}}
\put(77.00,19.33){\makebox(0,0)[cc]{$\pi_2$}}
\put(66.00,14.67){\circle*{1.33}}
\put(66.00,14.67){\vector(-3,4){2.33}}
\put(66.00,14.00){\vector(3,-4){1.67}}
\put(62.33,14.67){\makebox(0,0)[cc]{$x$}}
\put(68.67,14.00){\makebox(0,0)[cc]{$a$}}
\put(69.00,38.33){\line(-1,0){52.67}}
\put(69.00,28.00){\line(-1,0){52.33}}
\put(71.00,55.67){\line(-1,0){50.67}}
\put(70.67,10.67){\line(-1,0){46.33}}
\put(50.67,38.33){\circle*{1.33}}
\put(64.67,38.33){\vector(-1,0){9.33}}
\put(60.67,40.33){\makebox(0,0)[cc]{$x$}}
\put(64.33,28.00){\vector(-1,0){9.33}}
\put(60.33,24.67){\makebox(0,0)[cc]{$x$}}
\put(50.67,28.00){\circle*{1.33}}
\put(48.00,38.33){\vector(-1,0){11.33}}
\put(47.67,28.00){\vector(-1,0){11.00}}
\put(42.67,40.67){\makebox(0,0)[cc]{$a$}}
\put(42.67,24.67){\makebox(0,0)[cc]{$a$}}
\put(47.33,33.33){\makebox(0,0)[cc]{$\CC_0$}}
\bezier{196}(68.67,55.67)(33.67,52.33)(31.33,38.33)
\put(31.33,38.33){\line(0,-1){10.33}}
\bezier{188}(31.33,28.00)(40.67,11.00)(68.00,10.67)
\put(51.00,19.33){\makebox(0,0)[cc]{$\CC$}}
\put(79.00,7.00){\makebox(0,0)[cc]{$\bb$}}
\end{picture}

\end{center}
\begin{center}
\nopagebreak[4] Fig. \theppp.

\end{center}
\addtocounter{ppp}{1}

\endproof

\begin{lemma} \label{D}

Let $\Delta$ be a roll with base $L_j$, bounded by the inner
contour $q$ of an $L_j$-annulus, and the outer contour $p$ of an
$a$-annulus. Let $\ttt$ and $\ttt'$ be two consecutive (clockwise)
$\theta$-bands connecting $q$ and $p$. Let $\Gamma$ be the
subdiagram containing $\ttt$ and $\ttt'$ and all cells between
them. Let $V$ be the label of $\bott(\ttt)$. By Lemma \ref{-+},
$V\equiv V_1V_2$ where $(V_1)_a$ is positive and $(V_2)_a$ is
negative. We claim that either $|V_a|\le 2|p|$ or $|(V_1)_a|>|p|$.
\end{lemma}

\proof By contradiction, assume that $|(V_2)_a|
> |p|.$ Let $\CC_1,..., \CC_{|p|+1}$ be a series of consecutive
$a$-bands connecting $\ttt$ and $\ttt'$ enumerated from $q$ to
$p$, so $\CC_{|p|+1}$ is a subband of $\CC$. Then the length of
the $U_1,$ over $\xxx$ written on the bottom side of $\CC_1$ and
the length of the word $U'_1$ written on the top of $\CC_1$ are
related by the equality $|U'_1|=4|U_1|$ (see relations
(\ref{auxax})) because the word $(V_2)_a$ is negative.

\begin{center}
\unitlength=1.00mm \special{em:linewidth 0.4pt}
\linethickness{0.4pt}
\begin{picture}(108.00,72.33)
\bezier{272}(15.33,72.33)(20.00,35.67)(15.33,5.33)
\bezier{276}(98.67,72.33)(106.67,42.33)(98.67,5.33)
\put(13.33,38.00){\makebox(0,0)[cc]{$q$}}
\put(23.33,38.00){\makebox(0,0)[cc]{$\Gamma$}}
\put(108.00,38.33){\makebox(0,0)[cc]{$p$}}
\put(16.67,58.67){\line(1,0){85.00}}
\put(17.00,53.00){\line(1,0){85.00}}
\put(17.00,22.33){\line(1,0){84.67}}
\put(16.67,16.67){\line(1,0){84.00}}
\put(24.33,55.67){\makebox(0,0)[cc]{$\ttt$}}
\put(24.33,19.33){\makebox(0,0)[cc]{$\ttt'$}}
\put(32.00,53.00){\line(0,-1){30.67}}
\put(42.33,53.00){\line(0,-1){30.67}}
\put(52.67,53.00){\line(0,-1){30.67}}
\put(61.33,38.00){\makebox(0,0)[cc]{$\dots$}}
\put(89.67,53.00){\line(0,-1){30.67}}
\put(36.67,45.00){\makebox(0,0)[cc]{$\CC_1$}}
\put(47.67,45.00){\makebox(0,0)[cc]{$\CC_2$}}
\put(95.33,45.00){\makebox(0,0)[cc]{$\CC_{|p|+1}$}}
\put(28.33,27.67){\makebox(0,0)[cc]{$U_1$}}
\put(38.67,27.67){\makebox(0,0)[cc]{$U_1'$}}
\put(45.00,27.67){\makebox(0,0)[cc]{$U_2$}}
\put(55.00,27.67){\makebox(0,0)[cc]{$U_2'$}}
\end{picture}

\end{center}

\begin{center}
\nopagebreak[4] Fig. \theppp.

\end{center}
\addtocounter{ppp}{1}

Let us introduce similar notation $U_2$, $U_2'$ for words over
$\xxx$ on the bottom and top sides  $\CC_2$. Then as in Lemma
\ref{C} (4), we have $U_2 = x^{\pm 1}U'_1y^{\pm 1}$ for some
$x$-letters $x$ and $y.$ Notice that either $U_1$ is not empty or
$x^{\pm 1}\ne y^{\mp 1}$ because otherwise the $(\theta,a)$-cells
connected by $\CC_2$ would form a reducible pair of cells.
Therefore $|U_2|>|U_1|$. Similarly $|U_2|<|U_3|<...<|U_{|p|+1}|.$
But then $|U'_{|p|+1}|\ge 4|U_{|p|+1}|\ge 4(|p|)$, a contradiction
since the definition of $A_{|p|+1}$ implies that $|U'_{|p|+1}|\le
|p|.$ \endproof

\begin{lemma} \label{E}

In the notation of Lemma \ref{D}, let the label $X$ of the subpath
$\Gamma\cap q$ of $q$ be a power of a letter $x\in\xxx$, and $\CC$
be an arbitrary $a_i(L_j)$-band between $\ttt$ and $\ttt'$
starting on the ``positive" subpath of a side of $\ttt$ with label
$V_1$. Then the words in $\xxx$ on $\partial(\CC)$ are powers of
$x'=x(a_i(L_j), \tau)$, and all $a$-bands $\CC_1,\dots,\CC_r=\CC$
between the $L_j$-annulus $\bb$ and $\CC$, connecting $\ttt$ and
$\ttt'$ correspond to the same $a$-letter $a_i(L_j)$.
\end{lemma}

\proof Let $X'$ be a word in $\xxx$ written from $\ttt$ to $\ttt'$
on the top of the subband $\bar\bb$ of the $L_j$-annulus $\bb$
whose bottom side is labelled by $X$. Then relations (\ref{auxkx})
imply that $X'$ is a power of a letter from $\xxx$. let $Y_i$ and
$Y'_i$ are words in $\xxx$ written respectively on the bottom side
(the side that is closer to $\bb$) and the top side of the
$a$-band $\CC_i$. Then $x_1^{\pm 1} Y_1 x_2^{\pm 1} = X'$ for some
$x$-letters $x_1, x_2$ since the word $(V_1)_a$ is positive (see
Lemma \ref{Band}). Besides $Y_1$ is a product of fourth powers
according to the relations (\ref{auxax}). This immediately implies
that both $Y_1$ and  $Y_2$ are powers of a letter $x'$ (the
assumption that $Y_1$ is empty and $x_1^{\pm 1}$ and $x_2^{\pm 1}$
are mutually inverse implies that two $\theta$-cells in $\ttt$,
$\ttt'$ that are connected by $\CC_1$ cancel). The letters $x_1$
and $x_2$ also must coincide with $(x')^{\pm 1}$. Therefore $Y_1$
is an $a$-band where $a$ is determined by the letter $x'$ (see
relations (\ref{auxtheta})). Repeating this argument leads to the
same conclusion for $x$-subwords written on the sides of $a$-bands
$\CC_2,...,\CC_r=\CC$, so all these $a$-bands correspond to the
same letter from $\aaa$ determined by the $x$-letter $x'$.
\endproof

\begin{lemma} \label{F} Let a reduced annular diagram $\Delta$ consists of one $a$-annulus
with inner side $q$ and outer side $p$, and the edges with
positive labels $a=a_i(L_j)$ cut the diagram from $q$ to $p$ (that
is the edges with positive $a$-labels are oriented outward). Let
the word $w\equiv \Lab(q)$ have the form (\ref{*}) from Lemma
\ref{B} for $z=L_j$ and $s>0$. Suppose that the cyclic word $w$
admits two $(2,4)$-pairings, a $L_j$-best pairing $I$ and a
$\lel_j$-good pairing $II$, such that $I$ is a minus pairing, and
in $II$, all minus pairs correspond to occurrences of subwords
$\theta\iv X\theta$ where $X$ is a power of a letter
$x=x(a,\tau)$, $\theta=\theta(\tau,L_j)$. Suppose further that the
minus pairs of $II$ are normal pairs for $I$. Then there exists an
annular diagram $\Delta'$ with the labels of contours $w$, $w'$
where $w'$ also has the form (\ref{*}) but for $z=\lel_j$ and
possibly with different $x$-subwords than $w$; the diagram
$\Delta'$ has one $L_j$-annulus $\bb_0$ surrounding the hole of
$\Delta'$, other cells in $\Delta'$ correspond to relations
(\ref{auxtheta}) and (\ref{auxax}), and the number of
$\theta$-cells in $\Delta'$ is less than $s^2$.
\end{lemma}

\proof By Lemma \ref{B} (2), using the pairing $I$ we can replace
each pair of connected letters  $\dots \theta^{-1} \dots
\theta'\dots$ in the cyclic word $w$ by a pair of subwords
$\dots(a^{-1}x\theta^{-1})\dots (\theta'x'^{-1}a) \dots$ (starting
with the outermost pairs in $I$). The number of
$(\theta,a)$-relations used in the transition to the new word
$w_1$ is less than $s^2/2$. The form of the word $w_1$ given in
Lemma \ref{B} (2) and the form of relations (\ref{mainrel})
corresponding to rules from $\sss(2)\cup\sss(4)$, and relations
(\ref{auxtheta}) allows us to glue an $L_j$-annulus $\bb_0$ to the
annular diagram of conjugacy of $w$ and $w_1$ along the path with
label $w_1$ (the top side of $\bb_0$).

In the transition from $w$ to $w_1$, we do not change the
$\theta$-letters. Suppose that an occurrence of a subword of the
form $\theta^{-1}X\theta$ corresponds to a minus pair in $II$. By
Lemma \ref{A} this pair is connected in $I$ as well. This pair is
normal in $I$ by the assumption of the lemma. Therefore in the
subword $(a^{-1}x\theta^{-1})Y(\theta x^{-1}a)$ of the word $w_1$
obtained by Lemma \ref{B} (2) from the subword
$\theta^{-1}X\theta$ of $w$, the word $Y$ is a power of the same
letter $x$.

Let $w_2$ be the word written on the bottom side of our
$L_j$-annulus $\bb_0$. The subword of $w_1$ discussed in the
previous paragraph corresponds to the subword of the form
$(a_-^{-1}x_-\theta_-^{-1})Y_-^4(\theta_- x_-^{-1}a_-)$ in $w_2$
(see relations (\ref{auxkx}) and (\ref{mainrel})) where the minus
sign in the subscript means that we replace $L_j$ by $\lel_j$ in
$a_i(L_j)$, $x(a_i(L_j),\tau)$ and $\theta(\tau,L_j)$.  Since we
have a fourth power in the middle of this word, we can transform
this occurrence into $\theta_-^{-1}Y_-\theta_-$ without
$a$-letters applying relations (\ref{auxax}) and two relations
(\ref{auxtheta}).

     In the resulting word $w_3$,  $a$-letters occur only
in subwords of the form $(\theta x^{-1}a)^{\pm 1}$, where the
letter $\theta$ occurs in a plus pair of the pairing $II$. These
$a$-letters can be removed using relations (\ref{auxtheta}) and
(\ref{auxax}) by Lemma \ref{B} (1). As a result, we get a word
$w_4=w'$ and the diagram $\Delta'$ consists of a $L_j$-annulus
$\bb_0$ and other cells described above. This diagram satisfies
all the conditions of the lemma.
\endproof

\begin{lemma} \label{G} Let $\Delta$ be a roll with base $K_jL_j$ or $K_j^{-1}L_{j-1}$
($j$ is odd), whose history is not empty and does not contain
rules from $\bsss$. Suppose the outer contour $p$ of the roll is
the outer contour of some $a_i(L_j)$-annulus. Then either between
the $L_j$-annulus $\bb$ and $p$, there exists an $a$-annulus $\CC$
satisfying the conditions of Lemma \ref{F} and the length of every
$\theta$-band connecting $\bb$ and $\CC$ is at least $|p|$ or
$\bb$ and $p$ are connected by a $\theta$-band of length $\le
2|p|$.
\end{lemma}

\proof We can consider only the case when the base is $K_jL_j$.
Let $\Delta_1$ be the annular subdiagram in $\Delta$ bounded by
$\bb$ and $p$.  By Lemma \ref{D} we can assume that for every two
consequent $\theta$-bands $\ttt$ and $\ttt'$ connecting $\bb$ and
$p$, $|V_1|_a> |p|$ in the notation of Lemma \ref{D}. If we
enumerate all $a$-annuli between $\bb$ and $p$ from $\bb$ to $p$:
$\CC_1, \CC_2, \dots$, then we can let $\CC$ be the first
$a$-annulus which cannot be connected with $\bb$ by any
$\theta$-band of length $<|p|$.

Since $\Delta_1$ has an $a$-annulus by the lemma assumption, we
have that every maximal $a$-band $\CC'$ in $\Delta_1$ that starts
on $\bb$, ends on $\bb$, and we obtain a pairing in the word
$U'_a$ where $U'$ is written on the top side of $\bb$.

Since the $a$-letters in a relation (\ref{mainrel}) corresponding
to rules from $\sss(2)\cup\sss(4)$ determine the rule up to the
similarity relation, we have a $(2,4)$-pairing $I$ in the word
$U'$.

Any $a$-band $\CC'$ of $\Delta_1$, starting and ending on $\bb$
has a subband connecting a pair of consecutive $\theta$-bands
$\ttt$ and $\ttt'$. Moreover $\CC'$ must intersect the path
labelled by $V_1$ (in the above notation) because between $\CC'$
and $\ttt$ there are less than $|p|$ maximal $a$-bands and
$|V_1|_a\ge |p|$. Since $(V_1)_a$ is a positive word, $\CC'$
determines a minus pair of $a$-letters. Therefore $I$ is a
minus-pairing. Notice that this pairing is $L_j$-best. Indeed,
every maximal $\theta$-band in $\Delta$ intersects the annulus
$\CC$. Hence for each $\theta$-letter on the outer contour of
$\bb$, there exists a $(\theta,a)$-relation of the form
(\ref{auxtheta}) containing $\theta$. Thus the corresponding rule
from $\sss$ does not lock $L_jP_j$-sectors.

\begin{center}
\unitlength=1mm \special{em:linewidth 0.4pt} \linethickness{0.4pt}
\begin{picture}(103.00,74.33)
\put(12.33,6.33){\line(0,1){63.67}}
\put(17.00,6.33){\line(0,1){63.67}}
\put(50.33,6.33){\line(0,1){63.67}}
\put(55.00,6.33){\line(0,1){63.67}}
\put(90.67,6.33){\line(0,1){63.67}}
\put(95.33,6.33){\line(0,1){63.67}}
\put(50.33,39.00){\rule{52.67\unitlength}{2.33\unitlength}}
\put(50.33,32.33){\rule{52.67\unitlength}{2.33\unitlength}}
\bezier{164}(50.33,68.33)(31.67,62.00)(50.33,51.00)
\bezier{120}(50.33,65.67)(36.67,61.33)(50.33,54.00)
\bezier{296}(55.00,53.33)(88.00,41.67)(55.00,21.00)
\bezier{252}(55.00,51.33)(83.33,39.67)(55.00,24.33)
\bezier{136}(50.33,23.67)(35.67,14.33)(50.33,6.33)
\bezier{96}(50.33,20.33)(40.00,14.00)(50.33,8.67)
\bezier{176}(55.00,46.33)(75.00,38.67)(55.00,28.33)
\bezier{140}(55.00,44.00)(71.33,38.00)(55.00,30.33)
\bezier{212}(50.33,46.00)(25.67,37.00)(50.33,27.67)
\bezier{160}(50.33,44.00)(31.33,36.67)(50.33,30.33)
\put(83.67,44.67){\makebox(0,0)[cc]{$\ttt$}}
\put(83.67,27.67){\makebox(0,0)[cc]{$\ttt'$}}
\put(92.67,74.00){\makebox(0,0)[cc]{$\CC$}}
\put(53.00,74.33){\makebox(0,0)[cc]{$\bb$}}
\put(14.67,74.33){\makebox(0,0)[cc]{$\bb'$}}
\put(25.00,38.00){\makebox(0,0)[cc]{$\Delta_2$}}
\put(63.33,23.67){\makebox(0,0)[cc]{$\CC'$}}
\end{picture}

\end{center}

\begin{center}
\nopagebreak[4] Fig. \theppp.

\end{center}
\addtocounter{ppp}{1}

Now consider the subroll $\Delta_2$ bounded by $\bb$ and the
$K_j$-annulus $\bb'$. Lemma \ref{C} (2) gives a pairing in the
word $U_a$ where $U$ is the label of the bottom of $\bb$. This
automatically gives a $(2,4)$-pairing in the word  $U$ which we
shall denote by $II$. This pairing is $K_j$-good because an
$a$-band that connects two edges on the inner contour of $\bb$
cannot cross a $\theta$-band corresponding to a rule locking
$K_jL_j$-sectors (see relations (\ref{auxtheta})). By Lemma
\ref{C} (3,4) every minus pair of $\theta$-letters of $II$
corresponds to an occurrence of a subword $\theta^{-1}Y\theta$ in
$U$ where $Y$ is a power of an $x$-letter $x=x(a, \tau)$ ($a$ is
obviously determined by $\theta$ by the relations
(\ref{mainrel})). By Lemma \ref{A} the same pair of letters must
be connected in the pairing $I$, and by Lemma \ref{E} for
$\Delta_1$ the corresponding occurrence in $\Lab(p)$ also has the
form $\theta^{-1}X\theta$ where $X$ is a power of the letter
$x(a',\tau)$ where $a'$ is the ``brother" of $a$. The normality
condition of the bottom of $\CC$ follows from Lemma \ref{E} and
the choice of $\CC$. \endproof

\begin{lemma} \label{H} Let $\Delta$ be a roll with base $L_{j-1}^{-1}K_jL_j$ for some $j$
and history $h$, $|h|>0$. Let the outer contour $p$ of $\Delta$ be
the outer contour of some $a$-annulus, and $h$ does not contain
rules from $\bsss$. Then there exists a roll $\bar\Delta$ with the
same boundary labels and the same base as $\Delta$ where the
$L_j$-annulus can be connected with $p$ by a $\theta$-band of
length $\le 2|p|$.
\end{lemma}

\proof  By Lemma \ref{G} applied to the subroll with base $K_jL_j$
we can assume that there exists an $a$-annulus $\CC$ between the
$L_j$-annulus $\bb$ and $p$ satisfying the conditions of Lemmas
\ref{G} and \ref{F}. Hence the number of $(\theta,a)$-cells
between $\bb$ and $\CC$ is at least $|p|s$ which, in turn, at
least $s^2$ where $s$ is the length of the history of the roll.

     Consider the auxiliary diagram $\Delta'$
corresponding to $\CC$ as in the conclusion of Lemma \ref{F}. Let
$\Delta''$ be the diagram obtained from $\Delta'$ by first taking
the mirror image and then changing the indexes $L_j$ to $L_{j-1}$,
and $K_j$ to $K_{j}\iv$ in all labels (Lemma \ref{j-j'} allows us
to do that). Lemma \ref{F} and the form of relations
(\ref{mainrel}) corresponding to $K_j$ imply that there exists a
$K_j$-band $\bb'$ whose outer contour coincides with the inner
contour of $\Delta'$ and whose inner contour coincides with the
outer contour of $\Delta''$. The diagram $\Delta_0$ obtained by
gluing these three diagrams (($\Delta'$, $\bb'$ and $\Delta''$)
has outer contour $q_0$ labelled by $w$ (same as $\CC$) and inner
contour labelled by the word $w''$ obtained from $w$ by replacing
the index $L_j$ with $L_{j-1}$ in all letters.

Therefore we can do the following transformation with $\Delta$.
First cut $\Delta$ along $q_0$ (the inner contour of $\CC$), along
the outer contour $p_1$ of the $K_j$-annulus $\bb''$ of $\Delta$,
and along the inner contour $q_1$ of this annulus.  We get three
annular diagrams $\Delta_1, \Delta_2$ and $\Delta_3$ where
$\Delta_1$ is bounded by $p$ and $q_0$, $\Delta_2$ is bounded by
$q_0$ and $p_0$ and $\Delta_3$ is bounded by $q_1$ and $q$. Then
change index $K_j$ to $L_{j-1}\iv$, $L_j$ by $L_{j-1}$ in all
labels of edges of $\Delta_2$. Let $\Delta_2'$ be the resulting
diagram. The label of the outer contour of $\Delta_2'$ is
identical with the label of the inner contour  of $\Delta_0$ and
the label of the inner contour of $\Delta_2'$ coincides with the
label of $q_1$ as follows from the form of $K_j$-relations
(\ref{mainrel}) and (\ref{auxkx}). Let $\bar{\bar\Delta}$ be the
result of gluing $\Delta_1$, $\Delta_0$, $\Delta_2'$, and
$\Delta_3$ and reducing the resulting diagram. By Lemma \ref{F},
the number of $(\theta,a)$-cells between the $L_j$-annulus of
$\bar{\bar\Delta}$ and $p$ is smaller than the similar number for
$\Delta$. So our transformation moved the $L_j$-annulus ``closer"
to $p$.

The base of $\bar{\bar\Delta}$ before reducing was $L_{j-1}\iv
L_{j-1}L_{j-1}\iv K_jL_j$ (one $L_{j-1}$-band comes from
$\Delta_1$, and two from $\Delta_0$). After reducing two of the
$L_{j-1}$-bands can disappear.

If this does not happen then the history of $\bar{\bar\Delta}$ is
a word over $\sss(2)\cup\sss(4)$ (since other rules lock
$L_{j-1}\iv K_j$-sectors or $P_{j-1}\iv L_{j-1}\iv$-sectors).
Hence we can apply Lemma \ref{bezMivM} and remove two of the
$L_{j-1}$-annuli from $\bar{\bar \Delta}$.

Thus after the transformations described above we get a roll with
the same base, same labels of the contours as $\Delta$ but with
the $L_j$-annulus closer to $p$ than in $\Delta$. Moving the
$L_j$-annulus further toward $p$, we obtain the desired roll
$\bar\Delta$. \endproof

\subsection{Arbitrary rolls}
\label{ar}

\begin{lemma} \label{rolls}
Let $\Delta$ be a roll with contours $p$ and $q$. Then there
exists a roll $\Delta'$ with the same boundary labels as $\Delta$,
which contains a path $t$ connecting the two contours of
$\Delta'$, whose length is recursively bounded in terms of
$|p|+|q|$.
\end{lemma}

\proof By Lemma \ref{rollsprep}, we can assume that $\Delta$
satisfies (R1) and (R2).

If the history of $\Delta$ is empty then $\Delta$ contains only
$\gr$-cells and auxiliary cells corresponding to relations
(\ref{auxax}). Moreover, a $\gr$-cell cannot have common edges
with the $(\aaa,\xxx)$-cells. Thus $\Delta$ must be a union of
$\gr$-cells and diagrams over $\hhh_2$, and different components
of this union do not have common edges, i.e, $\Delta$ is the
diagram over the free product of groups. Hence the statement of
the lemma in that case follows from the solvability of the
conjugacy problem in $\gr$ and the conjugacy problem in $\hhh_2$
(Lemma \ref{Conjugacy}).

Hence we can assume that the history of $\Delta$ is not empty. Let
$\ttt_1,\ttt_2,...$ be the maximal $\theta$-bands in $\Delta$
connecting $q$ and $p$.

Denote by $\bb_1,\dots,\bb_r$ the basic annuli of the roll
$\Delta$ counted from $q$ to $p$. We are going to recursively
bound the number $r$ as function of $|p|, |q|$ in $\Delta$ or in a
roll with the same boundary labels and to use induction on $r$.
The case $r=0$ has been treated in Lemma \ref{bez bazy}.

Recall that there are no auxiliary $x$-cells in $P$- or $R$-bands.
Hence the length of every $P$- or $R$-annulus in $\Delta$ is at
most $\min(|p|,|q|)$. One may assume that there are no such
distinct bands with equal boundary labels, because this would
allow us to delete an annular subdiagram of $\Delta$ and to
reduce the number of $k$-annuli in $\Delta$. Since the number and
the lengths of $P$- and $R$-annuli are recursively bounded, our
task is reduced to rolls having nether $P$- nor $R$-annuli.

Assume that there is a $K_j$-annulus $\bb_s$ and a
$K_{j'}$-annulus $\bb_{s'}$ where $s'> s$. We choose $s'$ so that
$s'-s$ is minimal possible. Since there are neither $P$- nor $R$-
annuli among $\bb_s, \bb_{s+1},\dots, \bb_{s'}$, we have $j'=j$
by Lemma \ref{Base}, and the base of the subroll $\Delta'$ of
$\Delta$ bounded by $\bb_s$ and $\bb_{s'}$ has one of the forms:
\begin{equation}\label{jjj}
K_jL_jL_j\iv L_j...L_j\iv K_j\iv,\end{equation} or
\begin{equation}\label{j1j}
K_j\iv L_{j-1}L_{j-1}\iv...K_j, \end{equation} if $j$ is odd, or
\begin{equation}
K_jK_j\iv \hbox{ or } K_j\iv K_j
\end{equation}
if $j$ is even.

We shall assume that $j$ is odd, the other case is similar. By
Lemma \ref{Base} none of the rules in the history of the roll
locks $L_jP_j$-sectors if the base of $\Delta'$ has the form
(\ref{jjj}), and none of them locks $L_{j-1}P_{j-1}$-sectors if
the base has the form (\ref{j1j}). In particular, if the history
of the roll contains rules from $\bsss$, then $j\ne 1$ in the
(\ref{jjj}) case.  In addition, if there is a $\bar\Theta$-band in
$\Delta'$, then there are no $\Theta$-bands, because a common
admissible word for $\sss$ and $\bar\sss$ have no
$a(L_j)$-letters, and so, for any $j$, it has no subwords with
base $L_jL_j^{-1}$.

The subroll bounded by $\bb_s$ and $\bb_{s'}$ has no $\gr$-cells
by Lemma \ref{Band} and by Lemma \ref{No a-cells} applied to the
small subtrapezia of $\Delta'$. By Lemma \ref{bez1}, all $a$-edges
in the subroll between $\bb_t$ and $\bb_{t+1}$, $s\le t\le s'-1$
have labels from $\aba{z}$ where $z=z(t)$ depends on $t$ only ($z$
is either $K_j$ or $L_j$ if the base of $\Delta'$ has the form
(\ref{jjj}) and $z$ is either $K_j\iv$ or $L_{j-1}$ if the base
has the form (\ref{j1j}).

We are going to show how to eliminate the two $K_j$-annuli from
$\Delta$ without changing the boundary labels of $\Delta$. Suppose
that the base has the form (\ref{jjj}). Recall that by Lemma
\ref{Base}, there are no $\bar\Theta$-cells in $\Delta'$ if $j=1$.
Therefore every auxiliary relation involving letters from
$\aba{K_j}\cup\Theta(K_j)\cup\bar\Theta(K_j)\cup\xxx(K_j)$ and
involved in $\Delta'$, has a copy involving letters from
$\aba{L_{j-1}\iv}\cup\Theta(L_{j-1}\iv)\cup\bar\Theta(L_{j-1}\iv)\cup\xxx(L_{j-1}\iv)$,
every $L_j$-cell has a copy $L_{j-1}$-cell, and every auxiliary
relation involving letters from
$\aba{L_j}\cup\Theta(L_j)\cup\bar\Theta(L_j)\cup\xxx(L_j)$ (and
involved in $\Delta'$) has a copy involving letters from
$\aba{L_{j-1}}\cup\Theta(L_{j-1})\cup\bar\Theta(L_{j-1})\cup\xxx(L_{j-1})$.
Notice that the labels of the inner and outer side of a $K_j$-band
become identical after a substitution of each letter from
$\Theta(K_j)\cup\bar\Theta(K_j)\cup\xxx(K_j)$ by its ``brother" in
$\Theta(L_{j-1}\iv)\cup\bar\Theta(L_{j-1}\iv)\cup\xxx(L_{j-1}\iv)$
(see Lemma \ref{j-j'}).

Let us replace every label in the subroll of $\Delta$, bounded by
the outer side of $\bb_s$ and the inner side of $\bb_{s'}$, by its
copy in the sense of the previous paragraph. Then we obtain a roll
whose boundary labels coincide with the boundary labels of the
subroll $\Delta'$. Hence we can replace $\Delta'$ in $\Delta$ by a
subroll with fewer $K_j$-annuli, as desired. Let us call this
operation \label{oprem}{\em the operation of removing
$K_j$-annuli}.

In the case when the base has the form (\ref{j1j}) the procedure
of removing $K_j$-bands is similar. The only exception is the case
when $j=1$. In that case $(L_{j-1},\bar\Theta)$-cells do not have
$L_1$-copies. But then all $\theta$-cells of $\Delta'$ are
$\bar\Theta$-cells, as we noticed earlier. So we need to apply the
homomorphism that kills $\bar(a)K_j^{-1}$- and $\bar
a(L_{j-1})$-letters. This homomorphism sends relations to
relations, so the construction carries out without other changes.

Hence $\Delta$ contains at most one $K_j$-annulus. Suppose
$\Delta$ does not have $K_j$-annuli and the base of $\Delta$ has
at least 3 letters. Then the history of $\Delta$ is a words in
 $\sss(2)\cup\sss(4)\cup\bar\sss$ since other rules lock either
 $L_jP_j$- or $L_j^{-1}\lel_j^{-1}$-sectors. In this case the base
 can be shorten by Lemma \ref{bezMivM}. If there are $K_j$-annuli
 in $\Delta$ and its base has at most two letters, then we
can complete the proof by applying Lemma \ref{M-rolls}. So let us
assume that $\Delta$ contains exactly one $K_j$-annulus $\bb_s$.
Thus (up to an inside-out transformation of the roll) the base of
$\Delta$ has the form $...L_{j-1}L_{j-1}\iv K_jL_jL_j\iv...$ if
$j$ is odd or $K_j$ if $j$ is even (the $L_j$'s may be absent even
if $j$ is odd).

  If the base of $\Delta$ contains a
subwords $L_sL_s\iv$ and $L_s\iv L_s$ then none of the rules in
the history of $\Delta$ locks $L_s\iv z$- or $L_s$-sectors. Hence
the history of $\Delta$ is a word in
$\sss(2)\cup\sss(4)\cup\bar\sss$. By Lemma \ref{bezMivM}, every
subroll with base $L_s\iv L_s$ can be replaced by a subroll with
an empty base. Hence we can assume that the base of $\Delta$ does
not contain subwords $L_s\iv L_s$. Therefore the base of $\Delta$
has length at most $5$.

Thus $k$-annuli of $\Delta$ divide $\Delta$ into at most six
subrolls with empty bases. Let us number these subrolls
$\Delta_1,\Delta_2,..., \Delta_{r+1}$, $r\le 5$, counting from $q$
to $p$. As in Lemma \ref{M-rolls}, let us fix a ``big enough"
recursive function $f(x)$, it will be clear later, how big this
function should be.

If $\Delta_1$ or $\Delta_{r+1}$ has a path of length $\le
f(|p|+|q|)$ connecting its boundary components then using Lemma
\ref{tdist}, we can show that the number of cells in $\Delta_1$
or $\Delta_{r+1}$ is recursively bounded in terms of $|p|+|q|$,
and we can remove these subrolls and $\bb_1$ or $\bb_{t-1}$ from
$\Delta$ reducing the number of $k$-annuli in $\Delta$. Thus we
can assume that $\Delta_1$ and $\Delta_{r+1}$ have no short cuts.

Suppose that the base of $\Delta$ starts with $K_j$. By Lemma
\ref{a-annulus}, the roll $\Delta_1$ contains an $a$-annulus
$\CC$. It is easy to see that one can build a $K_j$-annulus $\bb'$
whose outer side label coincides with a inner side label of $\CC$.
Now cut the roll $\Delta$ along that side of $\CC$, and insert in
the hole the $K_j$-band $\bb'$ and its inverse (so that the two
$K_j$-bands cancel each other). As a result we get a non-reduced
annular diagram with three $K_j$-annuli. Using the operation of
removing $K_j$-annuli, we can now remove the bands $\bb_1$,
$\bb'$. After reducing the new annular diagram, we get a roll with
the same boundary labels and the same base as $\Delta$ but with
the $K_j$-band closer to $q$ than in $\Delta$. Let us call this
operation \label{opmov}{\em the operation of moving a $K_j$-band}.
We can use this operation to move the $K_j$-band towards $q$ until
there are no more $a$-annuli between it and $q$. By Lemma
\ref{a-annulus}, if that was the case then the number of cells
between this $K_j$-band and $q$ is recursively bounded, so we
would be able to remove the subroll bounded by the $K_j$-band and
$q$ from $\Delta$, and then apply Lemma \ref{M-rolls}. Similarly
we can treat the case when the base of the roll ends with $K_j$.

Suppose that the base of $\Delta$ starts with $L_{j-1}L_{j-1}\iv
K_j$. Then again by Lemma \ref{a-annulus} $\Delta_1$ contains an
$a$-annulus $\CC$ whose inner side label is the outer side label
of some $K_j$-band. Hence we can do the operation of removing a
$K_j$-band. As a result of this operation we obtain a roll whose
base starts with $K_j$ (and contains only one occurrence of
$K_j^{\pm 1}$), so we reduce the problem to the previous case.
Similarly we can argue when the base of $\Delta$ ends with
$L_jL_j\iv$.

Hence we can assume that the base of $\Delta$ has the form
$L_{j-1}\iv K_jL_j$. As before, we can assume that the number of
$a$-annuli in $\Delta$ between $p$ and the $L_j$-annulus, and
between $q$ and the $L_{j-1}$-annulus is big enough (in
particular, $>1$). This implies that the history of the roll
cannot contain rules that lock $L_jP_j$-sectors. The history of
$\Delta$ cannot contain both rules from $\sss$ and rules from
$\bar\sss$ because no $a$-letters are both of the form $a_i(L_j)$
and $\bar a(L_j)$. If the history is a word over $\sss$, we can
apply Lemma \ref{H} and move the $L_j$-annulus within a recursive
distance from $p$, and then eliminate this annulus from the
diagram reducing the base.

Hence we can assume that the history is a word in $\bsss$. Since
there exist no $(\bar a,x)$-relations, every $a$-annulus of
$\Delta$ does not contain $x$-edges, and the label of the outer
side of it consists of $\theta$-edges. Hence the lengths of all
$a$-annuli in $\Delta$ are bounded by $|p|+|q|$. Hence if the
number of $a$-annuli in $\Delta$ is big enough, there exist two
$a$-annuli with the same labels of outer contours. Then we can
reduce the number of $a$-annuli by removing the one of the two
$a$-annuli and all cells between them.
\endproof

\section{Arrangement of hubs}
\setcounter{equation}{0} \label{aod}

\begin{lemma} \label{LMqR} Let a reduced diagram $\Delta$ has exactly one hub $\Pi$. Assume
that a $K_j$-band $\bb_j$ and a $K_{j+1}$-band $\bb_{j+1}$ start
on $\Delta$ for $j\ne 1$, and a $\theta$-band $\ttt$ cross these
$K$-bands and does not cross any other $K$-bands. Assume further
that $\partial\Delta=p(j)q(j)t(j)q'(j)$ where $p(j)$ is a subpath
of $\partial\Pi$, $q(j)$ and $q'(j)$ are top or bottom paths of
$\bb_j$ and $\bb_{j+1}$, respectively, and $t(j)$ is the top of
$\ttt$.

Then one can attach a diagram over $\hhh_1$ to $\Delta$ along
$q'(j)t q(j)$ and obtain a diagram $\Delta'$ whose perimeter is
bounded from above by a linear function of the length of $\ttt$.
In addition, $\Lab(\partial\Delta')\equiv\Sigma$ if $\Lab(t(j))_a$
is empty.
\end{lemma}

\begin{center}

\unitlength=1mm \special{em:linewidth 0.4pt} \linethickness{0.4pt}
\begin{picture}(129.82,90.67)
\put(72.00,17.17){\oval(21.33,15.67)[]}
\put(61.33,14.67){\line(-2,1){53.33}}
\put(82.67,14.67){\line(5,3){45.67}}
\bezier{876}(6.67,41.33)(68.00,132.67)(128.33,41.67)
\put(61.33,19.33){\line(-2,1){51.67}}
\put(82.67,19.00){\line(5,3){43.33}}
\bezier{808}(12.00,39.00)(69.67,123.33)(123.67,39.33)
\bezier{508}(106.00,33.00)(70.67,84.67)(31.67,34.00)
\bezier{328}(46.00,27.00)(69.33,59.33)(96.00,27.00)
\put(56.29,22.08){\vector(-2,1){6.84}}
\put(42.90,28.63){\vector(-2,1){7.98}}
\put(50.02,20.37){\vector(-2,1){10.26}}
\put(110.72,31.77){\vector(-3,-2){10.83}}
\put(103.88,31.77){\vector(-3,-2){5.98}}
\put(93.34,25.50){\vector(-3,-2){6.55}}
\put(65.13,87.06){\vector(1,0){5.98}}
\put(75.38,9.26){\vector(-1,0){5.41}}
\put(20.38,54.85){\makebox(0,0)[cc]{$\ttt$}}
\put(6.42,47.73){\makebox(0,0)[cc]{$\bb_j$}}
\put(129.82,46.02){\makebox(0,0)[cc]{$\bb_{j+1}$}}
\put(40.05,48.87){\makebox(0,0)[cc]{$p_{22}(j)$}}
\put(70.82,50.29){\makebox(0,0)[cc]{$\Gamma_2(j)$}}
\put(41.19,33.48){\makebox(0,0)[cc]{$p_{21}(j)$}}
\put(99.32,32.62){\makebox(0,0)[cc]{$p_{23}(j)$}}
\put(70.54,31.20){\makebox(0,0)[cc]{$\Gamma_1(j)$}}
\put(70.54,44.59){\makebox(0,0)[cc]{$p_{21}(j)=p_{42}(j)\iv$}}
\put(56.29,26.07){\makebox(0,0)[cc]{$p_{11}(j)$}}
\put(87.64,25.50){\makebox(0,0)[cc]{$p_{31}(j)$}}
\put(106.45,73.38){\makebox(0,0)[cc]{$t(j)$}}
\put(67.98,69.67){\makebox(0,0)[cc]{$\dots$}}
\put(71.96,16.38){\makebox(0,0)[cc]{$\Pi$}}
\put(72.25,5.55){\makebox(0,0)[cc]{$p(j)$}}
\put(32.35,24.65){\makebox(0,0)[cc]{$q(j)$}}
\put(105.88,23.51){\makebox(0,0)[cc]{$q'(j)$}}
\put(51.45,34.05){\vector(4,3){3.70}}
\put(40.33,44.02){\vector(1,1){3.99}}
\end{picture}

\end{center}

\begin{center}
\nopagebreak[4] Fig. \theppp.

\end{center}
\addtocounter{ppp}{1}

\proof Let $\bar\Delta$ be a diagram obtained from $\Delta$ after
deleting of the hub $\Pi$. Obviously, $\bar\Delta$ is a trapezium
with base $\lel_jL_jP_jR_j\rer_j$. If we delete $\bb_j$ and
$\bb_{j+1}$ from $\bar\Delta$, we obtain a diagram $\Gamma(j)$
with contour
$p_{\Gamma(j)}q_{\Gamma(j)}t_{\Gamma(j)}q'_{\Gamma(j)}$ where
$q_{\Gamma(j)}$ ($q'_{\Gamma(j)}$) is a the inner top/bottom of
$\bb_j$ (of $\bb_{j+1}$), $p_{\Gamma(j)}$ is a subpath of
$\partial\Pi$ and $t_{\Gamma(j)}$ is a subpath of $t(j)$. Our goal
is to construct similar diagrams $\Gamma(j')$ for other indices
$j'$ so that the labels of $p_{\Gamma(j')}$, $q_{\Gamma(j')}$ and
$q'_{\Gamma(j')}$ have index $j'$ and they are the copies of the
labels of $p_{\Gamma(j)}$, $q_{\Gamma(j)}$ and $q'_{\Gamma(j)}$,
respectively, and $|t_{\Gamma(j')}|\le |t_{\Gamma(j)}|$. Then we
can attach every such $\Gamma(j')$ (or its mirror copy if $j'-j$
is odd) to $\Delta$ and insert copies of $\bb_{j}$ (if $j'-j$ is
even) or of $\bb_{j+1}$ (if $j'-j$ is odd) between the copies of
$\Gamma(j)$. Then we obtain the desired diagram $\Delta'$.

To construct $\Gamma(j')$ we subdivide the trapezium $\bar\Delta$
into subtrapezia $\Delta_1, \Delta_2,\dots$ of the first and the
second types by making cuts along the boundaries of the maximal
$\bar\Theta$-bands of $\bar\Delta$ so that the subtrapezia of the
first and the second types alternate in $\bar\Delta$. The same
cuts partition $\Gamma(j)$ into subdiagrams $\Gamma_1(j),
\Gamma_2(j),\dots,$ when counting from $\Pi$ to $t$. The boundary
of $\Gamma_s(j)$ is of the form
$p_{1s}(j)p_{2s}(j)p_{3s}(j)p_{4s}(j)$ where
$p_{2s}(j)=p_{4,s+1}(j)^{-1}$ and the paths $p_{1s}(j)$ (the paths
$p_{3s}(j)$) are subpaths of $q_{\Gamma(j)}$ (of
$q'_{\Gamma(j)}$).

Notice that $a(\lel_j)$-, $a(L_j)$- or $a(R_j)$-bands of a
subtrapezia $\Delta_s$ of the second type cannot start on the
boundary of $\Delta_s$ because the common letters of $\aaa$ and
$\bar\aaa$ belong to $\cup_j\aaa(P_{j})$. Hence such bands contain
only $\bar a$-edges, and consequently, $\Delta_s$ has no
$(a,x)$-cells.

Since the boundary label of every $\Gamma_s(j)$ obviously has
index $j$, we can apply Lemma \ref{j-j'} or Lemma \ref{j-j'(bar)},
(depending on the type of trapezia $\Delta_s)$ and obtain diagrams
$\Gamma_s(j')$ whose boundaries
$p_{1s}(j')p_{2s}(j')p_{3s}(j')p_{4s}(j')$ are labelled by words
having index $j'$; the labels of $p_{1s}(j')$ and $p_{3s}(j')$ are
copies of the labels of $p_{1s}(j)$ and $p_{3s}(j)$, respectively.
The lengths of $p_{2s}(j')$ and $p_{4s}(j')$ are not greater of
those for index $j$, but to construct the desired diagram
$\Gamma(j')$ from $\Gamma_s(j')$, $s=1,2,... $ (and from some
auxiliary cells) we have to compare the labels of the $p_{2s}(j')$
and $p_{4,s+1}(j')$ for $j'=1$. (They are equal for $j'\ne 1$.) We
will assume that $s=1$, and $\Delta_1$ is of the first type,
because other comparisons are similar.

The word $\Lab(p_{21}(1))$ is the copy of $V\equiv
\Lab(p_{21}(j))\equiv \Lab(p_{42}(j)^{-1})$ by the definition of
the mapping $\varepsilon_1$, and $\Lab(p_{42}(1))$ is obtained
from $V$ by deleting all $a$-letters. However, all $a$-letters of
$V$ are contained in the subword $U$ in $\bar a_1,\dots, \bar a_m$
written between the $P_j$- and the $R_j$-letters of $V$ (see lemma
\ref{Band}) because $\aaa$ and $\bar\aaa$ have no other common
letters with index $j$. We are going to prove that $U=1$ modulo
$\gr$-relations. This will allow us to construct an auxiliary
intermediate diagram consisting of $\gr$-cells and having the
boundary label $p_{21}(1)p_{42}(1)$. Thus the desired diagram
$\Gamma(1)$ will be built from $\Gamma_1(1),\Gamma_2(1),\dots$ and
the intermediate diagrams consisting of $\gr$-cells.

It remains to show that $U=1 (mod \gr ).$ For this purpose we
apply homomorphism $\delta$ from Lemma \ref{alpha} to the boundary
label $W_1$ of $\Gamma_1(j)$. Since the all $a$-letters of $W_1$
occur in the subword $U$, we have $U=1$ in $\bar\gr$. By Lemma
\ref{Emb}, $U=1$ modulo $\gr$-relations as desired. (When arguing
in this way on a diagram $\Gamma_s(j)$, $s>1$, we have to induct
on $s$ and take into account that, by the inductive hypothesis,
$\delta$ sends the complement of $U$ in the boundary label $W_s$
of $\Gamma_s(j)$ to 1.)

The last statement of the lemma clearly follows from the
construction of $\Delta'$. \endproof

\begin{lemma} \label{2 spokes}
Assume that a Van Kampen or an annular diagram $\Delta$ over the
group $\hhh$ contains two hubs, $\Pi_1$ and $\Pi_2$, and two
bands, a $K_j$-band $\bb_1$ and a $K_{j+1}$-band $\bb_2$
connecting the hubs  where $j\ne 1$. Suppose that the diagram
bounded by the two $k$-bands and the boundaries of $\Pi_1$,
$\Pi_2$ does not contain any other hubs and does not contain the
hole of the diagram (in the annular case). Then the diagram
$\Delta$ is not minimal: one can decrease the number of hubs by 2
while preserving the boundary label.
\end{lemma}

\begin{center}
\unitlength=1mm \special{em:linewidth 0.4pt} \linethickness{0.4pt}
\begin{picture}(117.85,61.12)
\put(30.07,53.57){\oval(15.96,15.10)[]}
\put(30.07,12.25){\oval(15.96,15.10)[]}
\put(29.50,53.43){\makebox(0,0)[cc]{$\Pi_2$}}
\put(29.50,11.82){\makebox(0,0)[cc]{$\Pi_1$}}
\bezier{156}(23.23,49.15)(12.69,35.47)(22.95,16.38)
\bezier{136}(24.66,47.44)(16.11,34.62)(24.37,18.09)
\bezier{132}(34.35,47.16)(42.33,35.47)(35.77,18.09)
\bezier{148}(36.63,48.87)(46.03,35.19)(37.20,16.38)
\bezier{172}(17.53,50.58)(30.07,33.48)(43.18,50.58)
\put(19.81,50.86){\makebox(0,0)[cc]{$\ttt$}}
\put(14.11,30.34){\makebox(0,0)[cc]{$\bb_1$}}
\put(44.61,30.34){\makebox(0,0)[cc]{$\bb_2$}}
\bezier{40}(51.16,32.05)(55.44,35.47)(58.57,32.05)
\bezier{40}(58.57,32.05)(62.56,28.92)(66.27,32.05)
\put(65.70,32.05){\vector(1,0){1.71}}
\put(89.63,53.57){\oval(15.96,15.10)[]}
\put(89.06,53.43){\makebox(0,0)[cc]{$\Pi_2$}}
\put(89.63,38.47){\oval(15.96,15.10)[]}
\put(89.06,38.04){\makebox(0,0)[cc]{$\Pi_1$}}
\bezier{172}(87.64,46.02)(61.99,51.72)(73.10,39.75)
\bezier{164}(89.06,46.02)(117.85,50.01)(112.43,39.75)
\bezier{60}(73.10,39.46)(76.81,35.76)(79.37,26.07)
\bezier{68}(112.43,39.75)(105.59,32.34)(107.87,26.07)
\bezier{216}(79.37,25.78)(93.62,2.99)(107.87,26.07)
\end{picture}

\end{center}

\begin{center}
\nopagebreak[4] Fig. \theppp.

\end{center}
\addtocounter{ppp}{1}

\proof
We can assume that there are no cells in $\Delta$ distinct from
the hubs, the cells of the two $k$-bands, and the cells entirely
lying between the two $k$-bands. In particular, $\Delta$ is simply
connected. We can assume also that the subdiagram $\Delta_0$ that
consists of the same cells with the exception of the hubs is
minimal.

It is clear that $\Delta_0$ is a trapezium with base
$\lel_jL_jP_jR_j\rer_j$ or it has no $\theta$-cells. In the last
case the $P_j$-edges of $\Pi_1$ and $\Pi_2$ must coincide. Thus
$\Pi_1$ and $\Pi_2$ cancel as desired.

Otherwise we consider the nearest $\theta$-band $\ttt$ to $\Pi_2$.
Its top has no $a$-edges because otherwise an $a$-band starting on
$\ttt$ must end on the hub $\Pi_2$, a contradiction. The we apply
Lemma \ref{LMqR} to the diagram bounded by $\Pi_1, \bb_1, \bb_2$
and $\ttt$. It says the, preserving the boundary label, we can
change $\Delta$ by a subdiagram with two hubs (see the last
statement of Lemma \ref{LMqR}) which have a common $P_j$-edge.
Then we cancel them as in the first case.


\endproof



     We extend  the  notion  of  \label{lgr}$l$-{\it graph}  as compared with
[Ol,97], \cite{SBR}, \cite{BORS}, \cite{talk}, \cite{OlSaBurns},
because here we need such planar graphs not only in a plane but
also in an annulus. In the annular case, we have to introduce 2
exterior vertices $v_0$ and $v^0$ lying respectively in the outer
and in the inner plane components of the complement of the
annulus. By definition of $l$-graph $\Gamma$, each of interior
vertices $v_1,\dots,v_n$ is of degree at least $l$, where $l\ge
6$, $n\ge 1$, and $\Gamma$ has neither loops nor 2-gons which both
vertices are interior.

    As in \cite{Ol97}, \cite{SBR}, \cite{BORS} one easily obtains
the following corollary of the Euler formula for planar graphs.

      \begin{lemma} \label{Euler}
      For every $l$-graph $\Gamma$ there exists an interior
vertex $o$ of degree $d\ge l$ connected with $v_0$ by $d_0$
successive edges (such that there are no vertices of $\Gamma$
between them) and with $v^0$ by $d^0$ successive edges where
$d_0+d^0\ge d-3$. The number of vertices $n$ does not exceed the
sum of degrees of $v_0$ and $v^0$ if $l\ge 7$.
\end{lemma}

     Lemma \ref{Euler} can be applied to the hub graph
$\Gamma(\Delta)$ of a minimal diagram $\Delta$ which contains
hubs. By definition each the interior vertices of $\Gamma(\Delta)$
are chosen inside all hubs, and the exterior vertex (vertices) are
chosen in the connected components of the complement of $\Delta$
(one vertex per component) on the plane. The medians of all the
$K_j$-bands of $\Delta$ ($j\ne 1$) starting on hubs serve as the
edges of $\Gamma(\Delta)$.

The next Lemma is similar to Lemma 11.4 in \cite{SBR}, Lemma 9.4
\cite{Ol97} and Lemma 4.23 in \cite{BORS}. It immediately follows
from Lemma \ref{2 spokes}.

\begin{lemma} \label{Graph}
If a minimal diagram over $\hhh$ contains hubs, then the graph
$\Gamma(\Delta)$ is an $(N-1)$-graph.
\end{lemma}

\begin{lemma} \label{Small disc}
Let $\Delta$ be an annular diagram over $\hhh$ with contours $p$
and $p'$ or a \vk diagram with contour $p$. Assume $\Delta$
contains a hub $\Pi$ with the following properties.

1. There exists a $K_j$-band $\bb$ and a $K_{j+1}$-band $\bb'$
connecting $\partial(\Pi)$ with $p$, $j\ne 1$.

2. The van Kampen subdiagram $\Gamma$, which is bounded by $p$,
$\Pi$ and $\bb, \bb'$ and includes $\bb, \bb'$ is reduced and
contains no hubs.

Then the lengths of the closest to $p$ $\theta$-band $\ttt$ that
crosses  both $\bb$ and $\bb'$ (if any exists) is $O(|p|^2)$, and
there is a path $t$ of length $O(|p|)$ in $\Gamma$ connecting the
$P_j$-edge of $\ttt$ (or the $P_j$-edge of $\Pi$ if the
$\theta$-band $\ttt$ does not exist) with $p$.
\end{lemma}

\proof Let $\Gamma'$ be the subdiagram of $\Gamma$ bounded by
$\topp(\ttt)$ (or by $\Pi$ if there is no such $\ttt$), by $\bb$,
$\bb'$ and $p$. By Lemma \ref{NoAnnul} and by conditions 2,3, all
maximal $\theta$-band of $\Gamma'$ must start or end on $p$. Hence
the number of such bands is at most $|p|.$ There exist at most
$|p|$ maximal $k$-bands in $\Gamma'$ by Lemma \ref{NoAnnul}. Since
there exist no $(k,\theta)$-annuli by Lemma \ref{NoAnnul}, the
number of $(k,\theta)$-cells in $\Gamma'$ is $O(|p|^2)$.  Since
$j\ne 1$, the band $\ttt$ has no $a(P_1)$-cells. Hence every
maximal $a$-band starting on $\ttt$ must end on a
$(k,\theta)$-cell or on $p$. Then Lemma \ref{Band} gives the
desired upper bound for the length of $\ttt$.

The maximal $P_j$-band $\bb_1$, that starts on the hub, must end
on $p$ and intersect $\ttt$ (if $\ttt$  exists). It consists of
$(P_j,\theta)$-cells only. Since any two of these cells cannot
belong to the same $\theta$-band (by Lemma \ref{NoAnnul}), the
length of $\bb_1$ (or of its part belonging to $\Gamma'$) is
bounded by the number ($O(|p|)$) of maximal $\theta$-bands in
$\Gamma'$. Hence a subpath of a side of $\bb_1$ is the path $t$
that connects $\ttt$ (or $\partial(\Pi)$) and $p$ and has length
$O(|p|)$.\endproof

\begin{center}

\unitlength=1mm \special{em:linewidth 0.4pt} \linethickness{0.4pt}
\begin{picture}(104.17,66.54)
\put(57.15,15.95){\oval(18.24,17.38)[]}
\bezier{72}(17.82,25.50)(17.53,34.33)(26.08,37.18)
\bezier{60}(26.08,37.18)(33.21,37.75)(35.49,45.16)
\bezier{112}(35.49,45.16)(38.34,59.41)(50.88,54.85)
\bezier{76}(50.88,54.85)(59.14,49.44)(66.27,54.57)
\bezier{108}(66.55,54.57)(77.09,61.69)(86.50,50.86)
\bezier{68}(86.78,50.86)(89.92,46.02)(95.62,36.90)
\bezier{80}(95.62,36.90)(103.31,24.08)(104.17,18.66)
\put(48.03,13.53){\line(-5,3){29.35}}
\put(48.03,16.67){\line(-5,3){27.93}}
\put(55.15,24.65){\line(0,1){28.21}}
\put(59.14,24.65){\line(0,1){27.36}}
\put(66.27,18.66){\line(1,1){24.79}}
\put(66.27,14.10){\line(1,1){26.79}}
\put(16.39,32.91){\makebox(0,0)[cc]{$\bb$}}
\put(56.86,55.71){\makebox(0,0)[cc]{$\bb_1$}}
\put(94.19,44.59){\makebox(0,0)[cc]{$\bb'$}}
\bezier{448}(26.65,23.22)(55.72,66.54)(95.62,22.37)
\bezier{428}(28.08,21.51)(57.15,62.83)(94.19,20.09)
\put(43.18,48.87){\makebox(0,0)[cc]{$\Gamma'$}}
\put(59.14,46.02){\vector(0,1){3.99}}
\put(59.14,44.85){\circle*{1.57}}
\put(61.42,48.01){\makebox(0,0)[cc]{$t$}}
\put(42.04,32.91){\makebox(0,0)[cc]{$\ttt$}}
\put(80.23,56.28){\vector(3,-2){3.13}}
\put(83.08,56.85){\makebox(0,0)[cc]{$p$}}
\put(57.15,16.10){\makebox(0,0)[cc]{$\Pi$}}
\end{picture}

\end{center}

\begin{center}
\nopagebreak[4] Fig. \theppp.

\end{center}
\addtocounter{ppp}{1}

\begin{lemma} \label{Disk off} Let $V$ and $V'$ be two words which are conjugate in the group
$\hhh$, and a reduced annular diagram $\Delta$ for this conjugacy,
contains hubs. Then there is a word $U$ of length
$O((|V|+|V'|)^2)$ such that (a) $U=1$ in $\hhh$, and a minimal
diagram for this equality contains one hub, (b) a minimal diagram
for conjugacy of $VU$ and $V'$, or $V$ and $V'U$, has fewer hubs
than $\Delta.$
\end{lemma}

\proof Let $p$ and $p'$ be contours of $\Delta$ labelled by $V$
and $V'.$ Assume $\Pi$ is a hub given by lemmas \ref{Euler} and
\ref{Graph}. Since $N\ge 8$, the hub is connected  with one of
contours $p$, $q$ by a pair of consecutive $K$-bands $\bb$ and
$\bb'$ which are not $K_1$-bands. Without loss of generality, we
assume that these two bands ends on $p$. Besides, there are no
other hubs between these $k$-bands. Denote by $\Gamma$ the
subdiagram bounded by $\Pi,$ by $p$ and by these two $k$-bands.

Let $\ttt$ be the closest to $p$ maximal $\theta$-band of
$\Gamma$. It is of length $O(|p|^2)$ (if $\ttt$ exists) by Lemma
\ref{Disk off} and the $P_j$-edge of $\ttt$ (or of $\Pi$) can be
connected with $p$ by a path $t$ of length $O(|p|)$. Lemma
\ref{LMqR} says that preserving the subdiagram $\Gamma'$ of
$\Gamma$ which is bounded by $\ttt$ and $p$, we can insert a
number of non-hub cells (and their mirror images) in $\Delta$, and
obtain a diagram $\bar\Delta$ (with the same boundary label as
$\Delta$) containing a subdiagram $\Delta'$ which (1) attached to
$\Gamma'$ along the top of $\ttt$; (2) contains exactly one hub;
(3) has linearly bounded perimeter as function of $|\ttt|$, and
therefore as function of $|V|$, since $|V|=|p|$.

Then we consider the path $t$, given by Lemma \ref{Small disc},
and the boundary path $y$ of $\Delta'$. One may suppose that
initial vertices of the paths $y$ and $t$ coincide. Let $p_1$ is a
cyclic permutation of $p$, starting with the terminal vertex of
$t$. Then the length of the null-homotopic path $t^{-1}yt$ is
$O(|p|^2)$, and the path $p_1ty^{\pm 1}t^{-1}$ bounds an annular
diagram with smaller number of hubs than $\Delta.$ Thus, to prove
the lemma, we just denote the label of $p_0ty^{\pm
1}t^{-1}p_0^{-1}$ by $U$ where $p_0$ connects the initial vertices
of $p$ and $t$.
\endproof

\section{The end of the proof}
\label{eop} \setcounter{equation}{0}

\begin{lemma} \label{conjK0}
The conjugacy problem is decidable for $\hhh_1$.
\end{lemma}

\proof Let $\Delta$ is any annular diagram over $\hhh_1$. We need
to find another diagram with the same boundary labels and
recursively bounded number of cells. By Lemma \ref{rollsprep2}, we
can assume that $\Delta$ is reduced and has minimal boundaries.

By Lemma \ref{Main1}, we can assume that the boundary of $\Delta$
contains $\theta$-edges. Now if $\partial(\Delta)$ contains no
$k$-edges then $\Delta$ is a roll and we can refer to Lemma
\ref{rolls}. Otherwise $\Delta$ is a union of a spiral and a
recursively bounded number of cells with recursively bounded
perimeters (by Lemma \ref{sp123}), and we can refer to Lemma
\ref{spiral}.
\endproof

\begin{lemma} \label{final}
The conjugacy problem is decidable for $\hhh$.
\end{lemma}

\proof It suffices to prove that there is an annular diagram over
$\hhh$ for the conjugacy of any two given conjugated words $w_1$
and $w_2$, whose number of cells and perimeters of cells are
recursively bounded as function of $|w_1|+|w_2|$. Consider a
minimal annular diagram $\Delta$ for such a conjugacy. The hub
graph $\Gamma(\Delta)$ is a $N-1$ graph by Lemma \ref{Graph}. By
Lemma \ref{Euler}, the number of hubs in $\Delta$ does not exceed
the number of edges of $\Gamma(\Delta)$, which are incident with
the two exterior vertices of $\Gamma(\Delta)$ because $N\ge 8$,
that is, the number of hubs is at most the sum of the lengths of
boundary paths $p$ and $q$ of $\Delta$.

Then Lemma \ref{Disk off} reduces the problem to the conjugacy
problem for the group $\hhh_1$ because the words $U$ from that
lemma are conjugate in $\hhh_1$ of the hub. Thus, the statement
follows from Lemma \ref{conjK0}.
\endproof

The next lemma completes the proof of Theorem \ref{Col}

\begin{lemma}\label{Frattini} The embedding of $\gr$ into $\hhh$
is a Frattini embedding.
\end{lemma}

\proof Let $u, v$ be reduced words over $a_1(P_1),...,a_m(P_1)$
(the generators of the copy of $\gr$ in $\hhh$ by Lemma
\ref{Emb}). Suppose that $u$ and $v$ are conjugate in $\hhh$.
Consider a minimal conjugacy diagram $\Delta$ with boundary labels
$u$ and $v$. Since $u$ and $v$ do not contain $k$-letters, by
Lemmas \ref{Euler} and \ref{Graph}, $\Delta$ contains no hubs, so
$\Delta$ is a diagram over $\hhh_1$.

Without loss of generality we can assume that there are no
$\gr$-cells which have common edges with the boundary of $\Delta$.
Indeed, otherwise we can remove these cells, replacing $u$ (resp.
$v$) by words that are conjugates of $u$ (resp. $v$) in $\gr$.

Suppose that $\Delta$ contains $k$-edges. Then by Lemma
\ref{NoAnnul}, $\Delta$ would contain $k$-annuli surrounding the
hole of $\Delta$. Since the boundary of $\Delta$ contains no
$\theta$-edges, the $k$-annuli contain no $\theta$-cells by Lemma
\ref{NoAnnul}. Therefore its boundary label is a word in $\xxx$.
By Lemma \ref{alpha} the $\delta$-images of $u$ and $v$ are
trivial in $\bar\gr$. Hence they are trivial in $\gr$ by Lemma
\ref{Emb}.


Suppose that $\Delta$ contains $\theta$-edges. Then $\Delta$
contains a $\theta$-annulus $\ttt$ surrounding the hole of
$\Delta$. Since $\Delta$ does not contain $k$-edges, all
$\theta$-edges belong to the same set $\Theta(z)$ or
$\bar\Theta(z)$. If $z\ne P_1$ then the $a$-bands starting on a
side of $\ttt$ cannot end on the boundary of $\Delta$ or on
$\gr$-cells. Therefore they must cross $\ttt$ twice, which
contradicts Lemma \ref{NoAnnul}. Therefore $z=P_1$. But then all
cells in $\ttt$ are commutativity cells corresponding to relations
(\ref{auxtheta}). Therefore the labels of the sides of $\ttt$
coincide and we can remove $\ttt$ from $\Delta$ reducing the
number of $\theta$-cells in $\Delta$. This contradicts to the
minimality of $\Delta$. Thus $\Delta$ contains no $\theta$-edges.

But this means that there are no cells in $\Delta$ that have
common edges with the contours of $\Delta$. Thus $\Delta$ has no
cells, and $u$ is a cyclic shift of $v$ (since $u, v$ are
reduced). Hence $u$ and $v$ are conjugate in $\gr$.\endproof

As we mentioned before, the {\bf proof of Theorem \ref{Col1}}
proceeds in almost the same way. We need to show that the
conjugacy problem in $\hhh$ is Turing reducible to the conjugacy
problem in $\gr$. Thus we have to expand the class of recursive
functions by the characteristic function of the set of pairs of
words $(u,v)$ which are conjugate in $\gr$ to the set of
elementary recursive functions, and then apply the usual operators
used to produce all recursive functions \cite{Mal}. Let us call
such functions $\gr$-recursive. For example, the word problem in
$\gr$ is $\gr$-recursive, by Clapham's theorem \cite{Cla} the word
problem in $\bar\gr$ is also $\gr$-recursive and so on.  The
reader can check that all the lemmas in this paper remain true if
we replace the word ``recursive" in their formulations by
``$\gr$-recursive". It is easy to see that this modification turns
the proof of Theorem \ref{Col} into a proof of Theorem \ref{Col1}.

\twocolumn

\noindent {\bf Subject index} \label{sind}
\bigskip

\noindent Accepted word \pr{acceptedw}\\ Admissible words of an
$S$-machines \pr{Admissible}\\ Admissible words of the $S$-machine
$\sss$ \pr{adm}\\ Admissible words of the $S$-machines $\bsss$ and
$\csss$ \pr{admbar}\\ Age of $\omega$ \pr{age}\\ Annular (Schupp)
diagram \pr{annular}\\ Annulus \pr{annulus}\\ Auxiliary relations
\pr{auxiliary}\\ Band \pr{Sband}\\ $k$-, $\theta$-, $a$-bands
\pr{ktabands}\\ Base of a $\theta$-band \pr{baseb}\\ Base of an
admissible word \pr{baseadm}\\ Base of a computation \pr{basec}\\
Base of a roll \pr{baser}\\ Base of a trapezium \pr{baset}\\ Basic
letters \pr{Basicletters}\\ $z$-best pairing \pr{zbest}\\ Bottom
path of a band \pr{bottt}\\ Bottom path of a trapezium \pr{bott}\\
Boundary of a band \pr{boundband}\\ Brief history of a computation
\pr{brh}\\ Cancellation pairing \pr{cancpar}\\ $(a,x)$-, $(k,x)$-,
$(\theta,a)$-, $(\theta,k)$-cell \pr{axcell}\\ $\gr$-cell
\pr{grcell}\\ Compressible diagram \pr{compressible}\\ Computation
associated with a trapezium \pr{compt}\\ $\gr$-computation
\pr{gcomp}\\ Conditions (R1)-(R4) \pr{RRRR}\\ Connected pair of
parentheses \pr{connectedp}\\ Connected pair of letters
\pr{connectedl}\\ Contracting (expanding) bands \pr{contracting}\\
$\bee$-coordinate \pr{beeco}\\ $\Omega$-coordinate \pr{omegaco}\\
Data associated with a trapezium \pr{datat}\\ $a$-, $k$-,
$\theta$-, $x$-edge (letter, cell) \pr{aedge}\\ Empty band
\pr{emptyband}\\ Fractional letters (words) \pr{fractional}\\
Frattini embedding \pr{Frattinif}\\ Free computation \pr{freec}\\
 Graded presentation \pr{gradedpres}\\
$z$-good pairing \pr{zgood}\\ $l$-graph \pr{lgr}\\ Height of a
trapezium \pr{heightt}\\ Historical period \pr{histper}\\ History
of a computation \pr{hiscomp}\\ History of a roll \pr{historyr}\\
History of a trapezium (ring) \pr{historyt}\\ Hub \pr{hub}\\ Inner
diagram of an annulus \pr{innerdiag}\\ Inner part of a sector
\pr{innerpart}\\ Inverse rule \pr{Inverserule}\\ Locked sectors
(by a rule) \pr{locked}\\
 $S$-machine
\pr{Smachines}\\ Main $k$-band of a spiral \pr{mainsp}\\ Main
relations \pr{mainrelations}\\ Maximal band \pr{maximalb}\\
Minimal boundary \pr{minimalb}\\ Minimal diagram \pr{minimald}\\
Minimal word \pr{minimalw}\\ Minus pair of letters
\pr{minuspair}\\ Minus pairing \pr{minusp}\\ Minus word
\pr{minusw}\\ Normal pair of letters \pr{normalp}\\ Operation of
moving $K_j$-annuli \pr{opmov}\\ Operation of removing
$K_j$-annuli \pr{oprem}\\ $(2,4)$-pairing \pr{24pairing}\\ Plus
pair of letters \pr{pluspair}\\ Positive rule \pr{pnr}\\
$(2,4)$-projection of a word \pr{24projection}\\ Quasiring
\pr{quasiringq}\\ Quasispiral \pr{qspirals}\\ Quasitrapezium
\pr{quasitrapezium}\\ Rank of a cell \pr{rankcell}\\ Reduced
computation \pr{reducedc}\\ Reduced or non-reduced diagram
\pr{reducedd}\\ Reducible pair of cells \pr{reduciblec}\\
$\tau$-regular word \pr{taur}\\ Related uniform words
\pr{relatedu}\\ $(a,x)$-, $(k,x)$-, $(\theta,a)$-,
$(\theta,k)$-relations \pr{axrel}\\ $\gr$-relations \pr{grrel}\\
Ring \pr{ringgg}\\ Ring computation \pr{ringcom}\\ Roll
\pr{roll}\\ $S$-rule \pr{Srule}\\ $S$-rule active with respect to
a sector \pr{activerule}\\ $S$-rule applicable to a word
\pr{applicable}\\ Sector of an admissible word \pr{sector}\\
$a$-similar rules \pr{asimilar}\\ Small trapezium \pr{smallt}\\
Spiral \pr{spirals}\\ Standard computation \pr{standardc}\\ State
letters \pr{State}\\ Tame computation \pr{tamec}\\ Tape letters
\pr{Tapeletters}\\ Top path of a band \pr{toppp}\\ Top path of a
trapezium \pr{topt}\\ Transition rules \pr{trrules}\\ Trapezium
\pr{trapeziumt}\\ Trapezium of the first (second, mixed) type
\pr{fsmtype}\\ Type of a diagram \pr{typed}\\ Type of a path
\pr{typep}\\ Uniform word \pr{uniformw}\\ \vk diagram \pr{vkd}\\
Wild computation \pr{wild}\\ \vskip 0.2 in

\noindent $\tool$ \pr{tool}\\ $\equiv$ \pr{equiv}\\
 $\aaa$ \pr{aaa}\\ $\aaa(z)$
\pr{aaaz}\\ $\bar\aaa$ \pr{bkta}\\ $\alpha_\tau$ \pr{alphatau}\\
$\beta$ \pr{beta}\\ $\bott(.)$ \pr{bottq}\\ $\br(.)$ \pr{brh1}\\
$\delta$ \pr{delta}\\ $\diff(.)$ \pr{diff}\\ $\eee$ \pr{eee}\\
$\bee$ \pr{beee}\\ $\gr$ \pr{grrr}\\ $\bar\gr$ \pr{bargr1}\\
$\gamma$ \pr{gamma}\\ $\hhh$ \pr{hhhh}\\ $\hhh_0$ \pr{hhh0}\\
$\hhh_1$ \pr{hhh1}\\ $\hhh_2$ \pr{hhh2}\\ $\kkk$ \pr{kkk}\\ $\tkk$
\pr{tkk}\\ $\bar\kkk$ \pr{bkta}\\ $\bar\tkk$ \pr{bkta}\\ $\lel_j$
\pr{lel}\\ $\rer_j$ \pr{rer}\\ $\Sigma$ \pr{sigmaaaa}\\ $\tsigma$
\pr{tsigmaaaa}\\ $\Sigma(w)$ \pr{sigmaw}\\
$\Sigma(w_1,w_2,w_3,w_4)$ \pr{sigmawwww}\\
$\bar\Sigma(w_1,w_2,w_3,w_4)$ \pr{barsigmawwww}\\
$\Sigma_{r,i}(w_1,w_2,w_3,w_4)$ \pr{sigmariwwww}\\
$\bar\Sigma_{r,i}(w_1,w_2,w_3,w_4)$ \pr{barsigmariwwww}\\
$\sss(\omega)$ \pr{somega}\\ $\Theta$ \pr{theta}\\ $\Theta(\tau)$
\pr{thetatau}\\ $\Theta(z)$ \pr{thetaz}\\ $\bar\Theta$
\pr{btheta}\\ $\topp(.)$ \pr{bottq}\\ $U(r,i)$ \pr{uri}\\$\xxx$
\pr{xxx}\\ $\xxx(\tau)$ \pr{xxxtau}\\ $\xxx(z)$ \pr{xxxz}\\
$\xxx(z,\tau)$ \pr{xxxztau}\\ $z_-$ \pr{zmin}\\ $z_+$ \pr{zplus}\\

 \onecolumn
\begin{minipage}[t]{2.9 in}
\noindent Alexander Yu. Ol'shanskii\\ Department of Mathematics\\
Vanderbilt University \\ alexander.olshanskiy@vanderbilt.edu\\
 and\\ Department of
Higher Algebra\\ MEHMAT\\
 Moscow State University\\
olshan@shabol.math.msu.su\\
\end{minipage}
\begin{minipage}[t]{2.6 in}
\noindent Mark V. Sapir\\ Department of Mathematics\\
Vanderbilt University\\
msapir@math.vanderbilt.edu\\
\end{minipage}

\addtocontents{toc}{\contentsline {section}{\numberline {
}{\noindent References} \hbox {}}{\pageref{bibl}}}

\addtocontents{toc}{\contentsline {section}{\numberline { }Subject
index \hbox {}}{\pageref{sind}}}


\begin{thebibliography}{99}

\label{bibl}
\newcommand{\bi}{\bibitem}
\newcommand{\nb}{\newblock}
\bi[AC]{AC} S. Aanderaa, D.E. Cohen. Modular machines and the
Higman-Clapham-Valiev embedding theorem. Word problems, II (Conf.
on Decision Problems in Algebra, Oxford, 1976), pp. 17--28, Stud.
Logic Foundations Math., 95, North-Holland, Amsterdam-New York,
1980.


\bi[BORS]{BORS} J. C. Birget, A.Yu. Ol'shanskii, E.Rips, M. V.
Sapir. \newblock Isoperimetric functions of groups and
computational complexity of the word problem. Annals of
Mathematics, 156, 2 (2002), 467-518


\bi[Bo]{Bo} W.~W. Boone.
\newblock Certain simple unsolvable problems in group theory .
\newblock {\it Proc. Kon. ned. akad. Wetensch. A}, (I) 57 (1954), 231--237,
(II) 57 (1954), 492--497, (III) 58 (1955), 252--256, (IV) 58
(1955), 571--577, (V) 60 (1957), 22--27, (VI) 60 (1957), 222--232.

\bi[Cla]{Cla} C. R. J. Clapham. \nb An embedding theorem for
finitely generated groups", Proc. London. Math. Soc. (3), 17,
1967, 419-430.

\bibitem[Col]{Collins} Donald J. Collins. \nb Conjugacy and the
Higman embedding theorem. Word problems, II (Conf. on Decision
Problems in Algebra, Oxford, 1976), pp. 81--85, Stud. Logic
Foundations Math., 95, North-Holland, Amsterdam-New York, 1980.


\bi[CM]{CM} D. J. Collins, C. F. Miller III. \nb  The conjugacy
problem and subgroups of finite index. Proc. London Math. Soc. (3)
34 (1977), no. 3, 535--556.

\bi[FW]{FW} N.J.Fine and H.S.Wilf, Uniqueness theorems for
periodic functions, Proc. AMS 16 (1965), 109-114

\bi[GK]{GK} A.V. Gorjaga, A.S.  Kirkinski\u\i.\nb The decidability
of the conjugacy problem cannot be transferred to finite
extensions of groups. (Russian) Algebra i Logika 14 (1975), no. 4,
393--406.

\bi[Hi]{Hi} G. Higman. Subgroups of finitely presented groups.
Proc. Roy. Soc. Ser. A, 262 (1961), 455--475.

\bi[Kal]{Kal} K.A. Kalorkoti. Decision problems in group theory.
Proc. London Math. Soc. (3) 44 (1982), no. 2, 312--332.

\bibitem[KS]{KharlSap}
O. G. Kharlampovich and M.V. Sapir.
\newblock Algorithmic problems in varieties.
\newblock Internat. J. Algebra Comput. 5 (1995), no. 4-5,
379--602.

\bibitem[KT]{KT} {\em Kourovka Notebook}. Unsolved Problems in Group Theory.
5th edition, Novosibirsk, 1976.

\bi[LS]{LS} Roger Lyndon and Paul Schupp.
\newblock {\it Combinatorial group theory}.
\newblock Springer-Verlag, 1977.


\bi[Mak]{Makanin} G. S. Makanin.
\nb Equations in a free group.
(Russian) Izv. Akad. Nauk SSSR Ser. Mat. 46 (1982), no. 6,
1199--1273, 1344.

\bi[Mal]{Mal} A.I. Mal'cev. Algorithms and recursive functions.
Translated from the first Russian edition by Leo F. Boron, with
the collaboration of Luis E. Sanchis, John Stillwell and Kiyoshi
Iséki Wolters-Noordhoff Publishing, Groningen 1970.

\bi[Ma]{Ma} Yu. I. Manin. The computable and the non-computable.
(Vychislimoe i nevychislimoe). (Russian) [B] Kibernetika. Moskva:
"Sovetskoe Radio".

\bi[Mil]{Mil} Charles F. Miller III. \nb On group-theoretic
decision problems and their classification. Annals of Mathematics
Studies, No. 68. Princeton University Press, Princeton, N.J.;
University of Tokyo Press, Tokyo, 1971.

\bi[Nov]{No} P.S. Novikov, On the algorithmic insolvability of the
word problem in group theory. American Mathematical Society
Translations, Ser 2, Vol. 9, pp. 1--122



\bi[Ol1]{Ol_book} A.~Yu. Ol'shanskii.
\newblock {\it The geometry of defining relations in groups},
\newblock Nauka, Moscow, 1989.

\bi[Ol2]{Ol97} A. Yu. Ol'shanskii.
\newblock On distortion of subgroups in finitely presented groups. Mat.
Sb., 1997, V.188, N 11, 51-98.

\bi[OlSa1]{OlSaBurns} A. Yu. Ol'shanskii, M. V. Sapir. Embeddings
of relatively free groups into finitely presented groups. Contemp.
Math., 264, 2000, 23-47.

\bi[OlSa2]{talk} A.Yu. Ol'shanskii, M. V. Sapir. Length and area
functions on groups and quasi-isometric Higman embeddings.
Internat. J. Algebra Comput. 11 (2001), no. 2, 137--170.

\bi[OlSa3]{OSamen} A.Yu. Olshanskii, M.V. Sapir. Non-amenable
finitely presented torsion-by cyclic groups, Publications of IHES,
\# 96, 2002.

\bi[OlSa4]{OlSa} A.Yu. Olshanskii, M.V. Sapir. Word, power, order
and conjugacy problems in groups. (in preparation)



\bi[Rot]{Rotman} J. Rotman.
\newblock {\it An Introduction to the Theory of Groups}.
\newblock Allyn \& Bacon, third edition, 1984.


\bi[SBR]{SBR} M. V. Sapir, J. C. Birget, E. Rips.
\newblock Isoperimetric and isodiametric functions of groups,
Annals of Mathematics, 157, 2 (2002), 345-466.


\bi[Tho]{Tho}R. J. Thompson.\nb Embeddings into finitely generated
simple groups which preserve the word problem. Word problems, II
(Conf. on Decision Problems in Algebra, Oxford, 1976), pp.
401--441, Stud. Logic Foundations Math., 95, North-Holland,
Amsterdam-New York, 1980.

\bi[Va]{Va} M.K.Valiev. On polynomial reducibility of the word
problem under embedding of recursively presented groups in
finitely presented groups. Mathematical foundations of computer
science 1975 (Fourth Sympos., Mariánské Lázn\v e, 1975), pp.
432--438. Lecture Notes in Comput. Sci., Vol. 32, Springer,
Berlin, 1975.

\end{thebibliography}
\end{document}